\documentclass[12pt]{article}
\usepackage[latin1]{inputenc}
\usepackage{amsfonts,amssymb,amsmath, epsfig}
\usepackage{color,graphicx,graphics}
\usepackage{amsthm}
\usepackage{amsmath,amstext,amssymb,amsfonts, amscd}
\usepackage[lofdepth,lotdepth]{subfig}
\usepackage{xcolor}
\usepackage{hyperref}          
\usepackage{stix}
\usepackage{empheq}
\usepackage{tikz-cd}
\usepackage{circledsteps}
\usepackage{accents}

\def\Box{\leavevmode\vbox{\hrule
     \hbox{\vrule\kern4pt\vbox{\kern4pt}%
           \vrule}\hrule}}
\def\blackbox{\leavevmode\vrule height 5pt width 4pt depth 0pt\relax}
\def\endproof{\null\hfill {$\blackbox$}\bigskip}

\def\paragraph#1{{\bf #1\ }}

\newtheorem{lemma}{Lemma}[section]  

\newtheorem{theorem}[lemma]{Theorem}

\newtheorem{corollary}[lemma]{Corollary}

\newtheorem{definition}[lemma]{Definition}

\newtheorem{proposition}[lemma]{Proposition}

\newtheorem{remark}{Remark}[section]

\newtheorem{example}{Example}[section]

\newtheorem{conjecture}{Conjecture}[section]

\title{Swarming by curvature control in arbitrary dimension} 
\author{P. Degond$^{(1)}$, A. Diez$^{(2)}$, A. Frouvelle$^{(3)}$} 
\date{} 
\begin{document}

\maketitle


\begin{center}
(1) Institut de Math\'ematiques de Toulouse ; UMR5219 \\
Universit\'e de Toulouse ; CNRS \\
UPS, F-31062 Toulouse Cedex 9, France\\
email: pierre.degond@math.univ-toulouse.fr

\bigskip

(2) RIKEN Center for Interdisciplinary Theoretical and Mathematical Sciences (iTHEMS), \\
RIKEN iTHEMS, Wako, Saitama 351-0198, Japan\\ 
email: antoine.diez@riken.jp

\bigskip

(3) CEREMADE, CNRS ; UMR 7534\\
Universit\'e Paris Dauphine -- PSL, 75016 Paris, France \\
email: frouvelle@ceremade.dauphine.fr

\end{center}

\vspace{0.5 cm}
\begin{abstract}
We consider an interacting particle system proposed in the literature to model fish behavior. In this model, the agents move at constant speed and control the curvature of their trajectory (i.e. the time-derivative of their velocity) so as to align their velocity with that of their neighbors, up to some noise. We provide a novel $n$-dimensional formulation of this model for any $n \geq 3$ and derive its mean-field kinetic formulation using bundle geometry concepts. The target of the paper is the derivation of a fluid model in the hydrodynamic limit. We show that this fluid model is the ``self-organized hydrodynamic'' (SOH) model already found in earlier work pertaining to the Vicsek model. The derivation is based on the introduction of appropriate ``generalized collision invariants'' (GCI). The action of the $n$-dimensional orthogonal group is used to reduce the expression of the GCI to a set of two functions satisfying a system of equations which is solved by means of a variational formulation. This leads to explicit formulas for the coefficients of the SOH model in terms of those of the original interacting particle system. 
\end{abstract}

\medskip
\noindent
{\bf Acknowledgements:} PD holds a visiting professor association with the Department of Mathematics, Imperial College London, UK.

\medskip
\noindent
{\bf Data statement:} No new data was generated in the course of this research. 

\medskip
\noindent
{\bf Conflict of interest:} The authors declare no conflict of interest. 

\medskip
\noindent
{\bf Key words: } Interacting particle system, kinetic equation, hydrodynamic model, tangent bundle,  generalized collision invariant, group action, orthogonal group, variational formulation

\medskip
\noindent
{\bf AMS Subject classification: } 35Q70, 35Q84, 36Q92, 82C31, 82C21, 82C22, 53B50, 57M60, 92D50
\vskip 0.4cm

\setcounter{equation}{0}
\section{Introduction}
\label{sec:intro}

Swarming is a feature shared by many living systems such as locust swarms \cite{sayin2025behavioral}, fish schools \cite{lin2025experimental} or semen \cite{creppy2016symmetry} to cite only a few. Systems of interacting motile individuals produce large scales structures such as density or orientation waves, cluster formation or oscillations. Such structures have large sizes in comparison with the interaction range of each agent. The mechanisms that generate them from local-scale agent interactions have been the object of intense scrutiny in the last three decades or so \cite{vicsek2012collective}. In spite of this, a lot of questions remain unanswered. Mathematical modelling provides a useful tool to unravel these mechanisms. 

An assembly of interacting particles or agents can be modelled at the microscopic scale by a system of coupled ordinary or stochastic differential equations (ODEs, or SDEs) for the agents' locations and velocities for instance. However, particle systems do not provide much qualitative insight into the formation of large-scale structures except by means of intensive numerical simulations \cite{chate2008collective, vicsek1995novel}. The situation is drastically different with continuum (or macroscopic, or fluid) models which encode average properties of the system such as the local particle density or mean velocity. Such fluid models models can be analyzed through e.g. stability analysis which provides orders of magnitudes of pattern sizes (see an application of this in e.g. \cite{barre2019modelling}). However, for the analysis to be relevant, a rigorous link between the particle system and the fluid model must be available, including a ``dictionary'' linking the physical constants of the particle system and those of the fluid model. Such a dictionary is seldom provided. 

One way to produce this dictionary is by means of the ``hydrodynamic limit'' (see e.g. \cite[Sect. 3.8]{cercignani2013mathematical} in the case of rarefied gas dynamics). This method is based on a scale separation assumption, assuming that the time and space scales of the evolution of the system as a whole are larger than those associated to particle interactions. These assumptions align well with the goal of understanding the emergence of the large-scale structures appearing in collective dynamics. With such scaling assumptions, the particle system, or rather, 
its statistical description through mean-field kinetic equations, becomes singularly perturbed by a small parameter and the hydrodynamic model is obtained in the limit of this parameter tending to zero. Other methods of derivation of hydrodynamic models exist, such as the moment method (see e.g. \cite{levermore1996moment}) but this method requires a closure assumption which is either arbitrary or phenomenological. By contrast, the hydrodynamic limit carries no assumption, but for the scaling assumptions made initially, which amount to select the applicability regime of the so-obtained hydrodynamic model. 

Among models of collective dynamics, the Vicsek model \cite{vicsek2012collective} occupies a prominent position with many follow-ups (see among others, \cite{bertin2006boltzmann, bertin2009hydrodynamic, chate2008collective, costanzo2018spontaneous, figalli2018global, gamba2016global, griette2019kinetic} and the review \cite{vicsek2012collective}). In this model, the particle speeds are supposed constant as a way to account for self-propulsion. Particles control their velocity (which is thus, up to scaling, a normalized vector) so as to align it with the average velocity of their neighbors, up to some noise reflecting the uncertainty in the control. A macroscopic version of the Vicsek model was proposed by Toner and Tu in \cite{toner1995long, toner1998flocks, toner2005hydrodynamics}. This phenomenological model does not provide any link between its coefficients and those of the Vicsek model, thus failing to establish the expected ``dictionary''. The hydrodynamic limit of the Vicek model was performed for the first time in \cite{degond2008continuum} (see also \cite{aceves2019hydrodynamic, degond2013macroscopic, degond2015phase}) and leads to a novel type of fluid models called ``Self-Organized Hydrodynamics (SOH)'' \cite{degond2011hydrodynamic} different from the Toner and Tu model \cite{toner1995long}. 

The swarming model by curvature-control proposed in the present paper is an elaboration of the Vicsek model. Instead of the velocity, the variable by which the particles control their motion is the time-derivative of the velocity (which is, up to scaling, equal to the curvature of the particles' trajectories, hence the terminology). Like in the Vicsek model, agents control their curvature so as to align their velocity with the average velocity of their neighbors. Again, noise is added to this curvature control to reflect uncertainties in the control. A preliminary version of this model for a single particle was originally proposed in 2 dimensions (2D) in~\cite{gautrais2009analyzing} based on the analysis of fish trajectories. This single particle model was shown to behave like a diffusion at large scales in \cite{degond2008large}. The complete, many-particle version of the curvature-control model was first proposed in 2D and its hydrodynamic limit, derived, in \cite{degond2011macroscopic}. Strikingly, its hydrodynamic limit appeared to be the same SOH model as for the Vicsek model, except for the expression of its coefficients which took into account the specifities of curvature-control. Later, A three-dimensional (3D) version of the particle curvature-control model was proposed in \cite{cavagna2015flocking}. 

In this paper, our goal is twofold. We first propose a $n$-dimensional (nD) version of the curvature-control model for arbitrary $n \geq 3$ (the case $n=2$ being already treated in \cite{degond2011macroscopic}). Compared with \cite{degond2011macroscopic}, this extension is, by far, non-trivial. Indeed, in 2D, the curvature is a scalar and the model is posed in standard flat space. In nD, ($n \geq 3$), the curvature is a normal vector to the velocity, so that the pair (velocity, curvature) lies on the tangent bundle to the sphere. It is thus necessary to use bundle geometry, which greatly complexifies the framework. In \cite{cavagna2015flocking}, this geometric complexity was overcome by the use of speficic concepts to 3D (namely cross product), not generalizable to higher dimensions. The second goal is to perform the hydrodynamic limit of the curvature-control model in nD, to show that, like in 2D, the corresponding fluid model is the SOH model (in its nD version), and to give an explicit link between its coefficients and the particle curvature-control model. A fluid model based on the curvature-control model of \cite{cavagna2015flocking} was proposed in \cite{yang2015hydrodynamics} using truncated moment expansions. It is significantly different from the SOH model (these differences will be outlined in Section \ref{subsec:hydro_main}). While, in the curvature-control model, particle interactions are of mean-field type, a 2D model built from data \cite{gautrais2012deciphering} has features reminiscent to Boltzmann-like interactions, but otherwise similar in spirit.  

Beside the Vicsek model and its variants, many other collective dynamics models have been proposed, such as the three-zones model \cite{aoki1982simulation, cao2020asymptotic} or the Cucker-Smale model \cite{aceves2019hydrodynamic, barbaro2016phase, carrillo2010asymptotic, cucker2007emergent, ha2009simple, ha2008particle, haskovec2021simple, motsch2011new} or the soft-core interaction model \cite{d2006self} to name only a few. Cucker-Smale and Vicsek type models can be connected via singular limits \cite{barbaro2012phase, bostan2013asymptotic, bostan2017reduced}.

As already outlined, the curvature-control model in nD requires a bundle geometric framework. Indeed, the particle system must be appropriately formulated in the tangent bundle to the sphere. The mean-field kinetic description of the particle system is a partial differential equation~(PDE) posed on this bundle, which also requires appropriate definitions of bundle differential operators. Classically, hydrodynamic limits strongly rely on appropriate conservation relations satisfied at the level of the mean-field kinetic model. However, here, there are not enough such conservation relations to conduct the hydrodynamic limit in this classical way. The  same problem was already found in the Vicsek case, and was overcome by the introduction of a weaker form of conservation relation called ``generalized collision invariant'' (GCI) \cite{aceves2019hydrodynamic, degond2013macroscopic, degond2015phase, degond2011hydrodynamic, degond2008continuum}. Here, the equations for the GCI are more complicated and do not lend themselves easily to analytic resolution. However, we show that the $n$-dimensional orthogonal group acts on the GCI and that the latter are invariant by one of its subgroups. This enables us to derive a reduced form of the GCI and to deduce explicit formulas for the coefficients of the resulting SOH model, hence producing the ``dictionary'' between the particle and fluid models as claimed. The GCI and their reduced form satisfy systems of PDE's which have novel variational formulations and which are shown to be uniquely solvable. 

\smallskip
To summarize, the main innovations contained in this work are 

\smallskip
\noindent
- a new particle curvature-control model in arbitrary dimension and its associated mean-field kinetic equations formulated in a bundle geometric framework,

\smallskip
\noindent
- the first establishment of its associated hydrodynamic limit, 

\smallskip
\noindent
- the use of group action to reduce the complexity of a key object involved in the hydrodynamic limit, the ``generalized collision invariant'' (GCI),

\smallskip
\noindent
-  the derivation and unique solvability of a novel kind of variational formulation for the reduced~GCI. 

\smallskip
The outline of the paper is as follows. The particle curvature control model is presented in Section \ref{sec:particle} and its mean-field kinetic formulation in Section \ref{sec:kinetic}. Then, the main result of the paper, the establishment of the hydrodynamic limit, is stated in Section \ref{sec:hydro}. The hydrodynamic model consists of two equations. The first one, the continuity equation, is easily   established in the same section. The second one, the velocity equation, requires the introduction of the GCI. This is the object of Section \ref{sec:gci}, where an existence and uniqueness theorem for the GCI is given. An action of the orthogonal group on the GCI is introduced in Section~\ref{sec:gract_gci}. Thanks to this action, a reduced form of the GCI consisting of two scalar functions (the reduced GCI pair) is derived in Section~\ref{sec:gci_redform}, and a variational formulation for the reduced GCI pair is introduced and studied. This finally enables us to derive the velocity equation in Section \ref{sec:macnorm_eq_j}, which completes the proof of the main theorem. A conclusion is drawn in section~\ref{sec:conclu}. Appendix \ref{sec:redGCIprf} provides the derivation of the reduced GCI pair, while Appendices \ref{sec:proofs_rem_formulas} and \ref{sec:macnorm_nD_vectGCI_eqs} give proofs of auxiliary results stated in the main text. Bundle geometry concepts are gradually introduced when needed, in separate subsections at the beginning of Sections \ref{sec:particle}, \ref{sec:kinetic}, \ref{sec:hydro}, \ref{sec:gract_gci} and \ref{sec:gci_redform}.

\setcounter{equation}{0}
\section{The particle curvature-control (PCC) model}
\label{sec:particle}

In this section, we derive the particle curvature-control (PCC) model which is the starting point of our study. We begin with an elementary presentation of the PCC model. Then, we introduce some geometric concepts which allow us to reformulate the PCC model in a geometric framework. The model is then scaled in dimensionless form.

\subsection{Elementary presentation of the PCC system}
\label{subsec:partelem}

We start by presenting the PCC model for a single particle before moving to a many-particle system. Then, we motivate the introduction of the geometric framework which will be fully developed in the forthcoming sections.

\subsubsection{Single-particle PCC model}
\label{subsubsec:partelem_single_noiseless}

We start with the single-particle noiseless PCC model and describe the motion of a particle under a ``force'' that acts on the curvature of its trajectory. Let a particle have position $X(t) \in {\mathbb R}^n$ at time $t \in [0,\infty)$, where $n \geq 3$ is the spatial dimension (we exclude the case $n=1$, for which the model does not make sense, and the case $n=2$ which has been already treated in \cite{degond2011macroscopic}). We assume that the particle has given constant speed $c_0$, i.e. its velocity $dX/dt$ satisfies $|(dX/dt)(t)|=c_0$ at all times $t$, where $|\cdot|$ denotes the Euclidean norm in ${\mathbb R}^n$. Thus, we can write $dX/dt = c_0 {\mathcal V}$ where ${\mathcal V}$ is the velocity direction (or equivalently, the unit vector tangent to the trajectory of the particle, oriented towards increasing times) satisfying
\begin{equation}
|{\mathcal V}| = 1. 
\label{eq:partelem:Vconst}
\end{equation}
In other words ${\mathcal V}(t) \in B$, $\forall t \in [0,\infty)$, where $B=:{\mathbb S}^{n-1}$ is the $n-1$-dimensional sphere. Now, we introduce the curvature of the particle trajectory $K(t)$ at time $t$. By definition, $K = d{\mathcal V}/ds$ is the derivative of the tangent vector ${\mathcal V}$ with respect to the arclength $s$ of the trajectory. Since the particle speed is $c_0$ and is constant, we have, up to a constant, $s = c_0 t$. Hence, 
\begin{equation}
K = \frac{1}{c_0} \, \frac{d{\mathcal V}}{dt}.
\label{eq:partelem:Kformula}
\end{equation} 
The vector $K(t)$ belongs to ${\mathbb R}^n$ but differentiating \eqref{eq:partelem:Vconst} with respect to $t$, we get 
\begin{equation}
{\mathcal V} \cdot K = 0, 
\label{eq:partelem:VdotK}
\end{equation} 
where $\cdot$ denotes the Euclidean inner-product in ${\mathbb R}^n$. Hence, $K(t)$ is at all times orthogonal to the velocity ${\mathcal V}(t)$. Differentiating this relation once more, we get 
\begin{equation} 
\frac{dK}{dt} \cdot {\mathcal V} = - c_0 |K|^2. 
\label{eq:partelem_dKdtdotV}
\end{equation}
This means that variations of the curvature cannot be arbitrary. Their projection along ${\mathcal V}$ must equal $- c_0 |K|^2$. In particular, only the component of $dK/dt$ normal to ${\mathcal V}$ can be controlled by external phenomena. 

Thus, let a time-dependent ``force'' $T(t)$ (where $T$ stands for ``torsion'' of the trajectory, i.e. the derivative of its curvature) be such that 
\begin{equation}
T \cdot {\mathcal V} = 0. 
\label{eq:partelem_TdotV}
\end{equation}
Then the single-particle motion is written 
\begin{empheq}[left=\empheqlbrace]{align}
&\frac{dX}{dt} = c_0 {\mathcal V}, \label{eq:partelem_dX} \\
& \displaystyle \frac{d {\mathcal V}}{dt} = c_0 K, \phantom{\frac{\frac{dX}{dt}}{\frac{dX}{dt}}} \label{eq:partelem_dV}\\ 
& \displaystyle \frac{d K}{dt} = T - c_0 |K|^2 {\mathcal V}. \label{eq:partelem_dK}
\end{empheq}
We note that the constraints \eqref{eq:partelem:Vconst} and \eqref{eq:partelem:VdotK} are satisfied at all times by the solutions of this system, as soon as they are satisfied at time $t=0$. 

\begin{figure}[htbp]
\centering

\subfloat[]{\includegraphics[trim={3.5cm 19.2cm 5.cm 4.5cm},clip,height= 3.5cm]{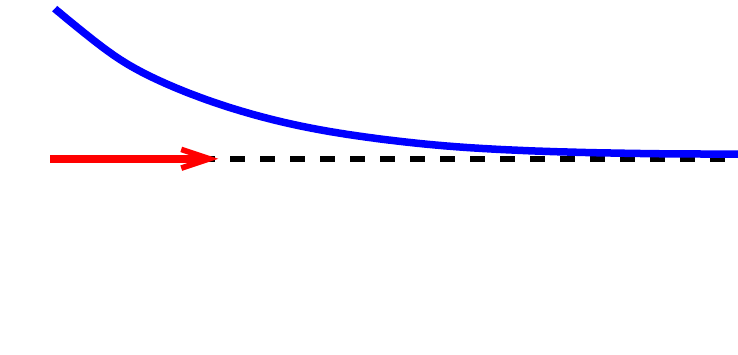}\label{subfig:Traject_strong}} \hspace{1.5cm}
\subfloat[]{\includegraphics[trim={3.5cm 18cm 5.cm 4.5cm},clip,height= 3.5cm]{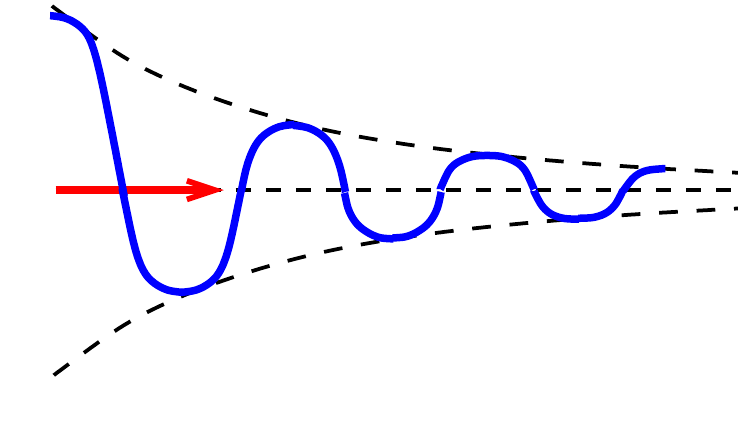}\label{subfig:Traject_weak}} 

\caption{Examples of single-particle trajectories in space in 2D (blue curve) and vector ${\mathcal J}$ (red arrow): (a) large alignment frequency $\nu > 4 c_0 \lambda |{\mathcal J}|$;  (b) small alignment frequency $\nu < 4 c_0 \lambda |{\mathcal J}|$. Note that the Vicsek case corresponds to the overdamped limit of \eqref{eq:pk_single_noise_2D} and that the Vicsek trajectories would thus resemble case (a).}
\label{fig:Trajectory}
\end{figure}

Now, we specify the choice of the force vector $T$ in the single-particle PCC model. We assume that we are given a target direction ${\mathcal J}$ (possibly depending on space and time) and that~$T$ is a relaxation force towards a fixed scalar multiple $\lambda >0$ of the projection $P_{{\mathcal V}^\bot} {\mathcal J}$ of ${\mathcal J}$ on the hyperplane orthogonal to ${\mathcal V}$, i.e. 
\begin{equation} 
T = - \nu \, \big( K - \lambda P_{{\mathcal V}^\bot} {\mathcal J} \big). 
\label{eq:partelem_Tformula}
\end{equation}
We recall that 
\begin{equation}
P_{{\mathcal V}^\bot} {\mathcal J} = {\mathcal J} - ({\mathcal J} \cdot {\mathcal V}) \, {\mathcal V}. 
\label{eq:partelem_projformula}
\end{equation}
The coefficient $\nu>0$ is the relaxation rate.  We notice that $T$ given by \eqref{eq:partelem_Tformula} satisfies \eqref{eq:partelem_TdotV} as it should. We claim that the PCC model for single-particle dynamics \eqref{eq:partelem_dX}-\eqref{eq:partelem_Tformula} models the alignment of the velocity ${\mathcal V}$ with the target vector ${\mathcal J}$. To show this, let us assume that ${\mathcal J}$ is constant and uniform (i.e. does not depend on $(x,t)$). Then, the equation for $({\mathcal V},K)$ decouples from that of $X$ and we focus on the former. We get 
$$ \frac{d{\mathcal V}}{dt} = c_0 K, \qquad \frac{dK}{dt} = - c_0 |K|^2 {\mathcal V} -\nu (K - \lambda P_{{\mathcal V}^\bot} {\mathcal J}) . $$
We further simplify the problem by considering a 2-dimensional motion (for illustration only, the core of the paper being concerned with dimensions $n \geq 3$). Hence, in an appropriate frame, we have 
$$ {\mathcal V}(t) = \left( \begin{array}{c} \cos \theta(t) \\ \sin \theta(t) \end{array} \right), \quad K(t) = \left( \begin{array}{c} - \sin \theta(t) \\ \cos \theta(t) \end{array} \right) \kappa(t), \quad {\mathcal J} = \left( \begin{array}{c} |{\mathcal J}| \\ 0 \end{array} \right), $$
thus defining $\theta(t)$ and $\kappa(t)$. Then, \eqref{eq:pk_dVK_single_noise} reduces to 
\begin{equation} 
\frac{d \theta}{dt} = c_0 \kappa, \qquad \frac{d \kappa}{dt} = - \nu (\kappa + \lambda |{\mathcal J}| \, \sin \theta), 
\label{eq:pk_single_noise_2D}
\end{equation}
or equivalently
$$ \frac{d^2 \theta}{dt^2} + \nu \frac{d \theta}{dt} + \nu c_0 \lambda |{\mathcal J}| \, \sin \theta = 0. $$
This is the equation of a damped pendulum. In particular, if $\theta$ is small (i.e. if ${\mathcal V}$ is almost aligned with ${\mathcal J}$), we can write 
$$ \frac{d^2 \theta}{dt^2} + \nu \frac{d \theta}{dt} + \nu c_0 \lambda |{\mathcal J}| \, \theta = 0, $$
whose solutions, according to the sign of $\nu - 4 c_0 \lambda |{\mathcal J}|$, are combinations of real exponentials or damped sinusoids both tending to $0$ as $t \to \infty$ (see Fig. \ref{fig:Trajectory} for illustrations). Hence, the large time behavior of the system shows that ${\mathcal V}$ aligns to ${\mathcal J}$, as claimed. It also shows that, when $\nu$ is not large, the particles are wiggling a long time around the target direction before reaching it asymptotically.

Finally, we introduce noise in the single-particle PCC model. We assume that some Brownian noise with mean-square variance $2D$, with $D$ a given constant positive diffusion coefficient, acts on the curvature $K$ only. The reason for not adding Brownian noise to the velocity ${\mathcal V}$ is that it would make ${\mathcal V}$ non time-differentiable hence making it impossible to define the curvature $K$. The precise definition of what Brownian noise in curvature only means (given that the curvature depends on ${\mathcal V}$ through the constraint \eqref{eq:partelem:VdotK}) depends on some geometric framework which will be described in forthcoming sections. At this stage, let us formally write the noisy single-particle PCC model as the following stochastic differential system: 
\begin{empheq}[left=\empheqlbrace]{align}
&dX = c_0 {\mathcal V} \, dt, \label{eq:partelem_dX_noisy} \\
& \displaystyle d {\mathcal V} = c_0 K \, dt, 
\phantom{\sqrt{2D}}
\label{eq:partelem_dV_noisy}\\ 
& \displaystyle dK  = \big( - \nu \, ( K - \lambda P_{{\mathcal V}^\bot} {\mathcal J}) - c_0 |K|^2 {\mathcal V} \big) \, dt + \sqrt{2D} \, dB_t, \label{eq:partelem_dK_noisy}
\end{empheq}
where $dB_t$ is a Brownian  noise.

\subsubsection{Many-particle PCC model}
\label{subsubsec:partelem_many_noisy}

We now introduce the many-particle PCC model. Let $N$ be the number of particles, and let the position, velocity direction and curvature of the $i$-th particle be given by the triple $(X_i(t), {\mathcal V}_i, K_i)$, where ${\mathcal V}_i$ and $K_i$ are constrained by \eqref{eq:partelem:Vconst}, \eqref{eq:partelem:VdotK}. In the many-particle PCC model (or in short, PCC model), each particle evolves independently through the single-particle PCC model introduced in the previous section, except that all the particle trajectories are coupled through the definition of the alignment vector ${\mathcal J}$. So, the trajectories of the particles are the solutions of the following system
\begin{empheq}[left=\empheqlbrace]{align}
&dX_i = c_0 {\mathcal V}_i \, dt, \label{eq:partelem_dXi} \\
&d {\mathcal V}_i = c_0 K_i \, dt, 
\phantom{\sqrt{2D}}
\label{eq:partelem_dVi} \\
&d K_i = \big( - \nu \, ( K_i - \lambda P_{{\mathcal V}_i^\bot} {\mathcal J}(X_i(t),t)) - c_0 |K_i|^2 {\mathcal V}_i \big) \, dt + \sqrt{2D} \, dB^i_t, \label{eq:partelem_dKi}
\end{empheq}
for $i = 1, \ldots, N$, where now $(dB^i_t)_{i=1, \ldots, N}$ are $N$ independent Brownian motions in $K_i$ only (the precise meaning of this will be given below). 

The vector ${\mathcal J}(x,t)$, which provides the coupling between the particles,  is defined as the mean of the velocity directions of the particles located in a neighborhood of point $x$ at time $t$. Hence, in~\eqref{eq:partelem_dKi}, ${\mathcal J}(X_i(t),t)$ will describe the mean of the velocity directions of the neighboring particles to Particle $i$. This feature turns the PCC model into an interacting particle system. There are two variants of the model, according to whether ${\mathcal J}$ is further normalized or not: 

\begin{itemize}
\item \textbf{Non-normalized PCC model.} In the non-normalized version, the vector ${\mathcal J}$ is simply the mean of the velocity directions ${\mathcal V}_j$ of the Particles $j$ located in a neighborhood of $x$ at time $t$. It is defined by 
\begin{equation}
{\mathcal J}(x,t) = \frac{1}{R^n} \sum_{\ell=1}^N {\mathcal K} \Big( \frac{|X_\ell(t) - x|}{R} \Big) \,{\mathcal V}_\ell(t), 
\label{eq:partelem_Ji_def}
\end{equation}
where ${\mathcal K}$: $[0,\infty) \to [0,\infty)$ is a sensing function satisfying 
\begin{equation} 
\int_{{\mathbb R}^n} {\mathcal K}(|x|) \, dx = 1, 
\label{eq:partelem_normalization_H}
\end{equation}
and $R>0$ is the sensing radius. This sensing function weights the importance of the particles in the neighboring of $x$. For instance, the simplest sensing function would be (up to normalization) the indicator function $\chi_{B(0,1)}$ of the ball $B(0,1)$ centered at $0$ with radius $1$: ${\mathcal K}(s) \sim \chi_{B(0,1)}$. In this case, the sensing function ${\mathcal K}$ in \eqref{eq:partelem_Ji_def} would just select those particles located in the ball centered at $X_i(t)$ with radius $R$. However, we could imagine more sophisticated sensing functions, where the importance of particle $j$ would be weighted by a function of the distance $|X_i(t) - X_j(t)|$. 

\item \textbf{Normalized PCC model.} In the normalized version, the vector ${\mathcal J}$ is the normalized mean of the velocity directions ${\mathcal V}_j$ of the Particles $j$ located in a neighborhood of $x$ at time $t$. In particular, this  assumes that this mean is nowhere zero. The vector ${\mathcal J}$ is given by: 
\begin{equation}
{\mathcal J}(x,t) = \frac{\displaystyle \sum_{\ell=1}^N {\mathcal K} \Big( \frac{|X_\ell(t) - x|}{R} \Big) \,{\mathcal V}_\ell(t)}{\displaystyle \bigg| \sum_{\ell=1}^N {\mathcal K} \Big( \frac{|X_\ell(t) - x)|}{R} \Big) \,{\mathcal V}_\ell(t) \bigg|}, 
\label{eq:partelem_Jinorm_def}
\end{equation}
with ${\mathcal K}$ satisfying the same assumptions as in the non-normalized case. Here, we must assume that the dynamics is defined only up to the first time (if it exists) where the denominator vanishes. 
\end{itemize}

The non-normalized PCC model consists of Eqs. \eqref{eq:partelem_dXi}-\eqref{eq:partelem_dKi}, \eqref{eq:partelem_Ji_def}, while the normalized PCC model, of Eqs. \eqref{eq:partelem_dXi}-\eqref{eq:partelem_dKi}, \eqref{eq:partelem_Jinorm_def}. The two systems are solved supposing that an initial condition $(X_{i0}, {\mathcal V}_{i0}, K_{i0})_{i=1, \ldots, N}$ is given. 

We now restrict to the three-dimensional case $n=3$ and establish that the non-normalized PCC model coincides with the model of \cite{cavagna2015flocking}. Introduce an arbitrary coefficient $\delta >0$ and define the new particle variables $S_i$ (for ``spin'') by 
\begin{equation} 
S_i = - \frac{1}{\delta} \, K_i \times {\mathcal V}_i, 
\label{eq:partelem_Si_def}
\end{equation}
where $\times$ denotes the vector product in three dimensions. Then, we have the relations
\begin{equation} 
K_i = \delta \, S_i \times {\mathcal V}_i, \qquad S_i \cdot {\mathcal V}_i = 0. 
\label{eq:partelem_rel_Ki_Si}
\end{equation}
From \eqref{eq:partelem_Si_def}, \eqref{eq:partelem_rel_Ki_Si}, \eqref{eq:partelem_dVi}, \eqref{eq:partelem_dKi} and \eqref{eq:partelem_Ji_def}, we deduce that the triples $(X_i, V_i, S_i)_{i=1}^N$ satisfy the system 

\begin{empheq}[left=\empheqlbrace]{align}
&\frac{dX_i}{dt} = c_0 {\mathcal V}_i, \label{eq:partelem_X_Cava} \\
&\frac{dV_i}{dt} = \frac{1}{\chi} S_i \times {\mathcal V}_i, \label{eq:partelem_V_Cava} \\
&\frac{dS_i}{dt} = {\mathcal V}_i \times \Big[ J \sum_{j=1}^N n_{ij} {\mathcal V}_j - \eta \frac{d {\mathcal V}_i}{dt} + \xi_i \Big] , \label{eq:partelem_S_Cava}
\end{empheq}
with 
\begin{equation}
\xi_i = \sqrt{(2d) \eta T} \, \frac{d B_t^i}{dt}, 
\label{eq:partelem_xi_def}
\end{equation}
and 
$$ \frac{1}{\chi} = c_0 \delta, \quad J = \frac{\nu \lambda}{\delta}, \quad n_{ij} = \frac{1}{R^3} {\mathcal K} \Big( \frac{|X_i - X_j|}{R} \Big), \quad \eta = \frac{\nu}{c_0 \delta}, \quad (2d) \eta T = \frac{2D}{\delta}. $$
Now, we can compare Equations \eqref{eq:partelem_V_Cava}-\eqref{eq:partelem_xi_def}  with Equations (1), (2) of \cite{cavagna2015flocking} and notice that they are identical, apart from the scaling factor $c_0$ in \eqref{eq:partelem_X_Cava}, which can be easily removed by introducing ${\mathcal W}_i = c_0 {\mathcal V}_i$. We conclude that the normalized PCC model is identical with the model of \cite{cavagna2015flocking} in dimension three. Hence, it can be viewed as an extension of \cite{cavagna2015flocking} to any dimensions, noting that the former was already an extension of the two-dimensional model of \cite{degond2011macroscopic}.

\subsubsection{A geometric view of the single-particle noiseless PCC model}
\label{subsubsec:partelem_single_geometric}

We return to System \eqref{eq:partelem_dX}-\eqref{eq:partelem_dK} and introduce some geometric concepts. Suppose that $t \mapsto K(t)$, $[0,\infty) \to {\mathbb R}^n$ satisfies \eqref{eq:partelem:VdotK} but is not necessarily given by \eqref{eq:partelem:Kformula}. Denote by $\Gamma = d{\mathcal V}/dt$ the particle acceleration and by ${\mathcal T} = d K/dt$ the time derivative of $K$. If we repeat the considerations of Section \ref{subsubsec:partelem_single_noiseless}, we find that 
\begin{equation} 
{\mathcal V} \cdot \Gamma = 0 \quad \textrm{ and } \quad {\mathcal T} \cdot {\mathcal V} + K \cdot \Gamma = 0. 
\label{eq:partelem:GammaTconst}
\end{equation}
Of course, if $K$ is given by \eqref{eq:partelem:Kformula}, we recover \eqref{eq:partelem_dKdtdotV}. We notice that, for a given pair $({\mathcal V},K)$ satisfying \eqref{eq:partelem:Vconst}, \eqref{eq:partelem:VdotK}, the set of pairs $(\Gamma, {\mathcal T}) \in {\mathbb R}^{2n}$ satisfying \eqref{eq:partelem:GammaTconst} is a $2(n-1)$-dimensional vector sub-space. Back to System
\eqref{eq:partelem_dX} - \eqref{eq:partelem_dK}, we can decompose
\begin{equation}
 \left( \begin{array}{c} c_0 K \\ T - c_0 |K|^2 {\mathcal V} \end{array} \right) = \left( \begin{array}{c} c_0 K \\ - c_0 |K|^2 {\mathcal V} \end{array} \right) + \left( \begin{array}{c} 0 \\ T \end{array} \right). 
\label{eq:partelem_decomp}
\end{equation}
We notice that both vectors at the right-hand side satisfy \eqref{eq:partelem:GammaTconst} (it is obvious for the former, while the latter requires \eqref{eq:partelem_TdotV}). 

We now introduce some elements of geometry. The $(n-1)$-dimensional sphere to which~${\mathcal V}$ belongs is a differentiable manifold and for each ${\mathcal V}$, a vector $K$ satisfying \eqref{eq:partelem:VdotK} is an element of the tangent space $T_{\mathcal V} B$ to $B$ at ${\mathcal V}$. The set of pairs $({\mathcal V}, K)$ where ${\mathcal V} \in B$ and $K \in T_{\mathcal V} B$ is called the tangent bundle of $B$ and will be denoted by $M$. The tangent space of $M$ at $({\mathcal V}, K)$, denoted by $T_{({\mathcal V}, K)} M$ consists of those pairs $(\Gamma, {\mathcal T})$ satisfying \eqref{eq:partelem:GammaTconst}. A vector field on $B$ is a map $k$: $B \to M$ such that, for any ${\mathcal V} \in B$, $k({\mathcal V}) \in T_{\mathcal V} B$. Like in classical analysis, vector fields can be interpreted as particle velocities which, after integration, give rise to trajectories on $B$. Likewise, a vector field on $M$ is a map $\xi$: $M \to N$ (where $N$ is the tangent bundle to $M$) such that for all $({\mathcal V}, K) \in M$, $\xi ({\mathcal V}, K) \in T_{({\mathcal V}, K)} M$. Hence, we see that the left-hand side of \eqref{eq:partelem_decomp} is a vector field on $M$ and that System \eqref{eq:partelem_dX}-\eqref{eq:partelem_dK} can be seen as a differential equation on ${\mathbb R}^n \times M$. 

In the present work, it will be important to distinguish motions on $M$ (i.e. vector fields on $M$) which are ``parallel to ${\mathcal V}$'' or ``parallel to $K$''. For this purpose, we notice that both vectors at the right-hand side of \eqref{eq:partelem_decomp} are vector fields on $M$. Furthermore, each of them depends on a vector field on $B$ (namely $K$ for the first one and $T$ for the second one). More generally, for any vector field $\Gamma$ on $B$, we can construct the following vector field on $M$: 
\begin{equation}
{\mathcal L}^H_{({\mathcal V}, K)} \Gamma = \left( \begin{array}{c} \Gamma({\mathcal V}) \\ - (\Gamma({\mathcal V}) \cdot K) {\mathcal V} \end{array} \right). 
\label{eq:partelem_LHdef}
\end{equation}
We note that the right-hand side of \eqref{eq:partelem_LHdef} satisfies \eqref{eq:partelem:GammaTconst} and defines a valid vector field on $M$. This vector field will be called the horizontal lift of $\Gamma$ at $({\mathcal V}, K)$. Loosely speaking, horizontal lifts define motions on $M$ in a direction ``parallel to ${\mathcal V}$''. We see that, moving parallel to ${\mathcal V}$ imposes $K$ to vary, as the second component of the vector \eqref{eq:partelem_LHdef} shows. This is because, as ${\mathcal V}$ varies, $K$ must vary too because of the constraint \eqref{eq:partelem:VdotK}. 

Likewise, for any vector field $T$ on $B$, we can construct the following vector field on $M$: 
\begin{equation}
{\mathcal L}^V_{({\mathcal V}, K)} T = \left( \begin{array}{c} 0 \\ T({\mathcal V}) \end{array} \right). 
\label{eq:partelem_LVdef}
\end{equation}
Again, we remark that the right-hand side of \eqref{eq:partelem_LVdef} satisfies \eqref{eq:partelem:GammaTconst} and defines a valid vector field on $M$. This vector field will be referred to as the vertical lift of $T$. Again, loosely speaking, vertical lifts define motions on $M$ in a direction ``parallel to $K$''. By contrast to motions in the ${\mathcal V}$-direction, motions in the $K$-direction do not involve any variation in ${\mathcal V}$, as the vanishing first component of the vector \eqref{eq:partelem_LVdef} shows. This is because the tangent space $T_{\mathcal V} B$ along which~$K$ varies is a vector space, and moving along the tangent space $T_{\mathcal V} B$ does not change its base-point ${\mathcal V}$. Hence, we conclude that motions along ${\mathcal V}$ and motions along $K$ do not play a symmetric role, and this will be one of the difficulties which requires the introduction of a rigorous geometric framework. 

Finally, with these definitions, \eqref{eq:partelem_dV}-\eqref{eq:partelem_dK} can be re-written
\begin{equation}
\frac{d}{dt} \left( \begin{array}{c} {\mathcal V} \\ K \end{array} \right) = c_0 {\mathcal L}^H_{({\mathcal V}, K)} K +  {\mathcal L}^V_{({\mathcal V}, K)} T. 
\label{eq:partelem_ddtVK}
\end{equation}

In the next subsection, we introduce this geometric framework more rigorously.

\subsection{Geometric preliminaries}
\label{subsec:particle_geom}

This subsection recalls the geometric concepts that were mentioned in the previous section. We refer to \cite{warner1983foundations} for an introduction to differential geometry, to \cite{docarmo1992riemannian, gallot1990riemannian} for Riemannian geometry and to \cite{greub1972connections} for the geometry of bundles. 
We refer to Figure \ref{fig:diagram} for a depiction of the various sets and maps involved in the forthcoming description.

\subsubsection{The tangent bundle to the sphere}
\label{subsubsec:particle_tangent}

We recall that we assume $n \geq 3$. We denote by $B = {\mathbb S}^{n-1}$ the $(n-1)$-dimensional sphere and by $F = {\mathbb R}^{n-1}$. We let $TB = (M, \pi, B, F)$ be the {\em tangent bundle} to $B$. In this quadruple, $M = \sqcup_{v \in B} T_v B$ is the support manifold to the vector bundle and consists of the disjoint union  of the tangent spaces $T_vB$ to $B$ at $v$, $B$ is the base, $F$ is the typical fiber and $\pi$: $M \to B$ is the base map. Since $B$ is embedded in ${\mathbb R}^n$, $M$ is embedded in ${\mathbb R}^{2n}$ and we can write: 
\begin{equation} 
\alpha \in M \quad \Longleftrightarrow \quad \big( \alpha=(v,\kappa) \in {\mathbb R}^{2n} \, \, \textrm{such that} \, \, |v|=1\, \, \mathrm{and} \, \, \kappa\cdot v = 0 \big), 
\label{eq:geom_M}
\end{equation}
where the symbols `` $\cdot$ '' and `` $| \, \, |$ '' are the Euclidean inner product and norm on ${\mathbb R}^n$. The notations~$v$ for elements of $B$ and $\kappa$ for elements of $T_v B$ stand for ``velocity'' and ``curvature'' respectively.  

Now, we introduce the tangent bundle $TM = (N, \bar{\pi},M,G)$ of $M$, where $N$ is the support manifold of $TM$, $G = F^2 = {\mathbb R}^{2n-2}$ is the typical fiber and $\bar{\pi}$: $N \to M$ is the base map. Now, $N$ is embedded in ${\mathbb R}^{4n}$. Let $\alpha = (v,\kappa) \in M$ (with $(v,\kappa)$ satisfying \eqref{eq:geom_M}). Then, we have 
\begin{eqnarray} 
\xi \in T_\alpha M & \Longleftrightarrow & \big( \xi = (a,\tau) \in {\mathbb R}^{2n} \, \, \textrm{such that} \, \, a \cdot v = 0 \, \, \mathrm{and} \, \, \tau \cdot v + a \cdot \kappa = 0 \big) \label{eq:geom_TM}\\
& \Longleftrightarrow &  \big( \xi =  (a,\bar \tau - (a \cdot \kappa) v) \, \, \textrm{such that} \, \, a, \, \bar \tau \in \{v\}^\bot \big). \nonumber
\end{eqnarray}
The notations $(a,\tau)$ for the elements of $T_\alpha M$ respectively stand for ``acceleration'' and ``torsion''.

The Euclidean structure of ${\mathbb R}^n$ induces a {\em Riemannian structure} on $B$: the metric on $T_v B$ is defined by 
\begin{equation} 
\big\langle \kappa_1 , \kappa_2 \big\rangle_v = \kappa_1 \cdot \kappa_2,
\label{eq:geom_innerprod_al}
\end{equation}
for all $\kappa_i \in T_v B = \{v\}^\bot$, $i=1, \, 2$. We will drop the subscript $v$ when the context is clear. The Riemannian structure on $B$ determines the {\em Levi-Civita connection}, i.e. the unique torsion-free linear connection on $B$ compatible with the Riemannian structure (here, torsion has its differential geometric meaning, and should not be confused with the torsion of a curve which has been used so far). Denote by ${\mathcal X}(B)$ the set of smooth vector fields on $B$ and let $k$, $h \in {\mathcal X}(B)$. The Levi-Civita connection generates a new vector field $\nabla^\mathrm{B}_k h$ (the gradient of $h$ along the field $k$). In general, we will assume that $k \in {\mathcal X}(B)$ can be extended into a smooth function $\bar k$: $U_B \to {\mathbb R}^n$, where $U_B$ is an open set of ${\mathbb R}^n$ containing~$B$. Then, thanks to \cite[Chapt. 2, Exercise 3]{docarmo1992riemannian}, the Levi-Civita connection is given by 
\begin{equation}
\nabla^\mathrm{B}_k h(v) = P_{v^\bot} \big((\bar k \cdot \boldsymbol{\nabla}_v) \bar h \big)(v),  
\label{eq:geom_LCconnect}
\end{equation}
We have denoted by $\boldsymbol{\nabla}_v$ the usual nabla operator operating on functions of $v \in {\mathbb R}^n$ and by $P_{v^\bot}$ the orthogonal projection ${\mathbb R}^n \to \{v\}^\bot$. Matrix-wise, $P_{v^\bot} = {\mathrm I}_n - v \otimes v$ (where ${\mathrm I}_n$ is the identity matrix of ${\mathbb R}^n$ and $\otimes$ stands for the tensor product).

\subsubsection{Vertical and horizontal bundle}
\label{subsubsec:particle_vert_hor}

Let $\alpha = (v, \kappa) \in M$. The space $T_\alpha M$ is decomposed into two components, which respectively represent motions ``in the direction of $v$'' and ``in the direction of $\kappa$''. We first start with motions in the direction of $\kappa$ and define the {\em vertical subspace} $V_\alpha = \mathrm{ker} \, (d \pi)_\alpha$, where $(d \pi)_\alpha$: $T_\alpha M \to T_v B$ is the differential of the base map $\pi$. For any $\xi = (a, \tau) \in T_\alpha M$ (with $(a,\tau)$ satisfying \eqref{eq:geom_TM}), we have 
\begin{equation}
(d \pi)_\alpha(\xi) = a, 
\label{eq:geom_Vertical_prf1}
\end{equation}
and thus 
\begin{equation} 
\xi \in V_\alpha \, \,  \Longleftrightarrow  \, \, \big( \xi = (0,\tau) \in {\mathbb R}^{2n} \, \, \textrm{such that} \, \, \tau \cdot v = 0 \big). 
\label{eq:geom_Vertical}
\end{equation} 
In particular, $\mathrm{dim} \, V_\alpha = n-1$. For any $\alpha \in M$, the map ${\mathcal L}^V_\alpha$: $T_v B \to V_\alpha$ defined by 
\begin{equation} 
{\mathcal L}^V_\alpha (\tau) = (0,\tau), \quad \forall \tau \in T_v B, 
\label{eq:geom_omega}
\end{equation}
is called the {\em vertical lift}. We note that 
\begin{equation}
(d \pi)_\alpha \circ {\mathcal L}_\alpha^V = 0.
\label{eq:geom_dpioL}
\end{equation}

Now, we define directions of motion ``in the $v$ direction''. This requires the  connection. We first define the horizontal lift. In general, a {\em horizontal lift} is a linear map ${\mathcal L}^H_\alpha$: $T_v B \to T_\alpha M$ such that 
\begin{equation} 
(d \pi)_\alpha \circ {\mathcal L}^H_\alpha = {\mathrm I}_{T_v B}, \quad \forall \alpha = (v,\kappa) \in M, 
\label{eq:geom_gam}
\end{equation}
where ${\mathrm I}_{T_v B}$ is the identity of $T_v B$. Then, the {\em horizontal subspace $H_\alpha$} associated with the horizontal lift ${\mathcal L}^H_\alpha$ is $H_\alpha = \textrm{Im} \, {\mathcal L}^H_\alpha$. The horizontal and vertical spaces are supplementary subspaces of $T_\alpha M$:  
\begin{equation} 
V_\alpha \oplus H_\alpha = T_\alpha M. 
\label{eq:geom_ValopHal}
\end{equation}
When $B$ is endowed with a linear connection $\nabla^{\mathrm{B}}$, there exists a {\em unique horizontal map} ${\mathcal L}^H$ associated with this connection \cite[Vol. 2, Chapt. VII, Sect. 6, Prop. IX]{greub1972connections}. It can be shown that, for the sphere endowed with the Levi-Civita connection \ref{eq:geom_LCconnect}, we have
\begin{equation} 
{\mathcal L}^H_\alpha(a) = \big(a,-(\kappa \cdot a) v \big), \quad \forall a \in T_v B. 
\label{eq:geom_horizmap}
\end{equation}
and
\begin{equation} 
\xi \in H_\alpha \, \,  \Longleftrightarrow  \, \, \big( \xi = (a,\tau) \in {\mathbb R}^{2n} \, \, \textrm{with} \, \, \tau = - (\kappa \cdot a) \, v \big) 
\label{eq:geom_Horizontal}
\end{equation}

\bigskip
\begin{figure}[htbp]
\centering
\includegraphics[clip,height= 7cm]{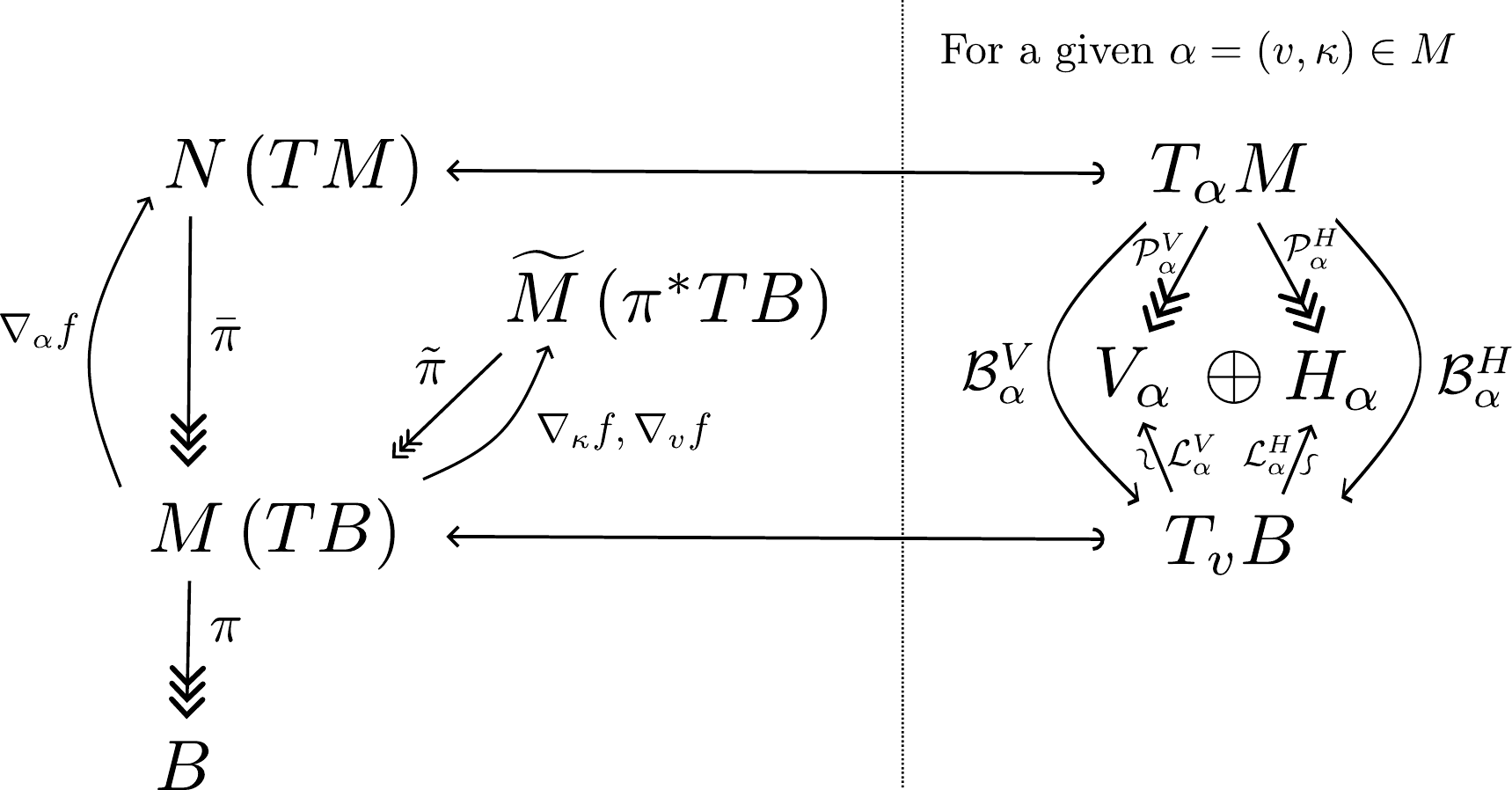}
\caption{Diagram of the sets and maps involved in the geometric framework. Arrows with a triple head mean surjections and arrows with an anchor are for injections.}
\label{fig:diagram}
\end{figure}

\subsection{Geometric presentation of the PCC system}
\label{subsec:particle_dynamics}

In this subsection, we revisit the derivation of the PCC model introduced in Subsection \ref{subsec:partelem} in light of the geometric concepts recalled in Subsection \ref{subsec:particle_geom}. Many different constructions of stochastic differential equations on manifolds and bundles have been proposed in the literature (see \cite{bucataru2011geometric} for the deterministic case and \cite{angst2015kinetic, baudoin2016hypocoercive, girolami2011riemann, grothaus2022hypocoercivity, joergensen1978construction, lelievre2012langevin, soloveitchik1995focker} for the stochastic case). Our construction relies on \cite{ikeda2014stochastic} (see also \cite{hsu2002stochastic}). 

We first consider the noiseless single-particle PCC model. Let a single particle be  characterized by its position $X \in {\mathbb R}^n$, and the pair ${\mathcal A} =:  ({\mathcal V}, K) \in M$ consisting of the direction of its velocity ${\mathcal V} \in {\mathbb S}^{n-1}$ and the curvature vector of its trajectory $K \in {\mathbb R}^n$ such that \eqref{eq:partelem:VdotK} holds. The assumptions of the model have been detailed in Section \ref{subsec:partelem}. The single-particle noiseless PCC model is written according to \eqref{eq:partelem_ddtVK}, which can be re-written in terms of ${\mathcal A} \in M$ according to
\begin{empheq}[left=\empheqlbrace]{align}
&\frac{dX}{dt} = c_0 {\mathcal V}, \label{eq:pk_dX_single_nonoise} \\
&\frac{d {\mathcal A}}{dt} = c_0 {\mathcal L}^H_{{\mathcal A}} K  - \nu {\mathcal L}^V_{{\mathcal A}} \big( K - \lambda P_{{\mathcal V}^\bot} ({\mathcal J}(X,t)) \big) 
\label{eq:pk_dVK_single_nonoise}
\end{empheq}

To pass to the noisy single-particle PCC model, we need to generate Brownian noise on~$M$ which acts on the curvature only. For this purpose, we generate Brownian noise on $B$ and then, lift it with the vertical lift to generate Brownian noise on the vertical subspace $V_{{\mathcal A}}$ of $T_{\mathcal A} M$. We know that Brownian noise on $B$ at ${\mathcal V}$ can be obtained by projecting Brownian noise on ${\mathbb R}^n$ onto $T_{\mathcal V} B$ using the Stratonovich construction \cite{hsu2002stochastic, ikeda2014stochastic}. Let $(\mathbf{e}_1, \ldots, \mathbf{e}_n)$ be the canonical basis of ${\mathbb R}^n$. Then, Brownian noise in the $K$ variable only with variance $2D$ is written $(2D)^{1/2} \sum_{j=1}^n {\mathcal L}^V_{{\mathcal A}} ( P_{{\mathcal V}^\bot} \mathbf{e}_j ) \circ dB_{j,t}$ where $(dB_{j,t})_{j=1, \ldots, n}$ are $n$ independent Brownian motions on ${\mathbb R}$ and the symbol $\circ$ indicates that the equation is meant in the Stratonovich sense. Hence, the noisy single-particle PCC model reads: 
\begin{empheq}[left=\empheqlbrace]{align}
&dX = c_0 {\mathcal V} \, dt, \label{eq:pk_dX_single_noise} \\
&d {\mathcal A} = c_0 {\mathcal L}^H_{{\mathcal A}} K \, dt  - \nu {\mathcal L}^V_{{\mathcal A}} \big( K - \lambda P_{{\mathcal V}^\bot} ({\mathcal J}(X,t)) \big) \, dt  + \sqrt{2D} \sum_{j=1}^n {\mathcal L}^V_{{\mathcal A}} \big( P_{{\mathcal V}^\bot} \mathbf{e}_j \big) \circ dB_{j,t}. 
\label{eq:pk_dVK_single_noise}
\end{empheq}

Finally, the noisy many-particle PCC model (or PCC model in short) is a direct generalization of the previous model. Using the notations of Section \ref{subsubsec:partelem_many_noisy}, this system reads
\begin{empheq}[left=\empheqlbrace]{align}
&dX_i = c_0 {\mathcal V}_i \, dt, \label{eq:pk_dXi} \\
&d {\mathcal A}_i = c_0 {\mathcal L}^H_{{\mathcal A}_i} K_i \, dt  - \nu {\mathcal L}^V_{{\mathcal A}_i} \big( K_i - \lambda P_{{\mathcal V_i}^\bot} ({\mathcal J}(X_i,t)) \big) \, dt + \sqrt{2D} \sum_{j=1}^n {\mathcal L}^V_{{\mathcal A}_i} \big( P_{{\mathcal V}_i^\bot} \mathbf{e}_j \big) \circ dB^i_{j,t}, \label{eq:pk_dVKi}
\end{empheq}
where now $(dB^i_{j,t})^{i=1, \ldots, N}_{j=1, \ldots, n}$ are $nN$ independent Brownian motions in ${\mathbb R}$. The vector ${\mathcal J}$ is defined by 
\eqref{eq:partelem_Ji_def} in the non-normalized case and by \eqref{eq:partelem_Jinorm_def} in the normalized case.

\subsection{Scaling of the PCC model}
\label{subsec:pk_scaling_particles}

\subsubsection{Scaling in the normalized case}
\label{subsubsec:pk_scaling_normalized}

In this case, ${\mathcal J}$ is dimensionless, so, the parameter $\lambda$ has the physical dimension of the inverse of a length. We will use $\lambda^{-1}$ as a spatial scale. Since $c_0$ is a speed, $(c_0 \lambda)^{-1}$ will be used as a time scale. We introduce the following changes of variables and unknowns: 
\begin{equation}
X_i = \lambda^{-1} \tilde X_i, \quad t = (c_0 \lambda)^{-1} \tilde t, \quad K_i = \lambda \tilde K_i,
\label{eq:pk_scaling_norm}
\end{equation} 
which defines dimensionless variables and unknowns $\tilde X_i$, $\tilde t$, $\tilde K_i$. We note that the unknown ${\mathcal V}_i$ is already dimensionless and does not require any scaling. The scaling on $K_i$ requires that we introduce a scaled variable $\tilde {\mathcal A}_i = (\tilde {\mathcal V}_i, \tilde K_i)$ which is obtained from the unscaled variable ${\mathcal A}_i = ({\mathcal V}_i, K_i)$ by setting 
\begin{equation}
{\mathcal A}_i = \Sigma_\lambda(\tilde {\mathcal A}_i) \quad \textrm{ or equivalently } \quad \tilde {\mathcal A}_i = \Sigma_{\lambda^{-1}} ({\mathcal A}_i). 
\label{eq:pk_scaling_Ai_norm}
\end{equation}
where for $\lambda >0$, the scaling map $\Sigma_\lambda$: $M \to M$, is defined by $\Sigma_\lambda(v,\kappa) = (v,\lambda \kappa)$, $\forall \alpha = (v,\kappa) \in M$. The scaling map satisfies 
\begin{equation}
(d \Sigma_\lambda)_\alpha \circ {\mathcal L}_\alpha^H = {\mathcal L}_{\Sigma_\lambda \alpha}^H, \qquad   
(d \Sigma_\lambda)_\alpha \circ {\mathcal L}_\alpha^V = \lambda \,  {\mathcal L}_{\Sigma_\lambda \alpha}^V. 
\label{eq:pk_scaling_LVtau}
\end{equation}
We now let $(\tilde X_i, \tilde {\mathcal A}_i)$ be functions of the rescaled time $\tilde t$ and define 
$$ \tilde {\mathcal J}( \tilde x,t) = {\mathcal J}(\lambda^{-1} \tilde x, (c_0 \lambda)^{-1} \tilde t). $$
We introduce three scaled constants 
\begin{equation}
\bar \nu = \frac{\nu}{c_0 \lambda} , \qquad  \bar D = \frac{D}{c_0 \lambda^3}, \qquad \bar R = \lambda R. 
\label{eq:pk_const_scaled}
\end{equation}
Then, using \eqref{eq:pk_scaling_LVtau}, we find that the normalized PCC model in the scaled variable $(\tilde X_i, \tilde {\mathcal A}_i)$ and the scaled time $\tilde t$ is given by (dropping the tildes and the bars for clarity) : 
\begin{empheq}[left=\empheqlbrace]{align}
&d X_i = {\mathcal V}_i \, d t, \label{eq:pk_dXi_scaled} \\
&d {\mathcal A}_i =  {\mathcal L}^H_{{\mathcal A}_i} K_i \, d t  - \nu {\mathcal L}^V_{{\mathcal A}_i} \big( K_i -  P_{{\mathcal V_i}^\bot} ({\mathcal J}( X_i, t)) \big) \, d t + \sqrt{2 D} \sum_{j=1}^n {\mathcal L}^V_{{\mathcal A}_i} \big( P_{{\mathcal V}_i^\bot} e_j \big) \circ dB^i_{j, t}, \label{eq:pk_dVKi_scaled}
\end{empheq}
with ${\mathcal J}$ given by 
\begin{equation}
{\mathcal J}( x, t) = \frac{\displaystyle \sum_{\ell=1}^N {\mathcal K} \Big( \frac{|X_\ell( t) - x|}{R} \Big) \,{\mathcal V}_\ell(t)}{\displaystyle \bigg| \sum_{\ell=1}^N {\mathcal K} \Big( \frac{|X_\ell( t) - x|}{R} \Big) \,{\mathcal V}_\ell(t) \bigg|}. 
\label{eq:pk_Jinorm_scaled}
\end{equation}

\subsubsection{Scaling in the non-normalized case}
\label{subsubsec:pk_scaling_nonnormalized}

We assume that ${\mathcal K}$ is dimensionless. Let $x_0$ be a spatial scale. ${\mathcal J}$ has the dimension of the inverse of a volume, i.e., $x_0^{-n}$ while $K_i$ has dimension $x_0^{-1}$. From the second term at the right-hand side of \eqref{eq:pk_dVKi}, $\lambda$ has the dimension of $x_0^{n-1}$. Thus, we can set the length scale $x_0 = \lambda^{1/(n-1)}$. So, we introduce the scaling 
\begin{equation}
X_i = \lambda^{\frac{1}{n-1}} \tilde X_i, \quad t = \lambda^{\frac{1}{n-1}}c_0^{-1} \tilde t, \quad K_i = \lambda^{-\frac{1}{n-1}} \tilde K_i, \quad {\mathcal J}=\lambda^{-\frac{n}{n-1}} \tilde {\mathcal J}, 
\label{eq:pk_scaling_nonnorm}
\end{equation} 
so defining the dimensionless variables and unknowns $\tilde X_i$, $\tilde t$, $\tilde K_i$, $\tilde {\mathcal J}$. We also define the following dimensionless constants: 
\begin{equation}
\bar \nu = \nu \lambda^{\frac{1}{n-1}} c_0^{-1} , \qquad  \bar D = D \lambda^{\frac{3}{n-1}} c_0^{-1}, \qquad \bar R = R \lambda^{-\frac{1}{n-1}}. 
\label{eq:pk_const_scaled_nonnorm}
\end{equation}
Then, the non-normalized PCC model in the scaled variables $(\tilde X_i, \tilde {\mathcal A}_i)$ and the scaled time $\tilde t$ is given by  \eqref{eq:pk_dXi_scaled}, \eqref{eq:pk_dVKi_scaled} (again dropping the tildes and the bars for clarity) with ${\mathcal J}$ given by 
\begin{equation}
{\mathcal J}( x, t) = \frac{1}{R^n} \sum_{\ell=1}^N {\mathcal K} \Big( \frac{|X_\ell( t) - x|}{R} \Big) \,{\mathcal V}_\ell(t). 
\label{eq:pk_Jinonnorm_scaled}
\end{equation}

\setcounter{equation}{0}
\section{The kinetic curvature-control (KCC) model}
\label{sec:kinetic}

In this section, we formally establish the kinetic mean-field model, below referred to as the kinetic curvature-control model (KCC model), associated with the PCC model. We first develop additional geometric concepts that are needed to achieve this goal. We then proceed to the derivation of the KCC model and conclude this section by presenting its hydrodynamic rescaling.

\subsection{Geometric preliminaries}
\label{subsec:kinetic_geom}

Again, we refer to Figure \ref{fig:diagram} for an illustration of the concepts described in this section and to \cite{docarmo1992riemannian, gallot1990riemannian, greub1972connections, warner1983foundations} for geometric background.

\subsubsection{Connection map}
\label{subsubsec:geom_horvertop}

In this section, we define additional operators linked to the horizontal and vertical lifts. The vertical (resp. horizontal) projection is the projection ${\mathcal P}^V_{\alpha}$ (resp. ${\mathcal P}^H_\alpha$) of $T_\alpha M$ on $V_\alpha$ parallel to $H_\alpha$ (resp. on $H_\alpha$ parallel to $V_\alpha$). They are given by 
\begin{eqnarray}
{\mathcal P}^V_\alpha \xi &=& \big( 0, \tau + (\kappa \cdot a) v \big) = \big( 0, \tau - (v \cdot \tau) v \big) = \big( 0, P_{v^\bot} \tau \big), \label{eq:geom_proj_Val} \\
{\mathcal P}^H_\alpha \xi &=& \big( a, - (\kappa \cdot a) v \big) , \label{eq:geom_proj_Hal}
\end{eqnarray}
for all $\xi = (a,\tau)$ in $T_\alpha M$ (i.e. satisfying \eqref{eq:geom_TM}). Next, we note that the vertical and horizontal lifts are isomorphism: ${\mathcal L}^V_\alpha$: $T_{\pi \alpha} B \to V_\alpha$ and ${\mathcal L}^H_\alpha$: $T_{\pi \alpha} B \to H_\alpha$. 
Their inverses are denoted by $({\mathcal L}^V_\alpha)^{-1}$ and $({\mathcal L}^H_\alpha)^{-1}$ and given by 
\begin{eqnarray} 
({\mathcal L}^V_\alpha)^{-1} \xi &=& \tau, \quad \forall \xi \in V_\alpha \quad \textrm{i.e.} \quad \forall \xi = (0,\tau) \quad \textrm{with} \quad \tau \cdot v = 0.
\label{eq:geom_inverse_LV} \\
({\mathcal L}^H_\alpha)^{-1} \xi &=& a, \quad \forall \xi \in H_\alpha \quad \textrm{ i.e. } \quad \forall \xi = (a, - (a \cdot \kappa) v) \quad \textrm{with} \quad a \cdot v = 0.
\label{eq:geom_inverse_LH}
\end{eqnarray}
Then, we define ${\mathcal B}^H_\alpha$ and ${\mathcal B}^V_\alpha$: $T_\alpha M \to T_v B$ by 
\begin{equation} 
{\mathcal B}^V_\alpha = ({\mathcal L}^V_\alpha)^{-1} \circ {\mathcal P}^V_\alpha, \qquad 
{\mathcal B}^H_\alpha = ({\mathcal L}^H_\alpha)^{-1} \circ {\mathcal P}^H_\alpha. 
\label{eq:geom_baseop_express}
\end{equation}
Thanks to \eqref{eq:geom_gam}, we have 
\begin{equation}
{\mathcal B}_\alpha^H = (d \pi)_\alpha, 
\label{eq:geom_baseH=dpi}
\end{equation}
and thus, ${\mathcal B}_\alpha^H$ is given by \eqref{eq:geom_Vertical_prf1}. We note that ${\mathcal B}_\alpha^H$ does not depend on the connection. On the other hand, the map ${\mathcal B}_\alpha^V$ does, and is called the connection map. We have 
\begin{equation} 
{\mathcal B}^V_\alpha \xi = P_{v^\bot} \tau, \qquad \forall \xi = (a,\tau) \quad \textrm{satisfying \eqref{eq:geom_TM}}. 
\label{eq:geom_BVBH}
\end{equation}

\subsubsection{Riemannian structure on $M$}
\label{subsubsec:geom_integ_Riemann_M}

We endow $M$ with the Riemannian structure induced by the {\em Sasaki metric} \cite{gudmundsson2002geometry, sasaki1958differential} and \cite[Chap. 3, Exercise 2]{docarmo1992riemannian}. For $\alpha \in M$ and $\xi_i \in T_\alpha M$, $i=1, \, 2$, the Sasaki metric on $T_\alpha M$ is defined by  
\begin{equation} 
\big\langle \hspace{-1.8mm} \big\langle \hspace{-0.1mm} \xi_1 , \xi_2 \big\rangle \hspace{-1.8mm} \big\rangle_{\alpha} = \big\langle {\mathcal B}^H_\alpha \xi_1, {\mathcal B}^H_\alpha \xi_2 \big \rangle_{\pi \alpha} + \big \langle {\mathcal B}^V_\alpha \xi_1, {\mathcal B}^V_\alpha \xi_2 \big \rangle_{\pi \alpha}. 
\label{eq:geom_Sasaki}
\end{equation}
Writing $\alpha = (v,\kappa) \in {\mathbb R}^{2n}$ and $\xi_i = (a_i, \tau_i) \in {\mathbb R}^{2n}$, with $(v,\kappa)$ satisfying \eqref{eq:geom_M} and $(a_i,\tau_i)$ satisfying \eqref{eq:geom_TM}, the Sasaki metric is written: 
\begin{equation} 
\big\langle \hspace{-1.8mm} \big\langle \hspace{-0.1mm} \xi_1 , \xi_2 \big\rangle \hspace{-1.8mm} \big\rangle_{\alpha} = a_1 \cdot a_2 + P_{v^\bot} \tau_1 \cdot P_{v^\bot} \tau_2.
\label{eq:geom_innerprod_M}
\end{equation}
We note that this metric is \textbf{not} the metric induced by the embedding $M \to {\mathbb R}^{2n}$. Indeed, the latter would add a term $(a_1 \cdot \kappa)(a_2 \cdot \kappa)$ to the right-hand side of \eqref{eq:geom_innerprod_M}. The vertical and horizontal subspaces are orthogonal for the Sasaki metric: 
$$ V_\alpha \stackrel{\bot}{\oplus} H_\alpha = T_\alpha M, \quad \forall \alpha \in M. $$
and thus, the projection operators ${\mathcal P}^V_\alpha$ and ${\mathcal P}^V_\alpha$ are self-adjoint with respect to the Sasaki metric. The Sasaki metric on $M$ is linked with the metric on $B$ by 
\begin{equation}
\big\langle \hspace{-1.8mm} \big\langle {\mathcal L}^H_\alpha \kappa, \xi \big\rangle \hspace{-1.8mm} \big\rangle_{\alpha}= \big \langle \kappa, {\mathcal B}^H_\alpha \xi \big \rangle_{\pi \alpha}, \qquad \big\langle \hspace{-1.8mm} \big\langle {\mathcal L}^V_\alpha \kappa, \xi \big\rangle \hspace{-1.8mm} \big\rangle_{\alpha} = \big \langle \kappa, {\mathcal B}^V_\alpha \xi \big \rangle_{\pi \alpha},
\label{eq:geom_SasvsEuc}
\end{equation}
for all $\alpha =(v,\kappa) \in M$, $\xi \in T_\alpha M$ and $\kappa \in T_v B$. Consequently, ${\mathcal B}^H_\alpha$ and ${\mathcal L}^H_\alpha$ are mutually adjoint and so are ${\mathcal B}^V_\alpha$ and ${\mathcal L}^V_\alpha$

\subsubsection{Integration on $M$}
\label{subsubsec:geom_orient_Sasaki}

We define the positive orientation of ${\mathbb R}^n$ as that given by the canonical basis $(\mathbf{e}_1, \ldots, \mathbf{e}_n)$, with $\mathbf{e}_i = (0, \ldots, 0,1,0, \ldots 0)$, where the $1$ is the $i$-th position. We choose the positive orientation of~$B$ such that an orthogonal basis $(\varepsilon_1, \ldots, \varepsilon_{n-1})$ of $T_v B$ is direct if and only if $(v,\varepsilon_1, \ldots, \varepsilon_{n-1})$ is a direct basis of ${\mathbb R}^n$. The volume form $\mathrm{Vol}_B$ associated with the Riemannian structure induced by the embedding of $B$ into ${\mathbb R}^n$ is given as follows: Let $(\varepsilon_1, \ldots, \varepsilon_{n-1})$ be a local, positively oriented orthonormal frame. Then, the volume form  $\mathrm{Vol}_B$ on $B$ has expression 
\begin{equation}
\mathrm{Vol}_B(v)(\kappa_1, \ldots, \kappa_{n-1}) = \mathrm{det} \big( \big \langle \varepsilon_i, \kappa_j \big \rangle_v \big), \quad \forall (\kappa_1, \ldots, \kappa_{n-1}) \in (T_v B)^{n-1}. 
\label{eq:geom_volformB}
\end{equation}
We define the following $n-1$-form on $M$ for all $(\xi_1, \ldots, \xi_{n-1}) \in (T_\alpha M)^{n-1}$ by
\begin{eqnarray}
\Omega(\alpha)(\xi_1, \ldots, \xi_{n-1}) &=& \mathrm{Vol}_B (\pi \alpha) \big( {\mathcal B}_\alpha^V(\xi_1), \ldots, {\mathcal B}_\alpha^V(\xi_{n-1}) \big).
\label{eq:geom_orient_form_bundle}
\end{eqnarray}
Then, 
\begin{equation}
\mathrm{Vol}_M = \pi^* \mathrm{Vol}_B \wedge \Omega. 
\label{eq:geom_sasvolform_compact}
\end{equation}
(where $\pi^* \mathrm{Vol}_B$ denotes the pullback of the form $\mathrm{Vol}_B$ of $B$ to a form on $M$ by the bundle base map $\pi$) is a volume form on $M$ associated with the Sasaki metric and we choose the orientation of $M$ provided by this form. We will denote by $|\mathrm{Vol}_B|$ and $|\mathrm{Vol}_M|$ the densities associated with the volume forms on $B$ and $M$ respectively. 

Following \cite[Vol. 1, Ch. 7]{greub1972connections}, for an integrable function $f$: $M \to {\mathbb R}$, we can integrate $f \Omega$ along the fibers ${\mathbb R}^{n-1}$ as follows 
\begin{equation} 
\Big( \int_{{\mathbb R}^{n-1}} f \, \Omega \Big)(v)  =: \int_{T_v B} f \, \Omega_v,   \label{eq:Fub_sigma(x)_def}
\end{equation}
where $\Omega_v$ is a $n-1$-form on $T_v B$ called the retrenchment of $\Omega$ on $T_v B$ and defined for all $\kappa$, $\eta_1, \ldots, \eta_{n-1}$  in $T_v B$ by 
$$ \Omega_v(\kappa) (\eta_1, \ldots, \eta_{n-1}) = \Omega(\alpha) \big( {\mathcal L}^V_\alpha \eta_1, \ldots , {\mathcal L}^V_\alpha \eta_{n-1}) = \mathrm{Vol}_B(v)(\eta_1, \ldots, \eta_n), $$
with $\alpha = (v,\kappa)$ and using \eqref{eq:geom_orient_form_bundle}. Then, the following Fubini theorem holds true 
\begin{equation} 
\int_M f \, \mathrm{Vol}_M = \int_B \Big( \int_{{\mathbb R}^{n-1}} f \, \Omega \Big) \, \mathrm{Vol}_B. 
\label{eq:Fub_Fubthm}
\end{equation}
The same relations are true replacing the volume forms by their densities. Let $(\varepsilon_1, \ldots, \varepsilon_{n-1})$ be a local, positively oriented orthonormal frame of $B$. Define a coordinate map $T_v B \to {\mathbb R}^{n-1}$, $\kappa \mapsto (\kappa^1, \ldots, \kappa^{n-1})$ with $\kappa^i = \langle \kappa , \varepsilon_i \rangle_v$. It maps isometrically $T_vB$ onto ${\mathbb R}^{n-1}$. Then, in the coordinates $(\kappa^1, \ldots, \kappa^{n-1})$, the retrenchment $\Omega_v$ on $T_v B$ has expression
\begin{equation} 
\Omega_v = d \kappa^1 \wedge \ldots \wedge d \kappa^{n-1}. 
\label{eq:von_Mises_norm_prf3}
\end{equation}

\subsubsection{Differential operators on $M$}
\label{subsubsec:geom_gradients}

Let $f$: $M \to {\mathbb R}$ be a smooth function. We recall that the gradient of $f$ is the vector field $\nabla_\alpha f \in {\mathcal X}(M)$ on $M$ defined by 
\begin{equation} 
\big\langle \hspace{-1.8mm} \big\langle \nabla_\alpha f(\alpha) , \xi \big\rangle \hspace{-1.8mm} \big\rangle_\alpha = (df)_\alpha (\xi), \quad \forall \alpha \in M, \quad \forall \xi \in T_\alpha M. 
\label{eq:geom_grad_def}
\end{equation}
Below, we will define partial gradients of functions on $M$ with respect to $v$ or to $\kappa$ separately, denoted by $\nabla_v f$ and $\nabla_\kappa f$ respectively. The index $\alpha$ permits to remember that $\nabla_\alpha f$ is the gradient of $f$ in the classical sense (i.e. with respect to the pair $\alpha = (v,\kappa)$). Let $\Xi \in {\mathcal X}(M)$ be a smooth vector field. The divergence of $\Xi$ is the scalar field $\nabla_\alpha \cdot \Xi$ defined by 
\begin{equation} 
\int_M \nabla_\alpha \cdot \Xi \, f \, |\mathrm{Vol}_M| = - \int_M  \big\langle \hspace{-1.8mm} \big\langle \Xi , \, \nabla_\alpha f  \big\rangle \hspace{-1.8mm} \big\rangle \, |\mathrm{Vol}_M|, 
\label{eq:geom_div_def}
\end{equation}
for all smooth functions $f$: $M \to {\mathbb R}$. In other words, gradient and divergence are formal anti-adjoint to each other with respect to the $L^2$ inner-product. The index $\alpha$ in $\nabla_\alpha \cdot \Xi$ serves the same purpose as for the gradient: it will help us to distinguish divergences with respect to the pair $\alpha = (v,\kappa)$ from partial divergences with respect to $v$ or to $\kappa$ separately, denoted by $\nabla_v \cdot \Xi$ and $\nabla_\kappa \cdot \Xi$ respectively. We recall an alternate formula for the divergence \cite[\S 4.8]{gallot1990riemannian}. Let $\Xi \in {\mathcal X}({\mathcal M})$ and $D \Xi$ be the linear map ${\mathcal X}({\mathcal M}) \to {\mathcal X}({\mathcal M})$, $\Upsilon \mapsto \nabla^{\mathrm{M}}_\Upsilon \Xi$, where $\nabla^{\mathrm{M}}$ is the Levi-Civita connection associated with the Sasaki metric. Then, 
\begin{equation}
\nabla_\alpha \cdot \Xi = \mathrm{Tr} D \Xi. 
\label{eq:pk_div_LCC}
\end{equation}

\medskip
We now define the horizontal (i.e. ``in the $v$ direction'') and vertical (i.e. ``in the $\kappa$ direction'') gradients and divergences which we will denote by using the symbols $\nabla_v$ and $\nabla_\kappa$ respectively. The gradients map smooth functions $f$: $M \to {\mathbb R}$ to smooth vector fields $\Xi$ such that $\Xi(\alpha) \in T_{\pi \alpha}B$. Such vector fields are smooth sections of the pullback bundle $\pi^* TB$. Specifically, the pullback bundle $\pi^* TB = (\tilde{M}, \tilde{\pi}, M, F)$ is such that $\tilde{M} = \sqcup_{\alpha \in M} T_{\pi \alpha} B$ and for any $(\alpha, \kappa) \in \tilde{M}$ with $\kappa \in T_{\pi \alpha} B$, we have $\tilde{\pi} (\alpha, \kappa) = \alpha$. The space of smooth sections of $\pi^* TB$ will be denoted by ${\mathcal Y}(M)$. The horizontal (resp. vertical)  gradients $\nabla_v f$ (resp. $\nabla_\kappa f$) are defined by 
\begin{equation} 
\nabla_v f(\alpha) = {\mathcal B}^H_\alpha \big( \nabla_\alpha f(\alpha) \big), \qquad \nabla_\kappa f(\alpha) = {\mathcal B}^V_\alpha \big( \nabla_\alpha f(\alpha) \big), 
\label{eq:geom_horvert_grad}
\end{equation}
and $\nabla_v f$, $\nabla_\kappa f \in {\mathcal Y}(M)$. The horizontal (resp. vertical) divergences $\nabla_v \cdot \Xi$ (resp. $\nabla_\kappa \cdot \Xi$) of a field $\Xi \in {\mathcal Y}(M)$ are smooth functions $M \to {\mathbb R}$ defined by 
\begin{equation}
\nabla_v \cdot \Xi  = \nabla_\alpha \cdot \big( {\mathcal L}^H \Xi \big), \qquad \nabla_\kappa \cdot \Xi = \nabla_\alpha \cdot \big({\mathcal L}^V \Xi  \big). 
\label{eq:geom_horvert_div}
\end{equation}
Thanks to \eqref{eq:geom_SasvsEuc}, gradient and divergence are formal $L^2$-anti-duals, i.e. 
\begin{equation}
\int_M \nabla_v \cdot \Xi \, f \, |\mathrm{Vol}_M| = - \int_M \big\langle \Xi  ,  \nabla_v f \big\rangle_{\pi \alpha} \, |\mathrm{Vol}_M|, 
\label{eq:geom_graddiv_dual}
\end{equation}
and similarly by replacing $\nabla_v$ by $\nabla_\kappa$. Classical vector calculus identities carry over to $\nabla_v$ and $\nabla_\kappa$, i.e. 
\begin{equation}
\nabla_v (fg) = f \nabla_v g + g \nabla_v f, \quad \nabla_v \cdot (f \Xi) = \big( \nabla_v f , \Xi \big)_v + f \nabla_v \cdot \Xi. 
\label{eq:geom_graddiv_calculus}
\end{equation}
and similarly by replacing $\nabla_v$ by $\nabla_\kappa$.

\subsubsection{Remarkable formulas}
\label{subsubsec:rem_formulas}

In this section, we give various useful lemmas whose proofs can be found in Appendix \ref{sec:proofs_rem_formulas}. 

\begin{lemma}[Divergences of a field which does not depend on $\kappa$]~

\noindent
Let $k \in {\mathcal X}(B)$ and denote by $\nabla \cdot k$ its divergence. We define an element of ${\mathcal Y}(M)$ by $ \tilde k = k \circ \pi$. Such elements of ${\mathcal Y}(M)$ can be viewed as ``not depending on $\kappa$''. Then, we have 
\begin{equation}
\nabla_v \cdot \tilde k = \nabla \cdot k, \qquad \nabla_\kappa \cdot \tilde k = 0. 
\label{eq:geom_divk_vertlift}
\end{equation}
\label{lem:geom_divk_vertlift}
\end{lemma}

\begin{lemma}[Gradient of a function which only depends on $v$]~

\noindent
(i) Let $\varphi$: $B \to {\mathbb R}$ be smooth and let $\nabla \varphi \in {\mathcal X}(B)$ be its gradient. Define $\tilde \varphi$: $M \to {\mathbb R}$ by 
$\tilde \varphi = \varphi \circ \pi$. 
Then, we have 
\begin{equation} 
\nabla_v \tilde \varphi (\alpha) = \nabla \varphi (\pi \alpha), \qquad \nabla_\kappa \tilde \varphi (\alpha) = 0, \qquad \forall \alpha \in M. 
\label{eq:geom_nab_fct_v_only}
\end{equation}

\noindent(ii)
Conversely, let $\psi$: $M \to {\mathbb R}$ be smooth and such that $\nabla_\kappa \psi \equiv 0$. Then, there exists $\varphi$: $B \to {\mathbb R}$, smooth, such that $\psi = \varphi \circ \pi$. 
\label{lem:geom_nab_fct_v_only}
\end{lemma}

\medskip
\begin{lemma}[Gradient of $\langle \kappa, \kappa \rangle_v$]~

\noindent
We consider the function $M \to {\mathbb R}$, $\alpha \mapsto \langle \kappa, \kappa \rangle_v$, with $\alpha = (v,\kappa)$. Then, we have 
\begin{equation}
\nabla_\kappa \langle \kappa, \kappa \rangle_v = 2 \kappa, \qquad \nabla_v \langle \kappa, \kappa \rangle_v = 0.
\label{eq:geom_na_kappa2}
\end{equation}
\label{lem:geom_na_kappa2}
\end{lemma}

\subsubsection{Symplectic structure, Hamiltonian and Poisson bracket}
\label{subsubsec:appnavkap_Poisson}

There is a symplectic structure on $M$ which combines the symplectic structure on the cotangent bundle of $B$ (see e.g. \cite{arnol2013mathematical}) with the metric. We define a one-form $\omega^{(1)}$ on $M$ as follows. Let $\alpha = (v,\kappa) \in M$. Let $\xi \in T_\alpha M$. Then, 
\begin{equation}
\omega^{(1)}(\alpha)(\xi) = \big\langle \kappa, (d \pi)_\alpha (\xi) \big\rangle_v. 
\label{eq:appnavkap_symp_om1}
\end{equation}
The symplectic form $\omega^{(2)}$ on $TM$ is given by 
\begin{equation}
\omega^{(2)} = d \omega^{(1)}. 
\label{eq:appnavkap_symp_om2}
\end{equation}
Thanks to the non-degeneracy of $\omega^{(2)}$, we can define a linear isomorphism $I_\alpha$: $T_\alpha M^* \to T_\alpha M$, $\sigma \to I_\alpha(\sigma)$ as follows: 
\begin{equation}
\sigma (\xi) = \omega^{(2)}(\alpha) \big( \xi, I_\alpha(\sigma) \big), \qquad \forall \xi \in T_\alpha M. 
\label{eq:appnavkap_symp_I_def}
\end{equation}
Let $\mathbf{H}$: $M \times (0,\infty) \to {\mathbb R}$, $(\alpha, t) \mapsto \mathbf{H}(\alpha, t)$ be a smooth function which will be called the Hamiltonian. The Hamiltonian vector field with Hamiltonian $\mathbf{H}$ is the vector field $I(d_\alpha \mathbf{H})$, where $d_\alpha$ stands for the derivative with respect to $\alpha$ for fixed $t$. A Hamiltonian flow is the flow of a Hamiltonian vector field, i.e. the flow associated with a differential equation of the form $ d \alpha / dt = I (d_\alpha \mathbf{H}(\cdot,t)) (\alpha(t))$.

As an example, let us consider the Hamiltonian $\mathbf{H}_{\mathcal J}$: $M \to {\mathbb R}$ given by 
\begin{equation}
\mathbf{H}_{\mathcal J}(\alpha) = - \nu v \cdot {\mathcal J} + \frac{1}{2} \langle \kappa, \kappa \rangle_v, \quad \alpha = (v,\kappa) \in M, 
\label{eq:equi_first_integral}
\end{equation}
where ${\mathcal J} \in {\mathbb R}^n$ is a given vector. 
Then, we have 
\begin{equation} 
{\mathcal L}^H_{{\mathcal A}_i} K_i  + \nu {\mathcal L}^V_{{\mathcal A}_i}   P_{{\mathcal V_i}^\bot} {\mathcal J} = I(d \mathbf{H}_{{\mathcal J}(X_i(t),t)}) \big( {\mathcal A}_i(t) \big), 
\label{eq:pk_dVKi_hamilton_0}
\end{equation}
where $d \mathbf{H}_{{\mathcal J}(X_i(t),t)}$ denotes the derivative of the function $M \to {\mathbb R}$, $\alpha \mapsto \mathbf{H}_{{\mathcal J}(X_i(t),t)}$, for fixed $t$. Thus, \eqref{eq:pk_dVKi_hamilton_0} is a Hamiltonian vector field with Hamiltonian $\mathbf{H}_{{\mathcal J}(X_i(t),t)}$. Hence, Equation \eqref{eq:pk_dVKi_scaled} for ${\mathcal A}_i$ in the PCC model can be written
\begin{equation}
d {\mathcal A}_i = I(d \mathbf{H}_{{\mathcal J}(X_i(t),t)}) \big( {\mathcal A}_i(t) \big) \,  d t - \nu {\mathcal L}^V_{{\mathcal A}_i} K_i \, dt  + \sqrt{2 \bar D} \sum_{j=1}^n {\mathcal L}^V_{{\mathcal A}_i} \big( P_{{\mathcal V}_i^\bot} e_j \big) \circ dB^i_{j, t}. 
\label{eq:pk_dVKi_hamilton}
\end{equation}
In this form, we notice that the first term at the right-hand side of \eqref{eq:pk_dVKi_hamilton} is Hamiltonian while the last two terms are dissipative. We will notice a similar structure on the kinetic model and this structure will give rise to several important properties of the model.  

For this purpose, it is convenient to introduce the Poisson bracket associated with the symplectic form $\omega^{(2)}$. It is defined for two smooth functions $f$ and~$h$: $M \to {\mathbb R}$ by
\begin{equation}
\{f,h\} = \omega^{(2)} \big( I(dh), I(df) \big) = df \big( I(dh) \big), 
\label{eq:appnavkap_Poibra_def}
\end{equation}
and is again a smooth function $M \to {\mathbb R}$. The Poisson bracket is an anticommutative bilinear form on the space of smooth functions that satisfies the Jacobi identity and the Leibnitz rule. We can show that
\begin{equation}
I(d \mathbf{H} ) = {\mathbb J} \nabla_\alpha \mathbf{H},  \qquad {\mathbb J} = {\mathcal L}^H \circ {\mathcal B}^V - {\mathcal L}^V \circ {\mathcal B}^H.
\label{eq:appnavkap_symp_I(dH)=JnaH}
\end{equation}
and consequently
\begin{equation}
\{f,h\} = \big \langle \nabla_v f, \nabla_\kappa h \big \rangle - \big \langle \nabla_\kappa f, \nabla_v h \big \rangle. 
\label{eq:appnavkap_Poibra_express}
\end{equation}
By Liouville's theorem, $\omega^{(2)}$ is preserved by a Hamiltonian flow (including the case where this Hamiltonian is time-dependent) and so, any exterior power of $\omega^{(2)}$ is preserved as well. The $n$-th exterior power $(\omega^{(2)})^{\wedge n}$ is equal to the Riemannian volume form $\mathrm{Vol}_M$ up to a constant, namely
\begin{equation}
(\omega^{(2)})^{\wedge n} = n! \, (-1)^{\frac{n(n+1)}{2}} \mathrm{Vol}_M. 
\label{eq:appnavkap_symp_vol_form}
\end{equation}
Thus, the flow of the vector field ${\mathbb J} \nabla_\alpha \mathbf{H}$ on $M$, where $\mathbf{H}$ is any smooth Hamiltonian: $M \times (0,\infty) \to {\mathbb R}$, $(\alpha, t) \mapsto \mathbf{H}(\alpha,t)$, preserves the Sasaki volume form. It implies that ${\mathbb J} \nabla_\alpha \mathbf{H}$ is divergence free. Hence, we have 
\begin{equation}
\nabla_v \cdot \nabla_\kappa \mathbf{H} - \nabla_\kappa \cdot \nabla_v \mathbf{H} = 0. 
\label{eq:equi_navnakH-nakanavH}
\end{equation}
It results alternate expressions of the bracket of two smooth functions $f$ and $h$: 
\begin{equation}
\{f,h\} = \nabla_v \cdot \big( f \,  \nabla_\kappa h \big) - \nabla_\kappa \cdot \big( f \,  \nabla_v h \big) = \nabla_\kappa \cdot \big( h \, \nabla_v f \big) - \nabla_v \cdot \big( h \, \nabla_\kappa f \big). 
\label{eq:equi_bracket_conserv}
\end{equation}
From this and \eqref{eq:geom_graddiv_dual} (and its analog when $\nabla_v$ is changed in $\nabla_\kappa$), we deduce that if $g$ is another smooth function $M \to {\mathbb R}$, we have 
\begin{equation}
\int_M \{f,h\} \, g  \, |\mathrm{Vol}_M| = \int_M f \, \{h,g\} \, |\mathrm{Vol}_M| . 
\label{eq:equi_bracket_integ_byparts}
\end{equation}
In particular, if $g=1$, we find 
\begin{equation}
\int_M \{f,h\} \, |\mathrm{Vol}_M| = 0. 
\label{eq:equi_bracket_integ_0}
\end{equation}

\subsection{Kinetic model}
\label{subsubsec:pk_kinetic}

We first derive the KCC model in the single-particle case. We then successively present the many-particle KCC models in the non-normalized and normalized cases.

\subsubsection{Single particle KCC model}
\label{subsubsec:pk_kinetic_single}

For the single-particle PCC model, applying \cite[Theorem V.1.2]{ikeda2014stochastic}, we get 

\begin{theorem}[Kinetic model for the single particle PCC model]~

\noindent
Suppose ${\mathcal J}$: ${\mathbb R}^n \times [0,\infty) \to {\mathbb R}^n$, $(x,t) \mapsto {\mathcal J}(x,t)$ is given and smooth. Consider the single-particle PCC model, given by \eqref{eq:pk_dXi_scaled}, \eqref{eq:pk_dVKi_scaled} with $i=N=1$ associated with this choice of ${\mathcal J}$. Let $(X(t),{\mathcal A}(t))$ be a trajectory of this stochastic differential system  with initial condition $(X_0, {\mathcal A}_0)$ drawn according to the probability distribution $f_0 \, |dx \wedge \mathrm{Vol}_M|$, where $f_0$: ${\mathbb R}^n \times M \to [0,\infty)$ is a given integrable function. Then, the probability distribution $f(\cdot,\cdot,t) \, |dx \wedge \mathrm{Vol}_M|$ of $(X(t),{\mathcal A}(t))$ at time $t$ (where $f(\cdot,\cdot,t)$: ${\mathbb R}^n \times M \to [0,\infty)$ is an integrable function) satisfies the following kinetic equation: 
\begin{equation}
\partial_t f + \nabla_x \cdot \big(v f \big) + \nabla_v \cdot ( \kappa f ) - \nu \nabla_\kappa \cdot  \big( (\kappa - P_{v^\bot} {\mathcal J}) \, f \big) - D \Delta_\kappa f = 0, 
\label{eq:pk_eq_f_single_noise}
\end{equation}
with $ \Delta_\kappa f =: \nabla_\kappa \cdot (\nabla_\kappa f)$ being the so-called vertical Laplacian of $f$ on $M$. Equation \eqref{eq:pk_eq_f_single_noise} will be referred to as the kinetic curvature-control (KCC) model. 
\label{thm:pk_eq_f_single_noise}
\end{theorem}

\begin{remark}
Eq. \eqref{eq:pk_eq_f_single_noise} has the same form as if $(v,\kappa)$ belonged to a standard linear space. However, we should keep in mind that the meaning of the operators $\nabla_v \cdot$ and $\nabla_\kappa \cdot$ is given by~\eqref{eq:geom_horvert_div}. 
\end{remark}

\medskip
\noindent
Let $\mathbf{H}_{\mathcal J}$ be the Hamiltonian given by \eqref{eq:equi_first_integral}. Thanks to \eqref{eq:geom_nab_fct_v_only} and \eqref{eq:geom_na_kappa2}, we have 
\begin{equation}
\nabla_v \mathbf{H}_{\mathcal J} = - \nu P_{v^\bot}{\mathcal J}, \qquad \nabla_\kappa \mathbf{H}_{\mathcal J} = \kappa. 
\label{eq:equi_gradients_first_int}
\end{equation}
Thus, in view of \eqref{eq:appnavkap_Poibra_express}, the kinetic equation \eqref{eq:pk_eq_f_single_noise} can be written 
\begin{equation}
\partial_t f + \nabla_x \cdot \big(v f \big) + \big\{ f , \mathbf{H}_{\mathcal J} \big\} - \nu \nabla_\kappa \cdot (\kappa \, f ) - D \Delta_\kappa f = 0. 
\label{eq:pk_eq_f_poibra}
\end{equation}

Eq. \eqref{eq:pk_eq_f_single_noise} is of Fokker-Planck type, i.e. the second-order operator is a partial diffusion which operates on the $\kappa$-variable only and not on the $x$ and $v$ variables. Fokker-Planck equations on manifolds have previously been considered e.g. in \cite{bismut2005hypoelliptic, calogero2012exponential, grothaus2022hypocoercivity, lebeau2005geometric, nier2024global}. However, the specific Fokker-Planck equation \eqref{eq:pk_eq_f_single_noise} posed on the cartesian product of ${\mathbb R}^n$ and the tangent bundle to the sphere ${\mathbb S}^{n-1}$ has not been studied previously, to the best of our knowledge.

\subsubsection{Many particle KCC (non-normalized case)}
\label{subsubsec:pk_kinetic_many_nonnormalized}

We consider the PCC model \eqref{eq:pk_dXi_scaled}, \eqref{eq:pk_dVKi_scaled}, \eqref{eq:pk_Jinonnorm_scaled} in the non-normalized case. We associate it with its empirical measure $\mu_t^N$ on ${\mathbb R}^n \times M$ such that 
$$ \mu_t^N(x,\alpha) = \frac{1}{N} \sum_{i=1}^N \delta_{(X_i(t), {\mathcal A}_i(t))}(x,\alpha), \quad x \in {\mathbb R}^N, \quad \alpha \in M, $$
and we denote by $\mu^N_0$ the empirical measure associated with the initial data $((X_{i0}, {\mathcal A}_{i0}))_{i=1, \ldots, N}$. We now consider $\mu_t^N$ and $\mu_0^N$ as random measures. An adaptation of classical results on mean-field limits (see e.g. the review \cite{chaintron2022propagation1, chaintron2021propagation2}) would show the following:

\begin{theorem}[Kinetic model for the many-particle PCC model (non-normalized case)]~

\noindent
Let $\mu^N_0$ converge as $N \to \infty$ in probability sense towards a deterministic absolutely continuous probability measure $f_0 \, |dx \wedge \mathrm{Vol}_M|$. Then, for any $t \in [0,\infty)$, $\mu_t^N$ converges in probability sense towards a deterministic absolutely continuous probability measure $f_t \, |dx \wedge \mathrm{Vol}_M|$ and $f(x,\alpha,t) \equiv f_t(x,\alpha)$ satisfies \eqref{eq:pk_eq_f_single_noise} with ${\mathcal J} = {\mathcal J}_f(x,t)$ given by: 
\begin{equation}
{\mathcal J}_f(x,t) = \frac{1}{R^n} \int_{{\mathbb R}^n \times M} {\mathcal K} \Big( \frac{|x-y|}{R} \Big) \, v \, f(y,\alpha,t) \, |dx \wedge \mathrm{Vol}_M|(y,\alpha). 
\label{eq:pk_eq_J}
\end{equation}
We do not precise the functional setting and refer to \cite{chaintron2022propagation1, chaintron2021propagation2} for details. 
\label{thm:pk_eq_f_many_nonorm}
\end{theorem}

\subsubsection{Many particle KCC (normalized case)}
\label{subsubsec:pk_kinetic_many_normalized}

We consider the particle system \eqref{eq:pk_dXi_scaled}, \eqref{eq:pk_dVKi_scaled}, \eqref{eq:pk_Jinorm_scaled} in the normalized case. The normalization brings a major difficulty to prove a rigorous convergence result from the particle dynamics to the kinetic model and has not been solved yet, even for the Vicsek model, except when the singularity brought by the normalization is regularized \cite{briant2022cauchy}. Here, we will assume that the following is true

\begin{conjecture}[Kinetic model for the many-particle PCC model (normalized case)]~

\noindent
With the same assumptions as for Theorem \ref{thm:pk_eq_f_many_nonorm}, and in the generic case where the denominator of~\eqref{eq:pk_Jinorm_scaled} does not vanish, $\mu_t^N$ converges towards a deterministic absolutely continuous probability measure $f_t \, |dx \wedge \mathrm{Vol}_M|$ and $f(x,\alpha,t) \equiv f_t(x,\alpha)$ satisfies \eqref{eq:pk_eq_f_single_noise} with ${\mathcal J} = {\mathcal J}_f(x,t)$ given by: 
\begin{equation}
{\mathcal J}_f(x,t) = \frac{\displaystyle \int_{{\mathbb R}^n \times M} {\mathcal K} \Big( \frac{|x-y|}{R} \Big) \, v \, f(y,\alpha,t) \, |dx \wedge \mathrm{Vol}_M|(y,\alpha)}{\displaystyle \bigg|\int_{{\mathbb R}^n \times M} {\mathcal K} \Big( \frac{|x-y|}{R} \Big) \, v \, f(y,\alpha,t) \, |dx \wedge \mathrm{Vol}_M|(y,\alpha) \bigg|}. 
\label{eq:pk_eq_Jnorm}
\end{equation}
\label{conject:pk_eq_f_many_norm}
\end{conjecture}

\subsection{Hydrodynamic scaling}
\label{subsec:pk_hydro_scaling}

In the remainder of this paper, \textbf{we focus on the normalized KCC model} and derive its hydrodynamic limit. Indeed, the non-normalized model leads to additional difficulties stemming from the existence of phase transitions which will be dealt with in future work. 

In the hydrodynamic scaling, temporal and spatial phenomena are supposed to occur at large scales compared with the dynamics of the $\alpha$ variable. To express this, as classical in kinetic theory \cite[Sect. 3.8]{cercignani2013mathematical}, we change the time and space scales and introduce the following change of variables:
$$ \hat x = \varepsilon x, \qquad \hat t = \varepsilon t, \quad \textrm{ with } \quad \varepsilon \ll 1.  $$
At the same time, we assume that the trajectory curvatures are small at the microscopic scale, so that they are order unity at the macroscopic scale. For that reason, we do not rescale $\kappa$. In the new variables $\hat x$, $\hat t$, the unknowns $f$ and ${\mathcal J}$ are parametrized by $\varepsilon$ and will be denoted by $f^\varepsilon$, ${\mathcal J}^\varepsilon$. Now, the equations for $f^\varepsilon$, ${\mathcal J}^\varepsilon$ are given in the normalized case by (dropping the hats on~$\hat x$,~$\hat t$): 
\begin{equation}
\partial_t f^\varepsilon +  \nabla_x \cdot (v f^\varepsilon ) + \frac{1}{\varepsilon} \big( \big\{ f^\varepsilon , \mathbf{H}_{{\mathcal J}^\varepsilon_{f^\varepsilon}} \big\} - \nu \nabla_\kappa \cdot (\kappa \, f^\varepsilon )  + D \nabla_\kappa f^\varepsilon \big) = 0, 
\label{eq:pk_eq_f_hydro_scaling}
\end{equation}
with ${\mathcal J}^\varepsilon_f$ given for any $f(x,\alpha,t)$ by 
\begin{equation}
{\mathcal J}^\varepsilon_f(x,t) = \frac{\displaystyle \int_{{\mathbb R}^n \times M} {\mathcal K} \Big( \frac{|x-y|}{\varepsilon R} \Big) \, v \,  f( y, \alpha, t) \, |dx \wedge \mathrm{Vol}_M|(y,\alpha)}{\displaystyle \bigg| \int_{{\mathbb R}^n \times M} {\mathcal K} \Big( \frac{|x-y|}{\varepsilon R} \Big) \, v \,  f( y, \alpha, t) \, |dx \wedge \mathrm{Vol}_M|(y,\alpha) \bigg|}, 
\label{eq:pk_eq_Jnorm_hydro_scaling}
\end{equation}
We expand \eqref{eq:pk_eq_Jnorm_hydro_scaling} in powers of $\varepsilon$. In view of \eqref{eq:partelem_normalization_H}, straightforward calculations give
\begin{equation}
{\mathcal J}^\varepsilon_f = \mathbf{u}_f + {\mathcal O}(\varepsilon^2), 
\label{eq:pk_expan_Jeps}
\end{equation}
with 
\begin{equation}
\mathbf{u}_f(x,t) = \frac{\displaystyle \int_M v \,  f( x, \alpha, t) \, |\mathrm{Vol}_M|(\alpha)}{\displaystyle \bigg| \int_M v \,  f( x, \alpha, t) \, |\mathrm{Vol}_M|(\alpha) \bigg|}, 
\label{eq:pk_eq_Jnorm_hydro_lim}
\end{equation}
We only keep the leading order in \eqref{eq:pk_expan_Jeps} as the following orders are not involved in the limit $\varepsilon \to 0$ which we are interested in. 

Finally, we are left with the following perturbation problem: 
\begin{equation}
\partial_t f^\varepsilon +  \nabla_x \cdot (v f^\varepsilon ) + \frac{1}{\varepsilon} \big( \big\{ f^\varepsilon , \mathbf{H}_{\mathbf{u}_{f^\varepsilon}} \big\} - \nu \nabla_\kappa \cdot (\kappa \, f^\varepsilon )  + D \nabla_\kappa f^\varepsilon \big) = 0, 
\label{eq:pk_eq_f_hydro_scal_simple}
\end{equation}
with $\mathbf{u}_f$ given by \eqref{eq:pk_eq_Jnorm_hydro_lim}. We will refer to this system as the scaled KCC model. Letting $\varepsilon \to 0$ in~\eqref{eq:pk_eq_f_hydro_scal_simple} is referred to as the hydrodynamic limit. The main result of this paper is the determination of the limit. It will be useful to introduce the collision operator
\begin{equation}
Q(f) = - \big\{ f , \mathbf{H}_{\mathbf{u}_f} \big\} + \nu \nabla_\kappa \cdot (\kappa \, f ) + D \Delta_\kappa f, 
\label{eq:pk_collision_operator}
\end{equation}
so that \eqref{eq:pk_eq_f_hydro_scal_simple} can be simply written
\begin{equation}
\partial_t f^\varepsilon +  \nabla_x \cdot (v f^\varepsilon ) = \frac{1}{\varepsilon} Q(f^\varepsilon). 
\label{eq:pk_eq_f_Q}
\end{equation}
We note that $Q$ is nonlinear due to the dependence of $\mathbf{u}$ on $f$.

\setcounter{equation}{0}
\section{The hydrodynamic model}
\label{sec:hydro}

In this section, we derive the equilibria of the collision operator. Then, we state the main theorem of this paper which establishes that the scaled KCC model \eqref{eq:pk_eq_f_Q} converges to a hydrodynamic-type model (referred to as the ``self-organized hydrodynamic'' model (SOH)) in the limit $\varepsilon \to 0$. Finally, we show that the first equation of the SOH model, the continuity equation is easily derived from the KCC model. But first of all, we introduce a useful geometric change of variables.

\subsection{Geometric preliminary: a change of variables in $B$}
\label{subsubsec:hydro_chg_var}

We denote by $B_\mathbf{u} = B \setminus \{ \pm \mathbf{u} \}$. We identify ${\mathbb S}^{n-2} \approx B \cap \{\mathbf{u}\}^\bot$. We introduce the diffeomorphism ${\mathbf V}_\mathbf{u}$: $(0,\pi) \times {\mathbb S}^{n-2} \to B_\mathbf{u}$, $(\theta, w) \mapsto \mathbf{V}_\mathbf{u}(\theta, w)$, defined by
\begin{equation} 
\mathbf{V}_\mathbf{u}(\theta, w) = v = \cos \theta \, \mathbf{u} + \sin \theta \, w,    
\label{eq:macnorm_V_def}
\end{equation}
where 
\begin{equation} 
\cos \theta = \mathbf{u} \cdot v, \qquad |w|=1, \qquad \mathbf{u} \cdot w = 0. 
\label{eq:macnorm_V_def_2}
\end{equation}
The diffeomorphism ${\mathbf V}_\mathbf{u}$ is depicted in dimension $n=3$ in Fig. \ref{fig:chgvar} of Section \ref{subsubsec:macnorm_chgvar_M} below. We orient ${\mathbb S}^{n-2}$ in such a way that $(\varepsilon_1, \ldots, \varepsilon_{n-2})$ is a direct orthonormal basis of $T_w {\mathbb S}^{n-2}$ if and only if $(e_\theta, \varepsilon_1, \ldots, \varepsilon_{n-2})$ is a direct orthnormal basis of $T_v {\mathbb S}^{n-1}$, where 
\begin{equation}
e_\theta = - \sin \theta \, \mathbf{u} + \cos \theta \, w. 
\label{eq:hydro_chgvar_ethetdef}
\end{equation}
The volume form $\mathrm{Vol}_{{\mathbb S}^{n-2}}$ on ${\mathbb S}^{n-2}$ is defined analogously to $\mathrm{Vol}_B$, i.e. it is the Riemannian volume form for the Riemannian structure induced by the embedding of ${\mathbb S}^{n-2}$ in ${\mathbb R}^n$. Now, we have (the proof is left to the reader)
\begin{equation}
|\mathrm{Vol}_B| = \sin^{n-2} \theta \, d \theta \, |\mathrm{Vol}_{{\mathbb S}^{n-2}}|. 
\label{eq:hydro_chgvar_v2thetaw}
\end{equation}

\subsection{Equilibria}
\label{subsec:hydro_equi}

Let $\mathbf{u} \in{\mathbb S}^{n-1}$. We define the (generalized) von Mises distribution $\mathbf{M}_{\mathbf{u}}$: $M \to {\mathbb R}$ by 
\begin{equation}
\mathbf{M}_{\mathbf{u}} = \frac{1}{Z} \exp \Big( - \frac{\nu}{D} \mathbf{H}_{\mathbf{u}} \Big),
\label{eq:equi_von_Mises}
\end{equation}
with $Z$ such that 
\begin{equation} 
\int_{M} \mathbf{M}_{\mathbf{u}} \, |\mathrm{Vol}_M| = 1. 
\label{eq:equi_von_Mises_normal}
\end{equation}

\begin{lemma}
The normalization factor $Z$ is given by 
\begin{equation}
Z = \Big( \frac{2 \pi D}{\nu} \Big)^{\frac{n-1}{2}} \, | {\mathbb S}^{n-2} | \int_0^\pi \exp \Big( \frac{\nu^2}{D} \cos \theta \Big) \, \sin^{n-2} \theta \, d \theta, 
\label{eq:von_Mises_norm}
\end{equation}
(where $| {\mathbb S}^{n-2} |$ denotes the $n-2$-dimensional measure of ${\mathbb S}^{n-2}$) and does not depend on $\mathbf{u}$. 
\label{lem:von_Mises_norm}
\end{lemma}

\noindent
\textbf{Proof.} Using Fubini's theorem \eqref{eq:Fub_Fubthm}, together with \eqref{eq:Fub_sigma(x)_def}, \eqref{eq:von_Mises_norm_prf3} and \eqref{eq:equi_first_integral}, we have  
\begin{eqnarray*}
Z &=& \int_M \exp \Big( - \frac{\nu}{D} \mathbf{H}_{\mathbf{u}}(\alpha) \Big) \, |\mathrm{Vol}_M|(\alpha) \\
&=& \int_B  \exp \Big(\frac{\nu^2}{D} (v \cdot \mathbf{u}) \Big) \, \bigg( 
\int_{{\mathbb R}^{n-1}} \exp \Big( - \frac{\nu}{2D} \sum_{i=1}^{n-1} (\kappa^i)^2 \Big) \, d \kappa^1 \ldots d \kappa^{n-1} \bigg) \, |\mathrm{Vol}_B|(v) \\
&=& \Big( \frac{2 \pi D}{\nu} \Big)^{\frac{n-1}{2}} \, \int_B \exp \Big( \frac{\nu^2}{D} v \cdot \mathbf{u} \Big) \, |\mathrm{Vol}_B|(v). 
\end{eqnarray*}
Now, we use the diffeomorphism described in Section \ref{subsubsec:hydro_chg_var} and thanks to \eqref{eq:hydro_chgvar_v2thetaw}, we get 
$$\int_B \exp \Big( \frac{\nu^2}{D} v \cdot \mathbf{u} \Big) \, |\mathrm{Vol}_B|(v) = | {\mathbb S}^{n-2} | \int_0^\pi \exp \Big( \frac{\nu^2}{D}\cos \theta \Big) \, 
\sin^{n-2} \theta \, d \theta, $$
which leads to \eqref{eq:von_Mises_norm}. \endproof

Thanks to the von Mises distribution, we can write the collision operator as follows
\begin{equation}
Q(f) = - \big\{f , \mathbf{H}_{\mathbf{u}_f} \big\} + D \nabla_\kappa \cdot \bigg( \mathbf{M}_{\mathbf{u}_f} \nabla_\kappa \Big( \frac{f}{\mathbf{M}_{\mathbf{u}_f}} \Big) \bigg). 
\label{eq:equi_coll_oper_alt}
\end{equation}

\begin{lemma}[Negativity of $Q$]~

\noindent
Let $f$: $M \to {\mathbb R}$ be smooth and decaying fast enough at infinity in $\kappa$. Then, we have
\begin{equation}
\int_M Q(f) \, \frac{f}{\mathbf{M}_{\mathbf{u}_f}} \, |\mathrm{Vol}_M| = - D \int_M \Big\langle \nabla_\kappa \Big( \frac{f}{\mathbf{M}_{\mathbf{u}_f}} \Big) \, , \, \nabla_\kappa \Big( \frac{f}{\mathbf{M}_{\mathbf{u}_f}} \Big)
\Big\rangle 
 \, \mathbf{M}_{\mathbf{u}_f} \, |\mathrm{Vol}_M| \leq 0. 
\label{eq:equi:coll_oper_negat}
\end{equation}
\label{lem:equi:coll_oper_negat}
\end{lemma}

\noindent
\textbf{Proof.} Using \eqref{eq:geom_graddiv_dual}, we have 
$$ \int_M \nabla_\kappa \cdot \bigg( \mathbf{M}_{\mathbf{u}_f} \nabla_\kappa \Big( \frac{f}{\mathbf{M}_{\mathbf{u}_f}} \Big) \bigg) \, \frac{f}{\mathbf{M}_{\mathbf{u}_f}} \, |\mathrm{Vol}_M| 
= - \int_M \Big\langle \nabla_\kappa \Big( \frac{f}{\mathbf{M}_{\mathbf{u}_f}} \Big) \, , \, \nabla_\kappa \Big( \frac{f}{\mathbf{M}_{\mathbf{u}_f}} \Big) \Big\rangle
 \, \mathbf{M}_{\mathbf{u}_f} \, |\mathrm{Vol}_M|. $$
On the other hand, using the fact that the Poisson bracket is a derivation as well as \eqref{eq:equi_von_Mises}, we can write 
\begin{eqnarray*} 
\int_M \big\{f , \mathbf{H}_{\mathbf{u}_f} \big\} \, \frac{f}{\mathbf{M}_{\mathbf{u}_f}} \, |\mathrm{Vol}_M| 
&=& \frac{Z}{2} \int_M \big\{f^2 , \mathbf{H}_{\mathbf{u}_f} \big\} \, \exp \Big( \frac{\nu}{D} \mathbf{H}_{\mathbf{u}_f} \Big) \, |\mathrm{Vol}_M| \\
&=& \frac{D Z}{2 \nu} \int_M \Big\{f^2 , \exp \Big( \frac{\nu}{D} \mathbf{H}_{\mathbf{u}_f} \Big) \Big\} \, |\mathrm{Vol}_M| = 0, 
\end{eqnarray*}
thanks to \eqref{eq:equi_bracket_integ_0}. Inserting the expression \eqref{eq:equi_coll_oper_alt} for the collision operator at the left-hand side of~\eqref{eq:equi:coll_oper_negat} and using the above computation leads to the right-hand side of \eqref{eq:equi:coll_oper_negat} and ends the proof.  \endproof

We introduce the coefficient $c_1$ called the ``order parameter'', defined by: 
\begin{equation}
c_1 = \frac{\displaystyle \int_B \exp \Big( \frac{\nu^2}{D} v \cdot \mathbf{u} \Big) \, (v \cdot \mathbf{u}) \, |\mathrm{Vol}_B|(v)}
{\displaystyle \int_B \exp \Big( \frac{\nu^2}{D} v \cdot \mathbf{u} \Big) \, |\mathrm{Vol}_B|(v)} = \frac{\displaystyle \int_0^\pi \exp \Big( \frac{\nu^2}{D} \cos \theta \Big) \, \cos \theta \, \sin^{n-2} \theta \, d \theta}
{\displaystyle \int_0^\pi \exp \Big( \frac{\nu^2}{D} \cos \theta \Big) \, \sin^{n-2} \theta \, d \theta},
\label{eq:equi_order_param}
\end{equation}
which does not depend on $\mathbf{u}$.

\begin{lemma}[Properties of the order parameter]~

\noindent
(i) $c_1 \in [0,1)$. 

\noindent
(ii) For $\mathbf{u} \in{\mathbb S}^{n-1}$, we have 
\begin{equation}
\int_M \mathbf{M}_{\mathbf{u}}(\alpha) \, v \, |\mathrm{Vol}_M|(\alpha) = 
c_1 \, \frac{\mathbf{u}}{|\mathbf{u}|} . 
\label{eq:equi_von_mises_flux}
\end{equation}
\label{lem:equi_von_mises_flux}
\end{lemma}

\noindent
\textbf{Proof.} This is a classical result which has been proved in the literature \cite{degond2013macroscopic}. \endproof

\begin{proposition}[Equilibria]~

\noindent
Let $f$: $M \to {\mathbb R}$ be smooth and decaying fast enough at infinity in $\kappa$. Then, we have 
\begin{equation}
Q(f) = 0, 
\label{eq:equi_collision_zero}
\end{equation}
if and only if \, $\exists \rho \geq 0$, \, $\exists \mathbf{u} \in {\mathbb S}^{n-1}$ such that 
\begin{equation}
f = \rho \mathbf{M}_{\mathbf{u}}, 
\label{eq:equi_equi}
\end{equation}
Distributions of the type \eqref{eq:equi_equi} are called ``equilibria''. 

\label{prop:equi_equi}
\end{proposition}

\noindent
\textbf{Proof.} Let $f$ be such that $Q(f)=0$. Then, by \eqref{eq:equi:coll_oper_negat}, we have
$$0 = \int_M Q(f) \, \frac{f}{\mathbf{M}_{\mathbf{u}_f}} \, |\mathrm{Vol}_M| = - D \int_M \Big\langle \nabla_\kappa \Big( \frac{f}{\mathbf{M}_{\mathbf{u}_f}} \Big) \, , \, \nabla_\kappa \Big( \frac{f}{\mathbf{M}_{\mathbf{u}_f}} \Big) \Big\rangle
 \, \mathbf{M}_{\mathbf{u}_f} \, |\mathrm{Vol}_M|. $$
Since the integrated quantity in the last integral is nonnegative and $\mathbf{M}_{\mathbf{u}_f}>0$, this implies $\nabla_\kappa ( f/\mathbf{M}_{\mathbf{u}_f}) = 0$, hence, by Lemma \ref{lem:geom_nab_fct_v_only} (ii), there exists a function $\varphi$: $B \to {\mathbb R}$ such that 
$$f  = \mathbf{M}_{\mathbf{u}_f} \, \varphi \circ \pi. $$
Then, we compute, thanks to \eqref{eq:geom_nab_fct_v_only}, 
\begin{eqnarray*}
\nabla_v f &=& \mathbf{M}_{\mathbf{u}_f} \, (\nabla_v \varphi) \circ \pi - \frac{\nu}{D} \, \mathbf{M}_{\mathbf{u}_f} \, \varphi \circ \pi \, \nabla_v \mathbf{H}_{\mathbf{u}_f}, \\
\nabla_\kappa f &=& - \frac{\nu}{D} \, \mathbf{M}_{\mathbf{u}_f} \, \varphi \circ \pi \, \nabla_\kappa \mathbf{H}_{\mathbf{u}_f}. 
\end{eqnarray*}
Thus, 
$$ Q(f) = - \{ f, \mathbf{H}_{\mathbf{u}_f} \} = - \mathbf{M}_{\mathbf{u}_f} \, \langle \kappa, (\nabla_v \varphi) \circ \pi \rangle.$$
Since $Q(f) = 0$, we have $\langle \kappa, (\nabla_v \varphi) \circ \pi \rangle = 0$, $\forall \alpha \in M$, which implies $\nabla_v \varphi = 0$, hence, $\varphi$ is a constant on $B$. Therefore, there exists a constant $\rho \in {\mathbb R}$ such that $f = \rho \mathbf{M}_{\mathbf{u}_f}$. Since $f \geq 0$, $\rho \in [0,\infty)$. Consequently, $f$ is of the form \eqref{eq:equi_equi}. Conversely, suppose $f$ is given by \eqref{eq:equi_equi} for arbitrary $\mathbf{u} \in {\mathbb S}^{n-1}$ and $\rho \in [0,\infty)$. Then, computing $\mathbf{u}_f$ using \eqref{eq:pk_eq_Jnorm_hydro_lim} with the help of \eqref{eq:equi_order_param}, \eqref{eq:equi_von_mises_flux} and of the fact that $c_1$ is non-negative leads to the fact that any $\mathbf{u}_f = \mathbf{u}$, which shows that such $f$ is a solution of \eqref{eq:equi_collision_zero}. This ends the proof. \endproof

\subsection{Hydrodynamic limit: main result}
\label{subsec:hydro_main}

The main result of this paper is the following 

\begin{theorem}[Hydrodynamic limit of the KCC model (formal result)]~

\noindent
Suppose $f^\varepsilon$ is a smooth solution of the scaled KCC model \eqref{eq:pk_eq_f_Q}, which converges as $\varepsilon \to 0$ towards $f$ as smoothly as needed. Then, there exist two functions $\rho$: ${\mathbb R}^n \times [0,\infty) \to (0,\infty)$, $(x,t) \mapsto \rho(x,t)$ and $\mathbf{u}$: ${\mathbb R}^n \times [0,\infty) \to {\mathbb S}^{n-1}$, $(x,t) \mapsto \mathbf{u}(x,t)$ such that 
\begin{equation}
f(x,\alpha,t) = \rho(x,t) \mathbf{M}_{\mathbf{u}(x,t)}(\alpha), \quad \forall (x,\alpha,t) \in {\mathbb R}^n \times M \times [0,\infty). 
\label{eq:macnorm_equi}
\end{equation}
Furthermore, $\rho$ and $\mathbf{u}$ satisfy the following ``self-organized hydrodynamic (SOH)'' system 
\begin{eqnarray}
&&\hspace{-1cm}
\partial_t \rho + \nabla_x \cdot \big( c_1 \rho \mathbf{u} \big) = 0. 
\label{eq:macnorm_continuity} \\
&&\hspace{-1cm}
\rho \big( \partial_t \mathbf{u} + c_2 (\mathbf{u} \cdot \nabla_x) \mathbf{u} \big) + c_3 P_{\mathbf{u}^\bot} \nabla_x \rho  = 0, 
\label{eq:macnorm_eq_evol_u}
\end{eqnarray}
with $c_1$ given by \eqref{eq:equi_order_param}, $c_3$, by 
\begin{equation}
c_3 = \frac{D}{\nu^2}.  
\label{eq:macnorm_eq_coeffs_c3}
\end{equation}
To define $c_2$, we first introduce $\mathbf{M}_\mathbf{u}^\natural$: $[0,\pi] \times {\mathbb R} \times [0,\infty) \to {\mathbb R}$, given by
\begin{equation}
 \mathbf{M}_\mathbf{u}^\natural = \mathbf{M}_\mathbf{u}^\natural (\theta, \kappa_\parallel, \kappa_\bot) = \frac{1}{Z} \exp \Big( - \frac{\nu}{D} \big( - \nu \cos \theta + \frac{1}{2} (\kappa_\parallel^2 +\sin^2 \theta \, \kappa_\bot^2) \big) \Big), 
\label{eq:macnorm_Munat}
\end{equation}
and the solution $(\zeta_w, \zeta_r)$, $[0,\pi] \times {\mathbb R} \times [0,\infty) \to {\mathbb R}^2$, of the following system: 
\begin{eqnarray}
&&\hspace{-1cm}
{\mathcal U} \zeta_w  - \kappa_\bot \, \zeta_r = \sin \theta, \label{eq:macnorm_nD_vectGCI_wcomp} \\
&&\hspace{-1cm}
{\mathcal U} \zeta_r - \frac{(n-3) D}{\sin^2 \theta \, \kappa_\bot^2} \zeta_r + \kappa_\bot \, \zeta_w = 0. \label{eq:macnorm_nD_vectGCI_kapTcomp}
\end{eqnarray}
where the operator ${\mathcal U}$ operates on any smooth function $\zeta$: $(0,\pi) \times {\mathbb R} \times [0,\infty) \to {\mathbb R}$, $(\theta, \kappa_\parallel, \kappa_\bot) \to \zeta$ as follows: 
\begin{eqnarray}
{\mathcal U} &=& {\mathcal U}_1 + {\mathcal U}_2,  \label{eq:macnorm_nD_Uoperator_def} \\
{\mathcal U}_1 \zeta &=& \kappa_\parallel \partial_\theta \zeta + \big(\cos \theta \, \sin \theta \, \kappa_\bot^2 - \nu \, \sin \theta \big) \partial_{\kappa_\parallel} \zeta -  2 \frac{\cos \theta}{\sin \theta} \, \kappa_\parallel \, \kappa_\bot \, \partial_{\kappa_\bot} \zeta, \label{eq:macnorm_nD_U1operator_def} \\
{\mathcal U}_2 \zeta &=&  D  \Big( \frac{1}{\mathbf{M}_\mathbf{u}^\natural} \partial_{\kappa_\parallel} ( \mathbf{M}_\mathbf{u}^\natural \partial_{\kappa_\parallel} \zeta ) + \frac{1}{\sin^2 \theta \, \mathbf{M}_\mathbf{u}^\natural \, \kappa_\bot^{n-3}} \, \partial_{\kappa_\bot} ( \mathbf{M}_\mathbf{u}^\natural \, \kappa_\bot^{n-3} \, \partial_{\kappa_\bot} \zeta ) \Big). \label{eq:macnorm_nD_U2operator_def}  
\end{eqnarray}
Then, $c_2$ is given by
\begin{equation}
c_2 = \frac{\displaystyle \int_{(0,\pi) \times {\mathbb R} \times (0,\infty)} \cos \theta \, (\mathbf{M}_\mathbf{u}^\natural \zeta_w)(\theta, \kappa_\parallel, \kappa_\bot) \, \sin^{2n-3} \theta \, \, \kappa_\bot^{n-3} \, d\theta \, d\kappa_\parallel \, d\kappa_\bot} 
{\displaystyle \int_{(0,\pi) \times {\mathbb R} \times (0,\infty)} (\mathbf{M}_\mathbf{u}^\natural \zeta_w)(\theta, \kappa_\parallel, \kappa_\bot) \, \sin^{2n-3} \theta \, \, \kappa_\bot^{n-3} \, d\theta \, d\kappa_\parallel \, d\kappa_\bot}. 
\label{eq:macnorm_nD_coef_c2}
\end{equation}
\label{thm:hydro_main_theorem}
\end{theorem}

\begin{remark}

\noindent
(i) System \eqref{eq:macnorm_nD_vectGCI_wcomp}, \eqref{eq:macnorm_nD_vectGCI_kapTcomp} will be given a variational formulation and shown to be well-posed below (see \eqref{eq:macnorm_redGCI_varform} for the variational formulation in dimension $n \geq 4$ and \eqref{eq:macnorm_redGCI_varform_n=3} in dimension $n=3$). Boundary conditions associated with this variational formulation will follow from the variational formulation. Only the component $\zeta_w$ is needed to compute the coefficient~$c_2$. However, $\zeta_w$ and $\zeta_r$ are coupled and $\zeta_w$ cannot be found independently of $\zeta_r$. 

\noindent
(ii) $\mathbf{M}_\mathbf{u}^\natural$ is nothing but the expression of the von Mises $M_\mathbf{u}$ in an appropriate change of variables which will be introduced below. 
\end{remark}

In dimension 2, it was already proved in \cite{degond2011macroscopic} that the hydrodynamic limit of the KCC model is the SOH System. Theorem \ref{thm:hydro_main_theorem} extends the validity of this result to any dimension. This extension is far from obvious as it brings geometric intricacies that are absent in dimension~2. Indeed, in dimension 2, the PCC and KCC models can be formulated in standard flat spaces and the function $\zeta_w$ is obtained by solving a single equation instead of a system like \eqref{eq:macnorm_nD_vectGCI_wcomp}, \eqref{eq:macnorm_nD_vectGCI_kapTcomp}. Previously, the SOH model was derived as the hydrodynamic limit of the kinetic Vicsek model in \cite{degond2008continuum} (in the normalized case) and in \cite{degond2013macroscopic, degond2015phase} (in the non-normalized case). It also appears as the hydrodynamic limit of a Cucker-Smale model in \cite{aceves2019hydrodynamic}. All the previous results are formal but rigorous local convergence results were proved in \cite{jiang2024kinetic, jiang2016hydrodynamic}. Local existence of smooth solutions of the SOH model was proved in \cite{degond2011hydrodynamic} and in \cite{zhang2017local} in the viscous case. Its numerical resolution has been performed in \cite{motsch2011numerical}. 

A fluid model has been derived from the curvature-control model of \cite{cavagna2015flocking} in \cite{yang2015hydrodynamics}. This model is quite different from the SOH model. It is an extension of the Toner-Tu model in which $\mathbf{u}$ is not normalized but instead, a relaxation term of the Ginzburg-Landau type is addded at the right-hand side of the momentum equation \eqref{eq:macnorm_eq_evol_u_2nd}. In addition, the model of \cite{yang2015hydrodynamics} involves another variable, the locally-averaged curvature, and thus takes the form of three coupled PDEs. The model of \cite{yang2015hydrodynamics} is derived through a moment expansion of the kinetic equation and adequate moment closure assumptions. We do not think that such a model could be derived from the kinetic equation \eqref{eq:pk_eq_f_Q} by a hydrodynamic limit, even modifying the underpinning scaling assumptions. It seems however possible to derive a model involving an additional equation for the average curvature at the price of slightly modifying the initial PCC model \eqref{eq:pk_dXi}, \eqref{eq:pk_dVKi} and its associated kinetic model \eqref{eq:pk_eq_f_Q}. This will be investigated in future work. 

The remainder of this paper is devoted to the proof of Theorem \ref{thm:hydro_main_theorem}. The establishment of the continuity equation~\eqref{eq:macnorm_continuity} is easy and is done in the next subsection. All the forthcoming developments will be formal except otherwise stated.

\subsection{The continuity equation}
\label{subsec:macnorm_prelim}

In this section, we start the proof of Theorem \ref{thm:hydro_main_theorem} by determining the limit of $f^\varepsilon$ and showing the validity of the continuity equation \eqref{eq:macnorm_continuity}. 

\begin{lemma}[Equilibrium in the spatially nonhomogeneous case]~

\noindent
Suppose $f^\varepsilon \to f$ as $\varepsilon \to 0$ as smoothly as needed. Then, $f$ is given by \eqref{eq:macnorm_equi}. 
\label{lem:macnorm_equi}
\end{lemma}

\noindent
\textbf{Proof.} From \eqref{eq:pk_eq_f_Q}, we have $Q(f^\varepsilon) = {\mathcal O}(\varepsilon)$, so, when $\varepsilon \to 0$, we get $Q(f) = 0$. Thus, we can apply Prop. \ref{prop:equi_equi} and obtain that for any given $(x,t)$, the function $\alpha \mapsto f(x,\alpha,t)$ satisfies \eqref{eq:equi_equi}, which shows \eqref{eq:macnorm_equi}. \endproof

\begin{proposition}[Continuity equation]~

\noindent
The functions $\rho$ and $\mathbf{u}$ satisfy the continuity equation \eqref{eq:macnorm_continuity}. 
\label{prop:macnorm_continuity}
\end{proposition}

\noindent
\textbf{Proof.} Integrating \eqref{eq:equi_coll_oper_alt} with respect to $\alpha \in M$, we have 
$$ \int_M Q(f) \, |\mathrm{Vol}_M| = 0. $$
Hence, integrating \eqref{eq:pk_eq_f_Q} with respect to $\alpha$ leads to 
\begin{equation} 
\partial_t \Big( \int_M f^\varepsilon \, |\mathrm{Vol}_M| \Big) + \nabla_x \cdot \Big( \int_M f^\varepsilon \, v \, |\mathrm{Vol}_M| \Big) =0. 
\label{eq:macnorm_continuity_prf1}
\end{equation}
As $\varepsilon \to 0$, we have 
\begin{eqnarray*}
\int_M f^\varepsilon \, |\mathrm{Vol}_M| & \longrightarrow & \rho \int_M \mathbf{M}_{\mathbf{u}} \, |\mathrm{Vol}_M|(\alpha) = \rho , \\
\int_M f^\varepsilon \, v \, |\mathrm{Vol}_M| & \longrightarrow & \rho \int_M \mathbf{M}_{\mathbf{u}} \, v \, |\mathrm{Vol}_M| = c_1 \rho \mathbf{u}, 
\end{eqnarray*}
by \eqref{eq:equi_von_Mises_normal} and \eqref{eq:equi_von_mises_flux} respectively. Hence, letting $\varepsilon \to 0$ in \eqref{eq:macnorm_continuity_prf1} leads to \eqref{eq:macnorm_continuity}. \endproof

\subsection{The velocity equation}
\label{subsec:macnorm_velocity_eq}

The establishment of the velocity equation \eqref{eq:macnorm_eq_evol_u} is far less straightforward. The main difficulty is the lack of non-trivial collision invariants, i.e. functions $\zeta$, $M \to {\mathbb R}$ other than the constants, such that 
$$ \int_M Q(f) \, \zeta \, |\mathrm{Vol}_M| = 0, \quad \forall f: \, M \to {\mathbb R}, $$ 
(note that the constant collision invariants give rise to the continuity equation established in the previous subsection). 

In rarefied gas dynamics \cite{cercignani2013mathematical}, non-trivial collision invariants are the microscopic momentum and energy. Integrating the Boltzmann equation against them gives rise to the momentum and energy conservation equations in the macroscopic model. Here, microscopic momentum or energy is not preserved, and there is no substitute that would allow us to easily obtain an evolution equation for $\mathbf{u}$. This problem was first encountered in the case of the  the Vicsek model in \cite{degond2008continuum}, where it was resolved by introducing the concept of generalized collision invariant (GCI), which is a weakened form of collision invariant. It was applied to the 2D curvature control model in~\cite{degond2011macroscopic}. However, the passage from 2D to arbitrary dimension requires crossing a conceptual gap, due to the fact that the GCI in the 2D curvature control model is a scalar, while it is a vector in higher dimensions. The components of the GCI in higher dimensions do not simply follow from a simple scaling up formula from the scalar GCI in the 2D model. 

Below, we split the derivation of the velocity equation in four sections. First, in Section \ref{sec:gci}, we introduce the concept of GCI (and more precisely, of vector GCI). We prove its existence and uniqueness thanks to the transcription of the vector GCI equation into a variational formulation in appropriate spaces. We show that this variational formulation falls in the framework developed in \cite{albritton2019variational}, from which existence and uniqueness are easily deduced. 

However, a more detailed form of the vector GCI is necessary to derive an explicit form of the velocity equation. A first step towards this goal is to exploit the action of the orthogonal group $\mathrm{O}_n$ on the vector GCI. Section \ref{sec:gract_gci} is devoted to showing that the vector GCI is invariant under the action of the isotropy subgroup of $\mathrm{O}_n$ which fixes $\mathbf{u}$ (which is classically isomorphic to $\mathrm{O}_{n-1}$). This suggests that the vector GCI is the solution of a variational formulation posed on spaces of functions that are invariant under $\mathrm{O}_{n-1}$. In this section, we derive such a variational formulation by averaging the previously found variational formulation over the Haar measure of~$\mathrm{O}_{n-1}$, and show the existence and uniqueness of its solution.

We further exploit the invariance of the vector GCI under the action of $\mathrm{O}_{n-1}$ in Section \ref{sec:gci_redform}, by showing that it is necessarily expressed in terms of two scalar functions called the reduced GCI pair, which are nothing but the pair $(\zeta_w, \zeta_r)$ introduced in Theorem \ref{thm:hydro_main_theorem}. We then derive a variational formulation for the reduced GCI pair from the variational formulation for the vector GCI posed on spaces of $\mathrm{O}_{n-1}$-invariant functions. Existence and uniqueness of the former is guaranteed by those of the latter. The variational formulation for the GCI pair can be interpreted formally as the PDE system \eqref{eq:macnorm_nD_vectGCI_wcomp}, \eqref{eq:macnorm_nD_vectGCI_kapTcomp}. 

The final step, performed in Section \ref{sec:macnorm_eq_j}, is to derive the velocity equation \eqref{eq:macnorm_eq_evol_u} from the knowledge of the GCI pair and to show the expression \eqref{eq:macnorm_nD_coef_c2} for $c_2$. The computations involved in this step heavily rely on the previous sections.

\setcounter{equation}{0}
\section{Generalized collision invariants (GCI)}
\label{sec:gci}

In this section, we introduce the concept of generalized collision invariant (GCI) and show the existence and uniqueness (up to constants) of the GCI.

\subsection{GCI: definition}
\label{subsec:macnorm_GCI}

We introduce a linear operator constructed from the collision operator by freezing the $\mathbf{u}$ variable: let $\mathbf{u} \in {\mathbb S}^{n-1}$ be given and fixed. For a smooth function $f$: $M \to {\mathbb R}$, we define 
\begin{equation}
{\mathcal Q}_{\mathbf{u}}(f) = - \big\{f , \mathbf{H}_{\mathbf{u}} \big\} + D \nabla_\kappa \cdot \bigg( \mathbf{M}_{\mathbf{u}} \nabla_\kappa \Big( \frac{f}{\mathbf{M}_{\mathbf{u}}} \Big) \bigg). 
\label{eq:macnorm_coll_oper_frozen}
\end{equation}
By \eqref{eq:equi_coll_oper_alt}, for any smooth function $f$: $M \to {\mathbb R}$, we have
\begin{equation}
Q(f) = {\mathcal Q}_{\mathbf{u}_f}(f). 
\label{eq:macnorm_coll_vs_coll_frozen}
\end{equation}
Since the bracket is a derivation, the formal $L^2$-adjoint of ${\mathcal Q}_{\mathbf{u}}$ is ${\mathcal Q}_{\mathbf{u}}^*$ given on any smooth function $\zeta$: $M \to {\mathbb R}$  by 
\begin{equation}
{\mathcal Q}_{\mathbf{u}}^*(\zeta) = \big\{\zeta , \mathbf{H}_{\mathbf{u}} \big\} + \frac{D}{\mathbf{M}_{\mathbf{u}}} \nabla_\kappa \cdot \big( \mathbf{M}_{\mathbf{u}} \nabla_\kappa \zeta \big). 
\label{eq:macnorm_coll_frozen_adjoint}
\end{equation}
Now, we define a GCI associated with $\mathbf{u}$ as follows: 

\begin{definition}[GCI]~

\noindent
Let $\mathbf{u} \in {\mathbb S}^{n-1}$. A GCI associated with $\mathbf{u}$ is a function $\zeta$: $M \to {\mathbb R}$ such that 
\begin{equation}
\int_M {\mathcal Q}_{\mathbf{u}}(f) \, \zeta \, |\mathrm{Vol}_M| = 0, \quad \forall f \, \, \textrm{ such that } \, \, P_{\mathbf{u}^\bot} \Big( \int_M f \, v \, |\mathrm{Vol}_M| \Big) = 0. 
\label{eq:macnorm_GCI_def}
\end{equation}
\label{def:macnorm_GCI_def}
\end{definition}

\begin{proposition}[Equation solved by the GCI]~

\noindent
Let $\mathbf{u} \in {\mathbb S}^{n-1}$. A smooth function $\zeta$: $M \to {\mathbb R}$ is a GCI  associated with $\mathbf{u}$ if and only if there exists a vector $A \in {\mathbb R}^n$ such that $A \cdot \mathbf{u} = 0$ and 
\begin{equation}
{\mathcal Q}_{\mathbf{u}}^*(\zeta)(\alpha) = v \cdot A, \quad \forall \alpha=(v,\kappa) \in M, 
\label{eq:macnorm_GCI_char}
\end{equation}
\label{prop:macnorm_GCI_char}
\end{proposition}

\noindent
\textbf{Proof.} In the space $L^2(M,{\mathbb R})$, consider the finite-dimensional linear subspace ${\mathcal F}_{\mathbf{u}}$ consisting of all functions of the form $\alpha \mapsto (v \cdot A)$ (with $\alpha = (v,\kappa)$), where $A \in {\mathbb R}^n$ is any vector such that $A \cdot \mathbf{u} = 0$. Eq. \eqref{eq:macnorm_GCI_def} can be equivalently written
$$ f \in {\mathcal F}_{\mathbf{u}}^\bot \, \, \Longrightarrow f \in \mathrm{Span} \{{\mathcal Q}_{\mathbf{u}}^*(\zeta) \}^\bot. $$
Since both ${\mathcal F}_{\mathbf{u}}$ and $\mathrm{Span} \{{\mathcal Q}_{\mathbf{u}}^*(\zeta) \}$ are finite dimensional, hence closed, this statement is again equivalent to 
$$ \mathrm{Span} \{ {\mathcal Q}_{\mathbf{u}}^*(\zeta) \} \subset {\mathcal F}_{\mathbf{u}}, $$ 
which is exactly what \eqref{eq:macnorm_GCI_char} means. \endproof
 
Now, we define a ``vector GCI'' as follows. Let $\mathbf{u} \in {\mathbb S}^{n-1}$ be given and let $(e_1, \ldots, e_{n-1})$ be an orthonormal basis of $\{ \mathbf{u} \}^\bot$. For each $k \in \{1, \ldots, n-1\}$, denote by $\zeta_k$ a GCI associated with $\mathbf{u}$ and the choice $A = e_k$ in Prop. \ref{prop:macnorm_GCI_char}. Define a vector GCI $\vec \zeta$ associated with $\mathbf{u}$ by
\begin{equation} 
\vec \zeta = \sum_{k=1}^{n-1} \zeta_k e_k. 
\label{eq:macnorm_vec_GCI_def}
\end{equation} 
Applying Prop. \ref{prop:macnorm_GCI_char}, we find that $\vec \zeta$ satisfies the vector equation 
\begin{equation}
{\mathcal Q}_{\mathbf{u}}^*(\vec \zeta)(\alpha) = P_{\mathbf{u}^\bot} v, \quad \forall \alpha=(v,\kappa) \in M, 
\label{eq:macnorm_vec_GCI_eq}
\end{equation}
componentwise. A vector GCI has the following important property, which justifies the introduction of the GCI concept: 

\begin{lemma}[A vector GCI cancels the collision operator]~

\noindent
For any smooth enough function $f$: $M \to {\mathbb R}$, let $\vec \zeta$ be a vector GCI associated with $\mathbf{u}_f$. Then, we have 
\begin{equation}
\int_M Q(f) \, \vec \zeta \, |\mathrm{Vol}_M| = 0. 
\label{eq:macnorm_vecGCI_cancel}
\end{equation}
 
\label{lem:macnorm_vecGCI_cancel}
\end{lemma}

\noindent
\textbf{Proof.}  $f$ satisfies 
$$ P_{\mathbf{u}_f^\bot} \Big( \int_M f \, v \, |\mathrm{Vol}_M| \Big) = \Big| \int_M f \, v \, |\mathrm{Vol}_M| \Big| \, P_{\mathbf{u}_f^\bot} \mathbf{u}_f = 0. $$
Since $\vec \zeta$ is a vector GCI associated with $\mathbf{u}_f$, according to Def. \ref{def:macnorm_GCI_def}, we deduce that 
$$ \int_M {\mathcal Q}_{\mathbf{u}_f} (f) \, \vec \zeta \, |\mathrm{Vol}_M| = 0. $$
With \eqref{eq:macnorm_coll_vs_coll_frozen}, we get \eqref{eq:macnorm_vecGCI_cancel}. 
\endproof

\subsection{Existence and uniqueness (up to constants) of the vector GCI}
\label{subsec:macnorm_GCI_exist}

We now prove the existence of solutions of the stationary linear Fokker-Planck equation \eqref{eq:macnorm_vec_GCI_eq}. In \cite{degond1986global, degond1987existence}, a variational framework for the linear time-dependent Fokker-Planck equation was proposed. It relied on \cite{baouendi1968equation} which made essential use of a result of Lions \cite{lions2013equations}. This method was used again in \cite{degond2011macroscopic}, this time in the case of the stationary Fokker-Planck equation, but there is a mistake in the proof. Recently, \cite{albritton2019variational} revisited this variational framework and showed a Poincar\'e estimate and the existence of solutions for both the stationary and time-dependent Fokker-Planck equations. In this section, we adapt the results of \cite{albritton2019variational} to the present geometric setting. Since the proofs involve only minor adaptations of \cite{albritton2019variational}, they are omitted. We first introduce the following functional spaces.

\subsubsection{Spaces of functions and distributions on $M$}
\label{subsubsec:macnorm_GCI_subspaces}

The following list of spaces is tedious but is needed for the sake of notations. We begin with definitions of spaces of scalar functions on $M$.

\medskip
\noindent
(i) The space $L^2_v(L^2_\kappa)$ is defined as the closure of the space of smooth functions $M \to {\mathbb R}$ for the norm 
\begin{equation} 
\| f \|_{L^2_v(L^2_\kappa)}^2 = \int_M |f|^2 \, \mathbf{M}_{\mathbf{u}} \, |\mathrm{Vol}_M|. 
\label{eq:macnorm_L2(L2)_norm}
\end{equation}
We denote by $( \hspace{-0.8mm} ( f , g) \hspace{-0.8mm} )_{L^2_v(L^2_\kappa)}$ the associated inner-product. We also define the closed hyperplane 
$$ L^2_v(L^2_\kappa)_0 = \Big\{ f \in L^2_v(L^2_\kappa) \quad \Big| \quad \int_M f  \, \mathbf{M}_{\mathbf{u}} \, |\mathrm{Vol}_M| = 0 \Big\}. $$
Both $L^2_v(L^2_\kappa)$ and $L^2_v(L^2_\kappa)_0$ are Hilbert spaces. 

\medskip
\noindent 
(ii) The space $\vec{L}^2_v(L^2_\kappa)$ is the space of sections $k$ of the pullback  bundle $\pi^* TB$ such that the map $M \to {\mathbb R}$, $(v,\kappa) \to \langle k , k \rangle_v$ belongs to $L^2_v(L^2_\kappa)$. It is a Hilbert space for the inner product 
\begin{equation} 
\big( \hspace{-1.2mm} \big(  k , h \big) \hspace{-1.2mm} \big)_{\vec{L}^2_v(L^2_\kappa)} = \int_M \langle k,h \rangle \, \mathbf{M}_{\mathbf{u}} \, |\mathrm{Vol}_M|. 
\label{eq:macnorm_vecL2(L2)_innerprod}
\end{equation}
and we denote by $\| k \|_{\vec{L}^2_v(L^2_\kappa)}^2$ the associated norm. 

\medskip
\noindent 
(iii) The space $L^2_v(H^1_\kappa)$ is the closure of the space of smooth functions $M \to {\mathbb R}$ for the norm 
\begin{equation} 
\| f \|_{L^2_v(H^1_\kappa)}^2 =  \| f \|_{L^2_v(L^2_\kappa)}^2 + \| \nabla_\kappa f \|_{\vec{L}^2_v(L^2_\kappa)}^2, 
\label{eq:macnorm_L2(H1)_norm}
\end{equation}
and we denote by $( \hspace{-0.8mm} ( f , g ) \hspace{-0.8mm} )_{L^2_v(H^1_\kappa)}$ the associated inner-product. It is a Hilbert space. 

\medskip
\noindent
(iv) The space ${\mathcal D}(M)$ is the space of smooth (infinitely differentiable) functions $M \to {\mathbb R}$ with compact support. Its dual ${\mathcal D}'(M)$ is the space of distributions on $M$. For $T \in {\mathcal D}'(M)$ and $\varphi \in {\mathcal D}(M)$, the duality bracket $ ( \hspace{-0.8mm} ( T , \varphi ) \hspace{-0.8mm} )_{{\mathcal D}',{\mathcal D}}$ is defined to extend the inner product of $L^2_v(L^2_\kappa)$. The weight $\mathbf{M}_{\mathbf{u}}$ in the duality bracket must be taken into account when defining the distributional derivative. Hence, for $T \in {\mathcal D}'$, we define $\nabla_v T$ and $\nabla_\kappa T$ by duality against an element $\varphi \in {\mathcal Y}(M)$ with compact support as follows
$$ \big( \hspace{-1.2mm} \big( \nabla_v T, \varphi \big) \hspace{-1.2mm} \big)_{{\mathcal D}',{\mathcal D}} = - 
\big( \hspace{-1.2mm} \big( T, \frac{1}{\mathbf{M}_{\mathbf{u}}} \nabla_v \cdot ( \mathbf{M}_{\mathbf{u}} \varphi ) \big) \hspace{-1.2mm} \big)_{{\mathcal D}',{\mathcal D}}, $$
and similarly replacing $v$ by $\kappa$. We note that 
$$ \frac{1}{\mathbf{M}_{\mathbf{u}}} \nabla_v \cdot ( \mathbf{M}_{\mathbf{u}} \varphi ) = \nabla_v \cdot \varphi + \frac{\nu^2}{D} P_{v^\bot} \mathbf{u} \cdot \varphi , \qquad \frac{1}{\mathbf{M}_{\mathbf{u}}} \nabla_\kappa \cdot ( \mathbf{M}_{\mathbf{u}} \varphi ) = \nabla_\kappa \cdot \varphi - \frac{\nu}{D} \kappa \cdot \varphi, $$ are smooth so that the above definitions make sense. Thanks to this definition, gradients in the usual and distributional senses coincide when applied to smooth functions. Since the bracket is a derivation, and $\mathbf{M}_\mathbf{u}$ is a function of $\mathbf{H}_\mathbf{u}$, we have 
\begin{equation} 
\big( \hspace{-1.2mm} \big( \{ T, \mathbf{H}_{\mathbf{u}} \}, \varphi \big) \hspace{-1.2mm} \big)_{{\mathcal D}',{\mathcal D}} = - 
\big( \hspace{-1.2mm} \big( T, \{ \varphi, \mathbf{H}_{\mathbf{u}} \} \big) \hspace{-1.2mm} \big)_{{\mathcal D}',{\mathcal D}}. 
\label{eq:macnorm_bracket_dual}
\end{equation}

\medskip
\noindent
(v) The dual space $L^2_v(H^{-1}_\kappa)$ of $L^2_v(H^1_\kappa)$ is a subspace of ${\mathcal D}'(M)$, namely, the subspace of continuous linear forms on ${\mathcal D}(M)$ that extend to elements of $L^2_v(H^1_\kappa)$. We denote by $ ( \hspace{-0.8mm} ( T , \varphi ) \hspace{-0.8mm} )_{L^2_v(H^{-1}_\kappa),L^2_v(H^1_\kappa)}$ the duality between $T \in L^2_v(H^{-1}_\kappa)$ and $\varphi \in L^2_v(H^1_\kappa)$. The space $L^2_v(H^{-1}_\kappa)$ is normed by the usual operator norm of linear forms, namely 
$$ \| T \|_{L^2_v(H^{-1}_\kappa)} = \sup_{\varphi \in L^2_v(H^1_\kappa), \, \| \varphi \|_{L^2_v(H^1_\kappa)} \leq 1}
 \big| ( \hspace{-0.8mm} ( T , \varphi ) \hspace{-0.8mm} )_{L^2_v(H^{-1}_\kappa),L^2_v(H^1_\kappa)} \big| . $$
$L^2_v(H^{-1}_\kappa)$ is a Hilbert space.

\medskip
Now, we extend these definitions to spaces of vector-valued functions and distributions. Let $(e_1, \ldots, e_{n-1})$ be an orthonormal basis of $\{ \mathbf{u} \}^\bot$. 

\medskip
\noindent 
(i) The space $L^2_v(L^2_\kappa)_\mathbf{u}$ is the set of vector valued functions $\vec{\psi}$: $M \to \{\mathbf{u}\}^\bot$ such that $\vec{\psi}_i = \vec{\psi} \cdot e_i \in L^2_v(L^2_\kappa)$, $\forall i \in \{1, \ldots, n-1 \}$, endowed with the norm
\begin{equation} 
\| \vec{\psi} \|_{L^2_v(L^2_\kappa)_\mathbf{u}}^2 = \sum_{i=1}^{n-1} \| \vec{\psi}_i\|_{L^2_v(L^2_\kappa)}^2,  
\label{eq:macnorm_nD_L2norm_vect_def}
\end{equation}
and associated inner product 
$$ ( \hspace{-0.8mm} ( \vec{\psi}, \vec{\varphi} ) \hspace{-0.8mm} )_{L^2_v(L^2_\kappa)_\mathbf{u}} = \sum_{i=1}^{n-1} ( \hspace{-0.8mm} ( \vec{\psi}_i, \vec{\varphi}_i ) \hspace{-0.8mm} )_{L^2_v(L^2_\kappa)} = \int_M \vec{\psi} \cdot \vec{\varphi}  \, \mathbf{M}_{\mathbf{u}} \, |\mathrm{Vol}_M|, $$
where $\vec{\psi} \cdot \vec{\varphi} = \sum_{i=1}^{n-1} \vec{\psi}_i \,  \vec{\varphi}_i$ is the usual inner product in $\{\mathbf{u}\}^\bot$. Both the norm and inner-product are independent of the choice of the orthonormal basis of $\{ u \}^\bot$. We also define the closed hyperplane 
$$ L^2_v(L^2_\kappa)_{\mathbf{u} 0} = \Big\{ \vec{\psi} \in L^2_v(L^2_\kappa)_\mathbf{u} \quad \Big| \quad \vec{\psi}_i \in L^2_v(L^2_\kappa)_0, \, \, \forall i \in \{1, \ldots, n-1 \} \Big\}. $$
Both $L^2_v(L^2_\kappa)_\mathbf{u}$ and $L^2_v(L^2_\kappa)_{\mathbf{u} 0}$ are Hilbert spaces. 

\medskip
\noindent 
(ii) Let $\vec{k}$ be a section of the bundle $\pi^*(TB) \otimes (M \times \{ \mathbf{u} \}^\bot)$. We can write $ \vec{k} = \sum_{i=1}^{n-1} \vec{k}_i \otimes e_i$, with $\vec{k}_i$ being a section of the pullback bundle $\pi^*(TB)$, i.e. $\vec{k}_i(\alpha) \in T_{\pi \alpha} B$, $\forall \alpha \in M$. The space $\vec{L}^2_v(L^2_\kappa)_\mathbf{u}$ is the space of such sections such that $\vec{k}_i \in \vec{L}^2_v(L^2_\kappa)$, $\forall i \in \{1, \ldots, n-1 \}$. This space is a Hilbert space when endowed with norm 
\begin{equation} 
\| \vec{k} \|_{\vec{L}^2_v(L^2_\kappa)_\mathbf{u}}^2 = \sum_{i=1}^{n-1} \| \vec{k}_i\|_{\vec{L}^2_v(L^2_\kappa)}^2,
\label{eq:macnorm_nD_L2norm_vect_vect_def}
\end{equation}
and associated inner product $( \hspace{-0.8mm} ( \vec{k}, \vec{h} ) \hspace{-0.8mm} )_{\vec{L}^2_v(L^2_\kappa)_\mathbf{u}}$. Again, both the norm and inner-product are independent of the choice of the orthonormal basis of $\{ u \}^\bot$. 

\medskip
\noindent 
(iii) Let $\vec{\psi}$: $M \to \{ \mathbf{u} \}^\bot$. Then, $\nabla_\kappa \vec{\psi}$, when it exists, is a section of $\pi^*(TB) \otimes (M \times \{ \mathbf{u} \}^\bot)$ and we can write $\nabla_\kappa \vec{\psi} = \sum_{i=1}^{n-1} \nabla_\kappa \vec{\psi}_i \otimes e_i$. The space $L^2_v(H^1_\kappa)_\mathbf{u}$ is given by 
$$ L^2_v(H^1_\kappa)_\mathbf{u} = \big\{ \vec{\psi} \in L^2_v(L^2_\kappa)_\mathbf{u} \, \, \big| \, \, \nabla_\kappa \vec{\psi} \in \vec{L}^2_v(L^2_\kappa)_\mathbf{u} \big\}. $$
It is a Hilbert space for the norm
$$ \| \vec{\psi} \|_{L^2_v(H^1_\kappa)_\mathbf{u}}^2 =  \| \vec{\psi} \|_{L^2_v(L^2_\kappa)_\mathbf{u}}^2 + \| \nabla_\kappa \vec{\psi} \|_{\vec{L}^2_v(L^2_\kappa)_\mathbf{u}}^2 =  \sum_{i=1}^{n-1} \| \vec{\psi}_i\|_{L^2_v(H^1_\kappa)}^2, $$
and associated inner-product $( \hspace{-0.8mm} ( \vec{\psi}, \vec{\varphi} ) \hspace{-0.8mm} )_{L^2_v(H^1_\kappa)_\mathbf{u}}$, which again, do not depend on the choice of the orthonormal basis of $\{ \mathbf{u} \}^\bot$.

\medskip
\noindent
(iv) The space ${\mathcal D}_\mathbf{u}(M)$ is the space of smooth (infinitely differentiable) functions $M \to \{ \mathbf{u} \}^\bot$ with compact support. Its dual ${\mathcal D}_\mathbf{u}'(M)$ is the space of distributions on $M$ with values in $\{ \mathbf{u} \}^\bot$. Let $\vec{T} \in {\mathcal D}_\mathbf{u}'(M)$. We can write $\vec{T} = \sum_{i=1}^{n-1} \vec{T}_i \, e_i$, with $\vec{T}_i \in {\mathcal D}'(M)$, $\forall i \in \{ 1, \ldots, n-1 \}$. For $\vec{T} \in {\mathcal D}_\mathbf{u}'(M)$ and $\vec{\psi} = \sum_{i=1}^{n-1} \vec{\psi}_i \, e_i \in {\mathcal D}_\mathbf{u}(M)$, the duality bracket $ ( \hspace{-0.8mm} ( \vec{T} , \vec{\psi} ) \hspace{-0.8mm} )_{{\mathcal D}_\mathbf{u}',{\mathcal D}_\mathbf{u}}$ is defined to extend the inner product of $L^2_v(L^2_\kappa)_\mathbf{u}$. We have 
\begin{equation} 
\big( \hspace{-1.2mm} \big( \vec{T} , \vec{\psi} \big) \hspace{-1.2mm} \big)_{{\mathcal D}_\mathbf{u}',{\mathcal D}_\mathbf{u}} = 
\sum_{i=1}^{n-1} \big( \hspace{-1.2mm} \big( \vec{T}_i , \vec{\psi}_i \big) \hspace{-1.2mm} \big)_{{\mathcal D}',{\mathcal D}}, 
\label{eq:macnorm_nD_duality_vector_distrib}
\end{equation}

\medskip
\noindent
(v) The dual space $L^2_v(H^{-1}_\kappa)_\mathbf{u}$ of $L^2_v(H^1_\kappa)_\mathbf{u}$ is a subspace of ${\mathcal D}_\mathbf{u}'(M)$, namely, the subspace of continuous linear forms on ${\mathcal D}_\mathbf{u}(M)$ that extend to elements of $L^2_v(H^1_\kappa)_\mathbf{u}$. We have 
$$ \vec{T} = \sum_{i=1}^{n-1} \vec{T}_i \, e_i \in L^2_v(H^{-1}_\kappa)_\mathbf{u} \quad \Longleftrightarrow \quad \vec{T}_i \in L^2_v(H^{-1}_\kappa), \quad 
\forall i \in \{ 1, \ldots, n-1 \}. $$
We denote by $ ( \hspace{-0.8mm} ( \vec{T} , \vec{\psi} ) \hspace{-0.8mm} )_{L^2_v(H^{-1}_\kappa)_\mathbf{u},L^2_v(H^1_\kappa)_\mathbf{u}}$ the duality between $\vec{T} \in L^2_v(H^{-1}_\kappa)_\mathbf{u}$ and $\vec{\psi} \in L^2_v(H^1_\kappa)_\mathbf{u}$. 
The identity \eqref{eq:macnorm_nD_duality_vector_distrib} is still valid for the duality between $L^2_v(H^{-1}_\kappa)_\mathbf{u}$ and $L^2_v(H^1_\kappa)_\mathbf{u}$ at the left-hand side and between $L^2_v(H^{-1}_\kappa)$ and $L^2_v(H^1_\kappa)$ at the right-hand side. The space $L^2_v(H^{-1}_\kappa)_\mathbf{u}$ is a Hilbert space when normed by the usual operator norm of linear forms, namely 
$$ \| \vec{T} \|_{L^2_v(H^{-1}_\kappa)_\mathbf{u}} = \sup_{\vec{\psi} \in L^2_v(H^1_\kappa)_\mathbf{u}, \, \| \vec{\psi} \|_{L^2_v(H^1_\kappa)_\mathbf{u}} \leq 1}
 \big| ( \hspace{-0.8mm} (  \vec{T} , \vec{\psi} ) \hspace{-0.8mm} )_{L^2_v(H^{-1}_\kappa)_\mathbf{u},L^2_v(H^1_\kappa)_\mathbf{u}} \big| . $$

\medskip
\noindent
(vi) For any smooth function $\mathbf{H}$: $M \to {\mathbb R}$ and any $\vec{T} = \sum_{i=1}^{n-1} \vec{T}_i \, e_i \in {\mathcal D}'_\mathbf{u}(M)$, we define $ \{ \vec{T}, \mathbf{H} \} = \sum_{i=1}^{n-1} \{ \vec{T}_i, \mathbf{H} \} \, e_i$. Like in the scalar case, for any $\vec{\psi} \in {\mathcal D}_\mathbf{u}(M)$, we have 
\begin{equation} 
\big( \hspace{-1.2mm} \big( \{ \vec{T}, \mathbf{H}_{\mathbf{u}} \}, \vec{\psi} \big) \hspace{-1.2mm} \big)_{{\mathcal D}_\mathbf{u}',{\mathcal D}_\mathbf{u}} = - 
\big( \hspace{-1.2mm} \big( \vec{T}, \{ \vec{\psi}, \mathbf{H}_{\mathbf{u}} \} \big) \hspace{-1.2mm} \big)_{{\mathcal D}_\mathbf{u}',{\mathcal D}_\mathbf{u}}. 
\label{eq:macnorm_bracket_dual_u}
\end{equation}
The space ${\mathbb Y}_\mathbf{u}$ is defined by 
$${\mathbb Y}_\mathbf{u} = \Big\{ \vec{\psi} \in L^2_v(H^1_\kappa)_\mathbf{u}, \quad \{ \vec{\psi}, \mathbf{H}_{\mathbf{u}} \} \in L^2_v(H^{-1}_\kappa)_\mathbf{u} \Big\}, $$
where $\{ \vec{\psi}, \mathbf{H}_{\mathbf{u}} \}$ is defined in the distributional sense by \eqref{eq:macnorm_bracket_dual_u}. The space
${\mathbb Y}_\mathbf{u}$ is normed by 
$$ \| \vec{\psi} \|_{{\mathbb Y}_\mathbf{u}}^2 = \| \vec{\psi} \|_{L^2_v(H^1_\kappa)_\mathbf{u}}^2 + \| \{ \vec{\psi}, \mathbf{H}_{\mathbf{u}} \} \|_{L^2_v(H^{-1}_\kappa)_\mathbf{u}}^2, $$
and is a Hilbert space. 

\bigskip
The following lemmas provide some properties of the space ${\mathbb Y}_\mathbf{u}$. Their proofs are fairly straightforward adaptations of corresponding results proved in \cite{albritton2019variational}, and are omitted. 

\begin{lemma}[Density of ${\mathcal D}_\mathbf{u}(M)$ in ${\mathbb Y}_\mathbf{u}$]~

\noindent
${\mathcal D}_\mathbf{u}(M)$ is dense in  ${\mathbb Y}_\mathbf{u}$. 
\label{lem:Lions_D_dense_in_Y}
\end{lemma}

\begin{lemma}[Poincar\'e inequality in ${\mathbb Y}_\mathbf{u}$]~

\noindent
There exists a constant $C>0$ such that for any $\vec{\psi} \in {\mathbb Y}_\mathbf{u} \cap L^2_v(L^2_\kappa)_{\mathbf{u} 0}$ we have 
\begin{equation}
\| \vec{\psi} \|_{L^2_v(L^2_\kappa)_\mathbf{u}} \leq C \big( \| \nabla_\kappa \vec{\psi} \|_{\vec{L}^2_v(L^2_\kappa)_\mathbf{u}} + \| \{ \vec{\psi}, \mathbf{H}_{\mathbf{u}} \} \|_{L^2_v(H^{-1}_\kappa)_\mathbf{u}} \big). 
\label{eq:Lions_poincare_Yu}
\end{equation}
\label{lem:Lions_poincare_Y}
\end{lemma}

\subsubsection{Proof of existence of the vector GCI}
\label{subsubsec:macnorm_GCI_exist}

Let $\vec{\chi} \in L^2_v(L^2_\kappa)_\mathbf{u}$. We are going to prove the existence of solutions to the following equation:
\begin{equation}
{\mathcal Q}_{\mathbf{u}}^*(\vec{\psi}) = \vec{\chi}, 
\label{eq:macnorm_GCI_eq_gene}
\end{equation}
provided that $\vec{\chi} \in L^2_v(L^2_\kappa)_{\mathbf{u} 0}$, and we will characterize the space of such solutions. For that purpose, we transform \eqref{eq:macnorm_GCI_eq_gene} into a variational formulation. Let $\vec{\varphi} \in {\mathcal D}_\mathbf{u}(M)$. Taking the inner-product of~\eqref{eq:macnorm_GCI_eq_gene} (in the space $\{ \mathbf{u} \}^\bot$) with ~$\mathbf{M}_{\mathbf{u}} \, \vec{\varphi}$, integrating over $M$ and using \eqref{eq:geom_graddiv_dual} supposing that $\vec{\psi}$ has all the requested regularity for this to be possible, we get 
$$ - \int_M \{ \vec{\psi}, \mathbf{H}_{\mathbf{u}} \} \cdot \vec{\varphi} \, \mathbf{M}_{\mathbf{u}} \, |\mathrm{Vol}_M| + D \int_M \sum_{i=1}^{n-1} \big \langle \nabla_\kappa {\vec \varphi}_i, \nabla_\kappa \vec{\psi}_i  \rangle \, \mathbf{M}_{\mathbf{u}} \, |\mathrm{Vol}_M| = - \int_M \vec{\chi} \cdot \vec{\varphi} \, \mathbf{M}_{\mathbf{u}} \, |\mathrm{Vol}_M|. $$
Now, we see that if both $\vec{\varphi}$ and $\vec{\psi}$ belong to ${\mathbb Y}_\mathbf{u}$, this formula has a meaning by interpreting the first integral as a duality bracket between $L^2_v(H^{-1}_\kappa)_\mathbf{u}$ and $L^2_v(H^1_\kappa)_\mathbf{u}$. This leads to the variational formulation 
\begin{equation}
\mbox{} \hspace{-0.5cm} \left\{ \begin{array}{l}
\vec{\psi} \in {\mathbb Y}_\mathbf{u}, \\
- \big( \hspace{-1.2mm} \big( \{ \vec{\psi}, \mathbf{H}_\mathbf{u} \} , \vec{\varphi} \big) \hspace{-1.2mm} \big)_{L^2_v(H^{-1}_\kappa)_\mathbf{u},L^2_v(H^1_\kappa)_\mathbf{u}} + D \big( \hspace{-1.2mm} \big( \nabla_\kappa \vec{\psi}, \nabla_\kappa \vec{\varphi} \big) \hspace{-1.2mm} \big)_{\vec{L}^2_v(L^2_\kappa)_\mathbf{u}} = - \big( \hspace{-1.2mm} \big( \vec{\chi} , \vec{\varphi} \big) \hspace{-1.2mm} \big)_{L^2_v(L^2_\kappa)_\mathbf{u}}, \quad \forall \vec{\varphi} \in {\mathbb Y}_\mathbf{u}, \end{array} \right. 
\label{eq:macnorm_GCI_gene_varform_u}
\end{equation}
By testing \eqref{eq:macnorm_GCI_gene_varform_u} with $\vec{\varphi} = e_i$ for $i=1, \ldots, n-1$ (which belong to ${\mathbb Y}_\mathbf{u}$) and noting that $\nabla_\kappa \vec{\varphi} = 0$ and $\{ \vec{\varphi}, \mathbf{H}_{\mathbf{u}} \} = 0$, we see that the existence of a solution to \eqref{eq:macnorm_GCI_gene_varform_u} requires that $\vec{\chi} \in L^2_v(L^2_\kappa)_{\mathbf{u} 0}$. The following theorem, whose proof is a fairly straightforward adaptation of \cite{albritton2019variational} and is omitted, shows that, for $\vec{\chi} \in L^2_v(L^2_\kappa)_\mathbf{u}$, $\vec{\chi}$ being in $L^2_v(L^2_\kappa)_{\mathbf{u} 0}$ is a necessary and sufficient condition for the existence of a solution to~\eqref{eq:macnorm_GCI_gene_varform_u}. 

\begin{theorem}[Existence of a solution to \eqref{eq:macnorm_GCI_gene_varform_u}]~

\noindent
Let $\vec{\chi} \in L^2_v(L^2_\kappa)_\mathbf{u}$. Then, the variational formulation \eqref{eq:macnorm_GCI_gene_varform_u} has a solution if and only if $\vec{\chi} \in L^2_v(L^2_\kappa)_{\mathbf{u} 0}$ and in that case, the solution is unique in ${\mathbb Y}_\mathbf{u} \cap L^2_v(L^2_\kappa)_{\mathbf{u} 0}$. 
\label{thm:macnorm_GCI_gene_varform_exist}
\end{theorem}

\noindent
From Theorem \ref{thm:macnorm_GCI_gene_varform_exist}, we immediately deduce the following 

\begin{corollary}[Existence and space of GCI]~

\noindent
(i) There exists a unique variational solution (in the sense of \eqref{eq:macnorm_GCI_gene_varform_u} with $\chi$ given by $\vec{\chi} (v,\kappa) = P_{\mathbf{u}^\bot} v $, $\forall \alpha = (v,\kappa) \in M$) in $\vec{\mathbb Y}_\mathbf{u} \cap L^2_v(L^2_\kappa)_{\mathbf{u} 0}$  to Equation \eqref{eq:macnorm_vec_GCI_eq}. This unique solution is denoted by~$\vec{\zeta}_\mathbf{u}$. The set of solutions to \eqref{eq:macnorm_vec_GCI_eq} in ${\mathbb Y}_\mathbf{u}$ is given by $\{ \vec{\zeta}_\mathbf{u} + \vec{C} \, | \, \vec{C} \in \{\mathbf{u}\}^\bot \}$. 
                                                                               
\medskip
\noindent
(ii) The space of (scalar) GCI associated with a given direction $\mathbf{u} \in {\mathbb S}^{n-1}$ is given by 
$$ {\mathcal G}_\mathbf{u} = \big\{ \vec{\zeta}_\mathbf{u} \cdot A + C \, \,  \big| \, \, A \in \{ \mathbf{u} \}^\bot  \, \,  \textrm{ and } \, \,  C \in {\mathbb R} \big\},$$ 
and is an $n$-dimensional vector space. 
\label{cor:macnorm_exist_GCI}
\end{corollary}

\noindent
\textbf{Proof.} (i) The function $M \to \{\mathbf{u}\}^\bot$, $(v,\kappa) \mapsto \vec{\chi} (v,\kappa) = P_{\mathbf{u}^\bot} v$ belongs to $L^2_v(L^2_\kappa)_{\mathbf{u} 0}$ because for each $i = 1, \ldots, n-1$, we have $P_{\mathbf{u}^\bot} v \cdot e_i  =v \cdot e_i$ and 
$$ \int_B (v \cdot e_i) \, \exp \Big( \frac{\nu}{D} v \cdot \mathbf{u} \Big) \, |\mathrm{Vol}_B| = 0, $$
by antisymmetry. Theorem \ref{thm:macnorm_GCI_gene_varform_exist} ensures the existence of $\vec{\zeta}_\mathbf{u}$, the unique solution in ${\mathbb Y}_\mathbf{u} \cap L^2_v(L^2_\kappa)_{\mathbf{u} 0}$ of \eqref{eq:macnorm_GCI_char} in the sense of the variational formulation \eqref{eq:macnorm_GCI_gene_varform_u}.

\medskip
\noindent
(ii) Follows immediately from Point (i) and Prop. \ref{prop:macnorm_GCI_char}. \endproof

\begin{remark}[Back to construction \eqref{eq:macnorm_vec_GCI_def}]~

\noindent
The uniqueness (up to constants) of the vector GCI shows that the construction \eqref{eq:macnorm_vec_GCI_def} does not depend on the choice of the orthonormal basis $\{ e_1, \ldots, e_{n-1} \}$ of $\{\mathbf{u}\}^\bot$. 
\end{remark}

\setcounter{equation}{0}
\section{Action of the orthogonal group on the vector GCI}
\label{sec:gract_gci}

In this section, we show that the orthogonal group $\mathrm{O}_n$ acts on the vector GCI $\vec{\zeta}_\mathbf{u}$ and that $\mathrm{O}_{n-1}$ (identified with the isotropy subgroup of $\mathrm{O}_n$ which fixes $\mathbf{u}$) fixes $\vec{\zeta}_\mathbf{u}$. Then, thanks to a suitably defined projection on subspaces of invariant functions under $\mathrm{O}_{n-1}$, we show that we can restrict the variational formulation satisfied by $\vec{\zeta}_\mathbf{u}$ to such spaces of invariant functions. We first review some geometric background about the action of $\mathrm{O}_n$ on the various objects involved in this study.

\subsection{Geometric preliminaries}
\label{subsec:macnorm_group_action}

In this section, we review how the orthogonal group $\mathrm{O}_n$ acts on $B$, $M$, the space ${\mathcal S}(M,{\mathbb R})$ of smooth functions $M \to {\mathbb R}$, the spaces ${\mathcal X}(M)$ and ${\mathcal Y}(M)$ of smooth vector fields on $M$ and of smooth sections of the pullback bundle $\pi^*(TB)$ respectively, and the various spaces of vector-valued functions and distributions defined in Section \ref{subsubsec:macnorm_GCI_subspaces}.

\subsubsection{Action of $\mathrm{O}_n$ on $B$ and $M$}
\label{subsubsec:macnorm_group_action_BM}

We recall that $\mathrm{O}_n$ acts on the left on ${\mathbb R}^n$ by the action $\mathrm{O}_n \times {\mathbb R}^n \to {\mathbb R}^n$, $(R,v) \mapsto Rv$ (meaning that the orthogonal matrix $R$ acts on the $n$-tuple of coordinates of $v$ in the canonical basis or ${\mathbb R}^n$). This action can be extended to a left action on ${\mathbb R}^{2n}$ by $(R,(v,\kappa)) \mapsto R(v,\kappa) = : (Rv,R\kappa)$ and to a left action on ${\mathbb R}^{4n}$ by $(R,(v,\kappa,a,\tau)) \mapsto R(v,\kappa,a,\tau) = : (Rv,R\kappa,Ra,R\tau)$. Now, since $|Rv| = |v|$, the action of $\mathrm{O}_n$ on ${\mathbb R}^n$ can be restricted to an action of $\mathrm{O}_n$ on $B$. Similarly, since $Rv \cdot R \kappa = v \cdot \kappa$, the action of $\mathrm{O}_n$ on ${\mathbb R}^{2n}$ can be restricted to an action of $\mathrm{O}_n$  on $M$. Finally, because $Ra \cdot Rv = a \cdot v$ and $Ra \cdot R \kappa + Rv \cdot R \tau = a \cdot \kappa + v \cdot \tau$, the action of $\mathrm{O}_n$ on ${\mathbb R}^{4n}$ can be restricted to an action of $\mathrm{O}_n$ on $N$, the support manifold of the tangent bundle $TM$. 

The action of $\mathrm{O}_n$ is compatible with the metrics on $B$: let $v \in B$ and $\kappa_1$, $\kappa_2 \in T_vB$. Then, 
\begin{equation} 
\big \langle R \kappa_1, R \kappa_2 \big \rangle_{Rv} = \big \langle \kappa_1, \kappa_2 \big \rangle_v, \qquad \forall R \in \mathrm{O}_n, 
\label{eq:macnorm_grac_metric_B}
\end{equation}
i.e. the action of $\mathrm{O}_n$ on $B$ leads to orthogonal transformations from $T_v B$ to $T_{Rv} B$. Similarly, it is compatible with the metric on $M$: let $\alpha \in M$ and $\xi_1$, $\xi_2 \in T_\alpha M$. Then, 
\begin{equation}  
\big\langle \hspace{-1.8mm} \big\langle R \xi_1, R \xi_2 \big\rangle \hspace{-1.8mm} \big\rangle_{R \alpha} = \big\langle \hspace{-1.8mm} \big\langle \xi_1, \xi_2 \big\rangle \hspace{-1.8mm} \big\rangle_\alpha, \qquad \forall R \in \mathrm{O}_n. 
\label{eq:macnorm_grac_metric_M}
\end{equation}
In other words, the action of $R$ on $M$ leads to orthogonal transformations from $T_\alpha M$ to $T_{R\alpha} M$. 

Now, the action of $\mathrm{O}_n$ commutes with the lifts and connection maps: let $\alpha \in M$ and $\xi \in T_\alpha M$. We have
\begin{equation}
{\mathcal B}_{R \alpha}^V (R \xi) = R {\mathcal B}_\alpha^V \xi, \qquad {\mathcal B}_{R \alpha}^H (R \xi) = R {\mathcal B}_\alpha^H \xi, \qquad \forall R \in \mathrm{O}_n, 
\label{eq:macnorm_gract_connmap}
\end{equation}
and likewise, let $\alpha \in M$, $\alpha = (v,\kappa)$ and let $a, \, \tau \in T_v B$. We have
\begin{equation}
{\mathcal L}_{R \alpha}^V (R \tau) = R {\mathcal L}_\alpha^V \tau, \qquad {\mathcal L}_{R \alpha}^H (R a) = R {\mathcal L}_\alpha^H a, \qquad \forall R \in \mathrm{O}_n. 
\label{eq:macnorm_gract_liftmap}
\end{equation}

\subsubsection{Action of $\mathrm{O}_n$ on functions and vector fields}
\label{subsec:macnorm_group_action_fct_vectfields}

The left action of $\mathrm{O}_n$ on $M$ induces a right action of $\mathrm{O}_n$ on ${\mathcal S}(M,{\mathbb R})$ by ${\mathcal S}(M,{\mathbb R}) \times \mathrm{O}_n \to {\mathcal S}(M,{\mathbb R})$, $(f,R) \to f^R$ where $f^R(\alpha) = f(R \alpha)$, $\forall \alpha \in M$. Similarly, there is a right action of $\mathrm{O}_n$ on ${\mathcal X} (M)$ by ${\mathcal X} (M) \times \mathrm{O}_n \to {\mathcal X} (M)$, $(\Xi,R) \to \Xi^R$ where $\Xi^R(\alpha) = R^T \Xi(R \alpha)$, $\forall \alpha \in M$. Finally, there is a right action of $\mathrm{O}_n$ on ${\mathcal Y} (M)$ by ${\mathcal Y} (M) \times \mathrm{O}_n \to {\mathcal Y} (M)$, $(\Xi,R) \to \Xi^R$ where $\Xi^R(\alpha) = R^T \Xi(R \alpha)$, $\forall \alpha \in M$. Simply we note that, in the case of ${\mathcal X} (M)$, the symbol $R^T \Xi$ means the action of $\mathrm{O}_n$ on $N$ while in the case of ${\mathcal Y} (M)$ the same symbol means the action of $\mathrm{O}_n$ on~$M$. 

For $R \in \mathrm{O}_n$, define $\sigma_R$: $M \to M$, $\alpha \mapsto R\alpha$. Then, $(d \sigma_R)_\alpha$: $T_\alpha M \to T_{R \alpha} M$ is given by $(d \sigma_R)_\alpha(\xi) = R \xi$. We note that $f^R = f \circ \sigma_R$. Now, the Sasaki volume form $|\mathrm{Vol}_M|$ is invariant by the action of $\mathrm{O}_n$ on $M$ in the sense that 
\begin{equation}
\sigma_R^* \mathrm{Vol}_M = \mathrm{Vol}_M, \qquad \forall R \in \mathrm{O}_n. 
\label{eq:macnorm_gract_Sasaki}
\end{equation}
In particular, the maps $\sigma_R$ preserve the orientation. Obviously, this implies that
$$ |\sigma_R^* \mathrm{Vol}_M| = |\mathrm{Vol}_M|, \qquad \forall R \in \mathrm{O}_n. $$
As a consequence, for any integrable function $f$: $M \to {\mathbb R}$ and for any $R \in \mathrm{O}_n$, we have 
\begin{equation}
\int_M f^R (\alpha) \, |\mathrm{Vol}_M|(\alpha) = \int_M f(\alpha) \, |\mathrm{Vol}_M|(\alpha). 
\label{eq:macnorm_gract_integral}
\end{equation}

Gradients and divergence transform under the action of $\mathrm{O}_n$ under the following rules.  Suppose $f \in {\mathcal S}(M,{\mathbb R})$ and $\Xi \in {\mathcal X}(M)$. Then, we have, for any $R \in \mathrm{O}_n$: 
\begin{equation}
\nabla_\alpha f^R = (\nabla_\alpha f)^R, \qquad \nabla_\alpha \cdot \Xi^R = (\nabla_\alpha \cdot \Xi)^R. 
\label{eq:macnorm_gract_grad_alpha}
\end{equation}
Now, suppose that $f \in {\mathcal S}(M,{\mathbb R})$ and $\Xi \in {\mathcal Y}(M)$. Then, we have, for any $R \in \mathrm{O}_n$: 
\begin{eqnarray}
\nabla_v f^R &=& (\nabla_v f)^R, \qquad \nabla_\kappa f^R = (\nabla_\kappa f)^R, \label{eq:macnorm_gract_grad_vkap} \\
\nabla_v \cdot \Xi^R &=& (\nabla_v \cdot \Xi)^R, \qquad \nabla_\kappa \cdot \Xi^R = (\nabla_\kappa \cdot \Xi)^R. \label{eq:macnorm_gract_div_vkap}
\end{eqnarray}
Suppose that $\Xi$, $\Upsilon \in {\mathcal Y}(M)$. Then, we have, for any $R \in \mathrm{O}_n$: 
\begin{equation}
\big \langle \Xi, \Upsilon \big \rangle^R = \big \langle \Xi^R, \Upsilon^R \big \rangle. 
\label{eq:macnorm_gract_innprod}
\end{equation}
Finally, Let $\mathbf{u} \in {\mathbb R}^n$. We can easily check that 
\begin{equation} 
(\mathbf{H}_\mathbf{u})^R = \mathbf{H}_{R^T \mathbf{u}}. 
\label{eq:macnorm_gract_Hu}
\end{equation}

\subsubsection{Action of $\mathrm{O}_n$ on vector-valued functions and distributions}
\label{subsubsec:macnorm_nD_gract_sob_vectval}

Let ${\mathcal S}(M,{\mathbb R}^n)$ be the space of smooth vector-valued function $\vec \varphi$: $M \to {\mathbb R}^n$. We define a right action of $\mathrm{O}_n$ on this space by:  ${\mathcal S}(M,{\mathbb R}^n) \times \mathrm{O}_n \to {\mathcal S}(M,{\mathbb R}^n)$, $(\vec \varphi, R) \mapsto \vec \varphi^R$ where $\vec \varphi^R(\alpha) = R^T \vec \varphi(R \alpha)$, $\forall \alpha \in M$. Let $\mathbf{u} \in {\mathbb S}^{n-1}$. We note that if $\vec \varphi(\alpha) \in \{\mathbf{u}\}^\bot$, then $\vec \varphi^R(\alpha) \in \{R^T\mathbf{u}\}^\bot$. In this section, we review how this action affects the various spaces of vector-valued functions and distributions defined in Section \ref{subsubsec:macnorm_GCI_subspaces}. The proofs are straightforward and are omitted.

\medskip
\noindent
(i) Let $R \in \mathrm{O}_n$. Then, $(L^2_v(L^2_\kappa)_\mathbf{u})^R = L^2_v(L^2_\kappa)_{R^T \mathbf{u}}$ and the inner-product is equivariant: 
\begin{equation} 
\big( \hspace{-1.2mm} \big( \vec{\psi}^R , \vec{\varphi}^R \big) \hspace{-1.2mm} \big)_{L^2_v(L^2_\kappa)_{R^T \mathbf{u}}} = \big( \hspace{-1.2mm} \big( \vec{\psi} , \vec{\varphi} \big) \hspace{-1.2mm} \big)_{L^2_v(L^2_\kappa)_\mathbf{u}}, \quad \forall \vec{\psi}, \, \vec{\varphi} \in L^2_v(L^2_\kappa)_\mathbf{u}, \quad \forall R \in \mathrm{O}_n. 
\label{eq:macnorm_nD_L2vL2ku_rot}
\end{equation}

\medskip
\noindent
(ii) We define a group action of $\mathrm{O}_n$ on the bundle $\pi^* (TB) \otimes (M \times {\mathbb R}^n)$ as follows: let $(e_1, \ldots, e_n)$ be an orthonormal basis of ${\mathbb R}^n$ and let $\vec{k}$ be a section of this bundle. We can write $\vec{k} = \sum_{i=1}^n \vec{k}_i \otimes e_i$, with $\vec{k}_i$ being a section of $\pi^* B$, for all $i \in \{1, \ldots, n\}$. We let $\vec{k}^R =  \sum_{i=1}^n \vec{k}_i^R \otimes R^T e_i$, where we recall that $\vec{k}_i^R(\alpha) = R^T  \vec{k}_i(R\alpha)$, $\forall \alpha \in M$. We easily verify that this group action does not depend on the choice of the orthonormal basis $(e_1, \ldots, e_n)$. Under this group action, we have $(\pi^* (TB) \otimes (M \times \{ \mathbf{u} \}^\bot))^R = (\pi^* (TB) \otimes (M \times \{ R^T \mathbf{u} \}^\bot))$ and the $\vec{L}^2_v(L^2_\kappa)_\mathbf{u}$ inner-product is equivariant by $\mathrm{O}_n$: 
\begin{equation} \big( \hspace{-1.2mm} \big( \vec{k}^R , \vec{h}^R \big) \hspace{-1.2mm} \big)_{\vec{L}^2_v(L^2_\kappa)_{R^T \mathbf{u}}} = \big( \hspace{-1.2mm} \big( \vec{k} , \vec{h} \big) \hspace{-1.2mm} \big)_{\vec{L}^2_v(L^2_\kappa)_\mathbf{u}}, \quad \forall \vec{k}, \, \vec{h} \in \vec{L}^2_v(L^2_\kappa)_\mathbf{u}, \quad \forall R \in \mathrm{O}_n, 
\label{eq:macnorm_nD_vecL2vL2ku_rot}
\end{equation}

\medskip
\noindent
(iii) Under the group action defined at Point (i), we have $(L^2_v(H^1_\kappa)_\mathbf{u})^R = L^2_v(H^1_\kappa)_{R^T \mathbf{u}}$ and for any $\vec{\psi} \in L^2_v(H^1_\kappa)_\mathbf{u}$, $\nabla_\kappa (\vec{\psi}^R) = (\nabla_\kappa \vec{\psi})^R$ (where at the right-hand side, the group action is that defined at Point (ii)). Furthermore, the $L^2_v(H^1_\kappa)_\mathbf{u}$ inner-product is equivariant by~$\mathrm{O}_n$: 
\begin{equation} \big( \hspace{-1.2mm} \big( \nabla_\kappa \vec{\psi}^R , \nabla_\kappa \vec{\varphi}^R \big) \hspace{-1.2mm} \big)_{\vec{L}^2_v(L^2_\kappa)_{R^T \mathbf{u}}} = \big( \hspace{-1.2mm} \big( \nabla_\kappa \vec{\psi} , \nabla_\kappa \vec{\varphi} \big) \hspace{-1.2mm} \big)_{\vec{L}^2_v(L^2_\kappa)_\mathbf{u}}, \quad \forall \vec{\psi}, \, \vec{\varphi} \in L^2_v(H^1_\kappa)_\mathbf{u}, \quad \forall R \in \mathrm{O}_n, 
\label{eq:macnorm_nD_L2vH1ku_rot}
\end{equation}
together with \eqref{eq:macnorm_nD_L2vL2ku_rot}.

\medskip
\noindent
(iv) We define an action of $\mathrm{O}_n$ on the space of ${\mathbb R}^n$-valued distribution ${\mathcal D}'(M, {\mathbb R}^n)$ (the dual of the space ${\mathcal D}(M, {\mathbb R}^n)$ of smooth functions $M \to {\mathbb R}^n$ with compact support) by 
\begin{equation}
\big( \hspace{-1.2mm} \big( \vec{T}^R , \vec{\psi} \big) \hspace{-1.2mm} \big)_{{\mathcal D}'(M, {\mathbb R}^n),{\mathcal D}(M, {\mathbb R}^n)} = \big( \hspace{-1.2mm} \big( \vec{T} , \vec{\psi}^{R^T} \big) \hspace{-1.2mm} \big)_{{\mathcal D}'(M, {\mathbb R}^n),{\mathcal D}(M, {\mathbb R}^n)}
, \quad \forall \vec{T} \in {\mathcal D}'(M, {\mathbb R}^n), \quad \forall \vec{\psi} \in {\mathcal D}(M, {\mathbb R}^n).  
\label{eq:macnorm_nD_gract_distrib_vect1} 
\end{equation}
More explicitely, let $(e_1, \ldots, e_n)$ be an orthonormal basis of ${\mathbb R}^n$ and let $\vec{T} \in {\mathcal D}'(M, {\mathbb R}^n)$. Write $\vec{T} = \sum_{i=1}^n \vec{T}_i \, e_i$ with $\vec{T}_i \in {\mathcal D}'(M)$ for all $i \in \{1, \ldots, n\}$. Then, we have 
\begin{equation}
\vec{T}^R = \sum_{i=1}^n \vec{T}_i^R \, R^T e_i.
\label{eq:macnorm_nD_grpact_distrib_vect2}
\end{equation}
Hence, we have $({\mathcal D}_\mathbf{u}'(M))^R = ({\mathcal D}_{R^T \mathbf{u}}'(M))$ and the duality bracket between ${\mathcal D}_\mathbf{u}'(M)$ and ${\mathcal D}_\mathbf{u}(M)$ is equivariant: 
\begin{equation}
\big( \hspace{-1.2mm} \big( \vec{T}^R , \vec{\psi}^R \big) \hspace{-1.2mm} \big)_{{\mathcal D}_{R^T \mathbf{u}}'(M),{\mathcal D}_{R^T \mathbf{u}}(M)} = \big( \hspace{-1.2mm} \big( \vec{T} , \vec{\psi} \big) \hspace{-1.2mm} \big)_{{\mathcal D}_\mathbf{u}'(M),{\mathcal D}_\mathbf{u}(M)}
, \quad \forall \vec{T} \in {\mathcal D}'_\mathbf{u}(M), \quad \forall \vec{\psi} \in {\mathcal D}_\mathbf{u}(M).  
\label{eq:macnorm_nD_gract_distrib_vect3} 
\end{equation}
We deduce that the duality bracket between $L^2_v(H^{-1}_\kappa)_\mathbf{u}$ and $L^2_v(H^1_\kappa)_\mathbf{u}$ is also equivariant: 
\begin{equation}
\big( \hspace{-1.2mm} \big( \vec{T}^R , \vec{\psi}^R \big) \hspace{-1.2mm} \big)_{L^2_v(H^{-1}_\kappa)_{R^T \mathbf{u}},L^2_v(H^1_\kappa)_{R^T \mathbf{u}}} = \big( \hspace{-1.2mm} \big( \vec{T} , \vec{\psi} \big) \hspace{-1.2mm} \big)_{L^2_v(H^{-1}_\kappa)_\mathbf{u},L^2_v(H^1_\kappa)_\mathbf{u}}
, \quad \forall \vec{T} \in L^2_v(H^{-1}_\kappa)_\mathbf{u}, \quad \forall \vec{\psi} \in L^2_v(H^1_\kappa)_\mathbf{u}.  
\label{eq:macnorm_nD_gract_distrib_vect4} 
\end{equation}
Then, thanks to \eqref{eq:macnorm_nD_duality_vector_distrib}, \eqref{eq:macnorm_gract_grad_vkap}, \eqref{eq:macnorm_gract_div_vkap}, \eqref{eq:macnorm_gract_Hu}, one easily verifies that for all $R \in \mathrm{O}_n$, we have 
\begin{equation}
\big\{ \vec{T}^R, \mathbf{H}_{R^T \mathbf{u}} \big\} = \big\{ \vec{T}, \mathbf{H}_\mathbf{u} \big\}^R, 
\quad \forall T \in {\mathcal D}'_\mathbf{u}(M). 
\label{eq:macnorm_nD_bracket_rot2_u} 
\end{equation}
We conclude that $({\mathbb Y}_\mathbf{u})^R = {\mathbb Y}_{R^T \mathbf{u}}$, $\forall R \in \mathrm{O}_n$.

\subsection{Action of the orthogonal group on the vector GCI}
\label{subsec:macnorm_vectorGCI}

We have the following

\begin{lemma}[Action of $\mathrm{O}_n$ on the vector GCI]~

\noindent
Let $\mathbf{u} \in {\mathbb S}^{n-1}$ and let $\vec \zeta_{\mathbf{u}}$ be the vector GCI associated with $\mathbf{u}$. Then, we have 
\begin{equation}
\vec \zeta_{\mathbf{u}}^R = \vec \zeta_{R^T \mathbf{u}}, \qquad \forall R \in \mathrm{O}_n. 
\label{eq:macnorm_gract_zetaj}
\end{equation}
\label{lem:macnorm_gract_zetaj}
\end{lemma}

\noindent
\textbf{Proof.} As seen in the proof of Corollary \ref{cor:macnorm_exist_GCI}, $\vec{\zeta}_\mathbf{u}$ is the unique solution of variational formulation~\eqref{eq:macnorm_GCI_gene_varform_u} with $\vec{\chi}$ given by $\vec{\chi}(\alpha) = P_{\mathbf{u}^\bot} v$, $\forall \alpha = (v,\kappa) \in M$. Then, the results of Section \ref{subsubsec:macnorm_nD_gract_sob_vectval} show that $\vec{\zeta}_\mathbf{u}^R$ is a solution of the same variational formulation but with $\vec{\chi}$ replaced by $\vec{\chi}^R$. We compute, for all $\alpha = (v,\kappa) \in M$:
$$ \vec{\chi}^R(\alpha) = R^T \vec{\chi} (R \alpha) = R^T P_{\mathbf{u}^\bot} (Rv) = R^T \big( Rv - (Rv \cdot \mathbf{u}) \mathbf{u} \big) = v - (v \cdot R^T \mathbf{u}) R^T \mathbf{u} = P_{R^T \mathbf{u}} v. $$
Hence, by uniqueness of the solution to the variational formulation \eqref{eq:macnorm_GCI_gene_varform_u}, we obtain \eqref{eq:macnorm_gract_zetaj}, which ends the proof. \endproof

With \eqref{eq:macnorm_gract_zetaj}, it is enough to determine $\vec \zeta_\mathbf{u}$ for a single vector $\mathbf{u}$, for instance the last vector $\mathbf{e}_n$ of the canonical basis of ${\mathbf R}^n$. Indeed, since $\textrm{O}_n$ acts transitively on $B$, there exists a (non-unique) orthogonal transformation  $R \in \textrm{O}_n$ such that $R^T \mathbf{e}_n = \mathbf{u}$. Then, $\vec \zeta_\mathbf{u} = \vec \zeta_{\mathbf{e}_n}^R$, which allows us to define $\vec \zeta_\mathbf{u}$ from $\vec \zeta_{\mathbf{e}_n}$. Thus, in all the developments about the determination of $\vec \zeta_\mathbf{u}$, we will specify $\mathbf{u} = \mathbf{e}_n$ and will use indifferently $\mathbf{u}$ or $\mathbf{e}_n$ to designate it. 

By \eqref{eq:macnorm_gract_zetaj}, $\vec \zeta_{\mathbf{u}}$ satisfies 
\begin{equation}
\vec \zeta_{\mathbf{u}}^R = \vec \zeta_{\mathbf{u}}, \qquad \forall R \in \mathrm{O}_n \, \, \textrm{ such that } \, \, R \mathbf{u} = \mathbf{u}.
 \label{eq:macnorm_gract_zeta_invar}
\end{equation}
The subgroup of $\mathrm{O}_n$ consisting of orthogonal transformations $R$ such that $R \mathbf{u} = \mathbf{u}$ is the isotropy subgroup of $\mathbf{u}$ and is isomorphic to $\mathrm{O}_{n-1}$. We will identify below these two groups. We will use this invariance relation to obtain a simplified form of $\vec \zeta_\mathbf{u}$. 
To start with, we first show that we can restrict variational formulation \eqref{eq:macnorm_GCI_gene_varform_u} to spaces of invariant functions under the action of $\mathrm{O}_{n-1}$

\subsection{Variational formulation in spaces of invariant functions under the action of $\mathrm{O}_{n-1}$}
\label{subsec:macnorm_nD_varform_invar}

We first define appropriate spaces of invariant functions and distributions under the action of~$\mathrm{O}_{n-1}$. Then, we define rotational averages, which are projections from functions or distribution spaces onto the corresponding subspaces of invariant objects. This eventually enables us to show that the variational formulation for the vector GCI \eqref{eq:macnorm_GCI_gene_varform_u} can be restricted to these invariant spaces.

\subsubsection{Spaces of invariant vector-valued functions and distributions on $M$}
\label{subsubsec:macnorm_nD_gract_sob_invar}

In this section, we define spaces of invariant vector-valued functions and distributions on $M$ and state some properties which are direct consequences of Section \ref{subsubsec:macnorm_nD_gract_sob_vectval}. 

\medskip
\noindent 
(i) We define the subspace of invariant functions $L^2_v(L^2_\kappa)_\mathrm{inv}$ of $L^2_v(L^2_\kappa)_\mathbf{u}$ by
$$ L^2_v(L^2_\kappa)_\mathrm{inv} = \big\{ \vec{\psi} \in L^2_v(L^2_\kappa)_\mathbf{u} \, \, | \, \, \vec{\psi}^R = \vec{\psi}, \, \, \forall R \in \mathrm{O}_{n-1} \big\}. $$
 The space $L^2_v(L^2_\kappa)_\mathrm{inv}$ is a closed subspace of $L^2_v(L^2_\kappa)_\mathbf{u}$ and hence, is a Hilbert space normed by~$\|~\cdot~\|_{L^2_v(L^2_\kappa)_\mathbf{u}}$.

\medskip
\noindent 
(ii)  We define 
$$ \vec{L}^2_v(L^2_\kappa)_\mathrm{inv} = \big\{ \vec{k} \in \vec{L}^2_v(L^2_\kappa)_\mathbf{u} \, \, | \, \, \vec{k}^R = \vec{k}, \, \, \forall R \in \mathrm{O}_{n-1} \big\}. $$
Again, $\vec{L}^2_v(L^2_\kappa)_\mathrm{inv}$ is a closed subspace of $\vec{L}^2_v(L^2_\kappa)_\mathbf{u}$ and hence, is a Hilbert space.

\medskip
\noindent 
(iii) We define
$$ L^2_v(H^1_\kappa)_\mathrm{inv} = \big\{ \vec{\psi} \in L^2_v(H^1_\kappa)_\mathbf{u} \, \, | \, \, \vec{\psi}^R = \vec{\psi}, \, \, \forall R \in \mathrm{O}_{n-1} \big\}. $$
We note that, for  $\vec{\psi} \in L^2_v(H^1_\kappa)_\mathrm{inv}$, $\nabla_\kappa \vec{\psi} \in \vec{L}^2_v(L^2_\kappa)_\mathrm{inv}$. The space $L^2_v(H^1_\kappa)_\mathrm{inv}$ is a closed subspace of $L^2_v(H^1_\kappa)_\mathbf{u}$ and hence, a Hilbert space, normed by $\| \cdot \|_{L^2_v(H^1_\kappa)_\mathbf{u}}$.

\medskip
\noindent
(iv) The space ${\mathcal D}_\mathrm{inv}(M)$ is the subspace of ${\mathcal D}_\mathbf{u}(M)$ of elements $\vec{\varphi}$ such that $\vec{\varphi}^R = \vec{\varphi}$. The space ${\mathcal D}'_\mathrm{inv}(M)$ is the subspace of ${\mathcal D}'_\mathbf{u}(M)$ of elements $\vec{T}$ such that $\vec{T}^R = \vec{T}$. We note that ${\mathcal D}'_\mathrm{inv}(M)$ is the dual of ${\mathcal D}_\mathrm{inv}(M)$.

\medskip
\noindent
(v) The subspace of invariant elements of $L^2_v(H^{-1}_\kappa)_\mathbf{u}$ by the group $\mathrm{O}_{n-1}$ is given by 
$$ L^2_v(H^{-1}_\kappa)_\mathrm{inv} = \big\{ \vec{T} \in L^2_v(H^{-1}_\kappa)_\mathbf{u} \, \, | \, \, \vec{T}^R = \vec{T}, \, \, \forall R \in \mathrm{O}_{n-1} \big\}. $$
It is the dual of $L^2_v(H^1_\kappa)_\mathrm{inv}$. It is again a closed subspace of $L^2_v(H^{-1}_\kappa)_\mathbf{u}$ and hence, a Hilbert space, still normed by $\| \cdot \|_{L^2_v(H^{-1}_\kappa)_\mathbf{u}}$.  

\medskip
\noindent
(vi) The subspace of invariant elements of ${\mathbb Y}_\mathbf{u}$ by the action of $\mathrm{O}_{n-1}$ is 
$$ {\mathbb Y}_\mathrm{inv} = \big\{ \vec{\psi} \in {\mathbb Y}_\mathbf{u} \, \, | \, \, \vec{\psi}^R = \vec{\psi}, \, \, \forall R \in \mathrm{O}_{n-1} \big\}. $$
For $\vec{\psi} \in {\mathbb Y}_\mathrm{inv}$, we have $\nabla_\kappa \vec{\psi} \in \vec{L}^2_v(L^2_\kappa)_\mathrm{inv}$ and $\{ \vec{\psi}, \mathbf{H}_\mathbf{u} \} \in L^2_v(H^{-1}_\kappa)_\mathrm{inv}$. The space ${\mathbb Y}_\mathrm{inv}$ is a closed subspace of ${\mathbb Y}_\mathbf{u}$ and hence, is a Hilbert space, still normed by $\| \cdot \|_{{\mathbb Y}_\mathbf{u}}$.

\begin{lemma}
We have
\begin{equation}
\vec{\psi} \in L^2_v(L^2_\kappa)_\mathrm{inv} \quad \Longrightarrow \quad \int_M \vec{\psi} \, M_\mathbf{u} \, | \mathrm{Vol}_M | = 0. 
\label{eq:macnorm_int_vecpsi=0}
\end{equation}
\label{lem:macnorm_int_vecpsi=0}
\end{lemma}

\noindent
\textbf{Proof.} Let $\vec{\psi} \in L^2_v(L^2_\kappa)_\mathrm{inv}$. Thanks to \eqref{eq:macnorm_gract_integral} and to the fact that $M_\mathbf{u}^R = M_\mathbf{u}$, $\forall R \in \mathrm{O}_{n-1}$, we have 
$$ \int_M \vec{\psi} \, M_\mathbf{u} \, | \mathrm{Vol}_M | = R \int_M \vec{\psi}^R \, M_\mathbf{u} \, | \mathrm{Vol}_M | = R \int_M \vec{\psi} \, M_\mathbf{u} \, | \mathrm{Vol}_M |, \quad \forall R \in \mathrm{O}_{n-1}. $$
Hence, $\int_M \vec{\psi} \, M_\mathbf{u} \, | \mathrm{Vol}_M |$ is a vector in $\{ \mathbf{u} \}^\bot$, which is invariant by $\mathrm{O}_{n-1}$. The only such vector is~$0$. \endproof.

\subsubsection{Rotational averages}
\label{subsubsec:macnorm_nD_rotaver}

We define projection maps from spaces of vector-valued functions and distributions on $M$ to corresponding spaces of invariant functions. These projection maps are nothing but rotational averages.  

\medskip
\noindent
(i) For $\vec{\varphi} \in {\mathcal D}_\mathbf{u}(M)$, we define its average $\vec{\varphi}_\mathrm{av}$ over $\mathrm{O}_{n-1}$ by 
\begin{equation}
\vec{\varphi}_\mathrm{av} = \int_{\mathrm{O}_{n-1}} \vec{\varphi}^R \, dR, 
\label{eq:macnorm_average_fcts}
\end{equation}
where $dR$ is the normalized Haar measure on $\mathrm{O}_{n-1}$. We have $\vec{\varphi}_\mathrm{av} \in {\mathcal D}_\mathrm{inv}(M)$. Furthermore, the map ${\mathcal D}_\mathbf{u}(M) \to {\mathcal D}_\mathbf{u}(M)$, $\vec{\varphi} \mapsto \vec{\varphi}_\mathrm{av}$ is continuous for the topology of ${\mathcal D}_\mathbf{u}(M)$ and 
$$ \nabla_v \vec{\varphi}_\mathrm{av} = (\nabla_v \vec{\varphi})_\mathrm{av}, \qquad 
\nabla_\kappa \vec{\varphi}_\mathrm{av} = (\nabla_\kappa \vec{\varphi})_\mathrm{av}. $$

\medskip
\noindent
(ii) Let $\vec{T} \in {\mathcal D}'_\mathbf{u}(M)$. Then, $\vec{T}_\mathrm{av}$ is the element of ${\mathcal D}'_\mathbf{u}(M)$ such that 
\begin{equation}
\big( \hspace{-1.2mm} \big( \vec{T}_\mathrm{av}, \vec{\varphi} \big) \hspace{-1.2mm} \big)_{{\mathcal D}'_\mathbf{u}, {\mathcal D}_\mathbf{u}} 
= \big( \hspace{-1.2mm} \big( \vec{T}, \vec{\varphi}_\mathrm{av} \big) \hspace{-1.2mm} \big)_{{\mathcal D}'_\mathbf{u}, {\mathcal D}_\mathbf{u}} 
= \int_{\mathrm{O}_{n-1}} \big( \hspace{-1.2mm} \big( \vec{T}^R, \vec{\varphi} \big) \hspace{-1.2mm} \big)_{{\mathcal D}'_\mathbf{u}, {\mathcal D}_\mathbf{u}} \, dR 
= \int_{\mathrm{O}_{n-1}} \big( \hspace{-1.2mm} \big( \vec{T}, \vec{\varphi}^R \big) \hspace{-1.2mm} \big)_{{\mathcal D}'_\mathbf{u}, {\mathcal D}_\mathbf{u}} \, dR, 
\label{eq:macnorm_average_distrib}
\end{equation}
where the last two equalities come from the linearity and continuity of $\vec{T}$. Thanks to the invariance of the Haar measure by matrix inversion, the definitions \eqref{eq:macnorm_average_distrib} and \eqref{eq:macnorm_average_fcts} are identical when $T \in {\mathcal D}_\mathbf{u}(M)$. The map ${\mathcal D}'_\mathbf{u}(M) \to {\mathcal D}'_\mathbf{u}(M)$, $\vec{T} \mapsto \vec{T}_\mathrm{av}$ is continuous for the topology of ${\mathcal D}_\mathbf{u}(M)$ and we have 
\begin{equation} 
\nabla_v \vec{T}_\mathrm{av} = (\nabla_v \vec{T})_\mathrm{av}, \qquad 
\nabla_\kappa \vec{T}_\mathrm{av} = (\nabla_\kappa \vec{T})_\mathrm{av}. 
\label{eq:macnorm_commut_av_deriv}
\end{equation}

\begin{lemma}[Properties of rotational averages]~

\noindent
(i) We have 
\begin{equation} 
\vec{\varphi} \in {\mathbb Y}_\mathbf{u} \quad \Longrightarrow \quad \vec{\varphi}_\mathrm{av} \in {\mathbb Y}_\mathrm{inv}, \qquad \textrm{and} \qquad 
\| \vec{\varphi}_\mathrm{av} \|_{{\mathbb Y}_\mathbf{u}} \leq \| \vec{\varphi} \|_{{\mathbb Y}_\mathbf{u}}. 
\label{eq:macnorm_aver_prop}
\end{equation}

\noindent
(ii) For all $\vec{\varphi} \in {\mathbb Y}_\mathbf{u}$ and $\vec{\psi} \in {\mathbb Y}_\mathrm{inv}$, we have 
\begin{eqnarray}
&&\hspace{-1cm}
\int_{\mathrm{O}_{n-1}} \big( \hspace{-1.2mm} \big( \vec{\psi} , \vec{\varphi}^R \big) \hspace{-1.2mm} \big)_{L^2_v(L^2_\kappa)_\mathbf{u}} \, dR = \big( \hspace{-1.2mm} \big( \vec{\psi} , \vec{\varphi}_\mathrm{av} \big) \hspace{-1.2mm} \big)_{L^2_v(L^2_\kappa)_\mathbf{u}}, \label{eq:macnorm_aver_prop_L2}
\\
&&\hspace{-1cm}
\int_{\mathrm{O}_{n-1}}  \big( \hspace{-1.2mm} \big( \nabla_\kappa \vec{\psi}, \nabla_\kappa \vec{\varphi}^R \big) \hspace{-1.2mm} \big)_{\vec{L}^2_v(L^2_\kappa)_\mathbf{u}} \, dR = \big( \hspace{-1.2mm} \big( \nabla_\kappa \vec{\psi}, \nabla_\kappa \vec{\varphi}_\mathrm{av} \big) \hspace{-1.2mm} \big)_{\vec{L}^2_v(L^2_\kappa)_\mathbf{u}},  \label{eq:macnorm_aver_prop_H1}
\\
&&\hspace{-1cm}
\int_{\mathrm{O}_{n-1}}  \big( \hspace{-1.2mm} \big( \{ \vec{\psi}, \mathbf{H}_{\mathbf{u}} \} , \vec{\varphi}^R \big) \hspace{-1.2mm} \big)_{L^2_v(H^{-1}_\kappa)_\mathbf{u},L^2_v(H^1_\kappa)_\mathbf{u}} \, dR = \big( \hspace{-1.2mm} \big( \{ \vec{\psi}, \mathbf{H}_{\mathbf{u}} \} , \vec{\varphi}_\mathrm{av} \big) \hspace{-1.2mm} \big)_{L^2_v(H^{-1}_\kappa)_\mathbf{u},L^2_v(H^1_\kappa)_\mathbf{u}}. \label{eq:macnorm_aver_prop_bracket}
\end{eqnarray}
\label{lem:macnorm_prop_averages}
\end{lemma}

\noindent
\textbf{Proof of Lemma \ref{lem:macnorm_prop_averages}.} (i) {\em Step 1.} By Cauchy-Schwarz inequality, Fubini's theorem and \eqref{eq:macnorm_nD_L2vL2ku_rot}, we have
\begin{eqnarray*} 
\| \vec{\varphi}_\mathrm{av} \|_{L^2_v(L^2_\kappa)_\mathbf{u}}^2 &=& \int_M \Big| \int_{\mathrm{O}_{n-1}} \vec{\varphi}^R \, dR \Big|^2 \, \mathbf{M}_\mathbf{u} \, |\mathrm{Vol}_M | \leq   \int_M \Big( \int_{\mathrm{O}_{n-1}} | \vec{\varphi}^R |^2 \, dR \Big) \, \mathbf{M}_\mathbf{u} \, |\mathrm{Vol}_M |, \\
&=& \int_{\mathrm{O}_{n-1}} \| \vec{\varphi}^R \|_{L^2_v(L^2_\kappa)_\mathbf{u}}^2 \, dR = \int_{\mathrm{O}_{n-1}} \| \vec{\varphi} \|_{L^2_v(L^2_\kappa)_\mathbf{u}}^2 \, dR = \| \vec{\varphi} \|_{L^2_v(L^2_\kappa)_\mathbf{u}}^2. 
\end{eqnarray*}
Here, for an element of $X$ of $\{\mathbf{u}\}^\bot$, we have denoted by $|X|^2 = \sum_{i=1}^{n-1} (X \cdot e_i)^2$. 

\medskip
\noindent
{\em Step 2.} Using Section \ref{subsubsec:macnorm_nD_gract_sob_vectval} (iii) and the same arguments as just above, we can also prove 
$$ \| \nabla_\kappa \vec{\varphi}_\mathrm{av} \|_{\vec{L}^2_v(L^2_\kappa)_\mathbf{u}}^2 \leq \| \nabla_\kappa \vec{\varphi} \|_{\vec{L}^2_v(L^2_\kappa)_\mathbf{u}}^2.$$

\medskip
\noindent
{\em Step 3.} We take $\vec{\varphi} \in L^2_v(H^{-1}_\kappa)_\mathbf{u}$ and show that $\vec{\varphi}_\mathrm{av} \in L^2_v(H^{-1}_\kappa)_\mathbf{u}$. Indeed, for all $\psi \in L^2_v(H^1_\kappa)_\mathbf{u}$, we have, thanks to the previous two steps
\begin{eqnarray*}
\big| \big( \hspace{-1.2mm} \big( \vec{\varphi}_\mathrm{av}, \vec{\psi} \big) \hspace{-1.2mm} \big)_{{\mathcal D}'_\mathbf{u}, {\mathcal D}_\mathbf{u}} \big| &=& \big| \big( \hspace{-1.2mm} \big( \vec{\varphi}, \vec{\psi}_\mathrm{av} \big) \hspace{-1.2mm} \big)_{L^2_v(H^{-1}_\kappa)_\mathbf{u}, L^2_v(H^1_\kappa)_\mathbf{u}} \big| \\
&\leq& 
\| \vec{\varphi} \|_{L^2_v(H^{-1}_\kappa)_\mathbf{u}} \, \| \vec{\psi}_\mathrm{av} \|_{L^2_v(H^1_\kappa)_\mathbf{u}} \leq \| \vec{\varphi} \|_{L^2_v(H^{-1}_\kappa)_\mathbf{u}} \, \| \vec{\psi} \|_{L^2_v(H^1_\kappa)_\mathbf{u}}, 
\end{eqnarray*}
which shows that $\vec{\varphi}_\mathrm{av} \in L^2_v(H^{-1}_\kappa)_\mathbf{u}$ and 
\begin{equation} 
\| \vec{\varphi}_\mathrm{av} \|_{L^2_v(H^{-1}_\kappa)_\mathbf{u}} \leq \| \vec{\varphi} \|_{L^2_v(H^{-1}_\kappa)_\mathbf{u}}.
\label{eq:macnorm_prop_averages_prf1}
\end{equation}

\medskip
\noindent
{\em Step 4.} For $\vec{\varphi} \in {\mathcal D}'_\mathbf{u}(M)$, we have, thanks to \eqref{eq:macnorm_commut_av_deriv},  
$$ \big\{ \vec{\varphi}_\mathrm{av} , \mathbf{H}_\mathbf{u} \big\} = \big\{ \vec{\varphi} , \mathbf{H}_\mathbf{u} \big\}_\mathrm{av}. $$
From \eqref{eq:macnorm_prop_averages_prf1}, it follows that if $\{ \vec{\varphi} , \mathbf{H}_\mathbf{u} \} \in L^2_v(H^{-1}_\kappa)_\mathbf{u}$, then, $\{ \vec{\varphi}_\mathrm{av} , \mathbf{H}_\mathbf{u} \} \in L^2_v(H^{-1}_\kappa)_\mathbf{u}$ and
$$ \| \big\{ \vec{\varphi}_\mathrm{av} , \mathbf{H}_\mathbf{u} \big\} \|_{L^2_v(H^{-1}_\kappa)_\mathbf{u}} \leq \| \big\{ \vec{\varphi} , \mathbf{H}_\mathbf{u} \big\} \|_{L^2_v(H^{-1}_\kappa)_\mathbf{u}}.$$

\medskip
\noindent
{\em Final step.} Clearly, $\vec{\varphi}_\mathrm{av}$ is invariant by $\mathrm{O}_{n-1}$ and collecting all the previous estimates, we get the estimate in the last statement of \eqref{eq:macnorm_aver_prop}. It follows that $\vec{\varphi}_\mathrm{av} \in {\mathbb Y}_\mathrm{inv}$ which completes the proof of item (i). 

\bigskip
\noindent
(ii) Eqs. \eqref{eq:macnorm_aver_prop_L2} and \eqref{eq:macnorm_aver_prop_H1} follow from Fubini's theorem and \eqref{eq:macnorm_aver_prop_bracket}, from \eqref{eq:macnorm_average_distrib}. \endproof

\subsubsection{Variational formulation for the vector GCI on invariant spaces}
\label{subsubsec:macnorm_nD_varform_invar}

In this section, we show that the variational formulation for the vector GCI \eqref{eq:macnorm_GCI_gene_varform_u} can be equivalently formulated within spaces of invariant functions.

\begin{lemma}[Variational formulation involving invariant spaces]~

\noindent
Let $\vec{\chi} \in L^2_v(L^2_\kappa)_\mathrm{inv}$. Then, $\vec{\psi}$ is the
unique solution of the variational formulation \eqref{eq:macnorm_GCI_gene_varform_u} if and only if it is the unique solution of the following variational formulation
\begin{equation}
\mbox{} \hspace{-0.5cm} \left\{ \begin{array}{l}
\vec{\psi} \in {\mathbb Y}_\mathrm{inv}, \\
- \big( \hspace{-1.2mm} \big( \{ \vec{\psi}, \mathbf{H}_{\mathbf{u}} \} , \vec{\varphi} \big) \hspace{-1.2mm} \big)_{L^2_v(H^{-1}_\kappa)_\mathbf{u},L^2_v(H^1_\kappa)_\mathbf{u}} + D \big( \hspace{-1.2mm} \big( \nabla_\kappa \vec{\psi}, \nabla_\kappa \vec{\varphi} \big) \hspace{-1.2mm} \big)_{\vec{L}^2_v(L^2_\kappa)_\mathbf{u}} \\
\hspace{8cm} = - \big( \hspace{-1.2mm} \big( \vec{\chi} , \vec{\varphi} \big) \hspace{-1.2mm} \big)_{L^2_v(L^2_\kappa)_\mathbf{u}}, \quad \forall \vec{\varphi} \in {\mathbb Y}_\mathrm{inv}. \end{array} \right. 
\label{eq:macnorm_GCI_gene_varform_inv}
\end{equation}
\label{lem:macnorm_varform_invar}
\end{lemma}

\noindent
\textbf{Proof.} We first note that the variational formulation \eqref{eq:macnorm_GCI_gene_varform_inv} has a unique solution because we can reproduce the proof of Theorem \ref{thm:macnorm_GCI_gene_varform_exist} using the invariant vector-valued version of the spaces used in the proof. We just note that, because of Lemma \ref{lem:macnorm_int_vecpsi=0}, an element of $\vec{\chi} \in L^2_v(L^2_\kappa)_\mathrm{inv}$ or ${\mathbb Y}_\mathrm{inv}$ automatically satisfies $\int_{M} \vec{\chi} \, M_\mathbf{u} \, |\mathrm{Vol}_M| = 0$. 

Now, let $\vec{\psi}$ be the unique solution of \eqref{eq:macnorm_GCI_gene_varform_u}. Then, $\vec{\psi} \in {\mathbb Y}_\mathrm{inv}$. Indeed, using $\vec{\varphi}^R$ as a test function in the variational formulation \eqref{eq:macnorm_GCI_gene_varform_u} and applying the results of Section \ref{subsubsec:macnorm_nD_gract_sob_vectval} and the fact that $\vec{\chi} \in L^2_v(L^2_\kappa)_\mathrm{inv}$, we find that $\vec{\psi}^R$ is another solution of the same variational formulation. Since this variational formulation has a unique solution, we deduce that $\vec{\psi}^R = \vec{\psi}$. Furthermore, \eqref{eq:macnorm_GCI_gene_varform_u} is still valid by testing with function $\vec{\varphi}$ in ${\mathbb Y}_\mathrm{inv}$. So, $\vec{\psi}$ is also the unique solution of \eqref{eq:macnorm_GCI_gene_varform_inv}. 

Conversely, Suppose $\vec{\psi}$ is the unique solution of \eqref{eq:macnorm_GCI_gene_varform_inv}. Then we obviously get $\vec{\psi} \in {\mathbb Y}_\mathbf{u}$. Now, let us take $\vec{\varphi} \in {\mathbb Y}_\mathbf{u}$ and use $\vec{\varphi}_\mathrm{av} \in {\mathbb Y}_\mathrm{inv}$ as a test function in \eqref{eq:macnorm_GCI_gene_varform_inv}. Then, by Lemma \ref{lem:macnorm_prop_averages}~(ii), we get 
$$ \int_{\mathrm{O}_{n-1}} \Big\{ - \big( \hspace{-1.2mm} \big( \{ \vec{\psi}, \mathbf{H}_{\mathbf{u}} \} , \vec{\varphi}^R \big) \hspace{-1.2mm} \big)_{L^2_v(H^{-1}_\kappa)_\mathbf{u},L^2_v(H^1_\kappa)_\mathbf{u}} 
+ D \big( \hspace{-1.2mm} \big( \nabla_\kappa \vec{\psi}, \nabla_\kappa \vec{\varphi}^R \big) \hspace{-1.2mm} \big)_{\vec{L}^2_v(L^2_\kappa)_\mathbf{u}} 
+ \big( \hspace{-1.2mm} \big( \vec{\chi} , \vec{\varphi}^R \big) \hspace{-1.2mm} \big)_{L^2_v(L^2_\kappa)_\mathbf{u}} \Big\} \, dR = 0. $$
But, using the results of Section \ref{subsubsec:macnorm_nD_gract_sob_vectval} to transfer the exponent $R$ onto $\vec{\psi}$ and the fact that $\vec{\psi}$ is invariant under the action of $\mathrm{O}_{n-1}$, we end up with $\vec{\psi}$ being a solution of \eqref{eq:macnorm_GCI_gene_varform_u}, which proves the claim. \endproof

\setcounter{equation}{0}
\section{Reduced GCI pair}
\label{sec:gci_redform}

In this section, we first show that, by the invariance property \eqref{eq:macnorm_gract_zeta_invar}, the vector GCI can be reduced to two scalar quantities named the reduced GCI pair. Then from the variational formulation for the vector GCI in invariant spaces \eqref{eq:macnorm_GCI_gene_varform_inv}, we derive a variational formulation for the reduced GCI pair in appropriate function spaces.

\subsection{From vector GCI to reduced GCI pair}
\label{subsection:gci_redform_passage}

In this subsection, we state the theorem establishing that the vector GCI can be expressed in terms of the reduced GCI pair. The proof of this theorem will be deferred to Appendix \ref{sec:redGCIprf}. We start by introducing a convenient change of variables.

\subsubsection{A change of variables in $M$}
\label{subsubsec:macnorm_chgvar_M}

\begin{figure}[htbp]
\centering
\includegraphics[trim={3.5cm 13cm 4.5cm 4.5cm},clip,height= 8cm]{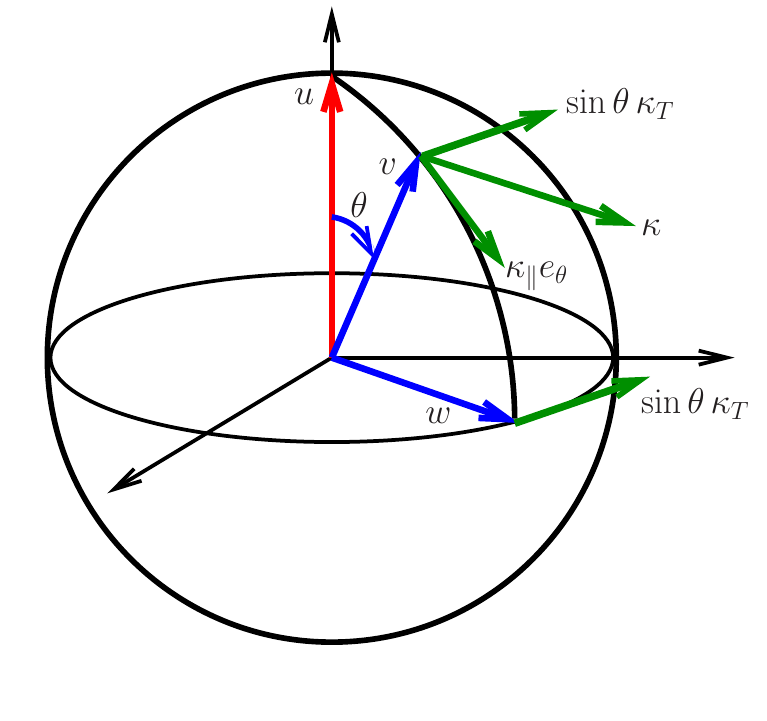}
\caption{The diffeomorphism ${\mathbf A}_\mathbf{u}$ depicted in dimension $n=3$. The velocity $v$ is changed into the pair $(\theta,w)$ where $\theta$ is the angle between $\mathbf{u}$ and $v$ and where $w$ belongs to the unit circle ${\mathbb S}^1$ in the plane orthogonal to $\mathbf{u}$ (diffeomorphism ${\mathbf V}_\mathbf{u}$ introduced in Section \ref{subsubsec:hydro_chg_var}). The tangent vector $\kappa$ to ${\mathbb S}^2$ at $v$ is changed into the pair $(\kappa_\parallel, \kappa_T)$ where $\kappa_\parallel$ is the projection of $\kappa$ on $e_\theta$, the director of the line normal to $v$ in the plane $\mathrm{Span}(\mathbf{u},v)$, and $\kappa_T$ is a tangent vector to ${\mathbb S}^1$ at $w$.}
\label{fig:chgvar}
\end{figure}

We consider the derivative of the diffeomorphism $\mathbf{V}_\mathbf{u}$ introduced in Section \ref{subsubsec:hydro_chg_var}. We recall that we denote by $B_\mathbf{u} = B \setminus \{ \pm \mathbf{u} \}$ and we denote by $M_\mathbf{u}$ the support manifold of $T B_\mathbf{u}$. We also recall that we identify ${\mathbb S}^{n-2} \approx B \cap \{\mathbf{u}\}^\bot$. We consider the tangent bundle $T(0,\pi)$ whose support manifold will be denoted by $M_V$ (we note that $M_V \approx (0,\pi) \times {\mathbb R}$) and its elements by $(\theta, \kappa_\parallel)$. Likewise, we consider the tangent bundle $T {\mathbb S}^{n-2}$ whose support manifold will be denoted by $M_T$ and its elements, by $(w, \kappa_T)$. We now introduce the bundle diffeomorphism ${\mathbf A}_\mathbf{u}$: $M_V \times M_T \to M_\mathbf{u}$ given by $ {\mathbf A}_\mathbf{u} = d \mathbf{V}_\mathbf{u}$,  where $\mathbf{V}_\mathbf{u}$ is given by \eqref{eq:macnorm_V_def}. Specifically, ${\mathbf A}_\mathbf{u}$ is given by 
\begin{equation} 
\mathbf{A}_\mathbf{u} \big( (\theta, \kappa_\parallel),(w, \kappa_T) \big) = \alpha = (v,\kappa) = \big( \cos \theta \, \mathbf{u} + \sin \theta \, w, \kappa_\parallel \, e_\theta + \sin \theta \, \kappa_T \big) , 
\label{eq:macnorm_dV_express}
\end{equation}
with $e_\theta$ given by \eqref{eq:hydro_chgvar_ethetdef} and where 
\begin{equation} 
\kappa_T \cdot w = \kappa_T \cdot \mathbf{u} = 0,
\label{eq:macnorm_dV_express_2}
\end{equation}
the latter being a consequence of the fact that $(w, \kappa_T) \in M_T$. We note that 
\begin{equation} 
w = \frac{1}{\sin \theta} (v - \cos \theta \, \mathbf{u}), \quad e_\theta = \frac{1}{\sin \theta} (\cos \theta \, v - \mathbf{u}), \quad \kappa_\parallel = \kappa \cdot e_\theta, \quad \kappa_T \cdot e_\theta = 0. 
\label{eq:macnorm_wethet_formulas}
\end{equation}
The diffeomorphism ${\mathbf A}_\mathbf{u}$ is depicted in dimension $n=3$ in Fig. \ref{fig:chgvar}. Now, thanks to the diffeomorphism ${\mathbf A}_\mathbf{u}$, for any smooth map $\varphi$: $M \to {\mathbb R}$, we can define $\varphi^\dagger$: $M_V \times M_T \to {\mathbb R}$ by 
\begin{equation}
\varphi^\dagger = \varphi \circ {\mathbf A}_\mathbf{u}.
\label{eq:macnorm_dagger_def}
\end{equation}

\medskip
To compute the jacobian of $\mathbf{A}_\mathbf{u}$, we first need to specify the volume form on $M_V \times M_T$ and its orientation. We start with $M_T$ and denote by $\alpha_T = (w, \kappa_T)$ a generic element of $M_T$ and by $\xi^T = (\sigma^w, \sigma^{\kappa_T})$ a generic element of $T_{\alpha_T} M_T$, i.e. $(\sigma^w, \sigma^{\kappa_T}) \in {\mathbb R}^{2n}$ with 
\begin{eqnarray}
&& \hspace{-1cm}
\sigma^w \cdot \mathbf{u} = 0, \qquad \sigma^w \cdot w = 0, \label{eq:Jacobian_Au_prf10} \\
&& \hspace{-1cm}
\sigma^{\kappa_T} \cdot \mathbf{u} = 0, \qquad \sigma^{\kappa_T} \cdot w + \sigma^w \cdot \kappa_T = 0. \label{eq:Jacobian_Au_prf11}
\end{eqnarray}
We define a $n-2$-form $\omega$ on $M_T$ in way similar to~$\Omega$ (see \eqref{eq:geom_orient_form_bundle}): we denote by $\pi_T$ the bundle projection $M_T \to {\mathbb S}^{n-2}$ and, for $\xi^T_i \in T_{\alpha_T} M_T$, $i \in \{1, \ldots, n-2\}$, we let 
\begin{equation} 
\omega(\alpha_T)(\ldots, \xi^T_i, \ldots) = \mathrm{Vol}_{{\mathbb S}^{n-2}} (\pi_T \alpha_T) \big( \ldots , \bar {\mathcal B}_{\alpha_T}^V (\xi^T_i), \ldots \big), 
\label{eq:Jacobian_Au_prf11.4}
\end{equation}
where $\bar {\mathcal B}_{\alpha_T}^V$ is the connection map $T_{\alpha_T} M_T \to T_{\pi_T \alpha_T} {\mathbb S}^{n-2}$ (the bar highlighting that it is different from the connection map for $M$). Like in the case of $B$, with have $\bar {\mathcal B}_{\alpha_T}^V (\sigma^w, \sigma^{\kappa_T}) = P_{w^\bot} \sigma^{\kappa_T}$. Again, 
\begin{equation} 
\mathrm{Vol}_{\mathrm{M}_\mathrm{T}} = \pi_T^* \mathrm{Vol}_{{\mathbb S}^{n-2}} \wedge \omega. 
\label{eq:Jacobian_Au_volMT}
\end{equation}
is a volume form on $M_T$ and we orient $M_T$ by this volume form. 

For $M_V \approx (0,\pi) \times {\mathbb R}$, we can proceed analogously, but since $M_V$ is included in the flat space~${\mathbb R}^2$, we can more simply define the volume form on $M_V$ as $\mathrm{Vol}_{\mathrm{M}_\mathrm{V}} = d \theta \wedge d \kappa_\parallel$, whose density is simply $d \theta \, d \kappa_\parallel$. We deduce that 
\begin{equation}
|\mathrm{Vol}_{\mathrm{M}_\mathrm{V} \times \mathrm{M}_\mathrm{T}} | =  | \mathrm{Vol}_{\mathrm{M}_\mathrm{T}}| \, d \theta \, d \kappa_\parallel. 
\label{eq:macnorm_nD_vol_MVXMT}
\end{equation}
A computation using spherical coordinates leads to 
\begin{equation}
|\mathbf{A}_\mathbf{u}^* \mathrm{Vol}_M| = \sin^{2(n-2)} \theta \, \, | \mathrm{Vol}_{\mathrm{M}_\mathrm{T}}| \, d \theta \, d \kappa_\parallel. 
\label{eq:macnorm_nD_AuvolM2}
\end{equation}

The retrenchment $\omega_w$ of $\omega$ on $T_w {\mathbb S}^{n-2}$ is computed similarly as that of $\Omega$ (see \eqref{eq:von_Mises_norm_prf3}). Let $(\varepsilon_1, \ldots, \varepsilon_{n-2})$ be a local, positively oriented orthonormal frame in an open subset of ${\mathbb S}^{n-2}$ containing $w$. We define a coordinate system on $T_w {\mathbb S}^{n-2}$ by $(\kappa_T^1, \ldots, \kappa_T^{n-2})$ where $\kappa_T^i = \kappa_T \cdot \varepsilon_i$. Thus, we can apply \eqref{eq:von_Mises_norm_prf3} except that here, we are concerned with ${\mathbb S}^{n-2}$ instead of ${\mathbb S}^{n-1}$. Hence, we have 
\begin{equation} 
\omega_w = d \kappa_T^1 \wedge \ldots \wedge d \kappa_T^{n-2}. 
\label{eq:macnorm_nD_coef_c2_prf1}
\end{equation}


\subsubsection{Reduction of the vector GCI and reduced GCI pair}
\label{subsubsec:gci_redform_statement}

We state the main theorem estblishing the expression of the vector GCI in terms of the reduced GCI pair. Its proof is given in Appendix \ref{sec:redGCIprf}.

\begin{theorem}[Reduction of the vector GCI]~

\noindent
Let $\vec{\psi}  \in {\mathcal S}(M, \{ \mathbf{u} \}^\bot)$ (the space of smooth functions $M \to \{ \mathbf{u} \}^\bot$) satisfy the invariance relation~\eqref{eq:macnorm_gract_zeta_invar} and let $\vec{\psi}^\dagger$ be obtained from $\vec{\psi}$ through the change of variables~\eqref{eq:macnorm_dagger_def}. Then, there exists a pair of functions: $(0,\pi) \times {\mathbb R} \times [0,\infty) \to {\mathbb R}^2$, $(\theta, \kappa_\parallel, \kappa_\bot) \to (\psi_w,\psi_r)$ such that 
\begin{eqnarray}
&&\hspace{-1cm}
\vec{\psi}^\dagger (\theta, \kappa_\parallel, w, \kappa_T) = \psi_w(\theta, \kappa_\parallel, |\kappa_T|) \, w + \psi_r(\theta, \kappa_\parallel, |\kappa_T|) \, \frac{\kappa_T}{|\kappa_T|}, \nonumber \\
&&\hspace{2cm}
\forall \big( (\theta, \kappa_\parallel), (w, \kappa_T) \big) \in M_v \times M_T  \, \, \textrm{ such that } \, \, \kappa_T \not = 0  , 
\label{eq:macnorm_nD_actionOn-1}
\end{eqnarray}
and 
\begin{equation}
\vec{\psi}^\dagger (\theta, \kappa_\parallel, w, 0) = \psi_w(\theta, \kappa_\parallel, 0) \, w , \quad \forall \big( (\theta, \kappa_\parallel), w \big) \in M_V \times {\mathbb S}^{n-2} .
\label{eq:macnorm_nD_actionOn-0.5}
\end{equation}
When $\vec{\psi} = \vec{\zeta}_\mathbf{u}$ is the vector GCI, the associated pair $(\zeta_w, \zeta_r)$ is called the reduced GCI pair. 
\label{thm:macnorm_nD_actionOn-1}
\end{theorem}

\begin{remark}[Consistency of \eqref{eq:macnorm_nD_actionOn-1}]~

\noindent
From Section \ref{subsec:macnorm_vectorGCI}, we know that $\vec{\psi}$, $w$, $\kappa_T$ all belong to $\{\mathbf{u}\}^\bot$, so that \eqref{eq:macnorm_nD_actionOn-1} makes sense. 
\end{remark}

\begin{example}[Reduced pair for the function $\chi(v,\kappa) = P_{\mathbf{u}^\bot} v$]~

\noindent
The function $\chi(v,\kappa) = P_{\mathbf{u}^\bot} v$ satisfies the invariance relation \eqref{eq:macnorm_gract_zeta_invar}. Its reduced pair is: 
\begin{equation} 
\chi_w (\theta, \kappa_\parallel, \kappa_\bot) = \sin \theta, \qquad \chi_w (\theta, \kappa_\parallel, \kappa_\bot) = 0. 
\label{ex:macnorm_nD_chi_reduced}
\end{equation}
\end{example}

\subsection{Variational formulation for the reduced GCI pair}
\label{subsec:macnorm_nD_varform_redGCI}

Our goal is now to derive a variational formulation for the reduced GCI pair from the variational formulation for $\vec{\zeta}_\mathbf{u}$ involving invariant spaces \eqref{eq:macnorm_GCI_gene_varform_inv}. We will have to distinguish two cases: the generic case $n \geq 4$ and the special case $n=3$. In the generic case, $(\zeta_w,\zeta_r)$ are functions of $(\theta, \kappa_\parallel, \kappa_\bot) \in {\mathcal N} =: (0,\pi) \times {\mathbb R} \times (0,\infty)$. So, the variational formulation for the reduced GCI pair will involve appropriate function or distribution spaces over ${\mathcal N}$ which we first introduce. Then, we establish this variational formulation and eventually derive the PDE and boundary conditions satisfied by the reduced GCI pair. The case $n=3$ will require a similar yet slightly different treatment which is summarized at the end.

\subsubsection{Spaces of functions on ${\mathcal N}$ and properties}
\label{subsubsec:macnorm_nD_weighted_subspaces}

In this subsection, we assume that $n \geq 4$. We start with the definition of appropriate function spaces on ${\mathcal N}$ and then, establish their link with spaces of invariant functions on $M$.

\bigskip
\noindent
(i) For $p$, $q \in [0,\infty)$, we define the space
$$ L^2_{p,q} = \Big\{ \psi: {\mathcal N} \to {\mathbb R}, \, \,  \big| \, \,   \| \psi \|_{L^2_{p,q}}^2 = : \int_{{\mathcal N}} |\psi|^2 \, \mathbf{M}_\mathbf{u}^\natural \, \kappa_\bot^q \, \sin^p \theta \, d \theta \, d \kappa_\parallel \, d\kappa_\bot < \infty \Big\}, $$
endowed with the Hilbert norm $\| \cdot \|_{L^2_{p,q}}$ and associated inner-product $( \hspace{-0.8mm} ( \cdot , \cdot) \hspace{-0.8mm} )_{L^2_{p,q}}$. We recall that $\mathbf{M}_\mathbf{u}^\natural$ is given by \eqref{eq:macnorm_Munat}. We denote by ${\mathbb L} = L^2_{2(n-2),n-3}$ and the associated norm and inner product by~$\| \cdot \|_{{\mathbb L}}$ and $( \hspace{-0.8mm} ( \cdot , \cdot) \hspace{-0.8mm} )_{{\mathbb L}}$ respectively.

\medskip
\noindent
(ii) We define the space 
$$ {\mathbb H} = \Big\{ \psi \in {\mathbb L}  \, \, \textrm{such that} \, \,
\| \psi \|_{{\mathbb H}}^2 =: \| \psi \|_{{\mathbb L}}^2 + \| \partial_{\kappa_\parallel} \psi \|_{{\mathbb L}}^2 + \| \partial_{\kappa_\bot} \psi \|_{L^2_{2(n-3),n-3}}^2 < \infty \big\}, $$
endowed with the Hilbert norm $\| \cdot \|_{{\mathbb H}}$ and associated inner-product $( \hspace{-0.8mm} ( \cdot , \cdot) \hspace{-0.8mm} )_{{\mathbb H}}$. We define the semi-norm 
$$ \| \psi \|_{\accentset{\textrm{\textbullet}}{{\mathbb H}}}^2 = \| \partial_{\kappa_\parallel} \psi \|_{{\mathbb L}}^2 + \| \partial_{\kappa_\bot} \psi \|_{L^2_{2(n-3),n-3}}^2, $$
so that 
$$  \| \psi \|_{{\mathbb H}}^2  =  \| \psi \|_{{\mathbb L}}^2  +  \| \psi \|_{\accentset{\textrm{\textbullet}}{{\mathbb H}}}^2. $$
We denote by  $( \hspace{-0.8mm} ( \cdot , \cdot) \hspace{-0.8mm} )_{\accentset{\textrm{\textbullet}}{{\mathbb H}}}$ the associated bilinear form.

\medskip
\noindent
We denote by ${\mathcal D}({\mathcal N})$ the space of $C^\infty$ functions with compact support on ${\mathcal N}$. We note the following important property, whose proof uses classical functional analytic arguments and is omitted: 

\begin{lemma}[Density of ${\mathcal D}({\mathcal N})$ in ${\mathbb H}$ when $n \geq 4$]~

\noindent
Suppose $n \geq 4$. Then, ${\mathcal D}({\mathcal N})$ is dense in ${\mathbb H}$. 
\label{lem:macnorm_nD_D(N)_densein_H1pq}
\end{lemma}

\begin{remark}[About Lemma \ref{lem:macnorm_nD_D(N)_densein_H1pq}]~

\noindent
The domain ${\mathcal N}$ has a boundary at $\kappa_\bot=0$ and the definition of ${\mathbb H}$ involves derivatives in $\kappa_\bot$. If $n \geq 4$ the inner product in the definition of ${\mathbb H}$ involves the weight $\kappa_\bot^{n-3}$ with an exponent $n-3 \geq 1$. Because of this weight, we can prove that that ${\mathcal D}({\mathcal N})$ is dense in ${\mathbb H}$. If $n=3$, this weight disappears and  ${\mathbb H}$ is akin to a classical Sobolev space~$H^1$ (i.e. without any weight in the inner-product) with respect to the $\kappa_\bot$ variable. Thus, in this case, ${\mathcal D}({\mathcal N})$ is not dense in~${\mathbb H}$. Whether ${\mathcal D}({\mathcal N})$ is dense in ${\mathbb H}$ or not impacts the structure of its dual space. Indeed, dual spaces of classical $H^1$ spaces cannot be identified with spaces of distributions when the domain involves boundaries. By contrast, duals of $H^1_0$-type spaces (which are the closures of ${\mathcal D}$ for the Sobolev norm) can (see \cite[Sections 3.7-3.14]{adams2003sobolev}). In the forthcoming sections, we will use the dual space of ${\mathbb H}$ and it will be of critical importance that it is a space of distributions, hence the restriction to $n \geq 4$. The case $n=3$ will require additional developments to circumvent the non-density of ${\mathcal D}({\mathcal N})$ in ${\mathbb H}$. There is no such problem with the other variables: the domain in $k_\parallel$ is ${\mathbb R}$ and has no boundary; the domain in $\theta$ is $(0, \pi)$ and has boundaries but the definition of ${\mathbb H}$ does not involve any derivative in $\theta$. 
\label{rem:macnorm_nD_D(N)_densein_H1pq}
\end{remark}

\medskip
\noindent
(iii) We define the space ${\mathbb K} = {\mathbb H} \cap  L^2_{2(n-3),n-5}$, endowed with the norm  
$$\| \psi \|_{{\mathbb K}}^2 = \| \psi \|_{{\mathbb H}}^2 + (n-3) \,  \| \psi \|_{L^2_{2(n-3),n-5}}^2. $$ 
The associated inner-product is denoted by $( \hspace{-0.8mm} ( \cdot , \cdot) \hspace{-0.8mm} )_{{\mathbb K}}$. We define the semi-norm 
$$ \| \psi \|_{\accentset{\textrm{\textbullet}}{{\mathbb K}}}^2 = \|  \psi \|_{\accentset{\textrm{\textbullet}}{{\mathbb H}}}^2 + (n-3) \, \| \psi \|_{L^2_{2(n-3),n-5}}^2, $$
so that 
$$  \| \psi \|_{{\mathbb K}}^2  =  \| \psi \|_{{\mathbb L}}^2  +  \| \psi \|_{\accentset{\textrm{\textbullet}}{{\mathbb K}}}^2 . $$
We denote by  $( \hspace{-0.8mm} ( \cdot , \cdot) \hspace{-0.8mm} )_{\accentset{\textrm{\textbullet}}{{\mathbb K}}}$ the associated bilinear form. Like for ${\mathbb H}$, we have the 

\begin{lemma}[Density of ${\mathcal D}({\mathcal N})$ in ${\mathbb K}$ when $n \geq 4$]~

\noindent
Suppose $n \geq 4$. Then, ${\mathcal D}({\mathcal N})$ is dense in ${\mathbb K}$. 
\label{lem:macnorm_nD_D(N)_densein_K1pq}
\end{lemma}

\medskip
\noindent
(iv) We define
$$ {\mathcal L}^2({\mathcal L}^2) = {\mathbb L}^2, $$
and denote its elements by $(\psi_w, \psi_r)$. It is endowed with the norm
$$ \| (\psi_w, \psi_r) \|_{{\mathcal L}^2({\mathcal L}^2)}^2 = \| \psi_w \|_{{\mathbb L}}^2 + \| \psi_r \|_{{\mathbb L}}^2. $$

\medskip
\noindent
(v) We define 
$$ {\mathcal L}^2({\mathcal H}^1) =  {\mathbb H} \times {\mathbb K}$$
endowed with the norm 
$$\| (\psi_w, \psi_r) \|_{{\mathcal L}^2({\mathcal H}^1)}^2 = \| \psi_w \|_{{\mathbb H}}^2 + \| \psi_r \|_{{\mathbb K}}^2. $$
and we have 
$$ \| (\psi_w, \psi_r) \|_{{\mathcal L}^2({\mathcal H}^1)}^2 = \| (\psi_w, \psi_r) \|_{{\mathcal L}^2({\mathcal L}^2)}^2 + \| \psi_w \|_{\accentset{\textrm{\textbullet}}{{\mathbb H}}}^2 + \| \psi_r \|_{\accentset{\textrm{\textbullet}}{{\mathbb K}}}^2. $$
In view of Lemmas \ref{lem:macnorm_nD_D(N)_densein_H1pq} and \ref{lem:macnorm_nD_D(N)_densein_K1pq}, since $n \geq 4$, the space ${\mathcal D}({\mathcal N})^2$ is dense in ${\mathcal L}^2({\mathcal H}^1)$.

\bigskip
We now establish the link between invariant function spaces on $M$ and function spaces on~${\mathcal N}$.

\begin{proposition}[Link between $L^2_v(L^2_\kappa)_\mathrm{inv}$ (resp. $L^2_v(H^1_\kappa)_\mathrm{inv}$) and ${\mathcal L}^2({\mathcal L}^2)$ (resp. ${\mathcal L}^2({\mathcal H}^1)$)]~

\noindent
(i) The function $\vec{\psi}$ belongs to $L^2_v(L^2_\kappa)_\mathrm{inv}$ if and only if $\vec{\psi}^\dagger$ (related to $\vec{\psi}$ through \eqref{eq:macnorm_dagger_def}) is written 
\begin{eqnarray}
&&\hspace{-1cm}
\vec{\psi}^\dagger (\theta, \kappa_\parallel, w, \kappa_T) = \psi_w(\theta, \kappa_\parallel, |\kappa_T|) \, w + \psi_r(\theta, \kappa_\parallel, |\kappa_T|) \, \frac{\kappa_T}{|\kappa_T|}, \nonumber \\
&&\hspace{6cm}
\textrm{for} \quad (\theta, \kappa_\parallel, w, \kappa_T) \in M_V \times M_T \, \, \textrm{a.e.} \, , 
\label{eq:macnorm_vecpsi_rot}
\end{eqnarray}
(where a.e. stands for ``almost everywhere'') and $(\psi_w,\psi_r) \in {\mathcal L}^2({\mathcal L}^2)$. Furthermore, defining $c$ as
\begin{equation} 
c = |{\mathbb S}^{n-2}| \, |{\mathbb S}^{n-3}|, 
\label{eq:macnorm_normal_coef}
\end{equation}
we have  
\begin{equation} 
\| \vec{\psi} \|_{L^2_v(L^2_\kappa)_\mathbf{u}}^2 = c \, \| (\psi_w, \psi_r) \|_{{\mathcal L}^2({\mathcal L}^2)}^2. 
\label{eq:macnorm_vecpsi_rot_L2L2}
\end{equation}

\noindent
(ii) The function $\vec{\psi}$ belongs to $L^2_v(H^1_\kappa)_\mathrm{inv}$ if and only if $\vec{\psi}^\dagger$ is given by \eqref{eq:macnorm_vecpsi_rot} and $(\psi_w,\psi_r) \in {\mathcal L}^2({\mathcal H}^1)$. Furthermore, 
\begin{equation} 
\| \vec{\psi} \|_{L^2_v(H^1_\kappa)_\mathbf{u}}^2 = c \, \| (\psi_w, \psi_r) \|_{{\mathcal L}^2({\mathcal H}^1)}^2, 
\label{eq:macnorm_vecpsi_rot_L2H1}
\end{equation}
and 
\begin{equation} 
\| \nabla_\kappa \vec{\psi} \|_{\vec{L}^2_v(L^2_\kappa)_\mathbf{u}}^2 = c \, \Big( 
\| \psi_w \|_{\accentset{\textrm{\textbullet}}{{\mathbb H}}}^2 + \| \psi_r \|_{\accentset{\textrm{\textbullet}}{{\mathbb K}}}^2 \Big).
\label{eq:macnorm_vecpsi_rot_L2H1_sn}
\end{equation}
\label{prop:macnorm_vecpsi_rot}
\end{proposition}

\noindent
\textbf{Proof.} (i) Let $\vec{\psi} \in {\mathcal D}_\mathrm{inv}(M)$. Then, Theorem \ref{thm:macnorm_nD_actionOn-1} applies and we get \eqref{eq:macnorm_vecpsi_rot}. We have 
$$ |\vec{\psi}^\dagger|^2 (\theta, \kappa_\parallel, w, \kappa_T) = \big( |\psi_w|^2 + |\psi_r|^2)(\theta, \kappa_\parallel, |\kappa_T|), \quad \textrm{for} \quad (\theta, \kappa_\parallel, w, \kappa_T) \in M_V \times M_T \, \, \textrm{a.e.} $$
Hence, using \eqref{eq:macnorm_nD_AuvolM2}, and the Fubini theorem for fiber bundles 
\eqref{eq:Fub_Fubthm}, we get 
\begin{eqnarray*}
&&\hspace{-1cm}
\int_M |\vec{\psi}|^2 \, \mathbf{M}_\mathbf{u} \, | \mathrm{Vol}_M | 
= \int_{M_V}  \Big( \int_{M_T} |\vec{\psi}^\dagger|^2 \, \mathbf{M}_\mathbf{u}^\dagger \, | \mathrm{Vol}_{\mathrm{M}_T} | \Big) \, \sin^{2(n-2)} \theta \, d \theta \, d \kappa_\parallel \\
&&\hspace{0.5cm}
= \int_{M_V}  \Big( \int_{M_T} \big( |\psi_w|^2 + |\psi_r|^2) \, \mathbf{M}_\mathbf{u}^\dagger \, | \mathrm{Vol}_{\mathrm{M}_T} | \Big) \, \sin^{2(n-2)} \theta \, d \theta \, d \kappa_\parallel \\
&&\hspace{0.5cm}
= \int_{M_V}  \Big[ \int_{{\mathbb S}^{n-2}} \Big( \int_{T_w {\mathbb S}^{n-2}} \big( |\psi_w|^2 + |\psi_r|^2) \, \mathbf{M}_\mathbf{u}^\dagger \, | \omega_w | \Big) \, |\mathrm{Vol}_B | \Big] \, \sin^{2(n-2)} \theta \, d \theta \, d \kappa_\parallel, 
\end{eqnarray*}
where $\omega_w$ is the retrenchment of $\omega$ on $T_w {\mathbb S}^{n-2}$, given by \eqref{eq:macnorm_nD_coef_c2_prf1}. Recalling the expression \eqref{eq:macnorm_Mudag} for $\mathbf{M}_\mathbf{u}^\dagger$, we see that the expression inside the inner integral only depends on $|\kappa_T|$. Hence, we can use spherical coordinates in $T_w {\mathbb S}^{n-2} \approx {\mathbb R}^{n-2}$ for the inner integral and find that the result is independent of $w$. Thus, the integral over ${\mathbb S}^{n-2}$ is just that of a constant and we get 
\begin{eqnarray*}
&&\hspace{-1cm}
\int_M |\vec{\psi}|^2 \, \mathbf{M}_\mathbf{u} \, | \mathrm{Vol}_M | = c \, \int_{M_V} \Big( \int_0^\infty \big( (|\psi_w|^2 + |\psi_r|^2) \, \mathbf{M}_\mathbf{u}^\natural \big)(\theta, \kappa_\parallel, \kappa_\bot) \, 
\\
&&\hspace{8cm}
\kappa_\bot^{n-3} \, d \kappa_\bot \Big) \, \sin^{2(n-2)} \theta \, d \theta \, d \kappa_\parallel, 
\end{eqnarray*}
with $\mathbf{M}_\mathbf{u}^\natural$ given by \eqref{eq:macnorm_Munat}. This is nothing but \eqref{eq:macnorm_vecpsi_rot_L2L2}. Now, \eqref{eq:macnorm_vecpsi_rot} and \eqref{eq:macnorm_vecpsi_rot_L2L2} extend to $\vec{\psi} \in L^2_v(L^2_\kappa)_\mathrm{inv}$ by density of ${\mathcal D}_\mathrm{inv}(M)$ in $L^2_v(L^2_\kappa)_\mathrm{inv}$. Conversely, if $(\psi_w,\psi_r) \in {\mathcal L}^2({\mathcal L}^2)$, then, $\vec{\psi}$ given in terms of $\vec{\psi}^\dagger$ by \eqref{eq:macnorm_vecpsi_rot} is invariant by $\mathrm{O}_{n-1}$ and belongs to $L^2_v(L^2_\kappa)_\mathbf{u}$ by virtue of \eqref{eq:macnorm_vecpsi_rot_L2L2}. Thus it is an element of $L^2_v(L^2_\kappa)_\mathrm{inv}$.

\medskip
\noindent
(ii) Let $\vec{\psi} \in {\mathcal D}_\mathrm{inv}(M)$. Then, the following formula is proven in Appendix \ref{sec:macnorm_nD_vectGCI_eqs}, Proposition \ref{prop:macnorm_nD_brack_zetu_express}(i)
\begin{equation} 
|(\nabla_\kappa \vec{\psi})^\dagger|^2 =  |\partial_{\kappa_\parallel} \psi_w |^2 + |\partial_{\kappa_\parallel} \psi_r |^2 + \frac{1}{\sin^2 \theta} \big( |\partial_{\kappa_\bot} \psi_w |^2 + |\partial_{\kappa_\bot} \psi_r |^2 \big) + \frac{n-3}{\sin^2 \theta \, |\kappa_T|^2} |\psi_r|^2. 
\label{eq:macnorm_nD_nakappsi_square}
\end{equation}
Now, inserting this into the integral over $M$, we can reproduce the calculation made in the proof of item (i) and we easily obtain \eqref{eq:macnorm_vecpsi_rot_L2H1} showing that $(\psi_w,\psi_r) \in {\mathcal L}^2({\mathcal H}^1)$. As in the proof of item~(i), the converse is obvious.  \endproof

\subsubsection{Spaces of distributions on ${\mathcal N}$ and properties}
\label{subsubsec:macnorm_nD_distrib_N}

In this subsection, we still assume that $n \geq 4$. We follow the same steps as in Subsection~\ref{subsubsec:macnorm_nD_weighted_subspaces}: we define distribution spaces on ${\mathcal N}$ and then, explore their link with spaces of invariant distributions on $M$ . 

\bigskip
\noindent
(i) We denote by ${\mathcal D}'({\mathcal N})$ the dual of ${\mathcal D}({\mathcal N})$, i.e. the space of distributions on ${\mathcal N}$. We define the duality bracket $( \hspace{-0.8mm} (T, \varphi ) \hspace{-0.8mm} )_{{\mathcal D}',{\mathcal D}}$ between $T \in {\mathcal D}'({\mathcal N})$ and $\varphi \in {\mathcal D}({\mathcal N})$ so as to extend the inner-product of~${\mathbb L}$. 

\smallskip
\noindent
(ii) We define the space ${\mathbb H}'$ as the dual of ${\mathbb H}$. By Lemma \ref{lem:macnorm_nD_D(N)_densein_H1pq}, ${\mathcal D}({\mathcal N})$ is dense in ${\mathbb H}$ and so, we can identify ${\mathbb H}'$ to the space of distributions $T \in  {\mathcal D}'({\mathcal N})$ such that 
$$ \| T \|_{{\mathbb H}'} =: \sup_{\varphi \in {\mathcal D}({\mathcal N}), \, \varphi \not = 0} \frac{\displaystyle |  \big( \hspace{-1.2mm} \big( T, \varphi \big) \hspace{-1.2mm} \big)_{{\mathcal D}',{\mathcal D}} |}{\displaystyle \| \varphi \|_{{\mathbb H}}} < \infty,  $$
endowed with $\| \cdot \|_{{\mathbb H}'}$ (see \cite[Sections 3.7-3.14]{adams2003sobolev}). 

\smallskip
\noindent
(iii) We define the space ${\mathbb K}'$ as the dual of ${\mathbb K}$. For the same reason as for ${\mathbb H}'$, we can identify ${\mathbb K}'$ with the space of distribution $T \in  {\mathcal D}'({\mathcal N})$ such that 
$$ \| T \|_{{\mathbb K}'} =: \sup_{\varphi \in {\mathcal D}({\mathcal N}), \, \varphi \not = 0} \frac{\displaystyle |  \big( \hspace{-1.2mm} \big( T, \varphi \big) \hspace{-1.2mm} \big)_{{\mathcal D}',{\mathcal D}} |}{\displaystyle \| \varphi \|_{{\mathbb K}}} < \infty,  $$
endowed with $\| \cdot \|_{{\mathbb K}'}$. 

\smallskip
\noindent
(iv) Let $\vec{T} \in {\mathcal D}'_\mathrm{inv}(M)$. We define two distributions $T_w$ and $T_r$ in ${\mathcal D}'({\mathcal N})$ such that for all $\varphi$ in ${\mathcal D}({\mathcal N})$, 
\begin{eqnarray} 
\big( \hspace{-1.2mm} \big( T_w, \varphi \big) \hspace{-1.2mm} \big)_{{\mathcal D}',{\mathcal D}} &=& c^{-1} \big( \hspace{-1.2mm} \big( \vec{T}, (w \varphi) \circ \mathbf{A}_\mathbf{u}^{-1} \big) \hspace{-1.2mm} \big)_{{\mathcal D}'_\mathbf{u},{\mathcal D}_\mathbf{u}}\, , \label{eq:macnorm_Tw_distrib} \\
\big( \hspace{-1.2mm} \big( T_r, \varphi \big) \hspace{-1.2mm} \big)_{{\mathcal D}',{\mathcal D}} &=& c^{-1} \big( \hspace{-1.2mm} \big( \vec{T}, (\frac{\kappa_T}{|\kappa_T|} \varphi) \circ \mathbf{A}_\mathbf{u}^{-1}  \big) \hspace{-1.2mm} \big)_{{\mathcal D}'_\mathbf{u},{\mathcal D}_\mathbf{u}}\, , \label{eq:macnorm_Tr_distrib}
\end{eqnarray}
where we refer to Section \ref{subsec:macnorm_vectorGCI} for the definition of $\mathbf{A}_\mathbf{u}$, and to Eq. \eqref{eq:macnorm_normal_coef} for that of $c$. At the right-hand sides of these formulas, $\varphi$ is considered as a function on $M_V \times M_T$ rather than just a function on ${\mathcal N}$. We note that since $\varphi$ has compact support on ${\mathcal N}$, $\kappa_T$ is never zero on the support of $\varphi$, so that $(w \varphi)\circ \mathbf{A}_\mathbf{u}^{-1}$ and $((\kappa_T/|\kappa_T|) \varphi) \circ \mathbf{A}_\mathbf{u}^{-1}$ are elements of ${\mathcal D}_\mathrm{inv}(M)$ and the definitions of $T_w$ and $T_r$ make sense. Furthermore, if $\vec{T} \in {\mathcal D}_\mathrm{inv}(M)$, it can be written $ \vec{T}^\dagger = T_w w + T_r (\kappa_T/|\kappa_T|)$ and, thanks to the definition of the duality bracket for ${\mathcal D}'({\mathcal N})$, this definition of $(T_w,T_r)$ coincides with \eqref{eq:macnorm_Tw_distrib}, \eqref{eq:macnorm_Tr_distrib}. Now, for $\vec{\varphi} \in {\mathcal D}_\mathrm{inv}(M)$, we have the characterization $ \vec{\varphi}^\dagger = \varphi_w w + \varphi_r (\kappa_T/|\kappa_T|)$. So, we can write 
\begin{equation} 
\big( \hspace{-1.2mm} \big( \vec{T}, \vec{\varphi} \big) \hspace{-1.2mm} \big)_{{\mathcal D}'_\mathbf{u},{\mathcal D}_\mathbf{u}} = c \Big( \big( \hspace{-1.2mm} \big( T_w, \varphi_w \big) \hspace{-1.2mm} \big)_{{\mathcal D}',{\mathcal D}} + \big( \hspace{-1.2mm} \big( T_r, \varphi_r \big) \hspace{-1.2mm} \big)_{{\mathcal D}',{\mathcal D}} \Big). 
\label{eq:macnorm_distrib_invar}
\end{equation}

\noindent
(v) We define 
$${\mathcal L}^2({\mathcal H}^{-1}) = {\mathbb H}' \times {\mathbb K}'. $$
It is the dual of ${\mathcal L}^2({\mathcal H}^1)$. Since $n \geq 4$, the space ${\mathcal D}({\mathcal N})^2$ is dense in ${\mathcal L}^2({\mathcal H}^1)$ and consequently, ${\mathcal L}^2({\mathcal H}^{-1})$ can be identified with the space of pairs of distributions $(T_w,T_r)$ such that 
\begin{eqnarray} 
&&\hspace{-1.5cm}
\| (T_w,T_r) \|_{{\mathcal L}^2({\mathcal H}^{-1})} 
\\
&&\hspace{-1cm}
=: \sup_{(\varphi_w,\varphi_r) \in {\mathcal D}({\mathcal N})^2 \, (\varphi_w,\varphi_r) \not = 0} \frac{\displaystyle \big|  \big( \hspace{-1.2mm} \big( T_w, \varphi_w \big) \hspace{-1.2mm} \big)_{{\mathbb H}',{\mathbb H}}  + \big( \hspace{-1.2mm} \big( T_r, \varphi_r \big) \hspace{-1.2mm} \big)_{{\mathbb K}',{\mathbb K}} \big|}{\displaystyle \big(\| \varphi_w \|_{{\mathbb H}}^2 + \| \varphi_r \|_{{\mathbb K}}^2 \big)^{1/2} } \\
&&\hspace{-1cm}
= \big(\| T_w \|_{{\mathbb H}'}^2 + \| T_r \|_{{\mathbb K}'}^2 \big)^{1/2} < \infty. 
\label{eq:macnorm_invar_distrib_TwTrnorm}
\end{eqnarray}
The space ${\mathcal L}^2({\mathcal H}^{-1})$ is normed by $\| \cdot \|_{{\mathcal L}^2({\mathcal H}^{-1})}$. We have 
\begin{equation} 
\vec{T} \in L^2_v(H^{-1}_\kappa)_\mathrm{inv} \quad \Longleftrightarrow \quad 
(T_w, T_r) \in {\mathcal L}^2({\mathcal H}^{-1}), 
\label{eq:macnorm_invar_distrib_charact}
\end{equation}
and 
\begin{equation} 
\| \vec{T} \|_{L^2_v(H^{-1}_\kappa)_\mathbf{u}}^2 =  c \, \| (T_w,T_r) \|_{{\mathcal L}^2({\mathcal H}^{-1})}^2, 
\label{eq:macnorm_norm_L2H-1_red}
\end{equation}

\medskip
\noindent
(vi) We refer to \eqref{eq:macnorm_nD_Uoperator_def} for the definition of operator ${\mathcal U}_1$. The space ${\mathbb Y}_{\mathcal N}$ is defined by 
\begin{eqnarray*}
&&\hspace{-1cm}
{\mathbb Y}_{\mathcal N} = \big\{ (\psi_w, \psi_r) \in {\mathcal L}^2({\mathcal H}^1) \quad \textrm{such that} 
\quad \big( {\mathcal U}_1 \psi_w - \kappa_\bot \psi_r, {\mathcal U}_1 \psi_r + \kappa_\bot \psi_w \big) \in {\mathcal L}^2({\mathcal H}^{-1}) \big\}, 
\end{eqnarray*}
where ${\mathcal U}_1 \psi_w$ and ${\mathcal U}_1 \psi_r$ are defined in the distributional sense (see \eqref{eq:macnorm_U1_skewadj_distrib} below), endowed with the norm
\begin{equation} 
\| (\psi_w, \psi_r) \|_{{\mathbb Y}_{\mathcal N}}^2 = \| (\psi_w, \psi_r) \|_{{\mathcal L}^2({\mathcal H}^1)}^2 + \| \big( {\mathcal U}_1 \psi_w - \kappa_\bot \psi_r, {\mathcal U}_1 \psi_r + \kappa_\bot \psi_w \big) \|_{{\mathcal L}^2({\mathcal H}^{-1})}^2. 
\label{eq:macnorm_norm_in_YN}
\end{equation}

\bigskip
We now expore the Link between ${\mathbb Y}_\mathrm{inv}$ and ${\mathbb Y}_{\mathcal N}$. We first explore some properties of the operator ${\mathcal U}_1$. 

\begin{lemma}[${\mathcal U}_1$ is formally skew-adjoint]~

\noindent
(i) For any $\varphi$, $\psi$ in ${\mathcal D}({\mathcal N})$, we have 
\begin{equation}
\big( \hspace{-1.2mm} \big( {\mathcal U}_1 \varphi, \psi \big) \hspace{-1.2mm} \big)_{{\mathbb L}} = - \big( \hspace{-1.2mm} \big( \varphi, {\mathcal U}_1 \psi \big) \hspace{-1.2mm} \big)_{{\mathbb L}}. 
\label{eq:macnorm_U1_skewadj}
\end{equation}

\noindent
(ii) Let $T \in {\mathcal D}'({\mathcal N})$. Then the operator ${\mathcal U}_1$ applied to $T$ in the distributional sense satisfies 
\begin{equation} 
\big( \hspace{-1.2mm} \big( {\mathcal U}_1 T, \varphi \big) \hspace{-1.2mm} \big)_{{\mathcal D}',{\mathcal D}} = - \big( \hspace{-1.2mm} \big( T, {\mathcal U}_1 \varphi \big) \hspace{-1.2mm} \big)_{{\mathcal D}',{\mathcal D}}, \quad \forall \varphi \in {\mathcal D}({\mathcal N}).
\label{eq:macnorm_U1_skewadj_distrib}
\end{equation}

\noindent
(iii) For $\vec{T} \in {\mathcal D}'_\mathrm{inv}(M)$, we have 
\begin{equation}
\big\{ \vec{T}, \mathbf{H}_\mathbf{u} \big\}^\dagger = \Big( w \big( {\mathcal U}_1 T_w - |\kappa_T| T_r \big) + \frac{\kappa_T}{|\kappa_T|} \big( {\mathcal U}_1 T_r + |\kappa_T| T_w \big) \Big), 
\label{eq:macnorm_bracket_distrib_red}
\end{equation}
where the operator ${\mathcal U}_1$ is applied to $T_w$ and $T_r$ in the distributional sense, by virtue of~\eqref{eq:macnorm_U1_skewadj_distrib}

\label{lem:macnorm_U1_skewadj}
\end{lemma}

\noindent
\textbf{Proof.} (ii) is an immediate consequence of (i). To show (i), we first remark that for $\vec{\psi} \in {\mathcal D}_\mathrm{inv}$, we have (see Appendix \ref{sec:macnorm_nD_vectGCI_eqs}, Proposition \ref{prop:macnorm_nD_brack_zetu_express}(ii))
\begin{equation}
\big\{ \vec{\psi} , \mathbf{H}_\mathbf{u} \big\}^\dagger = w \, \big( {\mathcal U}_1 \psi_w - |\kappa_T| \, \psi_r \big) + \frac{\kappa_{T}}{|\kappa_T|} \, \big( {\mathcal U}_1 \psi_r + |\kappa_T|  \, \psi_w \big), 
\label{eq:macnorm_nD_brack_zetu_express}
\end{equation}
Now, let $\varphi$, $\psi$ in ${\mathcal D}({\mathcal N})$ and define $\vec{\varphi} = (w \varphi) \circ \mathbf{A}_\mathbf{u}^{-1}$ and define $\vec{\psi}$ in a similar way. Then, $\vec{\varphi}$, $\vec{\psi}$ belong to ${\mathcal D}_\mathrm{inv}(M)$ and, thanks to \eqref{eq:macnorm_nD_brack_zetu_express}, we have 
$$ \big\{ \vec{\varphi}, \mathbf{H}_\mathbf{u} \big\} = \big( w {\mathcal U}_1 \varphi + \kappa_T \varphi \big) \circ \mathbf{A}_\mathbf{u}^{-1} , $$
and similarly for $\{ \vec{\psi}, \mathbf{H}_\mathbf{u} \}$. 
Now, from \eqref{eq:macnorm_bracket_dual_u} and given the regularity of $\varphi$ and $\psi$, we deduce that 
$$ \big( \hspace{-1.2mm} \big( \big( w {\mathcal U}_1 \varphi + \kappa_T \varphi \big) \circ \mathbf{A}_\mathbf{u}^{-1}, (w \psi) \circ \mathbf{A}_\mathbf{u}^{-1} \big) \hspace{-1.2mm} \big)_{L^2_v(L^2_\kappa)} = 
- \big( \hspace{-1.2mm} \big( (w \varphi) \circ \mathbf{A}_\mathbf{u}^{-1}, \big( w {\mathcal U}_1 \psi + \kappa_T \psi \big) \circ \mathbf{A}_\mathbf{u}^{-1} \big) \hspace{-1.2mm} \big)_{L^2_v(L^2_\kappa)}, $$
which leads to \eqref{eq:macnorm_U1_skewadj} since $w$ and $\kappa_T$ are orthogonal. 

\medskip
\noindent
(iii) Let $\vec{T} \in {\mathcal D}'_\mathrm{inv}(M)$ and $\vec{\psi} \in {\mathcal D}_\mathrm{inv}(M)$. Thanks to \eqref{eq:macnorm_bracket_dual_u}, \eqref{eq:macnorm_distrib_invar}, \eqref{eq:macnorm_nD_brack_zetu_express} and \eqref{eq:macnorm_U1_skewadj_distrib}, we have 
\begin{eqnarray*} 
&&\hspace{-1cm}
\big( \hspace{-1.2mm} \big( \{ \vec{T}, \mathbf{H}_{\mathbf{u}} \}, \vec{\psi} \big) \hspace{-1.2mm} \big)_{{\mathcal D}_\mathbf{u}',{\mathcal D}_\mathbf{u}} = - 
\big( \hspace{-1.2mm} \big( \vec{T}, \{ \vec{\psi}, \mathbf{H}_{\mathbf{u}} \} \big) \hspace{-1.2mm} \big)_{{\mathcal D}_\mathbf{u}',{\mathcal D}_\mathbf{u}} \\
&&\hspace{1cm}
= -c \Big[ \big( \hspace{-1.2mm} \big( T_w, {\mathcal U}_1 \psi_w - \kappa_\bot \psi_r \big) \hspace{-1.2mm} \big)_{{\mathcal D}',{\mathcal D}} + \big( \hspace{-1.2mm} \big( T_r, {\mathcal U}_1 \psi_r + \kappa_\bot \psi_w \big) \hspace{-1.2mm} \big)_{{\mathcal D}',{\mathcal D}} \Big] \\
&&\hspace{1cm}
= c \Big[ \big( \hspace{-1.2mm} \big( {\mathcal U}_1 T_w - \kappa_\bot T_r, \psi_w \big) \hspace{-1.2mm} \big)_{{\mathcal D}',{\mathcal D}} + \big( \hspace{-1.2mm} \big( {\mathcal U}_1 T_r + \kappa_\bot T_w, \psi_r \big) \hspace{-1.2mm} \big)_{{\mathcal D}',{\mathcal D}} \Big]. 
\end{eqnarray*}
Thanks to \eqref{eq:macnorm_distrib_invar}, this last expression is equivalent to \eqref{eq:macnorm_bracket_distrib_red}. \endproof

\begin{remark}[Direct computation]~

\noindent
Identity \eqref{eq:macnorm_U1_skewadj} can be cross-checked by a straightforward, albeit tedious, direct computation. 
\end{remark}

\medskip
With this lemma, we can show the 

\begin{proposition}[Link between ${\mathbb Y}_\mathrm{inv}$ and ${\mathbb Y}_{\mathcal N}$]~

\noindent
The function $\vec{\psi}$ belongs to ${\mathbb Y}_\mathrm{inv}$ if and only if $\vec{\psi}$ is written according to \eqref{eq:macnorm_vecpsi_rot} and $(\psi_w, \psi_r)$ belongs to ${\mathbb Y}_{\mathcal N}$. Furthermore, 
\begin{equation}
\| \vec{\psi} \|_{{\mathbb Y}_\mathrm{inv}}^2 = c \, \| (\psi_w, \psi_r) \|_{{\mathbb Y}_{\mathcal N}}^2. 
\label{eq:macnorm_vecpsi_rot_Y}
\end{equation}
\label{prop:macnorm_vecpsi_rot_Y}
\end{proposition}

\noindent
\textbf{Proof.} Let $\vec{\psi} \in {\mathbb Y}_\mathrm{inv}$. By Proposition \ref{prop:macnorm_vecpsi_rot}, $\vec{\psi}$ is given by \eqref{eq:macnorm_vecpsi_rot} with $(\psi_w, \psi_r) \in {\mathcal L}^2({\mathcal H}^1)$. Now, thanks to \eqref{eq:macnorm_invar_distrib_charact} and to \eqref{eq:macnorm_bracket_distrib_red}, $\{ \vec{\psi}, \mathbf{H}_\mathbf{u} \} \in L^2_v(H^{-1}_\kappa)$ if and only if $( {\mathcal U}_1 \psi_w - \kappa_\bot \psi_r, {\mathcal U}_1 \psi_r + \kappa_\bot \psi_w ) \in {\mathcal L}^2({\mathcal H}^{-1})$, which shows that $(\psi_w, \psi_r) \in {\mathbb Y}_{\mathcal N}$. We can move backwards this chain of implications, which shows the claimed equivalence. The equality \eqref{eq:macnorm_vecpsi_rot_Y} follows from \eqref{eq:macnorm_norm_in_YN}, \eqref{eq:macnorm_vecpsi_rot_L2L2} and \eqref{eq:macnorm_norm_L2H-1_red}. \endproof

\subsubsection{Expression of the variational formulation for the reduced GCI pair}
\label{subsubsubsec:macnorm_nD_red_GCI_varform}

Again, in this subsection, we assume that $n \geq 4$. We establish a variational formulation solved by the reduced GCI pair and derive the corresponding system of PDE \eqref{eq:macnorm_nD_vectGCI_wcomp}, \eqref{eq:macnorm_nD_vectGCI_kapTcomp} in strong form.

\begin{theorem}[Variational formulation]~

\noindent
Suppose $n \geq 4$. Let $\vec{\chi} \in L^2_v(L^2_\kappa)_\mathrm{inv}$ and define $(\chi_w, \chi_r) \in {\mathcal L}^2({\mathcal L}^2)$ as in \eqref{eq:macnorm_vecpsi_rot}. Then, $\vec{\psi}$ is the unique solution of the variational formulation \eqref{eq:macnorm_GCI_gene_varform_inv} if and only if $(\psi_w, \psi_r)$ linked to $\vec{\psi}$ by \eqref{eq:macnorm_vecpsi_rot} is the unique solution of the variational formulation
\begin{equation}
\mbox{} \hspace{-0.5cm} \left\{ \begin{array}{l}
(\psi_w, \psi_r) \in {\mathbb Y}_{\mathcal N}, \\
\textrm{ } \\
- \big( \hspace{-1.2mm} \big({\mathcal U}_1 \psi_w  - \kappa_\bot \psi_r \, ,  \, \varphi_w\big) \hspace{-1.2mm} \big)_{{\mathbb H}',{\mathbb H}} 
-   \big( \hspace{-1.2mm} \big( {\mathcal U}_1 \psi_r + \kappa_\bot \psi_w \, , \, \varphi_r \big) \hspace{-1.2mm} \big)_{{\mathbb K}',{\mathbb K}} 
+ D \big( \hspace{-1.2mm} \big( \psi_w \ , \, \varphi_w \big) \hspace{-1.2mm} \big)_{\accentset{\textrm{\textbullet}}{{\mathbb H}}} \\
\textrm{ } \\
\hspace{1.5cm} 
+ D \big( \hspace{-1.2mm} \big( \psi_r \, , \, \varphi_r \big) \hspace{-1.2mm} \big)_{\accentset{\textrm{\textbullet}}{{\mathbb K}}} = - \big( \hspace{-1.2mm} \big( \chi_w \, , \, \varphi_w \big) \hspace{-1.2mm} \big)_{{\mathbb L}} - \big( \hspace{-1.2mm} \big( \chi_r \, , \, \varphi_r \big) \hspace{-1.2mm} \big)_{{\mathbb L}},  
\qquad \forall (\varphi_w, \varphi_r) \in {\mathbb Y}_{\mathcal N}. \end{array} \right. 
\label{eq:macnorm_redGCI_varform}
\end{equation}
\label{thm:macnorm_redGCI_varform}
\end{theorem}

\noindent
\textbf{Proof.} This is a direct consequence of \eqref{eq:macnorm_bracket_distrib_red}, \eqref{eq:macnorm_distrib_invar} and the counterparts for bilinear forms of \eqref{eq:macnorm_vecpsi_rot_L2H1_sn} and \eqref{eq:macnorm_vecpsi_rot_L2L2}. The variational formulation \eqref{eq:macnorm_redGCI_varform} has a unique solution since \eqref{eq:macnorm_GCI_gene_varform_inv} has a unique solution. \endproof

\begin{remark}[About the variational formulation \eqref{eq:macnorm_redGCI_varform}]~

\noindent
The condition $(\varphi_w, \varphi_r) \in {\mathbb Y}_{{\mathcal N}}$ does not imply that $(\varphi_w, 0)$ nor $(0, \varphi_r)$ separately belong to ${\mathbb Y}_{{\mathcal N}}$. This is because $\varphi_w$ and $\varphi_r$ are entangled in the condition $( {\mathcal U}_1 \varphi_w - \kappa_\bot \varphi_r, {\mathcal U}_1 \varphi_r + \kappa_\bot \varphi_w) \in {\mathcal L}^2({\mathcal H}^{-1})$ of the definition of ${\mathbb Y}_{{\mathcal N}}$. Thus, if we separate \eqref{eq:macnorm_redGCI_varform} in two equations by taking either $\varphi_r = 0$ or $\varphi_w = 0$, we obtain a different variational formulation, which may be ill-posed.  
\end{remark}

From Theorem \ref{thm:macnorm_redGCI_varform}, we immediately get the following 

\begin{corollary}[Variational formulation satisfied by the reduced GCI pair]~

\noindent
Let $(\zeta_w, \zeta_r)$ be the reduced GCI pair. Then, $(\zeta_w, \zeta_r)$ satisfies the variational formulation \eqref{eq:macnorm_redGCI_varform}, with $(\psi_w, \psi_r)$ replaced by $(\zeta_w, \zeta_r)$ and $(\chi_w, \chi_r)$ given by \eqref{ex:macnorm_nD_chi_reduced}. \label{cor:macnorm_redGCI_varform}
\end{corollary}

In the following theorem, we provide a strong formulation of the problem associated with the variational formulation \eqref{eq:macnorm_redGCI_varform}. It involves the derivation of both the PDE and its boundary conditions. The latter are obtained formally assuming that the solution of the PDE is smooth, a property that is still open. 

\begin{theorem}[PDE system associated with the variational formulation \eqref{eq:macnorm_redGCI_varform} and boundary conditions]~

\noindent
We assume that $n \geq 4$. 

\smallskip
\noindent
(i) Let $(\zeta_w, \zeta_r) \in \mathbb Y_{{\mathcal N}}$ be the unique solution to the variational formulation \eqref{eq:macnorm_redGCI_varform}, with $(\psi_w, \psi_r)$ replaced by $(\zeta_w, \zeta_r)$ and $(\chi_w, \chi_r)$ given by \eqref{ex:macnorm_nD_chi_reduced}. Then, $(\zeta_w, \zeta_r)$ satisfies System \eqref{eq:macnorm_nD_vectGCI_wcomp}, \eqref{eq:macnorm_nD_vectGCI_kapTcomp} in the distributional sense. 

\smallskip
\noindent
(ii) (formal) Supposing that $(\zeta_w, \zeta_r)$ is smooth. Then the boundary conditions associated to the variational formulation \eqref{eq:macnorm_redGCI_varform} are: 
\begin{itemize}
\item[-] no boundary condition at $\theta = 0$ nor $\theta = \pi$ on neither $\zeta_w$ nor $\zeta_r$;
\item[-] dimension $n=4$: no boundary condition on $\zeta_w$ and homogeneous Dirichlet condition on $\zeta_r$ at $\kappa_\bot = 0$:
\begin{equation} 
 \zeta_r (\theta, \kappa_\parallel, 0) = 0; 
\label{eq:macnorm_nD_dirichlet}
\end{equation}
\item[-] dimension $n \geq 5$: no boundary condition at $\kappa_\bot=0$ on neither $\zeta_w$ nor $\zeta_r$. 
\end{itemize}
\label{thm:macnorm_nD_red_GCI_syst_and_bc}
\end{theorem}

\noindent
\textbf{Proof.} (i) Let $\varphi \in {\mathcal D}({\mathcal N})$. Then, $(\varphi,0) \in {\mathbb Y}_{{\mathcal N}}$ and $(0,\varphi) \in {\mathbb Y}_{{\mathcal N}}$. Thus, using these two pairs as test functions in \eqref{eq:macnorm_redGCI_varform} and interpreting the various bilinear forms as duality between distributions and test functions, we are led straightforwardly to \eqref{eq:macnorm_nD_vectGCI_wcomp}, \eqref{eq:macnorm_nD_vectGCI_kapTcomp}. 

\medskip
\noindent
(ii) Define ${\mathcal D}(\bar{\mathcal N})$ as the set of $C^\infty$ functions with compact support on an open set containing the closure $\bar{\mathcal N}$ of ${\mathcal N}$. We now investigate under which conditions on $\varphi \in {\mathcal D}(\bar{\mathcal N})$ we have $(\varphi, 0) \in {\mathbb Y}_{{\mathcal N}}$ or $(0,\varphi) \in {\mathbb Y}_{{\mathcal N}}$. We recall that we assume $n \geq 4$ and we will treat the case $n=3$ later on. 
\begin{itemize}
\item Since $\sin^2 \theta \leq 1$ and $\varphi$ is compactly supported in $\bar{\mathcal N}$, we rather straightforwardly have $\varphi$ and $\partial_{\kappa_\parallel} \varphi \in {\mathbb L}$ and $\partial_{\kappa_\bot} \varphi \in L^2_{2(n-3),n-3}$ , so that $\varphi \in {\mathbb H}$ for all values of $n \geq 4$. On the other hand, for $n=4$, the function $\kappa_\bot^{n-5} = 1/\kappa_\bot$ is not integrable and $\varphi \not \in L^2_{2(n-3),n-5}$. By contrast, if $n \geq 5$, $\varphi \in L^2_{2(n-3),n-5}$ and consequently, $\varphi \in {\mathbb K}$. 
\item Now, we clearly have ${\mathcal U}_1 \varphi = \tilde{\varphi}/\sin \theta$ with $\tilde \varphi \in {\mathcal D}(\bar{{\mathcal N}})$. So, $ \| {\mathcal U}_1  \varphi \|_{{\mathbb L}} = \| \tilde \varphi \|_{L^2_{2(n-3),n-3}} < \infty$ 
if $n \geq 4$. It follows that ${\mathcal U}_1 \varphi$ belongs to ${\mathbb L}$ and consequently, to ${\mathbb H}'$ and ${\mathbb K}'$. Similarly, since $\varphi$ is compactly spported, $\kappa_\bot \varphi$ belongs to ${\mathbb L}$ and so, to ${\mathbb H}'$ and ${\mathbb K}'$ as well. 
\end{itemize}
It follows that $(\varphi, 0) \in {\mathbb Y}_{{\mathcal N}}$ for $n \geq 4$ and that $(0,\varphi) \in {\mathbb Y}_{{\mathcal N}}$ for $n \geq 5$. 

Now, we take $\varphi \in {\mathcal D}(\bar{\mathcal N})$ and assume that the solution $(\varphi_w, \varphi_r)$ of the variational formulation~\eqref{eq:macnorm_redGCI_varform} is smooth. We use $(\varphi, 0) \in {\mathbb Y}_{{\mathcal N}}$ as a test function in \eqref{eq:macnorm_redGCI_varform} for all $n \geq 4$. By using Green's formula, and point (i) of the present lemma, we end up with
$$ \Big( \int_{(0,\pi) \times {\mathbb R}}  \big( \partial_{\kappa_\bot}\psi_w \, \varphi \, \mathbf{M}_\mathbf{u}^\dagger \,  \sin^{2(n-3)} \theta \, \kappa_\bot^{n-3} \big)  \, d \theta \, d \kappa_\parallel \Big)\big|_{\kappa_\bot = 0} = 0, $$
which is a tautology when $n \geq 4$. Performing the same computation using $(0,\varphi) \in {\mathbb Y}_{{\mathcal N}}$ as a test function in \eqref{eq:macnorm_redGCI_varform} for all $n \geq 5$, we get a tautology again. Now, in the case $n=4$, for $\varphi \in {\mathcal D}(\bar{\mathcal N})$ to be in ${\mathbb K}$, we need $\varphi(\theta, \kappa_\parallel, 0) = 0$ to compensate for the weight $1/\kappa_\bot$, which leads to the Dirichlet condition \eqref{eq:macnorm_nD_dirichlet}. This ends the proof. \endproof

\subsubsection{Case $n=3$}
\label{subsubsubsec:macnorm_nD_red_GCI_n=3}

In dimension $n=3$, we use spherical coordinates. Let $(e_1, e_2, \mathbf{u})$ be a direct orthonormal basis and let $(\theta, \varphi)$ be the spherical coordinates in this basis. We let $(v,e_\theta, e_\varphi)$ be the local spherical coordinates basis. Both $\theta$ and $e_\theta$ have been introduced in \eqref{eq:macnorm_V_def_2} and \eqref{eq:hydro_chgvar_ethetdef} respectively, and we have  
$$ w = \cos \varphi \, e_1 + \sin \varphi \, e_2, \qquad e_\varphi = - \sin \varphi e_1 + \cos \varphi \, e_2. $$
Because $T_w {\mathbb S}^{n-2}$ is one-dimensional and spanned by $e_\varphi$, we can write 
$$ \kappa_T = \bar{\kappa}_T \, e_\varphi, \quad \textrm{so that} \quad \bar{\kappa}_T = \left\{ \begin{array}{lll} |\kappa_T| & \textrm{if} & \bar{\kappa}_T \geq 0, \\  - |\kappa_T| & \textrm{if} & \bar{\kappa}_T \leq 0, \end{array} \right. $$
and $\bar{\kappa}_T \in {\mathbb R}$. 

For any function $\vec{\psi} \in {\mathcal S}(M, \{ \mathbf{u} \}^\bot)$ satisfying the invariance relation~\eqref{eq:macnorm_gract_zeta_invar}, we define the functions $\bar{\psi}_w$, $\bar{\psi}_r$: $\widetilde{{\mathcal N}} =: (0,\pi) \times {\mathbb R}^2 \to {\mathbb R}$, such that
\begin{equation} 
\bar{\psi}_w (\theta, \kappa_\parallel, \bar{\kappa}_T) = \psi_w (\theta, \kappa_\parallel, |\bar{\kappa}_T|) , \qquad \bar{\psi}_r (\theta, \kappa_\parallel, \bar{\kappa}_T) = \textrm{Sign}(\bar{\kappa}_T) \, \psi_r (\theta, \kappa_\parallel, |\bar{\kappa}_T|),  
\label{eq:macnorm_n=3_link_red_GCI}
\end{equation}
for all $(\theta, \kappa_\parallel, \bar{\kappa}_T) \in \widetilde{{\mathcal N}}$.
In view of \eqref{eq:macnorm_nD_actionOn-1}, any such function $\vec{\psi}$ can be written in the form 
\begin{equation} 
\vec{\psi}^\dagger(\theta, \kappa_\parallel, w, \kappa_T) = \bar{\psi}_w (\theta, \kappa_\parallel, \bar{\kappa}_T) \, w + \bar{\psi}_r (\theta, \kappa_\parallel, \bar{\kappa}_T) \, e_\varphi, \quad \forall \big( (\theta, \kappa_\parallel), (w, \kappa_T) \big) \in M_V \times M_T. 
\label{eq:macnorm_n=3_reduced_GCI}
\end{equation}
If $\vec{\psi} = \vec{\zeta}_\mathbf{u}$, the associated $(\bar{\zeta}_w, \bar{\zeta}_r)$ is the reduced GCI pair in dimension $n=3$. 

By analogy with the spaces defined in Subsections \ref{subsubsec:macnorm_nD_weighted_subspaces} and \ref{subsubsec:macnorm_nD_distrib_N}, we introduce the following spaces (notations that are a direct adaptation of those of the former sections are not repeated). 

\bigskip
\noindent
(i) For $p \in {\mathbb N}$, we define the Hilbert space:
$$ \widetilde{L^2_p} = \Big\{ \psi: \widetilde{{\mathcal N}} \to {\mathbb R}, \, \,  \big| \, \,   \| \psi \|_{\widetilde{L^2_p}}^2 = : \int_{\widetilde{{\mathcal N}}} |\psi|^2 \, \mathbf{M}_\mathbf{u}^\natural \, \sin^p \theta \, d \theta \, d \kappa_\parallel \, d\bar{\kappa}_T < \infty \Big\}. $$
We denote by $\widetilde{{\mathbb L}} = \widetilde{L^2_2}$.

\medskip
\noindent
(ii) We define
$$ \widetilde{{\mathbb H}} = \Big\{ \psi \in \widetilde{{\mathbb L}}  \, \, \big| \, \,  \| \psi \|_{\widetilde{{\mathbb H}}}^2 =: \| \psi \|_{\widetilde{{\mathbb L}}}^2 + \| \partial_{\kappa_\parallel} \psi \|_{\widetilde{{\mathbb L}}}^2 + \| \partial_{\kappa_\bot} \psi \|_{\widetilde{L^2_0}}^2 < \infty \big\}, $$
and recall that 
$$ \| \psi \|_{\widetilde{\accentset{\textrm{\textbullet}}{{\mathbb H}}}}^2 = \| \partial_{\kappa_\parallel} \psi \|_{\widetilde{{\mathbb L}}}^2 + \| \partial_{\kappa_\bot} \psi \|_{\widetilde{L^2_0}}^2.$$ 
Because $\widetilde{{\mathcal N}}$ has no boundary in the $\bar{\kappa}_T$ direction, the following density theorem holds true (this would not be the case if we used $\kappa_\bot$ instead of $\bar{\kappa}_T$, see Remark \ref{rem:macnorm_nD_D(N)_densein_H1pq}): 

\begin{lemma}[Density of ${\mathcal D}(\widetilde{{\mathcal N}})$ in $\widetilde{{\mathbb H}}$]~

\noindent
${\mathcal D} \big(\widetilde{{\mathcal N}}\big)$ is dense in $\widetilde{{\mathbb H}}$. 
\label{lem:macnorm_nD_D(N)_densein_H120}
\end{lemma}

\medskip
\noindent
(iii) We define $\widetilde{{\mathcal L}^2({\mathcal L}^2)} = \widetilde{{\mathbb L}}^2$ and denote its elements by $(\bar{\psi}_w, \bar{\psi}_r)$.

\medskip
\noindent
(iv) We define $\widetilde{{\mathcal L}^2({\mathcal H}^1)} = \widetilde{{\mathbb H}}^2$. Thanks to Lemma \ref{lem:macnorm_nD_D(N)_densein_H120}, ${\mathcal D}\big(\widetilde{{\mathcal N}}\big)^2$ is dense in $\widetilde{{\mathcal L}^2({\mathcal H}^1)}$. 

\begin{proposition}[Link between $L^2_v(L^2_\kappa)_\mathrm{inv}$ (resp. $L^2_v(H^1_\kappa)_\mathrm{inv}$) and $\widetilde{{\mathcal L}^2({\mathcal L}^2)}$ (resp. $\widetilde{{\mathcal L}^2({\mathcal H}^1)}$)]~

\noindent
We suppose $n = 3$. The function $\vec{\psi}$ belongs to $L^2_v(L^2_\kappa)_\mathrm{inv}$ (resp. $L^2_v(H^1_\kappa)_\mathrm{inv}$) if and only if $\vec{\psi}^\dagger$ is written according to \eqref{eq:macnorm_n=3_reduced_GCI} (in the a.e. sense) and $(\bar{\psi}_w,\bar{\psi}_r) \in \widetilde{{\mathcal L}^2({\mathcal L}^2)}$ (resp. $\widetilde{{\mathcal L}^2({\mathcal H}^1)}$). Furthermore, 
\begin{eqnarray} 
\| \vec{\psi} \|_{L^2_v(L^2_\kappa)_\mathbf{u}}^2 &=& 2 \pi \, \| (\bar{\psi}_w, \bar{\psi}_r) \|_{\widetilde{{\mathcal L}^2({\mathcal L}^2)}}^2, 
\label{eq:macnorm_vecpsi_rot_L2L2_n=3} \\
\| \vec{\psi} \|_{L^2_v(H^1_\kappa)_\mathbf{u}}^2 &=& 2 \pi \, \| (\bar{\psi}_w, \bar{\psi}_r) \|_{\widetilde{{\mathcal L}^2({\mathcal H}^1)}}^2, 
\label{eq:macnorm_vecpsi_rot_L2H1_n=3}
\end{eqnarray}
and 
\begin{equation} 
\| \nabla_\kappa \vec{\psi} \|_{\vec{L}^2_v(L^2_\kappa)_\mathbf{u}}^2 = 2 \pi \, \Big( \| \bar{\psi}_w \|_{\widetilde{\accentset{\textrm{\textbullet}}{{\mathbb H}}}}^2 + \| \bar{\psi}_r \|_{\widetilde{\accentset{\textrm{\textbullet}}{{\mathbb H}}}}^2 \Big).
\label{eq:macnorm_vecpsi_rot_L2H1_sn_n=3}
\end{equation}
\label{prop:macnorm_vecpsi_rot_n=3}
\end{proposition}

\begin{remark}[About the factor $2 \pi$ in \eqref{eq:macnorm_vecpsi_rot_L2L2_n=3} and \eqref{eq:macnorm_vecpsi_rot_L2H1_n=3}]~

\noindent 
From \eqref{eq:macnorm_normal_coef}, it seems at first glance that the factor appearing in \eqref{eq:macnorm_vecpsi_rot_L2L2_n=3} and \eqref{eq:macnorm_vecpsi_rot_L2H1_n=3} should be~$4 \pi$. However, here, we integrate over $\bar{\kappa}_T \in {\mathbb R}$ instead of over $\kappa_\bot \in [0,\infty)$. Hence a factor $2$ is absorbed in the duplication of the integration domain. 
\end{remark}

\noindent
\textbf{Proof.}
The proof of Prop. \ref{prop:macnorm_vecpsi_rot_n=3} is similar to that of Prop. \ref{prop:macnorm_vecpsi_rot} and is omitted. Let us just mention that the analog of \eqref{eq:macnorm_nD_nakappsi_square} is
\begin{equation} 
|(\nabla_\kappa \vec{\psi})^\dagger|^2 =  |\partial_{\kappa_\parallel} \bar{\psi}_w |^2 + |\partial_{\kappa_\parallel} \bar{\psi}_r |^2 + \frac{1}{\sin^2 \theta} \big( |\partial_{\bar{\kappa}_T} \bar{\psi}_w |^2 + |\partial_{\bar{\kappa}_T} \bar{\psi}_r |^2 \big). 
\label{eq:macnorm_nD_nakappsi_square_n=3}
\end{equation}
\endproof

\medskip
\noindent
(v) We define $\widetilde{{\mathbb H}}'$ as the dual of $\widetilde{{\mathbb H}}$. By Lemma \ref{lem:macnorm_nD_D(N)_densein_H120}, $\widetilde{{\mathbb H}}'$ is a distribution space. 

\medskip
\noindent
(vi) We define 
$$\widetilde{{\mathcal L}^2({\mathcal H}^{-1})} = \big( \widetilde{{\mathbb H}}' \big)^2. $$
Like in Subsection \ref{subsubsec:macnorm_nD_distrib_N}, we can identify this space with the space of pairs of distributions $(T_w,T_r)$ satisfying an analog of identities \eqref{eq:macnorm_invar_distrib_TwTrnorm} to \eqref{eq:macnorm_norm_L2H-1_red} (with $c = 2 \pi$), the precise writing of which are left to the reader. 

\medskip
\noindent
(vii) We define the operator $\widetilde{{\mathcal U}}_1$ acting on $\bar{\psi} \in \widetilde{{\mathbb H}}$ as follows: 
\begin{equation}
\widetilde{{\mathcal U}}_1 \bar{\psi} = \kappa_\parallel \partial_\theta \bar{\psi} + \big(\cos \theta \, \sin \theta \, \bar{\kappa}_T^2 - \nu \, \sin \theta \big) \partial_{\kappa_\parallel} \bar{\psi} -  2 \frac{\cos \theta}{\sin \theta} \, \kappa_\parallel \, \bar{\kappa}_T \, \partial_{\bar{\kappa}_T} \bar{\psi}. 
\label{eq:macnorm_nD_U1operator_def_n=3} 
\end{equation}
The space $\widetilde{{\mathbb Y}}_{\mathcal N}$ is defined by 
\begin{eqnarray*}
&&\hspace{-1cm}
\widetilde{{\mathbb Y}}_{\mathcal N} = \big\{ (\bar{\psi}_w, \bar{\psi}_r) \in \widetilde{{\mathcal L}^2({\mathcal H}^1)} \quad \textrm{such that} 
\quad \big( \widetilde{{\mathcal U}}_1 \bar{\psi}_w - \bar{\kappa}_T \bar{\psi}_r, \widetilde{{\mathcal U}}_1 \bar{\psi}_r + \bar{\kappa}_T \bar{\psi}_w \big) \in \widetilde{{\mathcal L}^2({\mathcal H}^{-1})} \big\}, 
\end{eqnarray*}
where $\widetilde{{\mathcal U}}_1 \psi_w$ and $\widetilde{{\mathcal U}}_1 \psi_r$ are defined in the distributional sense. Then, Lemma \ref{lem:macnorm_U1_skewadj} is satisfied, ``mutatis mutandis''. Let us just make the analog of \eqref{eq:macnorm_nD_brack_zetu_express} explicit: 
\begin{equation}
\big\{ \vec{\psi} , \mathbf{H}_\mathbf{u} \big\}^\dagger = w \, \big( \big( \widetilde{{\mathcal U}}_1 \bar{\psi}_w - \bar{\kappa}_T \bar{\psi}_r \big) + e_\varphi \, \big( \widetilde{{\mathcal U}}_1 \bar{\psi}_r + \bar{\kappa}_T \bar{\psi}_w \big).  
\label{eq:macnorm_nD_brack_zetu_express_n=3}
\end{equation}
This proof of this identity is analog to that of \eqref{eq:macnorm_nD_brack_zetu_express} and is omitted. Then, we get the analog of Prop. \ref{prop:macnorm_vecpsi_rot_Y}: 

\begin{proposition}[Link between ${\mathbb Y}_\mathrm{inv}$ and $\widetilde{{\mathbb Y}}_{\mathcal N}$]~

\noindent
We assume $n=3$. The function $\vec{\psi}$ belongs to ${\mathbb Y}_\mathrm{inv}$ if and only if $\vec{\psi}$ is written according to \eqref{eq:macnorm_n=3_reduced_GCI} and $(\bar{\psi}_w, \bar{\psi}_r)$ belongs to $\widetilde{{\mathbb Y}}_{\mathcal N}$. Furthermore, 
\begin{equation}
\| \vec{\psi} \|_{{\mathbb Y}_\mathrm{inv}}^2 = 2 \pi \, \| (\psi_w, \psi_r) \|_{\widetilde{{\mathbb Y}}_{\mathcal N}}^2. 
\label{eq:macnorm_vecpsi_rot_Y_n=3}
\end{equation}
\label{prop:macnorm_vecpsi_rot_Y_n=3}
\end{proposition}

We can now derive the variational formulation satisfied by $(\bar{\psi}_w, \bar{\psi}_r)$, i.e. perform the analog of Subsection \ref{subsubsubsec:macnorm_nD_red_GCI_varform}. We have the

\begin{theorem}[Variational formulation; case $n=3$]~

\noindent
We assume $n=3$. Let $\vec{\chi} \in L^2_v(L^2_\kappa)_\mathrm{inv}$ and define $(\bar{\chi}_w, \bar{\chi}_r) \in \widetilde{{\mathcal L}^2({\mathcal L}^2)}$ as in \eqref{eq:macnorm_n=3_reduced_GCI}. Then, $\vec{\psi}$ is the unique solution of the variational formulation \eqref{eq:macnorm_GCI_gene_varform_inv} if and only if $(\bar{\psi}_w, \bar{\psi}_r)$ linked to $\vec{\psi}$ by \eqref{eq:macnorm_n=3_reduced_GCI} is the unique solution of the variational formulation
\begin{equation}
\mbox{} \hspace{-0.5cm} \left\{ \begin{array}{l}
(\bar{\psi}_w, \bar{\psi}_r) \in \widetilde{{\mathbb Y}}_{\mathcal N}, \\
\textrm{ } \\
- \big( \hspace{-1.2mm} \big(\widetilde{{\mathcal U}}_1 \bar{\psi}_w  - \bar{\kappa}_T \bar{\psi}_r \, ,  \, \bar{\varphi}_w\big) \hspace{-1.2mm} \big)_{\widetilde{{\mathbb H}'},\widetilde{{\mathbb H}}} 
-  \big( \hspace{-1.2mm} \big( \widetilde{{\mathcal U}}_1 \bar{\psi}_r + \bar{\kappa}_T \bar{\psi}_w \, , \, \bar{\varphi}_r \big) \hspace{-1.2mm} \big)_{\widetilde{{\mathbb H}'},\widetilde{{\mathbb H}}} 
+ D  \big( \hspace{-1.2mm} \big( \bar{\psi}_w \ , \, \bar{\varphi}_w \big) \hspace{-1.2mm} \big)_{\widetilde{\accentset{\textrm{\textbullet}}{{\mathbb H}}}} 
\\
\textrm{ } \\
\hspace{1.5cm} 
+ D \big( \hspace{-1.2mm} \big( \bar{\psi}_r \, , \, \bar{\varphi}_r \big) \hspace{-1.2mm} \big)_{\widetilde{\accentset{\textrm{\textbullet}}{{\mathbb H}}}} = -  \big( \hspace{-1.2mm} \big( \bar{\chi}_w \, , \, \bar{\varphi}_w \big) \hspace{-1.2mm} \big)_{\widetilde{{\mathbb L}}} - \big( \hspace{-1.2mm} \big( \bar{\chi}_r \, , \, \bar{\varphi}_r \big) \hspace{-1.2mm} \big)_{\widetilde{{\mathbb L}}} ,  
\qquad \forall (\bar{\varphi}_w, \bar{\varphi}_r) \in \widetilde{{\mathbb Y}}_{\mathcal N}. \end{array} \right. 
\label{eq:macnorm_redGCI_varform_n=3}
\end{equation}
\label{thm:macnorm_redGCI_varform_n=3}
\end{theorem}

\noindent
followed by the

\begin{corollary}[Variational formulation satisfied by the reduced GCI pair; case $n=3$]~

\noindent
Let $(\bar{\zeta}_w, \bar{\zeta}_r)$ be the reduced GCI pair in dimension $n=3$. Then, $(\bar{\zeta}_w, \bar{\zeta}_r)$ satisfies the variational formulation \eqref{eq:macnorm_redGCI_varform_n=3}, with $(\bar{\psi}_w, \bar{\psi}_r)$ replaced by $(\bar{\zeta}_w, \bar{\zeta}_r)$ and $(\bar{\chi}_w, \bar{\chi}_r) = (\sin \theta, 0)$. 
\label{cor:macnorm_redGCI_varform_n=3}
\end{corollary}

The strong form associated with the variational formulation \eqref{eq:macnorm_redGCI_varform_n=3} is given in

\begin{theorem}[PDE system associated with the variational formulation \eqref{eq:macnorm_redGCI_varform_n=3} and boundary conditions, case $n=3$]~

\smallskip
\noindent
(i) Let $(\bar{\zeta}_w, \bar{\zeta}_r) \in \widetilde{{\mathbb Y}}_{{\mathcal N}}$ be the unique solution to the variational formulation \eqref{eq:macnorm_redGCI_varform_n=3} (with $(\bar{\psi}_w, \bar{\psi}_r)$ replaced by $(\bar{\zeta}_w, \bar{\zeta}_r)$ and $(\bar{\chi}_w, \bar{\chi}_r) = (\sin \theta, 0)$). Then, $(\zeta_w, \zeta_r)$ linked to $(\bar{\zeta}_w, \bar{\zeta}_r)$ through \eqref{eq:macnorm_n=3_link_red_GCI} satisfies System \eqref{eq:macnorm_nD_vectGCI_wcomp}, \eqref{eq:macnorm_nD_vectGCI_kapTcomp} in the distributional sense. 

\smallskip
\noindent
(ii) (formal) Supposing that $(\bar{\zeta}_w, \bar{\zeta}_r)$ is smooth. Then the boundary conditions implied by \eqref{eq:macnorm_n=3_link_red_GCI} and the variational formulation \eqref{eq:macnorm_redGCI_varform_n=3} are
\begin{itemize}
\item[-] no boundary condition at $\theta = 0$ nor $\theta = \pi$ on neither $\zeta_w$ nor $\zeta_r$;
\item[-] Neumann boundary conditions on $\zeta_w$ and homogeneous Dirichlet condition on $\zeta_r$ at~$\kappa_\bot = 0$:
\begin{equation} 
\partial_{\kappa_\bot} \zeta_w (\theta, \kappa_\parallel, 0) = 0, \qquad \zeta_r (\theta, \kappa_\parallel, 0) = 0; 
\label{eq:macnorm_nD_bc_n=3}
\end{equation}
\end{itemize}
\label{thm:macnorm_nD_red_GCI_syst_and_bc_n=3}
\end{theorem}

\noindent
\textbf{Proof.} The same proof as for Theorem \ref{thm:macnorm_nD_red_GCI_syst_and_bc} leads to a system of equation for $(\bar{\zeta}_w, \bar{\zeta}_r)$ posed on~$\widetilde{{\mathcal N}}$. By restricting it to ${\mathcal N}$, and using the fact that, on ${\mathcal N}$, $\bar{\kappa}_T = \kappa_\bot$, we recover System \eqref{eq:macnorm_nD_vectGCI_wcomp}, \eqref{eq:macnorm_nD_vectGCI_kapTcomp}. We get no boundary condition for the system posed on $\widetilde{{\mathcal N}}$. However, when restricting it to~${\mathcal N}$, the mere smoothness of $(\bar{\zeta}_w, \bar{\zeta}_r)$ implies Neumann boundary conditions on $\zeta_w$ and Dirichlet boundary conditions on $\zeta_r$ at $\kappa_\bot = 0$. This ends the proof. \endproof

\begin{remark}[About the variational formulation and the strong form]~

\noindent
Although the strong form can be restricted to a system of PDE posed on ${\mathcal N}$, this is not possible for the variational formulation. Indeed, the definition of the space $\widetilde{{\mathbb Y}}_{\mathcal N}$ involves the use of the dual space $ \widetilde{{\mathbb H}}'$ the definition of which requires a domain without boundaries, i.e. $\widetilde{{\mathcal N}}$. The variational formulation is the only formulation of the problem for which existence and uniqueness of solutions is proved. So, the reduction to a PDE posed on ${\mathcal N}$ as formulated in Theorem \ref{thm:macnorm_nD_red_GCI_syst_and_bc_n=3} is not useful for a rigorous mathematical use. Its only merit is to highlight the analogies and differences between  dimension 3 and higher dimensions. 
\end{remark}

\setcounter{equation}{0}
\section{Evolution equation for $\mathbf{u}$}
\label{sec:macnorm_eq_j}

We now use the reduced GCI pair to derive the equation for $\mathbf{u}$. This is the final step of the formal proof of Theorem \ref{thm:hydro_main_theorem}.

\begin{lemma}[Evolution equation for $\mathbf{u}$: first form]~

\noindent
The functions $\rho$ and $\mathbf{u}$ defined by \eqref{eq:macnorm_equi} satisfy the following equation
\begin{equation}
\int_M (\partial_t + v \cdot \nabla_x) \big( \rho \mathbf{M}_{\mathbf{u}} \big) \, \vec \zeta_{\mathbf{u}} \, |\mathrm{Vol}_M| = 0. 
\label{eq:macnorm_eq_evol_u_1st}
\end{equation}
\label{lem:macnorm_eq_evol_u_1st}
\end{lemma}

\noindent
\textbf{Proof.} We muliply \eqref{eq:pk_eq_f_Q} by $\vec \zeta_{\mathbf{u}_{f^\varepsilon}}$ and integrate over $M$. Thanks to \eqref{eq:macnorm_vecGCI_cancel}, we have 
$$\int_M Q(f^\varepsilon) \, \vec \zeta_{\mathbf{u}_{f^\varepsilon}} \, |\mathrm{Vol}_M| = 0. 
$$
Hence, we get 
\begin{equation} 
\int_M (\partial_t + v \cdot \nabla_x) f^\varepsilon \, \vec \zeta_{\mathbf{u}_{f^\varepsilon}} \, |\mathrm{Vol}_M|= 0. 
\label{eq:macnorm_eq_evol_u_1st_prf1}
\end{equation}
Now, when $\varepsilon \to 0$, we have formally $f^\varepsilon \to \rho \mathbf{M}_{\mathbf{u}}$ and $\mathbf{u}_{f^\varepsilon} \to \mathbf{u}_{\rho \mathbf{M}_{\mathbf{u}}} = \mathbf{u}$ by the proof of Prop.~\ref{prop:equi_equi}. Now, we show that, for any $\alpha \in M$, the map $\psi$: $\mathbf{u} \mapsto \zeta_\mathbf{u}(\alpha)$, ${\mathbb S}^{n-1} \to {\mathbb R}$ is smooth. We fix $\mathbf{u}_0 \in {\mathbb S}^{n-1}$. We first recall that ${\mathbb S}^{n-1}$ is diffeomorphic to the coset space $\mathrm{SO}_n/\mathrm{SO}_{n-1}$ where $\mathrm{SO}_{n-1}$ denotes the isotropy subgroup of $\mathbf{u}_0$. Hence, the projection $\pi$: $\mathrm{SO}_n \to {\mathbb S}^{n-1}$ has the structure of a quotient manifold. Now, thanks to \eqref{eq:macnorm_gract_zetaj}, the map $\varphi$, $R \mapsto \zeta_{R \mathbf{u}_0}(\alpha) = R \zeta_{\mathbf{u}_0}(R^T \alpha)$, $\mathrm{SO}_n \to {\mathbb R}$ is smooth. Since $\varphi = \psi \circ \pi$, it follows from \cite[Vol 1, Sect. 3.9, Prop. V]{greub1972connections} that $\psi$ is smooth. Hence, we have a pointwise convergence $\vec \zeta_{\mathbf{u}_{f^\varepsilon}} \to \vec \zeta_{\mathbf{u}}$ as $\varepsilon \to 0$. Thus letting $\varepsilon \to 0$ in \eqref{eq:macnorm_eq_evol_u_1st_prf1} and supposing that the convergence of $f^\varepsilon$ to $f$ is strong enough, we are led to \eqref{eq:macnorm_eq_evol_u_1st}, which ends the proof. \endproof

Now, we define the following tensors: 
\begin{eqnarray}
{\mathbb A} &=& \int_{M_V \times M_T} \sin^2 \theta \big( (\mathbf{M}_\mathbf{u} \vec \zeta_\mathbf{u})^\dagger(\theta, \kappa_\parallel, w, \kappa_T) \otimes w \otimes w \big) \,  \,  |\mathbf{A}_\mathbf{u}^* (\mathrm{Vol}_M)| , \label{eq:macnorm_AA_def} \\
{\mathbb B} &=&  \int_{M_V \times M_T} \sin \theta \, \big( (\mathbf{M}_\mathbf{u} \vec \zeta_\mathbf{u})^\dagger(\theta, \kappa_\parallel, w, \kappa_T) \otimes w \big) \,   |\mathbf{A}_\mathbf{u}^* (\mathrm{Vol}_M)|, \label{eq:macnorm_BB_def} \\
{\mathbb C} &=& \int_{M_V \times M_T} \sin \theta \, \cos \theta \, \big( (\mathbf{M}_\mathbf{u} \vec \zeta_\mathbf{u})^\dagger(\theta, \kappa_\parallel, w, \kappa_T) \otimes w \big) \,   |\mathbf{A}_\mathbf{u}^* (\mathrm{Vol}_M)|, \label{eq:macnorm_CC_def} \\
{\mathbb D} &=& \int_{M_V \times M_T} (\mathbf{M}_\mathbf{u} \vec \zeta_\mathbf{u})^\dagger(\theta, \kappa_\parallel, w, \kappa_T) \,  |\mathbf{A}_\mathbf{u}^* (\mathrm{Vol}_M)|, \label{eq:macnorm_DD_def} \\
{\mathbb E} &=& \int_{M_V \times M_T} \cos \theta \, (\mathbf{M}_\mathbf{u} \vec \zeta_\mathbf{u})^\dagger(\theta, \kappa_\parallel, w, \kappa_T)  \,  |\mathbf{A}_\mathbf{u}^* (\mathrm{Vol}_M)|, \label{eq:macnorm_EE_def}
\end{eqnarray}
where we recall that $|\mathbf{A}_\mathbf{u}^* (\mathrm{Vol}_M)|$ is given by \eqref{eq:macnorm_nD_AuvolM2}. Here, $w$ and $\vec \zeta_\mathbf{u}^\dagger$ are considered as vectors of ${\mathbb R}^n$ and ${\mathbb A}$ is a rank-3 tensor, ${\mathbb B}$ and ${\mathbb C}$ are rank-2 tensors (i.e. matrices) and ${\mathbb D}$ and ${\mathbb E}$ are vectors. These tensors are useful to transform the left-hand side of \eqref{eq:macnorm_eq_evol_u_1st_prf1} into explicit expressions, as shown next. We note that 
\begin{equation}
\mathbf{M}_\mathbf{u}^\dagger(\theta, \kappa_\parallel, w, \kappa_T) = \frac{1}{Z} \exp \bigg( - \frac{\nu}{D} \Big( - \nu \cos \theta + \frac{1}{2} (\kappa_\parallel^2 +\sin^2 \theta \, |\kappa_T|^2) \Big) \bigg). 
\label{eq:macnorm_Mudag}
\end{equation}
The right-hand side of \eqref{eq:macnorm_Mudag} does not depend on $w$  and defines the function $\mathbf{M}_\mathbf{u}^\natural = \mathbf{M}_\mathbf{u}^\natural(\theta, \kappa_\parallel, \kappa_T)$. 

\begin{lemma}[Evolution equation for $\mathbf{u}$: second form]~

\noindent
The functions $\rho$ and $\mathbf{u}$ defined by \eqref{eq:macnorm_equi} satisfy the following equation 
\begin{equation}
\frac{D^2}{\nu} \big( {\mathbb D} \partial_t \rho + {\mathbb E} (\mathbf{u} \cdot \nabla_x) \rho \big) +  {\mathbb B} \Big( \rho \partial_t \mathbf{u}  + \frac{D}{\nu^2} \nabla_x \rho \Big) + \rho {\mathbb C} (\mathbf{u} \cdot \nabla_x) \mathbf{u} + \rho {\mathbb A} : \nabla_x \mathbf{u} = 0,
\label{eq:macnorm_eq_evol_u_2nd}
\end{equation}
where the $:$ in the last term indicates contraction of the last two indices of ${\mathbb A}$ with the two indices of $\nabla_x \mathbf{u}$ while the second and third terms involve the matrix multiplication of ${\mathbb B}$ and ${\mathbb C}$ with the vector to their right. 
\label{lem:macnorm_eq_evol_u_2nd}
\end{lemma}

\noindent
\textbf{Proof.} We compute $(\partial_t + v \cdot \nabla_x) (\rho \mathbf{M}_{\mathbf{u}})$. This computation is straightforward but is reproduced here for the readers' convenience. We have 
$$(\partial_t + v \cdot \nabla_x) \big( \rho \mathbf{M}_{\mathbf{u}} \big) = \mathbf{M}_{\mathbf{u}} \Big( 
(\partial_t + v \cdot \nabla_x) \rho + \rho (d \log \mathbf{M}_{\mathbf{u}})_{\mathbf{u}} \big( (\partial_t + v \cdot \nabla_x) \mathbf{u} \big) \Big). $$
We have $ \log \mathbf{M}_{\mathbf{u}} = - (\nu/D) \mathbf{H}_{\mathbf{u}} - \log Z$. Hence, for $\sigma \in T_{\mathbf{u}} {\mathbb S}^{n-1}$, we have 
$$ (d \log \mathbf{M}_{\mathbf{u}})_{\mathbf{u}} (\sigma) = - (\nu/D) (d \mathbf{H}_{\mathbf{u}})_{\mathbf{u}}  (\sigma) = \frac{\nu^2}{D} P_{{\mathbf{u}}^\bot} v \cdot \sigma. $$
Thus, we have 
$$(\partial_t + v \cdot \nabla_x) \big( \rho \mathbf{M}_{\mathbf{u}} \big) = \mathbf{M}_{\mathbf{u}} \Big( 
(\partial_t + v \cdot \nabla_x) \rho + \frac{\nu^2}{D} \rho P_{{\mathbf{u}}^\bot} v \cdot \big( (\partial_t + v \cdot \nabla_x) \mathbf{u} \big) \Big). $$
We use the change of variables \eqref{eq:macnorm_dV_express} and get 
\begin{eqnarray}
&& \hspace{-1cm}
\frac{D}{\nu^2} \Big( (\partial_t + v \cdot \nabla_x) \big( \rho \mathbf{M}_{\mathbf{u}} \big) \Big)^\dagger = \mathbf{M}_{\mathbf{u}}^\dagger  \Big\{ \frac{D}{\nu^2} 
(\partial_t \rho + \cos \theta \, \mathbf{u} \cdot \nabla_x \rho) + \sin \theta \, w \cdot \Big(\rho \partial_t \mathbf{u} +  \frac{D}{\nu^2} \nabla_x \rho \Big) \nonumber \\
&& \hspace{1.5cm}
+  \rho \sin \theta \, \cos \theta \, w  \cdot \big((\mathbf{u} \cdot \nabla_x) \mathbf{u} \big)   + \rho \sin^2 \theta \, (w \otimes w): (\nabla_x \mathbf{u}) \Big\}. \label{eq:macnorm_eq_evol_u_2nd_prf1}
\end{eqnarray}
Now, tensorizing on the left by $\vec \zeta_\mathbf{u}^\dagger$ and integrating the result with respect to $|\mathbf{A}_\mathbf{u}^* (\mathrm{Vol}_M)|$, 
we get \eqref{eq:macnorm_eq_evol_u_2nd}, which ends the proof. \endproof

Now, thanks to the explicit expression \eqref{eq:macnorm_nD_actionOn-1} of $\vec{\zeta}_\mathbf{u}^\dagger$, we can compute the tensors ${\mathbb A}$ to ${\mathbb E}$ explicitly.

\begin{lemma}[Expressions of the tensors ${\mathbb A}$ to ${\mathbb E}$]~

\noindent
There exist two real numbers ${\mathbb b}$ and ${\mathbb c}$ such that
\begin{eqnarray}
{\mathbb A} &=& 0, \label{eq:macnorm_AA=0} \\
{\mathbb B} &=& {\mathbb b} P_{\mathbf{u}^\bot}, \quad {\mathbb C} = {\mathbb c} P_{\mathbf{u}^\bot}, \label{eq:macnorm_BB_CC} \\
{\mathbb D}&=&{\mathbb E}=0. \label{eq:macnorm_DD=EE=0}
\end{eqnarray}
Furthermore, we have
\begin{eqnarray}
&&\hspace{-1.5cm} 
{\mathbb b} = (-1)^n \, \frac{c}{n-1} 
\int_{{\mathcal N}} (\mathbf{M}_\mathbf{u}^\natural \zeta_w)(\theta, \kappa_\parallel, \kappa_\bot) \, \sin^{2n-3} \theta \, \, \kappa_\bot^{n-3} \, d\theta \, d\kappa_\parallel \, d\kappa_\bot, 
\label{eq:macnorm_b_express} \\
&&\hspace{-1.5cm} 
{\mathbb c} = (-1)^n \, \frac{c}{n-1} 
\int_{{\mathcal N}} (\mathbf{M}_\mathbf{u}^\natural \zeta_w)(\theta, \kappa_\parallel, \kappa_\bot) \, \cos \theta \, \sin^{2n-3} \theta \, \, \kappa_\bot^{n-3} \, d\theta \, d\kappa_\parallel \, d\kappa_\bot,  
\label{eq:macnorm_c_express}
\end{eqnarray}
where $c$ is given by \eqref{eq:macnorm_normal_coef}. 
\label{lem:macnorm_DD=EE=0}
\end{lemma}

\noindent
\textbf{Proof.} 
We use the same methodology as in the proof of Proposition \ref{prop:macnorm_vecpsi_rot}(i). We first consider ${\mathbb B}$ given by \eqref{eq:macnorm_BB_def}. Thanks to \eqref{eq:macnorm_nD_AuvolM2}, we have 
$$ {\mathbb B} = \int_{M_V \times M_T} \big( (\mathbf{M}_\mathbf{u} \zeta_\mathbf{u})^\dagger(\theta, \kappa_\parallel, w, |\kappa_T|) \otimes w \big) \, \sin^{2n-3} \theta \, \, |\mathrm{Vol}_{\mathrm{M}_T}| \, d \theta \, d \kappa_\parallel. $$
Then, we apply the Fubini theorem and get
\begin{eqnarray} {\mathbb B} &=& \int_{M_V} \Big( \int_{M_T }\big( (\mathbf{M}_\mathbf{u} \zeta_\mathbf{u})^\dagger(\theta, \kappa_\parallel, w, |\kappa_T|) \otimes w \big) \, |\mathrm{Vol}_{\mathrm{M}_T}| \Big) \, \sin^{2n-3} \theta \, \, d \theta \, d \kappa_\parallel \nonumber \\
&=& \int_{M_V} \, B(\theta, \kappa_\parallel) \, \sin^{2n-3} \theta \, \, d \theta \, d \kappa_\parallel. \label{eq:macnorm_DD=EE=0_prf1}
\end{eqnarray}
With \eqref{eq:macnorm_nD_actionOn-1}, we have
\begin{eqnarray*}
B &=& B_1 + B_2, \\
B_1 &=& \int_{M_T } (\mathbf{M}_\mathbf{u}^\natural \zeta_w) (\theta, \kappa_\parallel, |\kappa_T|) \, (w \otimes w) \, |\mathrm{Vol}_{\mathrm{M}_T}|, \\
B_2 &=& \int_{M_T } (\mathbf{M}_\mathbf{u}^\natural \zeta_r) (\theta, \kappa_\parallel, |\kappa_T|) \, \Big(\frac{\kappa_T}{|\kappa_T|} \otimes w \Big) \, |\mathrm{Vol}_{\mathrm{M}_T}|. 
\end{eqnarray*}
Thanks to \eqref{eq:Jacobian_Au_volMT}, \eqref{eq:Fub_sigma(x)_def} and to the Fubini theorem for fiber bundles \eqref{eq:Fub_Fubthm}, we can write
$$ B_1 = \int_{{\mathbb S}^{n-2}}  (w \otimes w) \, \Big( \int_{T_w {\mathbb S}^{n-2}} (\mathbf{M}_\mathbf{u}^\natural \zeta_w) (\theta, \kappa_\parallel, |\kappa_T|) \, |\omega_w| \Big) \, |\mathrm{Vol}_{{\mathbb S}^{n-2}}|, $$
where $\omega_w$ is the retrenchment of $\omega$ on $T_w {\mathbb S}^{n-2}$. Now, we can follow the same computation as in the proof of Proposition \ref{prop:macnorm_vecpsi_rot}(i) and get 
$$B_1 = |{\mathbb S}^{n-3}| \, \Big( \int_0^\infty  (\mathbf{M}_\mathbf{u}^\natural \zeta_w) (\theta, \kappa_\parallel, \kappa_\bot) \, \kappa_\bot^{n-3} \, d \kappa_\bot \Big) \, \Big( \int_{{\mathbb S}^{n-2}} (w \otimes w) \, |\mathrm{Vol}_{{\mathbb S}^{n-2}}| \Big). $$
It is readily seen that 
\begin{equation} 
\int_{{\mathbb S}^{n-2}} (w \otimes w) \, |\mathrm{Vol}_{{\mathbb S}^{n-2}}| = \frac{|{\mathbb S}^{n-2}|}{n-1} P_{\mathbf{u}^\bot}. 
\label{eq:macnorm_DD=EE=0_prf2}
\end{equation}
Hence, we have 
$$B_1 = \frac{|{\mathbb S}^{n-2}| \, |{\mathbb S}^{n-3}|}{n-1} \, \int_0^\infty  (\mathbf{M}_\mathbf{u}^\natural \zeta_w) (\theta, \kappa_\parallel, \kappa_\bot) \, \kappa_\bot^{n-3} \, d \kappa_\bot \, \, P_{\mathbf{u}^\bot}. $$
To compute $B_2$, we consider a direct orthonormal basis $(\mathbf{e}_1, \ldots, \mathbf{e}_{n-1}, \mathbf{u})$ of ${\mathbb R}^n$ and we compute the $(i,j)$-th entry of $B_2$ in this basis, i.e. 
$$ (B_2)_{ij} = \int_{M_T } (\mathbf{M}_\mathbf{u}^\natural \zeta_r) (\theta, \kappa_\parallel, |\kappa_T|) \, \frac{\kappa_T \cdot \mathbf{e}_i}{|\kappa_T|} \,  w_j \, |\mathrm{Vol}_{\mathrm{M}_T}|. $$
By the same manipulations as for $B_1$, we get 
$$ (B_2)_{ij} =, \int_{{\mathbb S}^{n-2}} w_j \,  \Big( \int_{T_w{\mathbb S}^{n-2}} (\mathbf{M}_\mathbf{u}^\natural \zeta_r) (\theta, \kappa_\parallel, |\kappa_T|) \, \frac{\kappa_T \cdot \mathbf{e}_i}{|\kappa_T|} \, |\omega_w| \Big) \, |\mathrm{Vol}_{{\mathbb S}^{n-2}}|. $$
Let $(\varepsilon_1, \ldots, \varepsilon_{n-2})$ be a local, positively oriented orthonormal frame of ${\mathbb S}^{n-2}$. Define a coordinate map $T_v B \to {\mathbb R}^{n-1}$, $\kappa \mapsto (\kappa^1, \ldots, \kappa^{n-2})$ with $\kappa^i = \langle \kappa , \varepsilon_i \rangle_v$.  Then, $|\omega_w|= d \kappa^1 \ldots d \kappa^{n-2}$ and we have 
\begin{eqnarray*}
&&\hspace{-1cm}
\int_{T_w{\mathbb S}^{n-2}} (\mathbf{M}_\mathbf{u}^\natural \zeta_r) (\theta, \kappa_\parallel, |\kappa_T|) \, \frac{\kappa_T \cdot \mathbf{e}_i}{|\kappa_T|} \, |\omega_w| \\
&& \hspace{1cm}
= \sum_{k=1}^{n-2} (\varepsilon_k \cdot \mathbf{e}_i) \, \int_{{\mathbb R}^{n-2}} (\mathbf{M}_\mathbf{u}^\natural \zeta_r) \Big(\theta, \kappa_\parallel, \big(\sum_{\ell=1}^{n-2} |\kappa_T^\ell|^2 \big)^{1/2} \Big) \, \frac{\kappa_T^k}{|\kappa_T|} \, d \kappa_T^1 \ldots d \kappa_T^{n-2}. 
\end{eqnarray*}
Now, by the change of variables $(\kappa_T^1, \ldots, \kappa_T^k, \ldots, \kappa_T^{n-2}) \to (\kappa_T^1, \ldots, - \kappa_T^k, \ldots, \kappa_T^{n-2})$, we realize that the integral in the $k$-th term of the sum at the right-hand side of the previous formula vanishes by antisymmetry. Hence, $B_2 = 0$. It follows that ${\mathbb B} = {\mathbb b} P_\mathbf{u}$ as it should, with ${\mathbb b}$ given by~\eqref{eq:macnorm_b_express}. Going back to the formula \eqref{eq:macnorm_CC_def}, we realize that the same formula holds for for ${\mathbb C}$, simply multiplying the integrated function by $\cos \theta$. 

Now, we consider ${\mathbb A}$. Reproducing the same computation as for ${\mathbb B}$, we would get an analog to \eqref{eq:macnorm_DD=EE=0_prf1} with $B$ replaced by $A$ such that  
\begin{eqnarray*}
A &=& A_1 + A_2, \\
A_1 &=& \int_{M_T } (\mathbf{M}_\mathbf{u}^\natural \zeta_w) (\theta, \kappa_\parallel, |\kappa_T|) \, (w \otimes w \otimes w) \, |\mathrm{Vol}_{\mathrm{M}_T}|, \\
A_2 &=& \int_{M_T } (\mathbf{M}_\mathbf{u}^\natural \zeta_r) (\theta, \kappa_\parallel, |\kappa_T|) \, \Big(\frac{\kappa_T}{|\kappa_T|} \otimes w \otimes w \Big) \, |\mathrm{Vol}_{\mathrm{M}_T}|. 
\end{eqnarray*}
Now, $A_2 = 0$ by the same argument as for $B_2$, and $A_1 = 0$ because we have an even number of $w$-factors inside the integral corresponding to \eqref{eq:macnorm_DD=EE=0_prf2} and the integral vanishes by antisymmetry. Hence, we get \eqref{eq:macnorm_AA=0}. To get \eqref{eq:macnorm_DD=EE=0}, the reasoning is the same, except that now there is only one $w$-factor inside the integral corresponding to \eqref{eq:macnorm_DD=EE=0_prf2} but the integral vanishes as well for the same reason. This ends the proof. \endproof

\begin{remark}[About the expressions of the tensors ${\mathbb A}$ to ${\mathbb E}$]~

\noindent
Expressions \eqref{eq:macnorm_AA=0}-\eqref{eq:macnorm_DD=EE=0} for the tensors ${\mathbb A}$ to ${\mathbb E}$ could be obtained by simple symmetry considerations using \eqref{eq:macnorm_gract_zeta_invar}. However, the expressions \eqref{eq:macnorm_b_express}, \eqref{eq:macnorm_c_express} of ${\mathbb b}$ and ${\mathbb c}$ require Theorem~\ref{thm:macnorm_nD_actionOn-1} and the use of the reduced GCI pair. 
\end{remark}

We are now able to state the final form of the equation for $\mathbf{u}$. The following lemma completes the proof of Theorem \ref{thm:hydro_main_theorem}.

\begin{lemma}[Evolution equation for $\mathbf{u}$: final form]~

\noindent
The functions $\rho$ and $\mathbf{u}$ defined by \eqref{eq:macnorm_equi} satisfy 
\eqref{eq:macnorm_eq_evol_u}, where $c_2$ is given by \eqref{eq:macnorm_nD_coef_c2} and $c_3$ by~\eqref{eq:macnorm_eq_coeffs_c3}. 
\label{lem:macnorm_eq_evol_u_final}
\end{lemma}

\noindent
\textbf{Proof.} We insert the expressions of the tensors ${\mathbb A}$ to ${\mathbb E}$ found in Lemma \ref{lem:macnorm_DD=EE=0} into \eqref{eq:macnorm_eq_evol_u_2nd} and get
$${\mathbb b} \, P_{\mathbf{u}^\bot} \Big( \rho \partial_t \mathbf{u}  - \frac{D}{\nu^2} \nabla_x \rho \Big) + {\mathbb c} \, \rho P_{\mathbf{u}^\bot}  \big( (\mathbf{u} \cdot \nabla_x) \mathbf{u} \big) = 0.$$
Noting that both vectors $\partial_t \mathbf{u}$ and $(\mathbf{u} \cdot \nabla_x) \mathbf{u}$ belong to $\{\mathbf{u}\}^\bot$, and dividing by ${\mathbb b}$, we find \eqref{eq:macnorm_eq_evol_u} with $c_2={\mathbb c}/{\mathbb b}$ and $c_3$ given by~\eqref{eq:macnorm_eq_coeffs_c3}. Now, thanks to \eqref{eq:macnorm_b_express}, \eqref{eq:macnorm_c_express}, we find that the expression of $c_2$ is prescisely given by \eqref{eq:macnorm_nD_coef_c2}, which ends the proof. \endproof

\setcounter{equation}{0}
\section{Conclusion and perspectives}
\label{sec:conclu}

In this paper, we have provided a new formulation of the particle curvature-control model proposed in \cite{cavagna2015flocking, degond2011macroscopic} which enables us to extend it to arbitrary dimensions $n \geq 3$ (previously the model was stated in dimensions 2 in \cite{degond2011macroscopic} and 3 in \cite{cavagna2015flocking}). We have also established the associated mean-field kinetic equations. Both the particle and kinetic models require a bundle geometric framework. Then, the main result of this paper is the derivation of the hydrodynamic limit of the kinetic model. We have shown that it coincides with the SOH model already found in the context of the Vicsek model. It requires the introduction of generalized collision invariants (GCI), which are solutions to a Fokker-Planck type equation. The action of the orthogonal group permits to express the GCI in terms of a pair of functions which satisfy a system of PDEs. Unique solvability of this system has been shown by means of a novel type of variational formulation. 

This work opens many directions of research. First of all, the derivation of the SOH model is formal and lacks a mathematically rigorous proof. Numerical simulations should be developped to assess the accuracy and applicability range of the SOH model. We would also like to extends the present work to the non-normalized case. In this case, the model exhibits phase transitions in a similar way to the Vicsek model. For the latter, the multiplicity and stability of the spatially-homogeneous equilbria has been uncovered in \cite{degond2013macroscopic, degond2015phase}. The corresponding results for the curvature-control model are still open. Some modifications of the model could be investigated, such as, for instance, restricting the curvature to a unit vector, to mimic the fact that animals have probably a range of preferences for their curvature. Another modification could be to replace the mean-field type of interaction by a Boltzmann like interaction, like in~\cite{gautrais2012deciphering}. Finally one could also try to modify the model such that its hydrodynamic limit involves an additional equation for the locally-averaged curvature, like in \cite{yang2015hydrodynamics}. This non-exhaustive list highlights the wealth of new research directions opened by the study of the curvature-control model of swarming behavior.


\newpage
\bibliographystyle{abbrv}
\bibliography{Biblio_PTWAn}

@book{warner1983foundations,
  title={Foundations of {D}ifferentiable {M}anifolds and {L}ie {G}roups},
  author={Warner, Frank W},
  year={1983},
  publisher={Springer}
}

@book{docarmo1992riemannian,
  title={Riemannian Geometry},
  author={do Carmo, Manfredo P.},
  year={1992},
  publisher={Birkh\"auser}
}

@book{greub1972connections,
  title={Connections, Curvature, and Cohomology},
  author={Greub, Werner and Halperin, Stephen and Vanstone, Ray},
  year={1972},
  publisher={Academic press}
}

@article{sasaki1958differential,
  title={On the differential geometry of tangent bundles of Riemannian manifolds},
  author={Sasaki, Shigeo},
  journal={Tohoku Math. J. (2)},
  volume={10},
  number={3},
  pages={338--354},
  year={1958},
  publisher={Mathematical Institute, Tohoku University}
}

@article{gudmundsson2002geometry,
  title={On the geometry of tangent bundles},
  author={Gudmundsson, Sigmundur and Kappos, Elias},
  journal={Expo. Math.},
  volume={20},
  number={1},
  pages={1--41},
  year={2002},
  publisher={Elsevier}
}

@book{ikeda2014stochastic,
  title={Stochastic differential equations and diffusion processes},
  author={Ikeda, Nobuyuki and Watanabe, Shinzo},
  year={2014},
  publisher={Elsevier}
}

@book{hsu2002stochastic,
  title={Stochastic analysis on manifolds},
  author={Hsu, Elton P},
  year={2002},
  publisher={American Mathematical Soc.}
}

@book{gallot1990riemannian,
  title={Riemannian geometry},
  author={Gallot, Sylvestre and Hulin, Dominique and Lafontaine, Jacques},
  volume={2},
  year={1990},
  publisher={Springer}
}

@article{nier2024global,
  title={Global subelliptic estimates for geometric {K}ramers-{F}okker-{P}lanck operators on closed manifolds},
  author={Nier, Francis and Sang, Xingfeng and White, Francis},
  journal={arXiv preprint arXiv:2402.07511},
  year={2024}
}

@book{arnol2013mathematical,
  title={Mathematical methods of classical mechanics},
  author={Arnol'd, Vladimir Igorevich},
  volume={60},
  year={2013},
  publisher={Springer Science \& Business Media}
}

@article{degond2013macroscopic,
  title={Macroscopic limits and phase transition in a system of self-propelled particles},
  author={Degond, Pierre and Frouvelle, Amic and Liu, Jian-Guo},
  journal={J. Nonlinear Sci.},
  volume={23},
  pages={427--456},
  year={2013},
  publisher={Springer}
}

@article{degond2015phase,
  title={Phase transitions, hysteresis, and hyperbolicity for self-organized alignment dynamics},
  author={Degond, Pierre and Frouvelle, Amic and Liu, Jian-Guo},
  journal={Arch. Ration. Mech. Anal.},
  volume={216},
  pages={63--115},
  year={2015},
  publisher={Springer}
}

@article{degond2008continuum,
  title={Continuum limit of self-driven particles with orientation interaction},
  author={Degond, Pierre and Motsch, S{\'e}bastien},
  journal={Math. Models Methods Appl. Sci.},
  volume={18},
  number={supp01},
  pages={1193--1215},
  year={2008},
  publisher={World Scientific}
}

@article{degond2011macroscopic,
  title={A macroscopic model for a system of swarming agents using curvature control},
  author={Degond, Pierre and Motsch, S{\'e}bastien},
  journal={J. Stat. Phys.},
  volume={143},
  pages={685--714},
  year={2011},
  publisher={Springer}
}

@book{lions2013equations,
  title={Equations diff{\'e}rentielles op{\'e}rationnelles et probl{\'e}mes aux limites},
  author={Lions, Jacques Louis},
  volume={111},
  year={2013},
  publisher={Springer-Verlag}
}

@article{albritton2019variational,
  title={Variational methods for the kinetic {F}okker-{P}lanck equation},
  author={Albritton, Dallas and Armstrong, Scott and Mourrat, J-C and Novack, Matthew},
  journal={Anal. PDE},
  volume={17},
  pages={1953--2010},
  year={2024}
}

@inproceedings{degond1986global,
  title={Global existence of smooth solutions for the {V}lasov-{F}okker-{P}lanck equation in $1 $ and $2 $ space dimensions},
  author={Degond, Pierre},
  booktitle={Ann. Sci. \'Ec. Norm. Super. (4)},
  volume={19},
  pages={519--542},
  year={1986}
}

@article{degond1987existence,
  title={Existence of solutions and diffusion approximation for a model {F}okker-{P}lanck equation},
  author={Degond, Pierre and Mas-Gallic, Sylvie},
  journal={Transport Theory and Statistical Physics},
  volume={16},
  number={4-6},
  pages={589--636},
  year={1987},
  publisher={Taylor \& Francis}
}

@article{baouendi1968equation,
  title={Sur une {\'e}quation d'{\'e}volution changeant de type},
  author={Baouendi, MS and Grisvard, Pierre},
  journal={J. Funct. Anal.},
  volume={2},
  number={3},
  pages={352--367},
  year={1968},
  publisher={Academic Press}
}

@article{chaintron2022propagation1,
  title={Propagation of chaos: a review of models, methods and applications. I. Models and methods},
  author={Chaintron, Louis-Pierre and Diez, Antoine},
  journal={Kinet. Relat. Models},
  volume={15},
  pages={895--1015},
  year={2022},
}

@article{chaintron2021propagation2,
  title={Propagation of chaos: a review of models, methods and applications. II. Applications},
  author={Chaintron, Louis-Pierre and Diez, Antoine},
  journal={Kinet. Relat. Models},
  volume={15},
  pages={1017--1173},
  year={2022},
}

@article{briant2022cauchy,
  title={Cauchy theory for general kinetic {V}icsek models in collective dynamics and mean-field limit approximations},
  author={Briant, Marc and Diez, Antoine and Merino-Aceituno, Sara},
  journal={SIAM J. Math. Anal.},
  volume={54},
  number={1},
  pages={1131--1168},
  year={2022},
  publisher={SIAM}
}

@article{bucataru2011geometric,
  title={A geometric setting for systems of ordinary differential equations},
  author={Bucataru, Ioan and Constantinescu, Oana and Dahl, Matias F},
  journal={Int. J. Geom. Methods Mod. Phys.},
  volume={8},
  number={06},
  pages={1291--1327},
  year={2011},
  publisher={World Scientific}
}

@article{angst2015kinetic,
  title={Kinetic Brownian motion on Riemannian manifolds},
  author={Angst, J{\"u}rgen and Bailleul, Isma{\"e}l and Tardif, Camille},
  journal={Electron. J. Probab.},
  volume={20},
  number={110},
  pages={1--40},
  year={2015},
}

@article{baudoin2016hypocoercive,
  title={Hypocoercive estimates on foliations and velocity spherical {B}rownian motion},
  author={Baudoin, Fabrice and Tardif, Camille},
  journal={Kinet. Relat. Models},
  volume={11},
  number={2},
  pages={1--23},
  year={2018},
}

@article{joergensen1978construction,
  title={Construction of the {B}rownian motion and the {O}rnstein-{U}hlenbeck process in a {R}iemannian manifold on basis of the {G}angolli-{M}c. {K}ean injection scheme},
  author={Joergensen, E},
  journal={Probab. Theory Related Fields},
  volume={44},
  number={1},
  pages={71--87},
  year={1978},
  publisher={Springer}
}

@article{soloveitchik1995focker,
  title={{F}ocker-{P}lanck equation on a manifold. Effective diffusion and spectrum},
  author={Soloveitchik, MR},
  journal={Potential Anal.},
  volume={4},
  pages={571--593},
  year={1995},
  publisher={Springer}
}

@article{grothaus2022hypocoercivity,
  title={Hypocoercivity of {L}angevin-type dynamics on abstract smooth manifolds},
  author={Grothaus, Martin and Mertin, Maximilian Constantin},
  journal={Stochastic Process. Appl.},
  volume={146},
  pages={22--59},
  year={2022},
  publisher={Elsevier}
}

@article{girolami2011riemann,
  title={Riemann manifold {L}angevin and {H}amiltonian {M}onte-{C}arlo methods},
  author={Girolami, Mark and Calderhead, Ben},
  journal={J. R. Stat. Soc. Ser. B. Stat. Methodol.},
  volume={73},
  number={2},
  pages={123--214},
  year={2011},
  publisher={Oxford University Press}
}

@article{lelievre2012langevin,
  title={Langevin dynamics with constraints and computation of free energy differences},
  author={Lelievre, Tony and Rousset, Mathias and Stoltz, Gabriel},
  journal={Math. Comp.},
  volume={81},
  number={280},
  pages={2071--2125},
  year={2012}
}

@article{bismut2005hypoelliptic,
  title={The hypoelliptic {L}aplacian on the cotangent bundle},
  author={Bismut, Jean-Michel},
  journal={J. Amer. Math. Soc.},
  volume={18},
  number={2},
  pages={379--476},
  year={2005}
}

@article{lebeau2005geometric,
  title={Geometric {F}okker-{P}lanck equations.},
  author={Lebeau, Gilles},
  journal={Port. Math.},
  volume={62},
  number={4},
  pages={469--530},
  year={2005},
  publisher={European Mathematical Society Publishing House}
}

@article{calogero2012exponential,
  title={Exponential convergence to equilibrium for kinetic {F}okker-{P}lanck equations},
  author={Calogero, Simone},
  journal={Comm. Partial Differential Equations},
  volume={37},
  number={8},
  pages={1357--1390},
  year={2012},
  publisher={Taylor \& Francis}
}

@book{cercignani2013mathematical,
  title={The mathematical theory of dilute gases},
  author={Cercignani, Carlo and Illner, Reinhard and Pulvirenti, Mario},
  volume={106},
  year={2013},
  publisher={Springer Science \& Business Media}
}

@article{degond2011hydrodynamic,
  title={Hydrodynamic models of self-organized dynamics: derivation and existence theory},
  author={Degond, Pierre and Liu, Jian-Guo and Motsch, Sebastien and Panferov, Vladislav},
  journal={Methods Appl. Anal.},
  volume={20},
  pages={089--114},
  year={2013},
}

@article{aceves2019hydrodynamic,
  title={Hydrodynamic limits for kinetic flocking models of {C}ucker-{S}male type},
  author={Aceves-S{\'a}nchez, Pedro and Bostan, Mihai and Carrillo, Jose-Antonio and Degond, Pierre},
  journal={Math. Biosci. Eng.},
  year={2019},
  volume={16},
  pages={7883--7910},
}

@article{jiang2016hydrodynamic,
  title={Hydrodynamic limits of the kinetic self-organized models},
  author={Jiang, Ning and Xiong, Linjie and Zhang, Teng-Fei},
  journal={SIAM J. Math. Anal.},
  volume={48},
  number={5},
  pages={3383--3411},
  year={2016},
  publisher={SIAM}
}

@article{jiang2024kinetic,
  title={From kinetic flocking model of {C}ucker--{S}male type to self-organized hydrodynamic model},
  author={Jiang, Ning and Luo, Yi-Long and Zhang, Teng-Fei},
  journal={Math. Models Methods Appl. Sci.},
  volume={34},
  number={13},
  pages={2395--2467},
  year={2024},
  publisher={World Scientific}
}

@article{zhang2017local,
  title={A local existence of viscous self-organized hydrodynamic model},
  author={Zhang, Teng-Fei and Jiang, Ning},
  journal={Non linear Anal. Real World Appl.},
  volume={34},
  pages={495--506},
  year={2017},
  publisher={Elsevier}
}

@article{motsch2011numerical,
  title={Numerical simulations of a nonconservative hyperbolic system with geometric constraints describing swarming behavior},
  author={Motsch, Sebastien and Navoret, Laurent},
  journal={Multiscale Model. Simul.},
  volume={9},
  number={3},
  pages={1253--1275},
  year={2011},
  publisher={SIAM}
}

@article{sayin2025behavioral,
  title={The behavioral mechanisms governing collective motion in swarming locusts},
  author={Sayin, Sercan and Couzin-Fuchs, Einat and Petelski, Inga and G{\"u}nzel, Yannick and Salahshour, Mohammad and Lee, Chi-Yu and Graving, Jacob M and Li, Liang and Deussen, Oliver and Sword, Gregory A and others},
  journal={Science},
  volume={387},
  number={6737},
  pages={995--1000},
  year={2025},
  publisher={American Association for the Advancement of Science}
}

@article{lin2025experimental,
  title={Experimental evidence of stress-induced critical state in schooling fish},
  author={Lin, Guozheng and Escobedo, Ram{\'o}n and Li, Xu and Xue, Tingting and Han, Zhangang and Sire, Cl{\'e}ment and Guttal, Vishwesha and Theraulaz, Guy},
  journal={bioRxiv},
  pages={2025--02},
  year={2025},
  publisher={Cold Spring Harbor Laboratory}
}

@article{creppy2016symmetry,
  title={Symmetry-breaking phase transitions in highly concentrated semen},
  author={Creppy, Adama and Plourabou{\'e}, Franck and Praud, Olivier and Druart, Xavier and Cazin, S{\'e}bastien and Yu, Hui and Degond, Pierre},
  journal={Journal of The Royal Society Interface},
  volume={13},
  number={123},
  pages={20160575},
  year={2016},
  publisher={The Royal Society}
}

@article{vicsek2012collective,
  title={Collective motion},
  author={Vicsek, Tam{\'a}s and Zafeiris, Anna},
  journal={Physics reports},
  volume={517},
  number={3-4},
  pages={71--140},
  year={2012},
  publisher={Elsevier}
}

@article{vicsek1995novel,
	author = {Vicsek, Tam{\'a}s and Czir{\'o}k, Andr{\'a}s and Ben-Jacob, Eshel and Cohen, Inon and Shochet, Ofer},
	journal = {Phys. Rev. Lett.},
	number = {6},
	pages = {1226},
	publisher = {APS},
	title = {Novel type of phase transition in a system of self-driven particles},
	volume = {75},
	year = {1995}}

@article{chate2008collective,
	author = {Chat{\'e}, Hugues and Ginelli, Francesco and Gr{\'e}goire, Guillaume and Raynaud, Franck},
	journal = {Phys. Rev. E},
	number = {4},
	pages = {046113},
	publisher = {APS},
	title = {Collective motion of self-propelled particles interacting without cohesion},
	volume = {77},
	year = {2008}}

@article{toner1998flocks,
	author = {Toner, John and Tu, Yuhai},
	journal = {Phys. Rev. E},
	number = {4},
	pages = {4828},
	publisher = {APS},
	title = {Flocks, herds, and schools: A quantitative theory of flocking},
	volume = {58},
	year = {1998}}

@article{barre2019modelling,
  title={Modelling pattern formation through differential repulsion},
  author={Barr{\'e}, Julien and Degond, Pierre and Peurichard, Diane and Zatorska, Ewelina},
  journal={Netw. Heterog. Media},
	pages = {307--352},
	volume = {15},
	year = {2020}
}

@article{gautrais2009analyzing,
  title={Analyzing fish movement as a persistent turning walker},
  author={Gautrais, Jacques and Jost, Christian and Soria, Marc and Campo, Alexandre and Motsch, S{\'e}bastien and Fournier, Richard and Blanco, St{\'e}phane and Theraulaz, Guy},
  journal={J. Math. Biol.},
  volume={58},
  pages={429--445},
  year={2009},
  publisher={Springer}
}

@article{degond2008large,
  title={Large scale dynamics of the persistent turning walker model of fish behavior},
  author={Degond, Pierre and Motsch, S{\'e}bastien},
  journal={J. Stat. Phys.},
  volume={131},
  pages={989--1021},
  year={2008},
  publisher={Springer}
}

@article{gautrais2012deciphering,
  title={Deciphering interactions in moving animal groups},
  author={Gautrais, Jacques and Ginelli, Francesco and Fournier, Richard and Blanco, St{\'e}phane and Soria, Marc and Chat{\'e}, Hugues and Theraulaz, Guy},
  journal={PLoS Comput. Biol.},
  volume={8},
  pages={e1002678},
  year={2012}
}

@article{cavagna2015flocking,
  title={Flocking and turning: a new model for self-organized collective motion},
  author={Cavagna, Andrea and Del Castello, Lorenzo and Giardina, Irene and Grigera, Tomas and Jelic, Asja and Melillo, Stefania and Mora, Thierry and Parisi, Leonardo and Silvestri, Edmondo and Viale, Massimiliano and others},
  journal={J. Stat. Phys.},
  volume={158},
  pages={601--627},
  year={2015},
  publisher={Springer}
}

@article{yang2015hydrodynamics,
  title={Hydrodynamics of turning flocks},
  author={Yang, Xingbo and Marchetti, M Cristina},
  journal={Phys. Rev. Lett.},
  volume={115},
  number={25},
  pages={258101},
  year={2015},
  publisher={APS}
}

@article{bertin2006boltzmann,
  title={Boltzmann and hydrodynamic description for self-propelled particles},
  author={Bertin, Eric and Droz, Michel and Gr{\'e}goire, Guillaume},
  journal={Phys. Rev. E},
  volume={74},
  number={2},
  pages={022101},
  year={2006},
  publisher={APS}
}

@article{cao2020asymptotic,
  title={Asymptotic flocking for the three-zone model},
  author={Cao, Fei and Motsch, Sebastien and Reamy, Alexander and Theisen, Ryan},
  journal={Math. Biosci. Eng},
  volume={17},
  number={6},
  pages={7692--7707},
  year={2020}
}

@article{costanzo2018spontaneous,
  title={Spontaneous emergence of milling (vortex state) in a {V}icsek-like model},
  author={Costanzo, A and Hemelrijk, CK},
  journal={Journal of Physics D: Applied Physics},
  volume={51},
  number={13},
  pages={134004},
  year={2018},
  publisher={IOP Publishing}
}

@article{griette2019kinetic,
  title={Kinetic equations and self-organized band formations},
  author={Griette, Quentin and Motsch, Sebastien},
  journal={Active Particles, Volume 2: Advances in Theory, Models, and Applications},
  pages={173--199},
  year={2019},
  publisher={Springer}
}

@article{barbaro2016phase,
  title={Phase transitions in a kinetic flocking model of {C}ucker--{S}male type},
  author={Barbaro, Alethea BT and Canizo, Jos{\'e} A and Carrillo, Jos{\'e} A and Degond, Pierre},
  journal={Multiscale Model. Simul.},
  volume={14},
  number={3},
  pages={1063--1088},
  year={2016},
  publisher={SIAM}
}

@article{barbaro2012phase,
  title={Phase transition and diffusion among socially interacting self-propelled agents},
  author={Barbaro, Alethea BT and Degond, Pierre},
  journal={Discrete Contin. Dyn. Syst. Ser. B},
  year={2014},
  volume={19},
  pages={1249--1278},
}

@article{carrillo2010asymptotic,
  title={Asymptotic flocking dynamics for the kinetic {C}ucker--{S}male model},
  author={Carrillo, Jos{\'e} A and Fornasier, Massimo and Rosado, Jes{\'u}s and Toscani, Giuseppe},
  journal={SIAM J. Math. Anal.},
  volume={42},
  number={1},
  pages={218--236},
  year={2010},
  publisher={SIAM}
}

@article{ha2009simple,
  title={A simple proof of the {C}ucker-{S}male flocking dynamics and mean-field limit},
  author={Ha, Seung-Yeal and Liu, Jian-Guo},
  journal={Commun. Math. Sci.},
  volume={7},
  pages={297--325},
  year={2009},
}

@article{ha2008particle,
  title={From particle to kinetic and hydrodynamic descriptions of flocking},
  author={Ha, Seung-Yeal and Tadmor, Eitan},
  journal={Kinet. Relat. Mod.},
  volume={1},
  pages={415--435},
  year={2008},
}

@article{haskovec2021simple,
  title={A simple proof of asymptotic consensus in the {H}egselmann--{K}rause and {C}ucker--{S}male models with normalization and delay},
  author={Haskovec, Jan},
  journal={SIAM J. Appl. Dyn. Syst.},
  volume={20},
  number={1},
  pages={130--148},
  year={2021},
  publisher={SIAM}
}

@article{motsch2011new,
  title={A new model for self-organized dynamics and its flocking behavior},
  author={Motsch, Sebastien and Tadmor, Eitan},
  journal={J. Stat. Phys.},
  volume={144},
  pages={923--947},
  year={2011},
  publisher={Springer}
}

@article{d2006self,
  title={Self-propelled particles with soft-core interactions: patterns, stability, and collapse},
  author={D'Orsogna, Maria R and Chuang, Yao-Li and Bertozzi, Andrea L and Chayes, Lincoln S},
  journal={Phys. Rev. Lett.},
  volume={96},
  number={10},
  pages={104302},
  year={2006},
  publisher={APS}
}

@article{bostan2013asymptotic,
  title={Asymptotic fixed-speed reduced dynamics for kinetic equations in swarming},
  author={Bostan, Mihai and Carrillo, Jose Antonio},
  journal={Math. Models Methods Appl. Sci.},
  volume={23},
  number={13},
  pages={2353--2393},
  year={2013},
  publisher={World Scientific}
}

@article{bostan2017reduced,
  title={Reduced fluid models for self-propelled particles interacting through alignment},
  author={Bostan, Mihai and Carrillo, Jose Antonio},
  journal={Math. Models Methods Appl. Sci.},
  volume={27},
  number={07},
  pages={1255--1299},
  year={2017},
  publisher={World Scientific}
}

@article{toner1995long,
  title={Long-range order in a two-dimensional dynamical {X}{Y} model: how birds fly together},
  author={Toner, John and Tu, Yuhai},
  journal={Phys. Rev. Lett.},
  volume={75},
  number={23},
  pages={4326},
  year={1995},
  publisher={APS}
}

@article{toner2005hydrodynamics,
  title={Hydrodynamics and phases of flocks},
  author={Toner, John and Tu, Yuhai and Ramaswamy, Sriram},
  journal={Annals of Physics},
  volume={318},
  number={1},
  pages={170--244},
  year={2005},
  publisher={Elsevier}
}

@article{Aoki1982simulation,
  title={A simulation study on the schooling mechanism in fish.},
  author={Aoki, I.},
  journal={Bull. Jpn. Soc. Sci. Fish},
  volume={48},
  number={8},
  pages={1081--1088},
  year={1982},
  publisher={The Japanese Society of Fisheries Science}
}

@article{cucker2007emergent,
  title={Emergent behavior in flocks},
  author={Cucker, Felipe and Smale, Steve},
  journal={IEEE Trans. Automat. Control},
  volume={52},
  number={5},
  pages={852--862},
  year={2007},
  publisher={IEEE}
}

@article{bertin2009hydrodynamic,
  title={Hydrodynamic equations for self-propelled particles: microscopic derivation and stability analysis},
  author={Bertin, Eric and Droz, Michel and Gr{\'e}goire, Guillaume},
  journal={Journal of Physics A: Mathematical and Theoretical},
  volume={42},
  number={44},
  pages={445001},
  year={2009},
  publisher={IOP Publishing}
}

@article{figalli2018global,
  title={Global well-posedness of the spatially homogeneous {K}olmogorov--{V}icsek model as a gradient flow},
  author={Figalli, Alessio and Kang, Moon-Jin and Morales, Javier},
  journal={Arch. Ration. Mech. Anal.},
  volume={227},
  pages={869--896},
  year={2018},
  publisher={Springer}
}

@article{gamba2016global,
  title={Global weak solutions for {K}olmogorov--{V}icsek type equations with orientational interactions},
  author={Gamba, Irene M and Kang, Moon-Jin},
  journal={Arch. Ration. Mech. Anal.},
  volume={222},
  pages={317--342},
  year={2016},
  publisher={Springer}
}

@article{levermore1996moment,
  title={Moment closure hierarchies for kinetic theories},
  author={Levermore, C David},
  journal={J. Stat. Phys.},
  volume={83},
  pages={1021--1065},
  year={1996},
  publisher={Springer}
}

@book{adams2003sobolev,
  title={Sobolev spaces},
  author={Adams, Robert A and Fournier, John JF},
  volume={140},
  year={2003},
  publisher={Elsevier}
}

\newpage

\appendix

\vspace{0.8cm}
\noindent
\textbf{\Large Appendix}

\vspace{-0.3cm}

\setcounter{equation}{0}
\section{Proofs of remarkable formulas (Section \ref{subsubsec:rem_formulas})}
\label{sec:proofs_rem_formulas}

\subsection{Proof of Lemma \ref{lem:geom_divk_vertlift}}
\label{subsec:proof_lem:geom_divk_vertlift}

{\color{red} Read from here}

For $k \in {\mathcal X}(B)$, we will denote ${\mathcal L}^V \tilde k$ by ${\mathcal L}^V k$ for simplicity, and similarly with $V$ replaced by $H$. Let $(\varepsilon_1, \ldots, \varepsilon_{n-1})$ be a local orthonormal frame on $B$. Then, it is readily checked that $({\mathcal L}^H \varepsilon_1, \ldots,$ ${\mathcal L}^H \varepsilon_{n-1}, {\mathcal L}^V \varepsilon_1, \ldots, {\mathcal L}^V \varepsilon_{n-1})$ is a local orthonormal frame on $M$. Now, using \eqref{eq:pk_div_LCC}, we can write 
\begin{equation}
 \nabla_v \cdot \tilde k = \nabla_\alpha \cdot \big( {\mathcal L}^H k \big) = \sum_{i=1}^{n-1} \Big( \big\langle \hspace{-1.8mm} \big\langle \nabla^{\mathrm{M}}_{{\mathcal L}^H \varepsilon_i} {\mathcal L}^H k, {\mathcal L}^H \varepsilon_i \big\rangle \hspace{-1.8mm} \big\rangle + \big\langle \hspace{-1.8mm} \big\langle \nabla^{\mathrm{M}}_{{\mathcal L}^V \varepsilon_i} {\mathcal L}^H k, {\mathcal L}^V \varepsilon_i \big\rangle \hspace{-1.8mm} \big\rangle
 \Big). 
\label{eq:geom_divk_vertlift_prf-1}
\end{equation}
By the torsonfreeness of the connection, we have 
$$ \nabla^{\mathrm{M}}_{{\mathcal L}^H \varepsilon_i} {\mathcal L}^H k = \nabla^{\mathrm{M}}_{{\mathcal L}^H k} {\mathcal L}^H \varepsilon_i + \big[ {\mathcal L}^H \varepsilon_i, {\mathcal L}^H k \big], $$
and since $ \langle \hspace{-1.2mm} \langle {\mathcal L}^H \varepsilon_i, {\mathcal L}^H \varepsilon_i \rangle \hspace{-1.2mm} \rangle = 1$ and the connection is compatible with the metric, we have 
$$  \big\langle \hspace{-1.8mm} \big\langle \nabla^{\mathrm{M}}_{{\mathcal L}^H k} {\mathcal L}^H \varepsilon_i, {\mathcal L}^H \varepsilon_i \big\rangle \hspace{-1.8mm} \big\rangle = 0, $$
and similarly with ${\mathcal L}^H \varepsilon_i$ replaced by ${\mathcal L}^V \varepsilon_i$. Thus, we get 
$$ \nabla_v \cdot \tilde k = \sum_{i=1}^{n-1} \Big( \big\langle \hspace{-1.8mm} \big\langle \big[ {\mathcal L}^H \varepsilon_i,  {\mathcal L}^H k \big], {\mathcal L}^H \varepsilon_i \big\rangle \hspace{-1.8mm} \big\rangle + \big\langle \hspace{-1.8mm} \big\langle \big[ {\mathcal L}^V \varepsilon_i,  {\mathcal L}^H k \big], {\mathcal L}^V \varepsilon_i \big\rangle \hspace{-1.8mm} \big\rangle
 \Big). 
$$
Now, for any two vector fields $k$, $h \in {\mathcal X}(B)$, we have \cite[Prop. 5.1]{gudmundsson2002geometry}: 
\begin{align}
\big[ {\mathcal L}^H h, {\mathcal L}^H k \big] &= {\mathcal L}^H \big( [h, k] \big) - {\mathcal L}^V \big( {\mathcal R}(h,k) \kappa \big), \label{eq:geom_gudkap_horhor} \\
\big[{\mathcal L}^H h, {\mathcal L}^V k \big] &= {\mathcal L}^V \big( \nabla^{\mathrm{B}}_h k \big), \label{eq:geom_gudkap_horvert} 
\end{align}
where the Levi-Civita connection $\nabla^{\mathrm{B}}$ on $B$ was used and ${\mathcal R}$ is its Riemann curvature tensor. Thus, we have 
\begin{eqnarray*} 
\nabla_v \cdot \tilde k &=& \sum_{i=1}^{n-1} \Big( \big\langle \hspace{-1.8mm} \big\langle {\mathcal L}^H \big( [\varepsilon_i, k] \big) - {\mathcal L}^V \big( {\mathcal R}(\varepsilon_i, k) \kappa \big), {\mathcal L}^H \varepsilon_i  \big \rangle \hspace{-1.8mm} \big\rangle - \big\langle \hspace{-1.8mm} \big\langle {\mathcal L}^V \big( \nabla_k^{\mathrm{B}} \varepsilon_i \big), {\mathcal L}^V \varepsilon_i \big\rangle \hspace{-1.8mm} \big\rangle \Big) \\
 &=& \sum_{i=1}^{n-1} \Big( \big\langle \hspace{-1.8mm} \big\langle {\mathcal L}^H \big( [\varepsilon_i, k] \big), {\mathcal L}^H \varepsilon_i  \big \rangle \hspace{-1.8mm} \big\rangle - \big\langle \hspace{-1.8mm} \big\langle {\mathcal L}^V \big( \nabla_k^{\mathrm{B}} \varepsilon_i \big), {\mathcal L}^V \varepsilon_i \big\rangle \hspace{-1.8mm} \big\rangle \Big), 
\end{eqnarray*}
where the term involving the Riemann curvature tensor has been cancelled due to orthogonality of the horizontal and vertical subspaces. Now, by the definition of the Sasaki metric, we can remove the operators ${\mathcal L}^H$ in the first term and ${\mathcal L}^V$ in the second term provided we replace the metric on $M$ by the metric on $B$. This leads to 
$$ \nabla_v \cdot \tilde k = \sum_{i=1}^{n-1} \Big( \big\langle [\varepsilon_i, k] , \varepsilon_i  \big \rangle - \big\langle \nabla_k^{\mathrm{B}} \varepsilon_i ,\varepsilon_i \big\rangle \Big). $$
We have $\langle \nabla_k^{\mathrm{B}} \varepsilon_i ,\varepsilon_i \rangle = 0$ again because $\langle \varepsilon_i ,\varepsilon_i \rangle = 1$. Using the same manipulations as previously, we also notice that 
$$ \sum_{i=1}^{n-1} \langle [\varepsilon_i, k] , \varepsilon_i \rangle = \sum_{i=1}^{n-1} \langle \nabla_{\varepsilon_i}^{\mathrm{B}} k , \varepsilon_i \rangle = \nabla \cdot k.$$ 
Hence, we have just proved the first formula~\eqref{eq:geom_divk_vertlift}.

\medskip
\noindent
For the second formula \eqref{eq:geom_divk_vertlift}, we proceed analogously. We have 
\begin{eqnarray}
 \nabla_\kappa \cdot \tilde k &=& \nabla_\alpha \cdot \big( {\mathcal L}^V k \big) = \sum_{i=1}^{n-1} \Big( \big\langle \hspace{-1.8mm} \big\langle \nabla^{\mathrm{M}}_{{\mathcal L}^H \varepsilon_i} {\mathcal L}^V k, {\mathcal L}^H \varepsilon_i \big\rangle \hspace{-1.8mm} \big\rangle + \big\langle \hspace{-1.8mm} \big\langle \nabla^{\mathrm{M}}_{{\mathcal L}^V \varepsilon_i} {\mathcal L}^V k, {\mathcal L}^V \varepsilon_i \big\rangle \hspace{-1.8mm} \big\rangle
 \Big). 
\label{eq:geom_divk_vertlift_prf1} \\
&=& \sum_{i=1}^{n-1} \Big( \big\langle \hspace{-1.8mm} \big\langle \big[ {\mathcal L}^H \varepsilon_i,  {\mathcal L}^V k \big], {\mathcal L}^H \varepsilon_i \big\rangle \hspace{-1.8mm} \big\rangle + \big\langle \hspace{-1.8mm} \big\langle \big[ {\mathcal L}^V \varepsilon_i,  {\mathcal L}^V k \big], {\mathcal L}^V \varepsilon_i \big\rangle \hspace{-1.8mm} \big\rangle \Big). \nonumber 
\end{eqnarray}
Now, for any two vector fields $k$, $h \in {\mathcal X}(B)$, we have \cite[Prop. 5.1]{gudmundsson2002geometry}: 
\begin{equation}
[{\mathcal L}^V h, {\mathcal L}^V k] = 0.  
\label{eq:geom_gudkap_vertvert}
\end{equation}
We also recall \eqref{eq:geom_gudkap_horvert} as well as the orthogonality of the horizontal and vertical subspaces. This leads to  
$$ \nabla_\kappa \cdot \tilde k = \sum_{i=1}^{n-1} \big\langle \hspace{-1.8mm} \big\langle {\mathcal L}^V \big( \nabla_{\varepsilon_i}^{\mathrm{B}} k \big) ,  {\mathcal L}^H \varepsilon_i \big\rangle \hspace{-1.8mm} \big\rangle   = 0. $$
Thus proving the second formula \eqref{eq:geom_divk_vertlift} and ending the proof. \endproof

\subsection{Proof of Lemma \ref{lem:geom_nab_fct_v_only}}
\label{subsec:proof_lem:geom_nab_fct_v_only}

\noindent
(i) Let $\alpha \in M$ and $v = \pi \alpha$. Let $a \in T_v B$. We have 
\begin{eqnarray*}
\big\langle \nabla_v \tilde \varphi(\alpha), a \big\rangle_v &=& \big\langle {\mathcal B}_\alpha^H \nabla_\alpha \tilde \varphi(\alpha), a \big\rangle_v = 
\big\langle \hspace{-1.8mm} \big\langle \nabla_\alpha \tilde \varphi(\alpha), {\mathcal L}_\alpha^H a \big\rangle \hspace{-1.8mm} \big\rangle_\alpha = (d \tilde \varphi)_\alpha \big( {\mathcal L}_\alpha^H a  \big) \\
&=& \big( (d \varphi)_{\pi \alpha} \circ (d \pi)_\alpha \big) \big( {\mathcal L}_\alpha^H a  \big) = (d \varphi)_v(a) = \big\langle \nabla_v \varphi(v), a \big\rangle_v, 
\end{eqnarray*}
where we have used \eqref{eq:geom_gam} to say that $(d \pi)_\alpha ({\mathcal L}_\alpha^H a) = a$. 

Similarly we get, for any $\tau \in T_v B$: 
\begin{eqnarray*}
\big\langle \nabla_\kappa \tilde \varphi(\alpha), \tau \big\rangle_v &=& \big\langle {\mathcal B}_\alpha^V \nabla_\alpha \tilde \varphi(\alpha), \tau \big\rangle_v = 
\big\langle \hspace{-1.8mm} \big\langle \nabla_\alpha \tilde \varphi(\alpha), {\mathcal L}_\alpha^V \tau \big\rangle \hspace{-1.8mm} \big\rangle_\alpha = (d \tilde \varphi)_\alpha \big( {\mathcal L}_\alpha^V \tau  \big) \\
&=& \big( (d \varphi)_{\pi \alpha} \circ (d \pi)_\alpha \big) \big( {\mathcal L}_\alpha^V \tau  \big) = 0, 
\end{eqnarray*}
thanks to \eqref{eq:geom_dpioL}. 

\medskip
\noindent
(ii) Let $\alpha = (v,\kappa)$ be given. We define $\theta$: ${\mathbb R} \to {\mathbb R}$ such that $ \theta (t) = \psi(v,\kappa t)$. Then, 
\begin{eqnarray*} \frac{d \theta}{dt}(t) &=& (d \psi)_{(v, \kappa t)} (0,\kappa) = (d \psi)_{(v, \kappa t)} \big( {\mathcal L}_{(v,\kappa t)}^V  \kappa \big) \\
&=& \big\langle \hspace{-1.8mm} \big\langle \nabla_\alpha \psi (v, \kappa t) , {\mathcal L}_{(v,\kappa t)}^V \kappa \big\rangle \hspace{-1.8mm} \big\rangle_{(v,\kappa t)} = \big \langle {\mathcal B}_{(v,\kappa t)}^V \nabla_\alpha \psi (v, \kappa t) , \kappa \big \rangle_v \\
&=& \big \langle \nabla_\kappa \psi (v, \kappa t) , \kappa \big \rangle_v = 0. 
\end{eqnarray*}
Hence, $\theta(0) = \theta(1)$ i.e. $\psi(v,0) = \psi(v,\kappa)$. Defining $\varphi$: $B \to {\mathbb R}$ by $\varphi(v) = \psi(v,0)$, we deduce that $\psi = \varphi \circ \pi$, which ends the proof. \endproof

\subsection{Proof of Lemma \ref{lem:geom_na_kappa2}}
\label{subsec:proof_lem:geom_na_kappa2}

\noindent
Let $\tau \in T_v B$. Then, we have 
\begin{eqnarray}
\big\langle \nabla_\kappa \langle \kappa, \kappa \rangle_v , \tau \big\rangle &=& 
\big\langle {\mathcal B}_\alpha^V \nabla_\alpha \langle \kappa, \kappa \rangle_v , \tau \big\rangle = 
\big\langle \hspace{-1.8mm} \big\langle \nabla_\alpha \langle \kappa, \kappa \rangle_v , {\mathcal L}_\alpha^V \tau \big\rangle \hspace{-1.8mm} \big\rangle_\alpha \noindent \\
&=& \big( d \langle \kappa, \kappa \rangle_v \big)_\alpha \big( {\mathcal L}_\alpha^V \tau \big) = \frac{d}{dt} \langle \tilde \kappa(t), \tilde \kappa(t) \rangle_{\tilde v(t)}\big|_{t=0}, \label{eq:geom_na_kappa2_prf1}
\end{eqnarray}
where $\tilde \alpha(t) = (\tilde v(t), \tilde \kappa(t))$ is the following path in $M$: $\tilde v(t) = v$, $\tilde \kappa(t) = \kappa + t \tau$. Indeed, we notice that $\tilde \alpha(0) = (v,\kappa) = \alpha$ and 
$$ \frac{d \tilde \alpha}{dt}\big|_{t=0} = (0, \tau) = {\mathcal L}_\alpha^V \tau. $$
Then, we straightforwardly get that 
$$ \frac{d}{dt} \langle \tilde \kappa(t), \tilde \kappa(t) \rangle_{\tilde v(t)}\big|_{t=0} = \langle 2 \kappa, \tau \rangle_v, $$
which, inserted into \eqref{eq:geom_na_kappa2_prf1}, leads to the first equation \eqref{eq:geom_na_kappa2}.  

Let now $a \in T_v B$. Similarly to the previous computation, we have 
\begin{eqnarray}
\big\langle \nabla_v \langle \kappa, \kappa \rangle_v , a \big\rangle &=& 
\big\langle {\mathcal B}_\alpha^H \nabla_\alpha \langle \kappa, \kappa \rangle_v , a \big\rangle = 
\big\langle \hspace{-1.8mm} \big\langle \nabla_\alpha \langle \kappa, \kappa \rangle_v , {\mathcal L}_\alpha^H a \big\rangle \hspace{-1.8mm} \big\rangle_\alpha \noindent \\
&=& \big( d \langle \kappa, \kappa \rangle_v \big)_\alpha \big( {\mathcal L}_\alpha^H a \big) = \frac{d}{dt} \langle \tilde \kappa(t), \tilde \kappa(t) \rangle_{\tilde v(t)}\big|_{t=0}, \label{eq:geom_na_kappa2_prf2}
\end{eqnarray}
where $\tilde \alpha(t) = (\tilde v(t), \tilde \kappa(t))$ is a path such that 
\begin{equation}
(\tilde v(0), \tilde \kappa(0)) = (v,\kappa), \qquad \frac{d}{dt}(\tilde v, \tilde \kappa)|_{t=0} = {\mathcal L}_\alpha^H a  = (a, - (a \cdot \kappa) v). 
\label{eq:geom_na_kappa2_prf3}
\end{equation}
Let $(e_1, \ldots, e_n)$ be an orthonormal frame of ${\mathbb R}^n$ such that $v = e_n$, $a = |a| e_{n-1}$ and $\kappa = \kappa_{n-2} e_{n-2} + \kappa_{n-1} e_{n-1}$. For given $\alpha$ and $a \in T_v B$, such a frame always exists. Now, define the path 
$$ \tilde v(t) = \sin (|a|t) e_{n-1} + \cos (|a|t) e_n, \qquad \tilde \kappa(t) = \kappa_{n-2} e_{n-2} + \kappa_{n-1} e_{n-1} - \kappa_{n-1} \tan(|a|t) e_n. $$
We easily see that $\tilde \alpha = (\tilde v, \tilde \kappa)$ is a path in $M$ (i.e. $|\tilde v|=1$, $\tilde \kappa \cdot \tilde v = 0$) defined in a neighborhood of $t=0$ and satisfying \eqref{eq:geom_na_kappa2_prf3}. 
Now, we have 
$$ \langle \tilde \kappa(t), \tilde \kappa(t) \rangle_{\tilde v(t)} = \kappa_{n-2}^2 + \kappa_{n-1}^2 \big( 1 + \tan^2(|a|t) \big). $$
Thus, 
$$ \frac{d}{dt} \langle \tilde \kappa(t), \tilde \kappa(t) \rangle_{\tilde v(t)} = 2 |a| \kappa_{n-1}^2 \tan(|a|t) \big( 1 + \tan^2(|a|t) \big), $$
which leads to 
$$  \frac{d}{dt} \langle \tilde \kappa(t), \tilde \kappa(t) \rangle_{\tilde v(t)}\big|_{t=0} = 0. $$
Inserting it in \eqref{eq:geom_na_kappa2_prf2}, we get the second equation of \eqref{eq:geom_na_kappa2}, which ends the proof. \endproof

\setcounter{equation}{0}
\section{Reduced GCI pair: proof of Theorem \ref{thm:macnorm_nD_actionOn-1}}
\label{sec:redGCIprf}

\subsection{Geometric preliminaries}
\label{subsec:redGCIprf_prelim}

We first remark that $\vec{\psi}^\dagger$  satisfies the following invariance under action of $\mathrm{O}_{n-1}$: 
\begin{equation}
\vec{\psi}^\dagger(\theta, \kappa_\parallel, Rw, R\kappa_T) = R \vec{\psi}^\dagger(\theta, \kappa_\parallel, w, \kappa_T), \quad \forall R \in \mathrm{O}_{n-1}. \label{eq:macnorm_invar_zeta} 
\end{equation}
We remark that $(Rw, R\kappa_T)$ is the image of $(w, \kappa_T)$ by the left group action of $\mathrm{O}_{n-1}$ on $M_T$. We decompose $\vec{\psi}^\dagger$ according to
\begin{equation}
\vec{\psi}^\dagger(\theta, \kappa_\parallel, w, \kappa_T)= \psi_w^\dagger(\theta, \kappa_\parallel, w, \kappa_T) \, w + \vec{\psi}_T^\dagger(\theta, \kappa_\parallel, w, \kappa_T), 
\label{eq:macnorm_zetadag_def1}
\end{equation}
with 
\begin{equation} 
\psi_w^\dagger(\theta, \kappa_\parallel, w, \kappa_T) \in {\mathbb R} \quad \textrm{and} \quad \vec{\psi}_T^\dagger(\theta, \kappa_\parallel, w, \kappa_T) \in T_w {\mathbb S}^{n-2}. 
\label{eq:macnorm_zetadag_def2}
\end{equation}
Straightforwardly, from \eqref{eq:macnorm_invar_zeta}, we get
\begin{eqnarray}
\psi_w^\dagger(\theta, \kappa_\parallel, Rw, R\kappa_T) &=& \psi_w^\dagger(\theta, \kappa_\parallel, w, \kappa_T), \quad \forall R \in \mathrm{O}_{n-1}, \label{eq:macnorm_invar_zetaR} \\
\vec{\psi}_T^\dagger(\theta, \kappa_\parallel, Rw, R\kappa_T) &=& R \vec{\psi}_T^\dagger(\theta, \kappa_\parallel, w, \kappa_T), \quad \forall R \in \mathrm{O}_{n-1}, \label{eq:macnorm_invar_zetaH}
\end{eqnarray}
Thus, we need to prove that 
\begin{equation}
\psi_w^\dagger(\theta, \kappa_\parallel, w, \kappa_T) = \psi_w(\theta, \kappa_\parallel, |\kappa_T|), \quad 
\vec{\psi}_T^\dagger(\theta, \kappa_\parallel, w, \kappa_T) = \psi_r(\theta, \kappa_\parallel, |\kappa_T|) \, \frac{\kappa_T}{|\kappa_T|}. 
\label{eq:macnorm_nD_actionOn-1_prf0}
\end{equation}

\medskip
We use the spherical coordinates on $B$ with pole $\mathbf{u}$. In this section (only) and for the sake of simplicity, we list the canonical basis of ${\mathbb R}^n$ in a different order, namely $(\mathbf{u}, \mathbf{e}_1, \ldots, \mathbf{e}_{n-1})$ with associated coordinates written $(x_u, x_1, \ldots, x_{n-1})$. Instead of the diffeomorphism \eqref{eq:macnorm_V_def}, we consider $\tilde {\mathbf V}_\mathbf{u}$: $(0,\pi)^{n-2} \times {\mathbb R} (\, \textrm{mod. } 2 \pi) \to \tilde B_\mathbf{u}$, \, $(\theta, \varphi_1, \varphi_2, \ldots \varphi_{n-2}) \mapsto v = \tilde {\mathbf V}_\mathbf{u}(\theta, \varphi_1, \varphi_2, \ldots \varphi_{n-2}) \in {\mathbb S}^{n-1}$ 
such that 
$$ \tilde {\mathbf V}_\mathbf{u} = \cos \theta \,  \mathbf{u} + \sin \theta \, w(\varphi_1, \varphi_2, \ldots \varphi_{n-2}), $$
with
$$ w(\varphi_1, \varphi_2, \ldots, \varphi_{n-2}) = \left( \begin{array}{c} 
0 \\
\cos \varphi_1 \\
\sin \varphi_1 \, \cos \varphi_2 \\
\sin \varphi_1 \, \sin \varphi_2 \, \cos \varphi_3\\
\vdots \\
\sin \varphi_1 \, \sin \varphi_2 \ldots \sin \varphi_{n-3} \, \cos \varphi_{n-2} \\ 
\sin \varphi_1 \, \sin \varphi_2 \ldots \sin \varphi_{n-3} \, \sin \varphi_{n-2}  
\end{array} \right) ,$$
and 
$$ \tilde B_\mathbf{u} = B \setminus \big\{ (x_u, x_1, x_2, \ldots, x_{n-3}, 0, 0) \, , | \, \, 
x_u^2 + x_1^2 + \ldots + x_{n-3}^2 = 1 \big\}, $$
being the excluded set corresponding to the points where $\tilde {\mathbf V}_\mathbf{u}$ is not uniquely invertible and thus, is not a diffeomorphism. We deduce that 
$$ \frac{\partial \tilde {\mathbf V}_\mathbf{u}}{\partial \theta} =  - \sin \theta \, \mathbf{u} + \cos \theta \, w =: e_\theta (\theta, \varphi_1, \ldots, \varphi_{n-2}),$$
and for $i \in \{1, \ldots, n-3 \}$, 
\begin{eqnarray*}
 \frac{\partial \tilde {\mathbf V}_\mathbf{u}}{\partial \varphi_i} &=&  \sin \theta \, \sin \varphi_1 \ldots \sin \varphi_{i-1} \, \left( \begin{array}{c} 
0 \\
\vdots \\
0 \\
- \sin \varphi_i \\
\cos \varphi_i \, \cos \varphi_{i+1} \\
\cos \varphi_i \, \sin \varphi_{i+1} \, \cos \varphi_{i+2} \\
\vdots \\
\cos \varphi_i \, \sin \varphi_{i+1}  \ldots \sin \varphi_{n-3} \, \cos \varphi_{n-2} \\
\cos \varphi_i \, \sin \varphi_{i+1}  \ldots \sin \varphi_{n-3} \, \sin \varphi_{n-2}
\end{array} \right) \\
&:=&  \sin \theta \, \sin \varphi_1 \ldots \sin \varphi_{i-1} \, \, e_{\varphi_i}, 
\end{eqnarray*}
where the last $0$ is on the $i$-th line, while, for $i=n-2$, 
$$ \frac{\partial \tilde {\mathbf V}_\mathbf{u}}{\partial \varphi_{n-2}} =  \sin \theta \, \sin \varphi_1 \ldots \sin \varphi_{n-3} \, \left( \begin{array}{c} 
0 \\
\vdots \\
0 \\
- \sin \varphi_{n-2} \\
\cos \varphi_{n-2}
\end{array} \right) :=  \sin \theta \, \sin \varphi_1 \ldots \sin \varphi_{n-3} \, e_{\varphi_{n-2}}, $$
The tuple $(e_{\varphi_1}, \ldots, e_{\varphi_{n-2}})$ form an orthonormal basis of $T_w {\mathbb S}^{n-2}$ while $(e_\theta, e_{\varphi_1}, \ldots, e_{\varphi_{n-2}})$ form an orthonormal basis of $T_v {\mathbb S}^{n-1}$. We note that $e_{\varphi_i}$ only depends on $(\varphi_i, \ldots, \varphi_{n-2})$ but we will still write $e_{\varphi_i}(\varphi_1, \ldots, \varphi_{n-2})$ or, by abuse of notation, $e_{\varphi_i}(w)$. In view of \eqref{eq:macnorm_dV_express}, a tangent vector $\kappa \in T_v B$ can be decomposed according to 
\begin{eqnarray}
&& \hspace{-1.5cm}
\kappa = \kappa_\parallel \, e_\theta(\theta, \varphi_1, \ldots, \varphi_{n-2}) + \sin \theta \, \kappa_T(\varphi_1, \ldots, \varphi_{n-2}, \kappa_{\varphi_1}, \ldots, \kappa_{\varphi_{n-2}}), \nonumber \\ 
&& \hspace{-1.5cm}
\kappa_T(\varphi_1, \ldots, \varphi_{n-2}, \kappa_{\varphi_1}, \ldots, \kappa_{\varphi_{n-2}}) = \sum_{i=1}^{n-2} \kappa_{\varphi_i} \,  e_{\varphi_i}(\varphi_1, \ldots, \varphi_{n-2}) \in T_w({\mathbb S}^2),
\label{eq:macnorm_nD_kappaT_decomp}
\end{eqnarray}
with $\kappa_{\varphi_i} \in {\mathbb R}$, $\forall i \in \{1, \ldots, n-2 \}$. We stress that $(\kappa_{\varphi_1}, \ldots, \kappa_{\varphi_{n-2}})$ denotes the components of $\kappa_T$ in $T_w {\mathbb S}^{n-2}$ and not those of $\kappa$. We  now define $\vec{\psi}^\flat$ such that 
\begin{equation}
\vec{\psi}^\dagger(\theta, \kappa_\parallel, w, \kappa_T) = \vec{\psi}^\flat (\theta, \kappa_\parallel, \varphi_1, \ldots, \varphi_{n-2}, \kappa_{\varphi_1}, \ldots, \kappa_{\varphi_{n-2}}).
\label{eq:macnorm_nD_zeta_coord}
\end{equation}
and define $\psi_w^\flat$, $\vec \psi_T^\flat$ accordingly. Furthermore, we decompose $\vec \psi_T^\flat$ as follows: 
\begin{eqnarray}
&&\hspace{-1cm}
\vec \psi_T^\flat(\theta, \kappa_\parallel, \varphi_1, \ldots, \varphi_{n-2}, \kappa_{\varphi_1}, \ldots, \kappa_{\varphi_{n-2}}) \nonumber \\
&&\hspace{2cm}
= \sum_{k=1}^{n-2} \psi_k^\flat(\theta, \kappa_\parallel, \varphi_1, \ldots, \varphi_{n-2}, \kappa_{\varphi_1}, \ldots, \kappa_{\varphi_{n-2}}) \, e_{\varphi_k}(\varphi_1, \ldots, \varphi_{n-2}), 
\label{eq:macnorm_nD_zetaTflat_decomp}
\end{eqnarray}
with $\psi_k^\flat(\theta, \kappa_\parallel, \varphi_1, \ldots, \varphi_{n-2}, \kappa_{\varphi_1}, \ldots, \kappa_{\varphi_{n-2}}) \in {\mathbb R}$.

\medskip
We will use \eqref{eq:macnorm_invar_zetaR}, \eqref{eq:macnorm_invar_zetaH}, with $R$ ranging in a generating set of $\mathrm{O}_{n-1}$. We remind that we identify $\mathrm{O}_{n-1}$ with the isotropy subgroup of $\mathbf{u}$ in $\mathrm{O}_n$, meaning that all orthogonal transformations of ${\mathbb R}^{n-1}$ are meant being orthogonal transformations of ${\mathbb R}^n$ that leave $\mathbf{u}$ invariant. We will mainly use three families of orthogonal transformation, which are successively presented in the following subsections.

\subsection{Action of orthogonal transformations}
\label{subsec:redGCIprf_actortho}

\subsubsection{First family of orthogonal transformations}
\label{subsubsec:gci_redform_first}

Fix $(\theta, \varphi_1, \ldots, \varphi_j, \ldots, \varphi_{n-2}) \in (0,\pi)^{n-2} \times {\mathbb R} (\, \textrm{mod. } 2 \pi)$. Define the rotation $R_{j,\beta}$ of ${\mathbb R}^{n-1}$ which maps $(w, e_{\varphi_1}, \ldots, e_{\varphi_{n-2}})(\theta, \varphi_1, \ldots, \varphi_j, \ldots, \varphi_{n-2})$ to $(w,$ $e_{\varphi_1}, \ldots, e_{\varphi_{n-2}})(\theta, \varphi_1, \ldots, \varphi_j + \beta, \ldots, \varphi_{n-2})$, for some $j \in \{1, \ldots, n-2 \}$ and some $\beta \in {\mathbb R}$ such that $(\theta, \varphi_1, \ldots, \varphi_j + \beta, \ldots, \varphi_{n-2}) \in (0,\pi)^{n-2} \times {\mathbb R} (\, \textrm{mod. } 2 \pi)$. This rotation, which maps a direct orthonormal basis of ${\mathbb R}^{n-1}$ to another one, is perfectly defined. For $\kappa \in T_v B$ written as in \eqref{eq:macnorm_nD_kappaT_decomp}, we have 
\begin{eqnarray*}
&& \hspace{-1cm}
R_{j,\beta} \kappa = \kappa_\parallel e_\theta(\theta, \varphi_1, \ldots, \varphi_j + \beta, \ldots, \varphi_{n-2}) + \sin \theta \, R_{j,\beta} \kappa_T, \\ 
&& \hspace{-1cm}
R_{j,\beta} \kappa_T = \sum_{i=1}^{n-2} \kappa_{\varphi_i} e_{\varphi_i}(\varphi_1, \ldots, \varphi_j + \beta, \ldots, \varphi_{n-2}).  
\end{eqnarray*}
Therefore, $(\kappa_\parallel, \kappa_{\varphi_1}, \ldots, \kappa_{\varphi_{n-2}})$ are unchanged through the action of $R_{j,\beta}$. 

\begin{itemize}

\item {\em Action of $R_{j,\beta}$ on $\psi_w^\flat$. }
We have
$$ \psi_w^\dagger(\theta, \kappa_\parallel, R_{j,\beta}w, R_{j,\beta}\kappa_T) =
\psi_w^\flat (\theta, \kappa_\parallel, \varphi_1, \ldots, \varphi_j + \beta, \ldots, \varphi_{n-2}, \kappa_{\varphi_1}, \ldots, \kappa_{\varphi_{n-2}}). $$ 
Thus, thanks to \eqref{eq:macnorm_invar_zetaR}, we get that $\psi_w^\flat$ does not depend on $\varphi_j$, $\forall j \in \{1, \ldots, n-2 \}$, i.e. 
\begin{equation}
\psi_w^\flat = \psi_w^\flat (\theta, \kappa_\parallel, \kappa_{\varphi_1}, \ldots, \kappa_{\varphi_{n-2}}). 
\label{eq:macnorm_nD_actionOn-1_prf01}
\end{equation}

\item {\em Action of $R_{j,\beta}$ on $\vec \psi_T^\dagger$.}
On the one hand, we have
\begin{eqnarray*}
&&\hspace{-1cm}
\vec \psi_T^\dagger(R_{j,\beta}v,R_{j,\beta}\kappa) =
\vec \psi_T^\flat (\theta, \kappa_\parallel, \varphi_1, \ldots, \varphi_j + \beta, \ldots, \varphi_{n-2}, \kappa_{\varphi_1}, \ldots, \kappa_{\varphi_{n-2}}) \\
&&\hspace{-1cm}
= \sum_{k=1}^{n-2} \psi_k^\flat(\theta, \kappa_\parallel, \varphi_1, \ldots, \varphi_j + \beta, \ldots, \varphi_{n-2}, \kappa_{\varphi_1}, \ldots, \kappa_{\varphi_{n-2}}) \, e_{\varphi_k}(\varphi_1, \ldots, \varphi_j + \beta, \ldots, \varphi_{n-2}),
\end{eqnarray*}
and on the other hand 
\begin{eqnarray*}
&&\hspace{-1cm}
R_{j,\beta} \vec \psi_T^\dagger(v,\kappa) = 
R_{j,\beta}  \vec \psi_T^\flat (\theta, \kappa_\parallel, \varphi_1, \ldots, \varphi_j, \ldots, \varphi_{n-2}, \kappa_{\varphi_1}, \ldots, \kappa_{\varphi_{n-2}}) \\
&&\hspace{0cm}
= \sum_{k=1}^{n-2} \psi_k^\flat(\theta, \kappa_\parallel, \varphi_1, \ldots, \varphi_j, \ldots, \varphi_{n-2}, \kappa_{\varphi_1}, \ldots, \kappa_{\varphi_{n-2}}) \, e_{\varphi_k}(\varphi_1, \ldots, \varphi_j + \beta, \ldots, \varphi_{n-2}).
\end{eqnarray*}
Thus, thanks to \eqref{eq:macnorm_invar_zetaH}, we get that for all $k$, $\psi_k^\flat$ is independent of $\varphi_j$, $\forall j \in \{1, \ldots, n-2 \}$, i.e. 
\begin{equation}
\psi_k^\flat = \psi_k^\flat (\theta, \kappa_\parallel, \kappa_{\varphi_1}, \ldots, \kappa_{\varphi_{n-2}}), \quad \forall k \in \{1, \ldots, n-2 \}. 
\label{eq:macnorm_nD_actionOn-1_prf02}
\end{equation}

\end{itemize}

\subsubsection{Second family of orthogonal transformations}
\label{subsubsec:gci_redform_second}

We return to $\vec \psi_T^\dagger$ and when $\kappa_T \not = 0$, we decompose 
\begin{equation} 
\vec \psi_T^\dagger (\theta, \kappa_\parallel, w, \kappa_T) = \psi_r^\dagger (\theta, \kappa_\parallel, w, \kappa_T) \frac{\kappa_T}{|\kappa_T|} + \vec \psi_C^\dagger (\theta, \kappa_\parallel, w, \kappa_T), \quad \vec \psi_C^\dagger (\theta, \kappa_\parallel, w, \kappa_T) \cdot \kappa_T = 0. 
\label{eq:macnorm_nD_actionOn-1_prf03}
\end{equation}
Taking the inner product of the first equation \eqref{eq:macnorm_nD_actionOn-1_prf03} with $\kappa_T$, we note that  
\begin{eqnarray} 
\psi_r^\dagger \Big(\theta, \kappa_\parallel, w, \sum_{m=1}^{n-2} \kappa_{\varphi_m} e_{\varphi_m}(w) \Big) &=& \frac{1}{|\kappa_T|}\sum_{k=1}^{n-2} \kappa_{\varphi_k} \, \psi_k^\flat (\theta, \kappa_\parallel, \kappa_{\varphi_1}, \ldots, \kappa_{\varphi_{n-2}}) \nonumber \\
&=:& \psi_r^\flat (\theta, \kappa_\parallel, \kappa_{\varphi_1}, \ldots, \kappa_{\varphi_{n-2}}) , 
\label{eq:macnorm_nD_actionOn-1_prf03.25}
\end{eqnarray}
thus defining $\psi_r^\flat$ and emphasizing that $\psi_r^\flat$ does not depend on $\varphi_1, \ldots, \varphi_{n-2}$. Now, we fix $\kappa_{T0} \not = 0$ given in $T_w {\mathbb S}^{n-2}$ and decompose $\kappa_{T0} = \sum_{k=1}^{n-2} (\kappa_0)_{\varphi_k} e_{\varphi_k}$. Suppose, without loss of generality,  that $(\kappa_0)_{\varphi_1} >0$ (if $(\kappa_0)_{\varphi_1} < 0$, we will just change the sign of $f_1$ below; if $(\kappa_0)_{\varphi_1} =0$, we just pick up a non-zero component $(\kappa_0)_{\varphi_j}$ as there always exists one since $\kappa_T \not = 0$). Then, define $f_1 = \kappa_{T0}/|\kappa_{T0}|$ and construct an orthonormal basis $(f_1, f_2, \ldots, f_{n-2})$ of $T_w {\mathbb S}^{n-2}$ by applying the Gram-Schmidt orthonormalization principle to the basis $(f_1, e_{\varphi_2}, \ldots, e_{\varphi_{n-2}})$. It is readily seen that $(w, f_1, \ldots, f_{n-2})$ is a direct orthonormal basis of ${\mathbb R}^{n-1}$. Then, any vector $\kappa_T \in T_w {\mathbb S}^{n-2}$ can be decomposed in this basis according to $\kappa_T = \sum_{k=1}^{n-2} \tilde \kappa_i f_i$ and $(\tilde \kappa_1, \ldots, \tilde \kappa_{n-2})$ depends smoothly of $(\kappa_{\varphi_1}, \ldots, \kappa_{\varphi_{n-2}})$. In particular, $\kappa_{T0}$ is decomposed into $\kappa_{T0} = \sum_{k=1}^{n-2} (\tilde  \kappa_0)_i f_i$ with $(\tilde  \kappa_0)_1= |\kappa_{T0}|$ and $(\tilde  \kappa_0)_i = 0$ for $i \in \{2, \ldots, n-2\}$. Then, we can define a change of variables of $\vec \psi_T^\dagger (\theta, \kappa_\parallel, w, \cdot)$ as follows (we stress that $(\theta, \kappa_\parallel, w)$ are kept fixed): 
$$ \vec \psi_T^\dagger (\theta, \kappa_\parallel, w, \kappa_T) = \vec \psi_T^\sharp (\tilde \kappa_1, \ldots, \tilde \kappa_{n-2}),  $$
and we define $\psi_r^\sharp$ and $\vec \psi_C^\sharp$ accordingly. Now, we decompose $\vec \psi_C^\sharp$ in the basis $(f_1, \ldots, f_{n-2})$: 
\begin{equation} 
\vec \psi_C^\sharp(\tilde \kappa_1, \ldots, \tilde \kappa_{n-2}) = \sum_{k=1}^{n-2} \psi_k^\sharp(\tilde \kappa_1, \ldots, \tilde \kappa_{n-2}) \, f_k.  
\label{eq:macnorm_nD_actionOn-1_prf03.5}
\end{equation}
We note that the second condition of \eqref{eq:macnorm_nD_actionOn-1_prf03} implies that 
\begin{equation} 
\sum_{k=1}^{n-2} \psi_k^\sharp(\tilde \kappa_1, \ldots, \tilde \kappa_{n-2}) \, \tilde \kappa_k = 0. 
\label{eq:macnorm_nD_actionOn-1_prf04}
\end{equation}
Hence, 
$$ \vec \psi_T^\sharp (\tilde \kappa_1, \ldots, \tilde \kappa_{n-2}) = \sum_{k=1}^{n-2} \Big[ \psi_r^\sharp (\tilde \kappa_1, \ldots, \tilde \kappa_{n-2}) \, \frac{\tilde \kappa_k}{|\kappa_T|} + \psi_k^\sharp (\tilde \kappa_1, \ldots, \tilde \kappa_{n-2}) \Big] \, f_k. $$

Now, we apply \eqref{eq:macnorm_invar_zetaH} with the orthogonal reflection $S_j$ in the hyperplane $\textrm{Span} \{w, f_1, \ldots, \hat f_j,$ $\ldots, f_{n-2} \}$ of ${\mathbb R}^{n-1}$ where $\hat f_j$ means that $f_j$ is omitted in the list and $j \in \{2, \ldots, n-2\}$. Since $S_j w = w$, the angles $(\varphi_1, \ldots, \varphi_{n-2})$ are unmodified by the operation of $S_j$. On the other hand, 
$$ S_j \kappa_T = \sum_{i \not = j} \tilde \kappa_i \, f_i - \tilde \kappa_j \, f_j, $$
meaning that $(\tilde \kappa_1, \ldots, \tilde  \kappa_{j-1}, \tilde \kappa_j, \tilde \kappa_{j+1}, \ldots, \tilde \kappa_{n-2})$ is changed into $(\tilde \kappa_1, \ldots, \tilde \kappa_{j-1}, - \tilde  \kappa_j, \tilde \kappa_{j+1}, \ldots, \tilde \kappa_{n-2})$. We have 
\begin{eqnarray*}
&&\hspace{-1cm}
\phantom{\sum_{k \not = j}} \vec \psi_T^\dagger(\theta, \kappa_\parallel, S_j w, S_j\kappa_T) = \vec \psi_T^\dagger(\theta, \kappa_\parallel, w, S_j\kappa_T) =  \vec \psi_T^\sharp (\tilde \kappa_1, \ldots, - \tilde  \kappa_j,  \ldots, \tilde \kappa_{n-2}) \\
&&\hspace{0.5cm}
= \sum_{k \not = j} \Big( \psi_r^\sharp(\tilde \kappa_1, \ldots, - \tilde  \kappa_j,  \ldots, \tilde \kappa_{n-2}) \, \frac{\tilde \kappa_k}{|\kappa_T|} + \psi_k^\sharp(\tilde \kappa_1, \ldots, - \tilde  \kappa_j,  \ldots, \tilde \kappa_{n-2}) \Big) f_k \\
&&\hspace{2cm} 
+ \Big( - \psi_r^\sharp(\tilde \kappa_1, \ldots, - \tilde  \kappa_j,  \ldots, \tilde \kappa_{n-2}) \, \frac{\tilde \kappa_j}{|\kappa_T|} + \psi_j^\sharp(\tilde \kappa_1, \ldots, - \tilde  \kappa_j,  \ldots, \tilde \kappa_{n-2}) \Big) f_j. 
\end{eqnarray*}
On the other hand, we have
\begin{eqnarray*}
&&\hspace{-1cm}
\phantom{\sum_{k \not = j}} S_j \vec \psi_T^\dagger(\theta, \kappa_\parallel, w, \kappa_T) = S_j \vec \psi_T^\sharp (\tilde \kappa_1, \ldots, \tilde  \kappa_j,  \ldots, \tilde \kappa_{n-2}) \\
&&\hspace{1cm}
 =\sum_{k\not = j} \Big( \psi_r^\sharp(\tilde \kappa_1, \ldots, \tilde \kappa_{n-2}) \, \frac{\tilde \kappa_k}{|\kappa_T|} + \psi_k^\sharp(\tilde \kappa_1, \ldots, \tilde \kappa_{n-2}) \Big) f_k \\
&&\hspace{3cm}
- \Big( \psi_r^\sharp(\tilde \kappa_1, \ldots, \kappa_{n-2}) \, \frac{\tilde \kappa_j}{|\kappa_T|} + \psi_j^\sharp(\tilde \kappa_1, \ldots, \tilde \kappa_{n-2}) \Big) f_j. 
\end{eqnarray*}
Now, we evaluate these two expressions at $\kappa_{T0}$. By construction (see \eqref{eq:macnorm_nD_actionOn-1_prf04}), we have 
$$ \kappa_{T0} = (\tilde  \kappa_0)_1 \,  f_1, \qquad \psi_1^\sharp((\tilde  \kappa_0)_1, 0 , \ldots, 0) = 0, $$
with $(\tilde  \kappa_0)_1 = |\kappa_{T0}|$, and we remind that $j \in \{2, \ldots, n-2 \}$. Hence, thanks to \eqref{eq:macnorm_invar_zetaH}, we get 
\begin{eqnarray*}
&&\hspace{-1cm}
\psi_r^\sharp((\tilde  \kappa_0)_1, 0,  \ldots, 0)  \, \frac{(\tilde  \kappa_0)_1}{|\kappa_T|} +  \sum_{k=1}^{n-2} \psi_k^\sharp((\tilde  \kappa_0)_1, 0,  \ldots, 0) \, f_k 
= \psi_r^\sharp((\tilde  \kappa_0)_1, 0,  \ldots, 0)  \, \frac{(\tilde  \kappa_0)_1}{|\kappa_T|} \\
&&\hspace{4cm}
+  \sum_{k \not = j} \psi_k^\sharp((\tilde  \kappa_0)_1, 0,  \ldots, 0) \, f_k - \psi_j^\sharp((\tilde  \kappa_0)_1, 0,  \ldots, 0) \, f_j. 
\end{eqnarray*}
Hence, $ \psi_j^\sharp((\tilde  \kappa_0)_1, 0,  \ldots, 0) = 0$, $\forall j \in \{2, \ldots, n-2\}$. With \eqref{eq:macnorm_nD_actionOn-1_prf03.5}, \eqref{eq:macnorm_nD_actionOn-1_prf04}, it follows that $\psi_C^\sharp((\tilde  \kappa_0)_1, 0,$  $\ldots, 0) = 0$, i.e. $ \vec \psi_C^\dagger(\theta, \kappa_\parallel, w, \kappa_{T0}) = 0$. Since, this equality is validy for all $\kappa_{T0} \not = 0$, we conclude that 
$$ \vec \psi_C^\dagger (\theta, \kappa_\parallel, w, \kappa_T) = 0, \quad \forall \kappa_T \in T_w {\mathbb S}^{n-2}, \quad \kappa_T \not = 0. $$
Hence, we can write 
\begin{equation}
\vec \psi_T^\dagger (\theta, \kappa_\parallel, w, \kappa_T) = \psi_r^\dagger (\theta, \kappa_\parallel, w, \kappa_T) \, \frac{\kappa_T}{|\kappa_T|}, \quad \forall \kappa_T \in T_w {\mathbb S}^{n-2}, \quad \kappa_T \not = 0,  
\label{eq:macnorm_nD_actionOn-1_prf05}
\end{equation}
with $\psi_r^\flat$ not depending on the angles $(\varphi_1, \ldots, \varphi_{n-2})$ as shown by \eqref{eq:macnorm_nD_actionOn-1_prf03.25}. 

In the case where $\kappa_T = 0$, then, for all element of $\mathrm{O}_{n-1}$ such that $R w = w$, \eqref{eq:macnorm_invar_zetaH} gives
$$ \vec \psi_T^\dagger (\theta, \kappa_\parallel, w, 0) = R \vec \psi_T^\dagger (\theta, \kappa_\parallel, w, 0). $$
Such orthogonal transformations span $\mathrm{O}_{n-2}$, the set of orthogonal transformations of $T_w {\mathbb S}^{n-2}$ and $\vec \psi_T^\dagger (\theta, \kappa_\parallel, w, 0)$ is an element of $T_w {\mathbb S}^{n-2}$ which is a fixed point of all elements of $\mathrm{O}_{n-2}$. Such element must be equal to zero. Hence we get 
\begin{equation}
\vec \psi_T^\dagger (\theta, \kappa_\parallel, w, 0) = 0. 
\label{eq:macnorm_nD_actionOn-1_prf055}
\end{equation}

\subsubsection{Third family of orthogonal transformations}
\label{subsubsec:gci_redform_third}

We now use the planar rotation $R_{j ,\ell ,\beta}$ of ${\mathbb R}^{n-1}$ of angle $\beta$ in the plane $\textrm{Span}\{e_{\varphi_j},e_{\varphi_\ell}\}$. In this rotation, no angle is changed because $w$ is unchanged, being orthogonal to both $e_{\varphi_j}$ and $e_{\varphi_\ell}$. For this reason, we drop the dependences of $e_\theta$, $w$ and $e_{\varphi_i}$ on $(\theta, \varphi_1, \ldots, \varphi_{n-2})$ since no confusion is possible. Now, we get 
$$ R_{j, \ell ,\beta} e_{\varphi_j}  = \cos \beta \, e_{\varphi_j} + \sin \beta \, e_{\varphi_\ell}, \quad  R_{j, \ell, \beta} e_{\varphi_\ell}  = - \sin \beta \, e_{\varphi_j} + \cos \beta \, e_{\varphi_\ell}. $$
Hence , For $\kappa \in T_v B$ written as in \eqref{eq:macnorm_nD_kappaT_decomp}, we have 
\begin{eqnarray*}
&& \hspace{-1cm}
R_{j,\ell,\beta} \kappa = \kappa_\parallel e_\theta + \sin \theta \, R_{j,\ell,\beta} \kappa_T, \\ 
&& \hspace{-1cm}
R_{j,\ell,\beta} \kappa_T = \sum_{i\not \in \{j,\ell\}} \kappa_{\varphi_i} e_{\varphi_i} + (\kappa_{\varphi_j} \, \cos \beta - \kappa_{\varphi_\ell} \, \sin \beta) \, e_{\varphi_j} + (\kappa_{\varphi_j} \, \sin \beta + \kappa_{\varphi_\ell} \, \cos \beta) \, e_{\varphi_\ell}.  
\end{eqnarray*}
Therefore, $(\kappa_{\varphi_1}, \ldots , \kappa_{\varphi_{n-2}})$ is changed into $(\kappa_{\varphi_1}, \ldots , \kappa_{\varphi_j} \, \cos \beta - \kappa_{\varphi_\ell} \, \sin \beta, \ldots, \kappa_{\varphi_j} \, \sin \beta + \kappa_{\varphi_\ell} \, \cos \beta,$ $\ldots, \kappa_{\varphi_{n-2}})$. We return to the coordinates \eqref{eq:macnorm_nD_zeta_coord} which define $\psi_w^\flat(\theta, \kappa_\parallel, \kappa_{\varphi_1}, \ldots, \kappa_{\varphi_{n-2}})$ and $\psi_r^\flat(\theta, \kappa_\parallel,$ $\kappa_{\varphi_1}, \ldots, \kappa_{\varphi_{n-2}})$.

\begin{itemize}

\item {\em Action of $R_{j,\ell\beta}$ on $\psi_w^\flat$.} Eq. \eqref{eq:macnorm_invar_zetaR} gives
\begin{equation}
\psi_w^\dagger(\theta, \kappa_\parallel, R_{j,\ell,\beta}w, R_{j,\ell,\beta}\kappa_T) = \psi_w^\dagger(\theta, \kappa_\parallel, w, \kappa_T). 
\label{eq:macnorm_nD_actionOn-1_prf06}
\end{equation}
We have
\begin{eqnarray*} 
&&\hspace{-1cm}
\psi_w^\dagger(\theta, \kappa_\parallel, R_{j,\ell,\beta}w, R_{j,\ell,\beta}\kappa_T)\\
&&\hspace{1cm}
= \psi_w^\flat (\theta, \kappa_\parallel, \kappa_{\varphi_1}, \ldots , \kappa_{\varphi_j} \, \cos \beta - \kappa_{\varphi_\ell} \, \sin \beta, \ldots, \kappa_{\varphi_j} \, \sin \beta + \kappa_{\varphi_\ell} \, \cos \beta, \ldots, \kappa_{\varphi_{n-2}}). 
\end{eqnarray*} 
Thus, \eqref{eq:macnorm_nD_actionOn-1_prf06} leads to 
\begin{eqnarray} 
&&\hspace{-2cm}
\psi_w^\flat (\theta, \kappa_\parallel, \kappa_{\varphi_1}, \ldots , \kappa_{\varphi_j} \, \cos \beta - \kappa_{\varphi_\ell} \, \sin \beta, \ldots, \kappa_{\varphi_j} \, \sin \beta + \kappa_{\varphi_\ell} \, \cos \beta, \ldots, \kappa_{\varphi_{n-2}}) \nonumber \\
&&\hspace{2cm}
=\psi_w^\flat (\theta, \kappa_\parallel, \kappa_{\varphi_1}, \ldots , \kappa_{\varphi_j}, \ldots, \kappa_{\varphi_\ell}, \ldots, \kappa_{\varphi_{n-2}}), \quad \forall \beta \in {\mathbb R}. \label{eq:macnorm_nD_actionOn-1_prf1}
\end{eqnarray} 
Define 
$$ \kappa_{j \ell} = (\kappa_{\varphi_j}^2 + \kappa_{\varphi_\ell}^2)^{1/2}, \quad \omega_j = \frac{\kappa_{\varphi_j}}{\kappa_{j \ell}} = : \cos \gamma, \quad \omega_\ell = \frac{\kappa_{\varphi_\ell}}{\kappa_{j \ell}} =: \sin \gamma, $$
with $\gamma \in {\mathbb R}$. Then, \eqref{eq:macnorm_nD_actionOn-1_prf1} can be written
\begin{eqnarray*} 
&&\hspace{-2cm} \psi_w^\flat (\theta, \kappa_\parallel, \kappa_{\varphi_1}, \ldots , \kappa_{j \ell} \, \cos (\gamma + \beta), \ldots, \kappa_{j \ell} \, \sin (\gamma + \beta), \ldots, \kappa_{\varphi_{n-2}}) \\
&&\hspace{2cm}
= \psi_w^\flat (\theta, \kappa_\parallel, \kappa_{\varphi_1}, \ldots , \kappa_{\varphi_j}, \ldots, \kappa_{\varphi_\ell}, \ldots, \kappa_{\varphi_{n-2}}), \quad \forall \beta \in {\mathbb R}. 
\end{eqnarray*}
Taking $\beta = - \gamma$, we get 
\begin{eqnarray*} 
&&\hspace{-2cm} \psi_w^\flat (\theta, \kappa_\parallel, \kappa_{\varphi_1}, \ldots , \kappa_{\varphi_j}, \ldots, \kappa_{\varphi_\ell}, \ldots, \kappa_{\varphi_{n-2}}) \\
&&\hspace{3cm}
= \psi_w^\flat (\theta, \kappa_\parallel, \kappa_{\varphi_1}, \ldots , \kappa_{j \ell}, \ldots,0, \ldots, \kappa_{\varphi_{n-2}}). 
\end{eqnarray*}
Proceeding successively with $(j,\ell) = (n-3,n-2)$, $(n-4, n-3)$, \ldots, $(2,3)$, $(1,2)$, we get 
$$ \psi_w^\flat (\theta, \kappa_\parallel, \kappa_{\varphi_1}, \ldots , \kappa_{\varphi_j}, \ldots, \kappa_{\varphi_\ell}, \ldots, \kappa_{\varphi_{n-2}}) 
= \psi_w^\flat (\theta, \kappa_\parallel, |\kappa_T|, 0, \ldots , 0). $$
Therefore, defining $\psi_w$ by,  
$$\psi_w (\theta, \kappa_\parallel, |\kappa_T|) = \psi_w^\flat (\theta, \kappa_\parallel, |\kappa_T|, 0, \ldots , 0), $$
we get the first equality \eqref{eq:macnorm_nD_actionOn-1_prf0}.

\item {\em Action of $R_{j,\ell, \beta}$ on $\psi_r^\flat$.}

Using \eqref{eq:macnorm_nD_actionOn-1_prf05},  we have 
$$ \vec \psi_T^\dagger(\theta, \kappa_\parallel, R_{j,\ell,\beta}w, R_{j,\ell,\beta} \kappa_T) 
= \psi_r^\dagger(\theta, \kappa_\parallel, R_{j,\ell,\beta}w, R_{j,\ell,\beta} \kappa_T)
\, \frac{R_{j,\ell,\beta} \kappa_T}{|R_{j,\ell,\beta} \kappa_T|} . $$
On the other hand we have  
$$ R_{j,\ell,\beta} \vec \psi_T^\dagger(\theta, \kappa_\parallel, w, \kappa_T) = 
\psi_r^\dagger(\theta, \kappa_\parallel, w, \kappa_T) \frac{R_{j,\ell,\beta} \kappa_T}{|\kappa_T|}. $$
Thus, thanks to \eqref{eq:macnorm_invar_zetaH} and the fact that $|R_{j,\ell,\beta} \kappa_T| = |\kappa_T|$, we get 
$$ \psi_r^\dagger(\theta, \kappa_\parallel, R_{j,\ell,\beta}w, R_{j,\ell,\beta} \kappa_T) = \psi_r^\dagger(\theta, \kappa_\parallel, w, \kappa_T), $$
which is the same equation as \eqref{eq:macnorm_nD_actionOn-1_prf06} for $\psi_w^\dagger$. Therefore, the same conclusion applies. Hence, with 
$$\psi_r (\theta, \kappa_\parallel, |\kappa_T|) = \psi_r^\flat (\theta, \kappa_\parallel, |\kappa_T|, 0, \ldots , 0), $$
we get the second equality \eqref{eq:macnorm_nD_actionOn-1_prf0}. 
\end{itemize}

Conversely, we easily check that \eqref{eq:macnorm_nD_actionOn-1} satisfies the invariance relation \eqref{eq:macnorm_gract_zeta_invar}. This ends the proof of Theorem \ref{thm:macnorm_nD_actionOn-1}. \endproof

\setcounter{equation}{0}
\section{Auxilliary computations for Subsection \ref{subsec:macnorm_nD_varform_redGCI}}
\label{sec:macnorm_nD_vectGCI_eqs}

\textbf{In all this section, we assume that $\vec{\psi} \in {\mathcal D}_\mathrm{inv}(M)$}. 

\subsection{Main statements}
\label{subsec:auxil_main}

We let $(\mathbf{e}_1, \ldots, \mathbf{e}_{n-1}, \mathbf{u})$ be the canonical basis of ${\mathbb R}^n$ ordered in the natural order (which is different from the order used in the proof of Theorem \ref{thm:macnorm_nD_actionOn-1}). Thus, $\{\mathbf{u} \}^\bot = \textrm{Span}(\mathbf{e}_1, \ldots, \mathbf{e}_{n-1})$ and we denote by $\psi_i = \vec{\psi}  \cdot \mathbf{e}_i$. From Section \ref{subsec:macnorm_vectorGCI}, we know that $\kappa_T \in \{\mathbf{u} \}^\bot$ and we define $\kappa_{T i} = \kappa_T \cdot \mathbf{e}_i$, $i=1, \ldots, n-1$ (note that $\kappa_{T i}$ is different from $\kappa_T^i$ defined in Subsection \ref{subsubsec:macnorm_chgvar_M}). Finally, we obviously have $\mathbf{u}_i =: \mathbf{u} \cdot \mathbf{e}_i = 0$. Now, we decompose $\vec{\psi}$ according to \eqref{eq:macnorm_vecpsi_rot}. We first state the following auxilliary two lemmas, which will be proved in Subsections \ref{subsec:proof_lem:macnorm_nD_nathkapthkapt} and \ref{subsec:proof_lem:macnorm_nD_nawikapi} respectively.

\begin{lemma}[gradients of $\theta$, $\kappa_\parallel$ and $|\kappa_T|$]~

\noindent
We have 
\begin{eqnarray}
&&\hspace{-1.5cm}
\nabla_v \theta = e_\theta,  
\quad \nabla_v \kappa_\parallel = \cos \theta \, \kappa_T, 
\quad \nabla_v |\kappa_T| = - \frac{\cos \theta}{\sin^2 \theta} \, \kappa_\parallel \, \frac{\kappa_T}{|\kappa_T|} - \frac{\cos \theta}{\sin \theta} \, |\kappa_T| \, e_\theta,
\label{eq:macnorm_nD_navthkapthkapt} \\
&&\hspace{-1.5cm}
\nabla_\kappa \theta = 0, \quad \nabla_\kappa \kappa_\parallel = e_\theta, \quad \nabla_\kappa |\kappa_T| = \frac{1}{\sin \theta} \, \frac{\kappa_T}{|\kappa_T|} . \label{eq:macnorm_nD_nakapthkapthkapt} 
\end{eqnarray}
\label{lem:macnorm_nD_nathkapthkapt}
\end{lemma}

\begin{lemma}[Gradients of $w_i$ and $\kappa_{Ti}/|\kappa_T|$]~

\noindent
For all $i \in \{1, \ldots, n-1\}$, we have 
\begin{eqnarray}
&& \hspace{-1.5cm}
\nabla_v w_i = \frac{1}{\sin \theta} \, \big( P_{v^\bot} \mathbf{e}_i - \cos \theta  \, w_i  \, e_\theta \big) , \qquad \nabla_\kappa w_i = 0, \label{eq:macnorm_nD_nawi} \\
&& \hspace{-1.5cm}
\nabla_v \frac{\kappa_{Ti}}{|\kappa_T|} = \frac{\cos^2 \theta}{\sin^2 \theta} \frac{\kappa_\parallel}{|\kappa_T|} \, w_i \, e_\theta + \Big( - \frac{w_i}{\sin \theta} + \frac{\cos \theta}{\sin^2 \theta} \frac{\kappa_\parallel \, \kappa_{Ti}}{|\kappa_T|^2} \Big) \frac{\kappa_T}{|\kappa_T|} - \frac{\cos \theta}{\sin^2 \theta} \frac{\kappa_\parallel}{|\kappa_T|} \, P_{v^\bot} \mathbf{e}_i, \label{eq:macnorm_nD_navkapi} \\
&& \hspace{-1.5cm}
\nabla_\kappa \frac{\kappa_{Ti}}{|\kappa_T|} = - \frac{\cos \theta}{\sin \theta} \frac{w_i}{|\kappa_T|} \, e_\theta  - \frac{1}{\sin \theta} \frac{\kappa_{Ti}}{|\kappa_T|^2} \frac{\kappa_T}{|\kappa_T|} + \frac{1}{\sin \theta \, |\kappa_T|} \, P_{v^\bot} \mathbf{e}_i. \label{eq:macnorm_nD_nakapkapi}
\end{eqnarray}
\label{lem:macnorm_nD_nawikapi}
\end{lemma}

\medskip
Thanks to these lemmas, we can state the  

\begin{proposition}[gradients of $\psi_i$ in terms of $\psi_w$ and $\psi_r$]
We have 
\begin{eqnarray}
&& \hspace{-1cm}
(\nabla_v \psi_i)^\dagger = e_\theta \bigg( w_i \, \Big( \partial_\theta \psi_w 
- \frac{\cos \theta}{\sin \theta} |\kappa_T| \, \partial_{\kappa_\bot} \psi_w - \frac{\cos \theta}{\sin \theta} \, \psi_w + \frac{\cos^2 \theta}{\sin^2 \theta} \, \frac{\kappa_\parallel}{|\kappa_T|} \, \psi_r \Big) \nonumber \\
&& \hspace{1cm} 
+ \frac{\kappa_{Ti}}{|\kappa_T|} \, \Big( \partial_\theta \psi_r - \frac{\cos \theta}{\sin \theta} \, |\kappa_T| \, \partial_{\kappa_\bot} \psi_r  \Big) \bigg) \nonumber \\
&& \hspace{0cm} 
+ \frac{\kappa_T}{|\kappa_T|} \bigg( w_i \Big( \cos \theta \, |\kappa_T| \, \partial_{\kappa_\parallel} \psi_w - \frac{\cos \theta}{\sin^2 \theta} \, \kappa_\parallel \, \partial_{\kappa_\bot} \psi_w - \frac{1}{\sin \theta} \, \psi_r \Big) \nonumber \\
&& \hspace{1cm} 
+ \frac{\kappa_{Ti}}{|\kappa_T|} \, \Big( \cos \theta \, |\kappa_T| \, \partial_{\kappa_\parallel} \psi_r - \frac{\cos \theta}{\sin^2 \theta} \, \kappa_\parallel \, \partial_{\kappa_\bot} \psi_r + \frac{\cos \theta}{\sin^2 \theta} \, \frac{\kappa_\parallel}{|\kappa_T|} \, \psi_r \Big) \bigg) \nonumber  \\
&& \hspace{0cm} 
+ P_{v^\bot} \mathbf{e}_i \Big( \frac{1}{\sin \theta} \, \psi_w - \frac{\cos \theta}{\sin^2 \theta} \, \frac{\kappa_\parallel}{|\kappa_T|} \, \psi_r \Big), 
\label{eq:macnorm_nD_nabla_v_psi_i}\\
&& \hspace{-1cm}
(\nabla_\kappa \psi_i)^\dagger = e_\theta \bigg( w_i \Big( \partial_{\kappa_\parallel} \psi_w - \frac{\cos \theta}{\sin \theta} \frac{1}{|\kappa_T|} \, \psi_r \Big) + \frac{\kappa_{Ti}}{|\kappa_T|} \, \partial_{\kappa_\parallel} \psi_r \bigg) \nonumber \\
&& \hspace{0cm} 
+ \frac{\kappa_T}{|\kappa_T|} \, \frac{1}{\sin \theta} \bigg(  w_i \, \partial_{\kappa_\bot} \psi_w + \frac{\kappa_{Ti}}{|\kappa_T|} \Big( \partial_{\kappa_\bot} \psi_r -  \frac{1}{|\kappa_T|} \psi_r \Big) \bigg) \nonumber \\
&& \hspace{0cm} 
+ P_{v^\bot} \mathbf{e}_i \,  \frac{1}{\sin \theta \, |\kappa_T|} \, \psi_r. \label{eq:macnorm_nD_diff_zetu_express_prf-1}
\end{eqnarray}
\label{prop_app:gradpsi}
\end{proposition}

\noindent
\textbf{Proof.} By the chain rule, we get 
\begin{eqnarray}
\nabla_v \psi_i 
&=& w_i \Big( \partial_\theta \psi_w \, \nabla_v \theta + \partial_{\kappa_\parallel} \psi_w \, \nabla_v \kappa_\parallel + \partial_{\kappa_\bot} \psi_w \, \nabla_v (|\kappa_T|) \Big) + \psi_w \, \nabla_v \, w_i \nonumber \\
&+& \frac{\kappa_{Ti}}{|\kappa_T|} \Big( \partial_\theta \psi_r \, \nabla_v \theta + \partial_{\kappa_\parallel} \psi_r \, \nabla_v \kappa_\parallel + \partial_{\kappa_\bot} \psi_r \, \nabla_v (|\kappa_T|) \Big) + \psi_r \, \nabla_v \, \Big( \frac{\kappa_{Ti}}{|\kappa_T|} \Big), \label{eq:macnorm_nD_vectGCI_eqs_prf2_1} 
\end{eqnarray}
and similarly for $\nabla_\kappa \psi_i$ by changing $v$ into $\kappa$ everywhere. Now, inserting the formulas found in Lemmas~\ref{lem:macnorm_nD_nathkapthkapt} and~\ref{lem:macnorm_nD_nawikapi}, in these expressions, we get the result after some straightforward computations. \endproof

As a consequence of this proposition, we get

\begin{proposition}[Expression of $|\nabla_\kappa \vec{\psi}|^2$ and $\{ \vec{\psi} , \mathbf{H}_\mathbf{u} \}$]~

\noindent
(i) $|(\nabla_\kappa \vec{\psi})^\dagger|^2$ is given by \eqref{eq:macnorm_nD_nakappsi_square}. 

\noindent
(ii) $\big\{ \vec{\psi} , \mathbf{H}_\mathbf{u} \big\}^\dagger$ is given by \eqref{eq:macnorm_nD_brack_zetu_express}. 
\label{prop:macnorm_nD_brack_zetu_express}
\end{proposition}

\medskip
\noindent
\textbf{Proof.} (i) We take the squared norm of  \eqref{eq:macnorm_nD_diff_zetu_express_prf-1}. We note that both vectors~$e_\theta$ and $\kappa_T/|\kappa_T|$ are normed and they are orthogonal to each other. On the other hand, we have 
\begin{eqnarray*}
|P_{v^\bot} e_i|^2 &=& |e_i|^2 - (e_i \cdot v)^2 = 1 - \sin^2 \, \theta w_i^2, \\
e_\theta \cdot P_{v^\bot} e_i &=& e_\theta \cdot e_i = \cos \theta \, w_i, \\
\frac{\kappa_T}{|\kappa_T|} \cdot P_{v^\bot} e_i &=& \frac{\kappa_T}{|\kappa_T|} \cdot e_i = \frac{\kappa_{Ti}}{|\kappa_T|}. 
\end{eqnarray*}
It results, after some cancellations:
\begin{eqnarray*}
&&\hspace{-1cm} 
|(\nabla_\kappa \vec{\psi}_i)^\dagger|^2 
= \Big( w_i \, \partial_{\kappa_\parallel} \psi_w  + \frac{\kappa_{Ti}}{|\kappa_T|} \, \partial_{\kappa_\parallel} \psi_r \Big)^2 
+ \frac{1}{\sin^2 \theta} \Big(  w_i \, \partial_{\kappa_\bot} \psi_w + \frac{\kappa_{Ti}}{|\kappa_T|} \, \partial_{\kappa_\bot} \psi_r \Big)^2 \\
&&\hspace{7cm} 
+ \Big( 1 - w_i^2 - \frac{\kappa_{Ti}}{|\kappa_T|} \Big) \frac{1}{\sin^2 \theta \, |\kappa_T|^2} \, \psi_r^2. 
\end{eqnarray*}
Then, using that $w$ and $\kappa_T/|\kappa_T|$ are normed and orthogonal to each other, summing over $i \in \{1, \ldots, n-1\}$, we get \eqref{eq:macnorm_nD_nakappsi_square}.

\medskip
\noindent
(ii) We have 
\begin{equation} 
\big\{ \vec{\psi} , \mathbf{H}_\mathbf{u} \big\}  = \big\langle \kappa, \nabla_v \vec{\psi} \big \rangle + \nu \big\langle P_{v^\bot} \mathbf{u}, \nabla_\kappa \vec{\psi} \big \rangle.  
\label{eq:macnorm_nD_vectGCI_eqs_prf15}
\end{equation}
We compute:
\begin{eqnarray*}
\big \langle \kappa, e_\theta \big \rangle &=& \kappa_\parallel,  \\
\big \langle \kappa, \frac{\kappa_T}{|\kappa_T|} \big \rangle &=& \sin \theta \, |\kappa_T|, \\
\big \langle \kappa, P_{v^\bot} \mathbf{e}_i \big \rangle &=& \big \langle \kappa, \mathbf{e}_i \big \rangle = \kappa_\parallel e_\theta \cdot \mathbf{e}_i + \sin \theta \, \kappa_{Ti} = \cos \theta \, \kappa_\parallel \, w_i + \sin \theta \, |\kappa_T| \, \frac{\kappa_{Ti}}{|\kappa_T|}. 
\end{eqnarray*}
and
\begin{eqnarray*}
\big \langle P_{v^\bot} \mathbf{u} , e_\theta \big \rangle &=& \mathbf{u} \cdot e_\theta = - \sin \theta,  \\
\big \langle P_{v^\bot} \mathbf{u} , \frac{\kappa_T}{|\kappa_T|} \big \rangle &=& 
\mathbf{u} \cdot \frac{\kappa_T}{|\kappa_T|} = 0, \\
\big \langle P_{v^\bot} \mathbf{u}, P_{v^\bot} \mathbf{e}_i \big \rangle &=& - (\mathbf{u} \cdot v) \, (\mathbf{e}_i \cdot v) = - \cos \theta \, \sin \theta \, w_i. 
\end{eqnarray*}
Inserting the first set of expressions into \eqref{eq:macnorm_nD_nabla_v_psi_i}, multiplying the second set by $\nu$ and insrting it into \eqref{eq:macnorm_nD_diff_zetu_express_prf-1}, and adding the two resulting expressions, 
after some cancellations, we get
\begin{equation} 
\big\{ \psi_i , \mathbf{H}_\mathbf{u} \big\}^\dagger = w_i \big( {\mathcal U}_1 \psi_w - |\kappa_T| \psi_r \big) + \frac{\kappa_{Ti}}{|\kappa_T|} \big( {\mathcal U}_1 \psi_r + |\kappa_T| \psi_r \big). 
\label{eq:prop:macnorm_nD_brack_zetu_express_prf1}
\end{equation}
Then, multiplying this expression by $\mathbf{e}_i$ and summing over $i$,  we get \eqref{eq:macnorm_nD_brack_zetu_express}. \endproof

\subsection{Proof of Lemma \ref{lem:macnorm_nD_nathkapthkapt}}
\label{subsec:proof_lem:macnorm_nD_nathkapthkapt}

 We first show the following: 
\begin{equation}
\theta = \cos^{-1} (v \cdot \mathbf{u}), \qquad \kappa_\parallel = - \frac{\kappa \cdot \mathbf{u}}{\sin \theta}, \qquad |\kappa_T|= \frac{( |\kappa|^2 - \kappa_\parallel^2 )^{1/2}}{\sin \theta} . \label{eq:macnorm_nD_thet_kapthet_|kap|2}
\end{equation}
Indeed, the first formula \eqref{eq:macnorm_nD_thet_kapthet_|kap|2} directly comes from the first formula \eqref{eq:macnorm_V_def_2}. The second formula~\eqref{eq:macnorm_nD_thet_kapthet_|kap|2} comes from the second and third formulas \eqref{eq:macnorm_wethet_formulas} and from the fact that $\kappa \cdot v = 0$. Finally, the third formula \eqref{eq:macnorm_nD_thet_kapthet_|kap|2} comes from \eqref{eq:macnorm_dV_express} and the fourth formula \eqref{eq:macnorm_wethet_formulas}.

\medskip 
\noindent
Then, taking the gradient with respect to $v$ of the first expression in \eqref{eq:macnorm_nD_thet_kapthet_|kap|2} and using the first expression \eqref{eq:geom_nab_fct_v_only} and the second equation \eqref{eq:macnorm_wethet_formulas}, we get 
$$ \nabla_v \theta = - \big( 1 - (v \cdot \mathbf{u})^2 \big)^{-1/2} \nabla_v (v \cdot \mathbf{u}) = - \frac{1}{\sin \theta} \, P_{v^\bot} \mathbf{u} = e_\theta, $$
which is the first formula \eqref{eq:macnorm_nD_navthkapthkapt}.

\medskip 
\noindent
Now, by the second expression in \eqref{eq:macnorm_nD_thet_kapthet_|kap|2}, and thanks to the previous identity, we get 
\begin{equation}
 \nabla_v \kappa_\parallel = \frac{\cos \theta}{\sin^2 \theta} \, (\kappa \cdot \mathbf{u}) \, \nabla_v \theta - \frac{1}{\sin \theta} \, \nabla_v (\kappa \cdot \mathbf{u}) = - \frac{\cos \theta}{\sin \theta} \, \kappa_\parallel \, e_\theta - \frac{1}{\sin \theta} \, \nabla_v (\kappa \cdot \mathbf{u}). 
\label{eq:macnorm_nD_nathkapthkapt_prf1}
\end{equation}
Now, let $a \in T_v B$ and denote by $\alpha$ the pair $(v, \kappa)$ as usual. We compute 
\begin{eqnarray*}  
\big \langle \nabla_v (\kappa \cdot \mathbf{u}), a \big\rangle_v &=& \big \langle {\mathcal B}^H_\alpha \nabla_\alpha (\kappa \cdot \mathbf{u}) , a \big\rangle_v = \big\langle \hspace{-1.8mm} \big\langle \nabla_\alpha (\kappa \cdot \mathbf{u}) , {\mathcal L}^H_\alpha a \big\rangle \hspace{-1.8mm} \big\rangle_\alpha \\
&=& \big( d (\kappa \cdot \mathbf{u}) \big)_\alpha ( {\mathcal L}^H_\alpha a ) = \frac{d}{dt} (\tilde \kappa(t) \cdot \mathbf{u}) \big|_{t=0}, 
\end{eqnarray*}
where $\tilde \alpha(t) = (\tilde v(t), \tilde \kappa(t))$ is a path in $M$ such that 
$\tilde \alpha(0) = (v,\kappa)$  and $(d/dt) (\tilde v(t), \tilde \kappa(t))|_{t=0} = {\mathcal L}^H_\alpha a = (a, - (a \cdot \kappa) v)$, by \eqref{eq:geom_horizmap}. Let $(f_1, \ldots, f_n)$ be an orthonormal basis such that $v=f_n$, $a = |a| f_{n-1}$, $\kappa = \tilde \kappa_{n-2} f_{n-2} + \tilde \kappa_{n-1} f_{n-1}$ with $(\tilde \kappa_{n-2} , \tilde \kappa_{n-1}) \in {\mathbb R}^2$. Then define the path $\tilde \alpha(t)$ as follows: 
$$ \tilde v(t) = \sin (|a|t) \, f_{n-1} + \cos (|a|t) \, f_n, \qquad 
\tilde \kappa(t) = \tilde \kappa_{n-2} f_{n-2} + \tilde \kappa_{n-1} f_{n-1} - \tilde \kappa_{n-1} \, \tan (|a|t) \, f_n. $$
Then, this path satisfies the above requirements and we have 
$$  \tilde \kappa(t) \cdot \mathbf{u} = \tilde \kappa_{n-2} f_{n-2} \cdot \mathbf{u} + \tilde \kappa_{n-1} f_{n-1}\cdot \mathbf{u} - \tilde \kappa_{n-1} \, \tan (|a|t) \, f_n \cdot \mathbf{u}. $$
So, 
$$ \frac{d}{dt} (\tilde \kappa(t) \cdot \mathbf{u}) \big|_{t=0} = - |a| \, \tilde \kappa_{n-1} \, f_n \cdot \mathbf{u} = - (a \cdot \kappa) \, (v \cdot \mathbf{u}) = \big \langle - \cos \theta \, \kappa , a \big \rangle_v. $$
Thus, 
$$ \nabla_v (\kappa \cdot \mathbf{u}) = - \cos \theta \, \kappa.  $$ 
Inserting this into \eqref{eq:macnorm_nD_nathkapthkapt_prf1}, we get 
$$  \nabla_v \kappa_\parallel = \frac{\cos \theta}{\sin \theta} \, ( - \kappa_\parallel \, e_\theta + \kappa) = \cos \theta \, \kappa_T, $$ 
which is the second formula \eqref{eq:macnorm_nD_navthkapthkapt}.

\medskip 
\noindent
Thanks to the third equation \eqref{eq:macnorm_nD_thet_kapthet_|kap|2}, we have 
\begin{equation} 
\nabla_v |\kappa_T| = \frac{1}{2(|\kappa|^2 - \kappa_\parallel^2 )^{1/2} \, \sin \theta} \big( \nabla_v |\kappa|^2 - 2 \kappa_\parallel \nabla_v \kappa_\parallel \big) 
- (|\kappa|^2 - \kappa_\parallel^2 )^{1/2} \frac{\cos \theta}{\sin^2 \theta} \, \nabla_v \theta. 
\label{eq:macnorm_nD_nathkapthkapt_prf1.5}
\end{equation}
Now, in view of the second formula in \eqref{eq:geom_na_kappa2}, we have $\nabla_v |\kappa|^2 = 0$. Inserting the first two equations of \eqref{eq:macnorm_nD_navthkapthkapt} into \eqref{eq:macnorm_nD_nathkapthkapt_prf1.5}, we get the third equation \eqref{eq:macnorm_nD_navthkapthkapt}.

\medskip 
\noindent
The first equation \eqref{eq:macnorm_nD_nakapthkapthkapt} follows directly from the second equation \eqref{eq:geom_nab_fct_v_only}.

\medskip 
\noindent
By the second equation \eqref{eq:geom_nab_fct_v_only} and the second equation \eqref{eq:macnorm_nD_thet_kapthet_|kap|2}, we get 
\begin{equation}
\nabla_\kappa \kappa_\parallel = - \frac{1}{\sin \theta} \nabla_\kappa (\kappa \cdot \mathbf{u} ). 
\label{eq:macnorm_nD_nathkapthkapt_prf2}
\end{equation}
Now, let $\tau \in T_v B$. We compute 
\begin{eqnarray*} 
\big \langle \nabla_\kappa (\kappa \cdot \mathbf{u}) , \tau \big\rangle_v &=& \big \langle {\mathcal B}^V_\alpha \nabla_\alpha (\kappa \cdot \mathbf{u}) , \tau \big\rangle_v = \big\langle \hspace{-1.8mm} \big\langle \nabla_\alpha (\kappa \cdot \mathbf{u}) , {\mathcal L}^V_\alpha \tau \big\rangle \hspace{-1.8mm} \big\rangle_\alpha \\
&=& \big( d (\kappa \cdot \mathbf{u}) \big)_\alpha ( {\mathcal L}^V_\alpha \tau ) = \frac{d}{dt} (\tilde \kappa(t) \cdot \mathbf{u}) \big|_{t=0}, 
\end{eqnarray*}
where $\tilde \alpha(t) = (\tilde v(t), \tilde \kappa(t))$ is a path in $M$ such that 
$\tilde \alpha(0) = (v,\kappa)$  and $(d/dt) (\tilde v(t), \tilde \kappa(t))|_{t=0} = {\mathcal L}^V_\alpha \tau = (0, \tau)$, by \eqref{eq:geom_omega}. The path $\tilde \alpha(t) = (v, \kappa + t \tau)$ satisfies these requirements. We have $\tilde \kappa(t) \cdot \mathbf{u} = (\kappa \cdot \mathbf{u}) + t (\tau \cdot \mathbf{u})$, hence $(d/dt) (\tilde \kappa(t) \cdot \mathbf{u})|_{t=0} = (\tau \cdot \mathbf{u}) = \langle \tau, P_{v^\bot} \mathbf{u} \rangle_v$, showing that 
$$ \nabla_\kappa (\kappa \cdot \mathbf{u}) = P_{v^\bot} \mathbf{u} = - \sin \theta \, e_\theta. $$
Inserting this formula into \eqref{eq:macnorm_nD_nathkapthkapt_prf2} leads to the  second formula \eqref{eq:macnorm_nD_nakapthkapthkapt}.

\medskip 
\noindent
Finally, thanks to the first formula in \eqref{eq:geom_na_kappa2}, we get 
$$ \nabla_\kappa |\kappa_T| = \frac{1}{\sin \theta \, (|\kappa|^2 - \kappa_\parallel^2 )^{1/2}} \big( \kappa - \kappa_\parallel \, \nabla_\kappa \kappa_\parallel \big), $$
and thanks to the second equation \eqref{eq:macnorm_nD_nakapthkapthkapt}, this leads to the third equation \eqref{eq:macnorm_nD_nakapthkapthkapt}. This ends the proof of Lemma \ref{lem:macnorm_nD_nathkapthkapt}. \endproof

\subsection{Proof of Lemma \ref{lem:macnorm_nD_nawikapi}}
\label{subsec:proof_lem:macnorm_nD_nawikapi}

With the first equation \eqref{eq:macnorm_wethet_formulas}, we have 
$ w_i = (1 - (v \cdot \mathbf{u})^2)^{-1/2} v_i$. So, 
\begin{eqnarray*}
\nabla_v w_i &=& \big(1 - (v \cdot \mathbf{u})^2 \big)^{-1/2} \, \nabla_v (v \cdot \mathbf{e}_i) + \big(1 - (v \cdot \mathbf{u})^2 \big)^{-3/2} \, v_i \, (v \cdot \mathbf{u}) \, \nabla_v (v \cdot \mathbf{u}) , \\
&=& \frac{1}{\sin \theta} \, P_{v^\bot} \mathbf{e}_i + \frac{\cos \theta}{\sin^3 \theta}  \, v_i \, P_{v^\bot} \mathbf{u} 
= \frac{1}{\sin \theta} \, P_{v^\bot} \mathbf{e}_i - \frac{\cos \theta}{\sin \theta}  \, w_i  \, e_\theta, 
\end{eqnarray*}
which is the first formula \eqref{eq:macnorm_nD_nawi}. 

\medskip
\noindent
The second formula \eqref{eq:macnorm_nD_nawi} follows directly from the second formula of \eqref{eq:geom_nab_fct_v_only}.

\medskip 
\noindent
With the first equation \eqref{eq:macnorm_dV_express_2}, we get 
\begin{equation}
e_\theta \cdot \mathbf{e}_i = \cos \theta \, w_i. 
\label{eq:macnorm_nD_nawikapi_prf0}
\end{equation}
Then, thanks to \eqref{eq:macnorm_dV_express} and the third equation \eqref{eq:macnorm_nD_thet_kapthet_|kap|2}, we have 
\begin{equation}
\frac{\kappa_{Ti}}{|\kappa_T|} = \frac{\kappa_i - \kappa_\parallel \, \cos \theta \, w_i}{(|\kappa|^2 - \kappa_\parallel^2)^{1/2}}. 
\label{eq:macnorm_nD_nawikapi_prf1}
\end{equation}
We deduce that 
\begin{eqnarray*}
&&\hspace{-1cm}
\nabla_v \frac{\kappa_{Ti}}{|\kappa_T|} = \frac{1}{(|\kappa|^2 - \kappa_\parallel^2)^{1/2}} \big( \nabla_v \kappa_i - \cos \theta \, w_i \, \nabla_v \kappa_\parallel +  \kappa_\parallel \, w_i \, \sin \theta \, \nabla_v \theta - \kappa_\parallel \, \cos \theta \, \nabla_v w_i \big) \\
&&\hspace{3.5cm}
- \frac{1}{2} \frac{1}{(|\kappa|^2 - \kappa_\parallel^2)^{3/2}} \big( \nabla_v |\kappa|^2 - 2 \kappa_\parallel \, \nabla_v \kappa_\parallel \big) \big( \kappa_i - \kappa_\parallel \, \cos \theta \, w_i \big). 
\end{eqnarray*}
By the same proof as for $\nabla_v (\kappa \cdot \mathbf{u})$ in the previous lemma, we have 
$$ \nabla_v \kappa_i = \nabla_v (\kappa \cdot \mathbf{e}_i) = - (v \cdot \mathbf{e}_i) \, \kappa = - \sin \theta \, w_i \, (\kappa_\parallel \, e_\theta + \sin \theta \, \kappa_T). $$
Then, the second equation \eqref{eq:geom_na_kappa2}, the first equation \eqref{eq:macnorm_nD_nawi} and Lemma \ref{lem:macnorm_nD_nathkapthkapt} lead to
\begin{eqnarray*}
&&\hspace{-1cm}
\nabla_v \frac{\kappa_{Ti}}{|\kappa_T|} = \frac{1}{\sin \theta \, |\kappa_T|} \Big( - \sin \theta \, w_i \, (\kappa_\parallel \, e_\theta + \sin \theta \, \kappa_T) - \cos \theta^2 \, w_i \, \kappa_T +   \kappa_\parallel \, w_i \, \sin \theta \, e_\theta \\
&&\hspace{1cm}
- \kappa_\parallel \, \frac{\cos \theta}{\sin \theta} \big( P_{v^\bot} \mathbf{e}_i - \cos \theta  \, w_i  \, e_\theta  \big) \Big) + \frac{\cos \theta}{\sin^3 \theta} \, \frac{\kappa_\parallel}{|\kappa_T|^2} \big( \kappa_i - \kappa_\parallel \, \cos \theta \, w_i \big) \, \frac{\kappa_T}{|\kappa_T|}, 
\end{eqnarray*}
which, after rearrangement, leads to \eqref{eq:macnorm_nD_navkapi}.

\medskip 
\noindent
Similarly, thanks to \eqref{eq:macnorm_nD_nawikapi_prf1}, we obtain 
\begin{eqnarray*}
&&\hspace{-1cm}
\nabla_\kappa \frac{\kappa_{Ti}}{|\kappa_T|} = \frac{1}{(|\kappa|^2 - \kappa_\parallel^2)^{1/2}} \big( \nabla_\kappa \kappa_i - \cos \theta \, w_i \, \nabla_\kappa \kappa_\parallel \big) \\
&&\hspace{3.5cm}
- \frac{1}{2} \frac{1}{(|\kappa|^2 - \kappa_\parallel^2)^{3/2}} \big( \nabla_\kappa |\kappa|^2 - 2 \kappa_\parallel \, \nabla_\kappa \kappa_\parallel \big) \big( \kappa_i - \kappa_\parallel \, \cos \theta \, w_i \big)
\end{eqnarray*}
By the same proof as for $\nabla_\kappa (\kappa \cdot \mathbf{u})$ in the previous lemma, we have 
$$ \nabla_\kappa \kappa_i = \nabla_\kappa (\kappa \cdot \mathbf{e}_i) = P_{v^\bot} \mathbf{e}_i. $$
Then, the first equation \eqref{eq:geom_na_kappa2} and Lemma \ref{lem:macnorm_nD_nathkapthkapt} lead to
$$ \nabla_\kappa \frac{\kappa_{Ti}}{|\kappa_T|} = \frac{1}{\sin \theta \, |\kappa_T|} \big( P_{v^\bot} \mathbf{e}_i -  \cos \theta \, w_i \, e_\theta \big) - \frac{1}{\sin^3 \theta \, |\kappa_T|^3} \big( \kappa - \kappa_\parallel \, e_\theta \big) \big( \kappa_i - \kappa_\parallel \, \cos \theta \, w_i \big), $$
which, after rearrangement, yields \eqref{eq:macnorm_nD_nakapkapi}. \endproof

\end{document}